\documentclass[a4paper,11pt,twoside]{book}
\RequirePackage[a4paper,head=1cm]{geometry} 

\usepackage{fancyhdr} 
\usepackage{verbatim}   
\fancypagestyle{pagevide}{%
	\fancyhead{}   
	\fancyfoot[C]{}  
} 
\fancypagestyle{plain}{%
	\fancyhead{} 
	\fancyfoot[C]{\thepage}
}
\usepackage{pstcol,pst-fill,pst-grad}
\usepackage{textcomp}
\usepackage[french]{babel}
\usepackage[utf8]{inputenc}
\usepackage{graphicx}
\usepackage{amsmath}
\usepackage{amsfonts}
\usepackage{amssymb}
\usepackage{array}
\usepackage{psfrag}
\usepackage{dsfont}

\usepackage{dsfont}
\newtheorem{Definition}{Définition}[section]
\newtheorem{Theoreme}{Théorème}[section]

\newtheorem{Lemme}{Lemme}[section]
\newtheorem{Corollaire}{Corollaire}[section]
\newtheorem{Proposition}{Proposition}[section]
\newtheorem{Remarque}{\bf Remarque}[chapter]

\def \A{\vec{A}}
\def \B{\vec{B}} 
\def \fe{\vec{f}} 

\def \vu{\vec{u}}

\def \vv{\vec{v}}

\def \P{\mathbb{P}}

\def \Rt{\mathbb{R}^{3}}
\def \Zt{\mathbb{Z}^3}
\def \finpv{\hfill $\blacksquare$  \\ \newline }
\def \pv{{\bf{Preuve.}}~}
\def \dm{{\bf{D\'emonstration.}}~}

\def \vg{\vec{g}}

\def \U{\vec{U}} 
\def \W{\vec{W}} 
\def \V{\vec{V}}

\def \ds{\displaystyle}

\def \Zt{\mathbb{Z}^3}
\def \R{\mathbb{R}}

\begin{document}    
\thispagestyle{empty}
\vspace{-2cm}

\voffset-10pt

\noindent
\hspace*{-1cm}\hbox{\includegraphics[width=8.6cm]{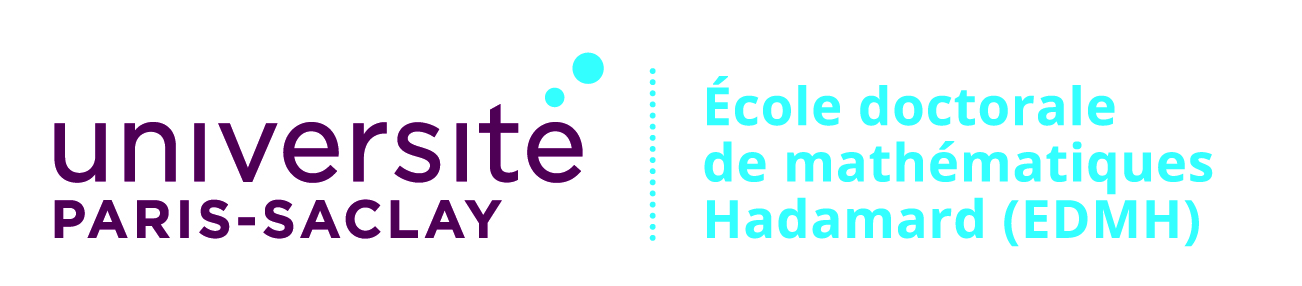}} 
\hfill
\hbox{\includegraphics[width=2cm]{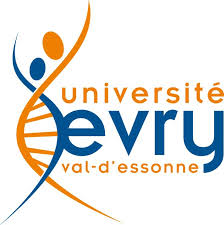}}
\hfill
\hbox{\includegraphics[width=1.6cm]{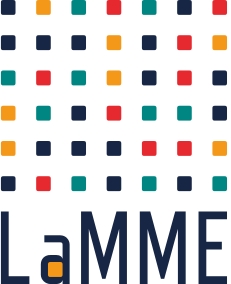}}
\vspace{7mm}

\begin{center}
	{\Large\bf TH\`ESE DE DOCTORAT}
\end{center}
\begin{center}
	{de }
\end{center}
\begin{center}
	{\Large\sc l'Universit\'e Paris-Saclay}\\
	\vspace*{0.4cm}
	\'Ecole doctorale de math\'ematiques Hadamard (EDMH, ED 574)\\  
	\vspace*{0.4cm} 
	{\small \it \'Etablissement d'inscription : } 
	    Universit\'e d'\'Evry-Val d'Essonne\\
	\vspace*{0.2cm} 
	\vspace*{0.2cm} 
	{\small \it Laboratoire d'accueil :    
	  Laboratoire de math\'ematiques et mod\'elisation d'\'Evry, UMR 8071 CNRS-INRA\\ }
	\vspace*{0.2cm}
\end{center}

\begin{center}
	{\it Sp\'ecialit\'e de doctorat : } 
	{\large Math\'ematiques appliquées}
\end{center}

\vspace{5mm}

\begin{center}
	{\large\bf Oscar JARR\'IN}
\end{center}

\vspace{3mm}

\begin{center}
	{\Large \textbf{Descriptions déterministes de la turbulence dans les équations de Navier-Stokes}}
\end{center}

\vspace{10mm}

\noindent{\small \it Date de soutenance~: } 20 Juin 2018

\vspace{5mm}

\noindent
{\small \it Apr\`es avis des rapporteurs~: }
\begin{tabular}{l} {\sc Isabelle GALLAGHER} (Université-Diderot)\vspace{1mm}  \\
	{\sc  Taoufik  HMIDI} (Université de Rennes 1)\\
\end{tabular}

\vspace{8mm}

\noindent
{\small \it Jury de soutenance~:  \\
\\}
\begin{tabular}{ll}
	{\sc   Lorenzo BRANDOLESE}&(Université Claude-Bernard-Lyon 1) {\small Examinateur}\vspace{1mm}\\
	{\sc   Diego  CHAMORRO}&(Université d'\'Evry-Val d'Essonne) {\small Codirecteur de th\`ese}\vspace{1mm}\\
	{\sc   Isabelle GALLAGHER}&(Université-Diderot) {\small Rapporteur}\vspace{1mm}\\
	{\sc   Pierre Gilles LEMARI\'E-RIEUSSET}&(Université d'\'Evry-Val d'Essonne) {\small Directeur de thèse}\vspace{1mm}\\
	{\sc   Roger  LEWANDOWSKI}&(Université de Rennes 1) {\small Président de jury}\vspace{1mm}\\
\end{tabular}

\vfill
\noindent
\hbox{\includegraphics[width=2.2cm]{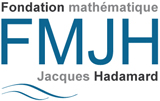}}
\hfill 
\hfill
\hbox{{\small {\bf NNT : 2018SACLE010}}}
\hfill
\hfill \includegraphics[width=1cm]{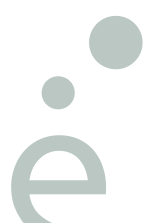}
\newpage
\thispagestyle{empty}

\emph{"$\cdots$  Il était permis à un étudiant en mathématiques d'être incroyablement ignorant. C'était mon cas, et il m'en reste une grande naïveté dans beaucoup de domaines importants. Mais j'avais et j'ai conservé le goût de résoudre des problèmes $\cdots$ "}

\begin{flushright}
Jean-Pierre Kahane 
\end{flushright} 
\newpage 
\thispagestyle{empty}
\textbf{{\huge Remerciements}}\\ 

Toute cette aventure pour venir en France et continuer ma formation scientifique a été pleine de personnes avec qui j'ai partagé de nombreux moments.  Dans ces lignes je ne fais que les  mentionner mais ces mots expriment un fort et sincère remerciement. \\  

Je tiens tout d'abord à remercier mes directeurs de thèse, Diego Chamorro et Pierre-Gilles Lemarié-Rieusset, pour toutes ces années de travail  et pour m'avoir initié à la recherche. Je les remercie non seulement pour leur qualité scientifique mais aussi pour leur qualité humaine.  Je poursuis leur bon exemple de chercheurs sérieux en envisageant toujours des résultats de la plus haute qualité  possible. \\

Je remercie également mes rapporteurs de thèse Mme. Isabelle Gallagher et M. Taoufik Hmidi qui m'ont fait le grand honneur de rapporter ce travail de recherche  et pour les commentaires précieux qu'ils ont effectués. De la même façon,  toute ma gratitude va à  M. Lorenzo Brandolese et M. Roger Lewandowski pour avoir accepté de faire partie de mon jury de soutenance.\\

Je veux remercier les membres du Laboratoire de Mathématiques et Modélisation d'\'Evry qui m'ont accompagné  pendant ma thèse. Je remercie particulièrement Arnaud Gloter,  directeur du laboratoire, et aussi les autres membres avec qui j'ai partagé de nombreux  déjeuners très agréables: Stéphane Menozzi, Christophe Profeta, Abass Sagna et Vincent Torri. \\

Un grand remerciement à Valérie Picot, secrétaire du Laboratoire de  Mathématiques d'\'Evry, pour sa grande efficacité et aussi pour sa grande gentillesse et bienveillance. Je remercie également El Maouloud Ould Baba, ingénieur informatique du laboratoire, pour toute son aide aux divers problèmes informatiques. \\

Je remercie aussi  les autres doctorants du laboratoire: Igor, Mohammed et Chiara, en leur souhaitant une bonne continuation dans leurs travaux de recherche. Un grand remerciement  à Kawther Mayoufi, avec qui nous avons  partagé beaucoup d'expériences comme étudiants en thèse,  les hautes et les basses de la recherche, ainsi que de nombreuses conférences partout en France.\\ 

Pendant ma dernière année de thèse j'ai fait aussi partie du département de mathématiques dans le cadre d'un poste d'ATER à l'Université-Dauphine . Je voudrais remercier également les membres de ce département, particulièrement Daniela Tonon et  Alexandre Afgoustidis, avec qui j'ai eu le  plaisir de travailler  en assurant leurs cours de travaux dirigés.\\

Por otro lado, quiero agradecer a las personas que han estado  presentes en todo este camino.  A Ricardo por su gran amistad y su ayuda en todo este tiempo. De igual manera a Yandira e Israel  quienes me acompañaron durante los primeros años de tesis. \\

A Joan por estar  siempre presente a pesar de la distancia y por toda la alegría que trajo en sus visitas a Par\'is. \\

A Elizabeth y Mar\'ia Jos\'e por su amistad y su apoyo. En particular por la deliciosa comida ecuatoriana de Majo.
\newpage
\thispagestyle{empty}
A Pedro, Randy y Ana Julia quienes han hecho que la fuerza ecuatoriana est\'e presente. \\

A todos ellos les deseo la m\'as exitosa carrera profesional, esperando verlos pronto convertidos en colegas con los cuales podamos trabajar juntos por la escuela matemática ecuatoriana. \\

Agradezco igualmente a Marco Calahorrano y Juan Carlos Trujillo, profesores de la Escuela Politécnica Nacional, quienes me guiaron en los primeros años de estudio de  la matemática. \\

En este punto no puedo dejar de agradecer a la Asociación AMARUN, gracias a la cual toda esta aventura de venir a Francia fue posible; y gracias a la cual tuve la oportunidad de desarrollar otras competencias tan importantes en el qu\'ehacer  de un investigador. \\

Agradezco también a los demás amigos y personas que me han acompañado de una manera u otra: Fernando, Jessi, Pato, Esteban y a mis primos Dany y Mart\'in. \\

Finalmente, un gran agradecimiento a mi familia por toda su ayuda y apoyo incondicional. En particular a mi hermano y mis padres: la culminación de esta etapa es también el resultado de todos los esfuerzos cotidianos que hicieron por mi hermano y por m\'i.   
 	
\tableofcontents
	
	
	
	
\chapter{L'étude déterministe de la loi de dissipation d'énergie}\label{Chap.1}
\section{Introduction}   

Dans ce chapitre nous allons étudier d'un point de vue déterministe la loi de dissipation d'énergie proposée par la théorie de la turbulence de Kolmogorov  également appelée \emph{``théorie K41''}. Cette loi de dissipation d'énergie fut développée par A.N. Kolmogorov dans les articles \cite{Kolm1}, \cite{Kolm2} et \cite{Kolm3} et correspond à une étude expérimentale de la quantité d'énergie cinétique dissipée (sous forme de chaleur)  dans un fluide  qui se trouve en régime turbulent.\\
\\
\'Etant donné que  les équations de Navier-Stokes donnent un modèle mathématique pour étudier le comportement des fluides, le but de ce chapitre est d'étudier aussi rigoureusement que possible  la loi de dissipation d'énergie dans le cadre des équations de Navier-Stokes déterministes que nous introduisons rapidement: nous nous intéressons donc à l'évolution  au cours du temps d'un fluide, via l'étude de son champ de vitesse en tout point de l'espace et à chaque instant. Nous supposons que le fluide est de densité constante (sa masse volumique est constante),  qu'il est incompressible (l'espace occupé par une certaine quantité de fluide à chaque instant peut changer de forme, mais pas de volume), et l'on supposera  qu'il est visqueux. Les équations qui décrivent un tel fluide sont les équations de Navier-Stokes qui s'écrivent de la façon suivante:
\begin{equation}\label{N-S-physique} 
\left\lbrace \begin{array}{lc}\vspace{3mm}
\rho \left( \partial_t\vu+(\vu\cdot \vec{\nabla}) \vu \right)=  \mu\Delta \vu  - \vec{\nabla}p_r+\vec{f}_{ext}, \quad \rho\, div(\vu)=0,\quad \mu>0, &\\
\rho \vu(0,\cdot)= \rho \vu_0,\quad \rho\,  div(\vu_0)=0.&\\
\end{array}\right. 
\end{equation}   
Dans ce système d'équations aux dérivées partielles non linéaires, la fonction vectorielle $\vu:[0,+\infty[\times \Rt\longrightarrow \Rt$ représente le champ de vitesse du fluide tandis que la fonction scalaire  $p_r:[0,+\infty[\times \Rt \longrightarrow \mathbb{R}$ représente sa pression, ces deux fonctions $\vu$ et $p_r$ sont les inconnues, tandis que la constante de densité du fluide $\rho>0$ (qui représente sa masse volumique et qui est constante pour les fluides incompressibles), la constante de viscosité dynamique du fluide $\mu>0$ (qui caractérise la résistance à l'écoulement d'un fluide incompressible),   la fonction vectorielle $\vec{f}_{ext}: [0,+\infty[\times \Rt\longrightarrow \Rt$ (qui correspond à l'ensemble des forces extérieures agissant sur le fluide) et la fonction vectorielle $\vu_0: \Rt \longrightarrow \Rt$ (qui représente la donnée initiale)  sont les données du problème.  L'équation $\rho\, div(\vu_0)=0$ et $\rho\, div(\vu)=0$ représente la condition d'incompressibilité du fluide au cours du temps. \\
\\
 Observons maintenant  que comme nous avons supposé que le fluide est de densité constante alors nous pouvons diviser les équations ci-dessus par la constante de densité $\rho>0$ et en écrivant maintenant  la force  $\fe=\frac{\vec{f}_{ext}}{\rho}$, qui correspond à la densité massique des forces agissant sur le fluide; la pression cinétique  $p=\frac{p_r}{\rho}$ et la constante de viscosité  cinétique du fluide $\nu=\frac{\mu}{\rho}$, qui est le quotient entre la viscosité dynamique et la masse volumique et qui représente la capacité de cohésion entre les particules du fluide;  nous obtenons ainsi le système de Navier-Stokes avec lequel l'on travaillera tout au long de cette thèse: 
\begin{equation}\label{N-S}  
\left\lbrace \begin{array}{lc}\vspace{3mm}
 \partial_t\vu =  \nu\Delta \vu -(\vu\cdot \vec{\nabla}) \vu - \vec{\nabla}p+\fe, \quad  div(\vu)=0,\quad \nu>0, &\\
\vu(0,\cdot)=\vu_0,\quad div(\vu_0)=0.&\\
\end{array}\right.
\end{equation} 
\`A ce stade, il convient de préciser que nous avons  tout d'abord présenté le système de Navier-Stokes (\ref{N-S-physique}) car dans la Section \ref{Sec:cadre-periodique} nous aurons besoin de revenir à ce système  pour introduire  quelques quantités physiques associées à ces équations. Néanmoins, la densité $\rho>0$ étant toujours constante  nous avons alors que les équations (\ref{N-S-physique}) sont  équivalentes aux équations (\ref{N-S}) et donc nous continuerons notre exposé en considérant ces dernières équations.\\ 
\\
L'existence des solutions faibles globales en temps du problème de Cauchy (\ref{N-S}) est un résultat classique qui a été développé par J. Leray en $1934$ \cite{Leray}: pour $\vu_0\in L^2(\Rt)$ une donnée initiale et $\fe \in L^2([0,+\infty[,\dot{H}^{-1}(\Rt))$ une force extérieure, les solutions faibles $(\vu, p)$ du problème (\ref{N-S}) vérifient:
 $$ \vu \in L^{\infty}([0,+\infty[,L^{2}(\Rt))\cap L^{2}([0,+\infty[,\dot{H}^{1}(\Rt)),$$
 $$ p \in L^{2}([0,+\infty[,\dot{H}^{-\frac{1}{2}}(\Rt)),$$ et de plus elles  vérifient certaines propriétés fondamentales comme par exemple l'inégalité d'énergie: pour tout $T>0$, on a
$$ \Vert \vu(T,\cdot)\Vert^{2}_{L^2}+2\nu \int_{0}^{T}\Vert \vec{\nabla} \otimes \vu(t,\cdot)\Vert^{2}_{L^2}dt\leq \Vert \vu_0\Vert^{2}_{L^2}+2\int_{0}^{T}\int_{\Rt}\fe(t,x)\cdot \vu(t,x)\,dxdt,$$
cette inégalité sera un outil qui sera largement utilisé par la suite.  \\ 
\\ 
Pour étudier  la loi de dissipation d'énergie de Kolmogorov dans le modèle déterministe des équations de Navier-Stokes, nous allons considérer les solutions faibles de Leray de ces équations car, comme nous expliquerons plus en détail dans la Section \ref{Sec:cadre-periodique} ci-dessous, nous aurons besoin de travailler avec des solutions globales en temps. \\
\\    
Une fois que nous avons introduit rapidement les équations de Navier-Stokes, nous nous concentrons maintenant sur la loi de dissipation d'énergie proposée par la théorie K41 et pour énoncer cette loi, tout d'abord, nous avons besoin de faire  une courte introduction sur l'idée phénoménolo\-gique sous-jacente.  \\
\\
Le but de la théorie K41 est de décrire  le comportement du champ de vitesse $\vu(t,x)$  dans un certain intervalle de temps $t$ et aux échelles de longueur $\ell=\vert x \vert$  où le fluide se trouve dans un état de turbulence pleine. Indiquons rapidement que les différentes échelles de longueur qui interviennent dans l'étude de la turbulence seront expliquées par le biais du modèle de cascade d'énergie et nous y reviendrons plus tard.\\
\\
 Maintenant il convient de faire  une discussion importante sur l'intervalle du temps où l'on s'attend à étudier la turbulence. \\
\\  
\`A  partir d'une donnée initiale $\vu_0\in L^2(\Rt)$ qui représente le champs de vitesse à l'instant $t=0$ et d'une force $\fe \in L^2([0,+\infty[, \dot{H}^{-1})$ nous nous intéressons à étudier  l'évolution du fluide au cours du temps au moyen de son champ de vitesse $\vu(t,\cdot)$ qui est une solution de Leray des équations de Navier-Stokes. Nous allons supposer que la force $\fe$ agit suffisamment fort sur le fluide de sorte qu'à partir d'un certain temps $t_1>0$ le fluide est en état de turbulence pleine.   De plus, si nous supposons que cette force vérifie $\fe \in L^2([0,+\infty[,  L^2(\Rt))$ alors on sait que la solution de Leray vérifie 
$ \ds{\lim_{t\longrightarrow +\infty} \Vert \vu(t,\cdot)\Vert_{L^2}=0}$ (voir le livre  \cite{PGLR1}, Corollaire $12,1$ page 357 pour une preuve de ce résultat) et alors nous observons qu'à partir d'un certain temps $t_2>t_1$  le fluide quitte son état de turbulence pleine; et donc  nous sommes censés étudier la turbulence dans l'intervalle de temps $]t_1,t_2[$. \\
\\
Néanmoins,  il n'est pas totalement trivial de donner une estimation précise de cet intervalle  et dans l'état actuel de nos connaissances cette question n'a pas de réponse satisfaisante, que ce soit du point de vue mathématique ou du point de vue physique (voir les articles \cite{FMRT2} et \cite{Vigneron}).  Nous pouvons alors remarquer  que  l'étude déterministe de la turbulence  dans le cadre d'une force extérieure qui dépend de la variable temporelle est un problème très compliqué car  nous sommes obligés de trouver un certain  intervalle $]t_1,t_2[$  où l'on puisse assurer que le fluide est en turbulence pleine.  Ce problème ne sera donc pas considéré ici et de cette façon nous considérons dorénavant $\fe$ une force extérieure  \emph{stationnaire.} \\ 
\\
Cette condition d'une force  stationnaire pour l'étude déterministe de la turbulence est largement considérée dans la littérature  (voir par exemple l'article \cite{Childress} de S. Childress, les notes du cours  \cite{Const} de P. Constantin, les articles \cite{DoerFoias,FMRT1,FMRT3} de F. Foias,  R. Temam \emph{et al.},  ainsi que le livre \cite{FMRTbook} de ces derniers auteurs) et se base sur l'idée suivante: comme $\fe$ est une fonction qui ne dépend que de la variable d'espace alors cette force agit sur le fluide de façon \emph{indépendante} du temps et donc une fois que le fluide est en régime turbulent nous pouvons supposer que ce fluide restera dans ce régime turbulent au cours du temps.  Dans le monde réel nous pouvons observer aussi  des fluides qui sont constamment en état turbulent: si l'on regarde, par exemple, une cascade d'eau de grande hauteur alors on observe que l'eau qui tombe par cette cascade est toujours en état turbulent.\\
\\
Ainsi, la force $\fe$ étant une fonction stationnaire  nous allons alors faire une étude déterministe de la théorie K41 dans le régime asymptotique lorsque le temps tend vers l'infini  et ce régime  asymptotique sera caractérisé par  des moyennes en temps long sur le champ de vitesse du fluide  $\vu$ (voir la Définition \ref{Def_U_varepsilon_per}, page \pageref{Def_U_varepsilon_per}  pour une définition précise). Ce régime en temps long  nous permettra donc d'étudier la théorie  K41 en prenant en compte seulement les différentes  échelles de longueur où le fluide se trouve en état turbulent. \\
\\
Expliquons maintenant  comment la turbulence peut être visualisée lorsqu'on regarde ces échelles de longueur. Cette théorie se base sur le modèle de \emph {cascade d'énergie}, introduite par L. F. Richardson en $1922$ dans son livre \cite{Richardson} et formalisée par A.N. Kolmogorov dans son article \cite{Kolm1} en $1941$, qui postule que si l'énergie cinétique est introduite dans le fluide par l'action d'une force extérieure $\fe$, alors, dans un régime turbulent  le mécanisme  de dissipation d'énergie sous forme de chaleur  (dû aux  forces de viscosité du fluide) n'est pas effectif et l'action de l'énergie cinétique est expliquée par une ``cascade d'énergie". \\
\\
Cette  cascade d'énergie explique tout simplement que les grands tourbillons se \emph{cassent} en des tourbillons plus petits. Soyons un peu plus précis: si la force extérieure agit (ou se fait ressentir) à une échelle $\ell_0>0$, échelle qui sera  nommée  ``l'échelle d'injection d'énergie" et si $\varepsilon_0>0$ est le taux d'injection d'énergie (qui provient de divers facteurs) à l'échelle $\ell_0$, alors l'énergie cinétique, introduite à cette échelle $\ell_0$ au taux $\varepsilon_0>0$, est transférée aux échelles de longueur plus petites $\ell>0$ (avec $\ell_0>>\ell$) à un certain taux $\varepsilon_T>0$.\\ 
\\
Ainsi, cette cascade d'énergie continue jusqu'à l'échelle de longueur $\ell_D$, nommée  \emph{l'échelle de dissipation de Kolmogorov} (voir l'expression (\ref{Def_l_D})  pour une définition plus précise) et à cette échelle  l'énergie cinétique provenant de l'échelle  de longueur supérieure $\ell$ (avec $\ell_D<\ell<<\ell_0$) est finalement dissipée sous forme de chaleur (au taux $\varepsilon_D>0$) par l'action directe des forces de viscosité du fluide.  \\
\\
Cette cascade d'énergie  a lieu pour les échelles de longueur $\ell$ qui appartiennent à l'intervalle $ ]\ell_D,\ell_0[$,   nommé \emph{l'intervalle d'inertie}, et  pour faire une exposition plus complète du modèle cascade d'énergie nous nous concentrons  sur les échelles $\ell_0$ et $\ell_D$ qui définissent cet intervalle d'inertie.  \\
\\
Nous avons déjà mentionné que l'échelle d'injection d'énergie $\ell_0>0$ est l'échelle de longueur à laquelle la force extérieure $\fe$ introduit l'énergie cinétique dans le fluide. On peut penser, par exemple, à un agitateur de taille $\ell_0$ qui en agitant un fluide visqueux et incompressible lui communique constamment de l'énergie cinétique.  Cette échelle $\ell_0$ est un paramètre du modèle qui sera donc fixé par la force $\fe$.  \\
\\
Une fois que l'énergie est injectée à l'échelle $\ell_0$, cette énergie est transférée par les forces d'inertie du fluide vers des échelles de plus en plus petites et ce processus s'arrête lorsqu'on arrive à l'échelle $\ell_D$ qui est suffisamment  petite ($\ell_D<<\ell_0$) pour que l'énergie y soit dissipée sous forme de chaleur.  Dans sa théorie K41,  Kolmogorov  introduit l'idée que cette échelle  $\ell_D$  ne dépend que du taux de dissipation d'énergie $\varepsilon_D$ et de la constante de viscosité du fluide $\nu$ et  Kolmogorov suggère de définir  l'échelle de dissipation d'énergie par  la relation
\begin{equation}\label{l_D}
\ell_D=(\varepsilon_D)^{\alpha}\nu^{\beta},
\end{equation} où l'on cherche a déterminer  les exposants $\alpha,\beta\in \mathbb{R}$ qui ne dépendent pas du fluide (voir les articles \cite{Kolm1}, \cite{Kolm2} et \cite{Kolm3} ou le livre \cite{McDonough} pour une discussion à ce sujet et pour plus de détails). \\
\\
De cette façon, pour déterminer les exposants $\alpha,\beta$ ci-dessus, Kolmogorov réalise une analyse des dimensions physiques: en effet,  observons tout d'abord que  le taux de dissipation $\varepsilon_D$ a une dimension physique $\frac{longueur^2}{temps^3}$. En effet, si $U>0$ est la vitesse ``caractéristique" du fluide, qui représente la vitesse moyenne du mouvement (voir la Définition \ref{Def_U_varepsilon_per}, page \pageref{Def_U_varepsilon_per}), nous savons que la quantité  $U^2$ est la quantité moyenne d'énergie cinétique du fluide (voir les livres \cite{JacTab} et \cite{McDonough}) et comme la vitesse caractéristique $U$ a une une dimension physique  $\frac{longueur}{temps}$ alors la quantité moyenne d'énergie $U^2$ a une dimension physique   $\frac{longueur^2}{temps^2}$. Ensuite, étant donné que le taux de dissipation d'énergie $\varepsilon_D$ mesure  la variation de l'énergie cinétique par rapport au temps (voir toujours les livres \cite{JacTab} et \cite{McDonough}) nous savons que $\varepsilon_D$ a une dimension physique  $\frac{\text{énergie}}{temps}$ et comme l'énergie est mesurée en $\frac{longueur^2}{temps^2}$ nous obtenons ainsi que $\varepsilon_D$ a une dimension  $\frac{longueur^2}{temps^3}$. D'autre part la constante de viscosité $\nu$ caractérise la résistance du milieu à un écoulement écoulement uniforme et cette constante  a  une  dimension physique  $\frac{longueur^2}{temps}$  (voir le livre \cite{McDonough} pour plus de détails).   \\ 
\\
Ainsi,  le terme à droite de (\ref{l_D}) a une dimension physique  $\frac{longueur^{2\alpha +2\beta}}{temps^{3\alpha+\beta}}$, mais, comme ce terme doit avoir une dimension de $longueur$ alors il faut que l'on ait les identités $2\alpha +2\beta=1$ et $3\alpha+\beta=0$ d'où nous avons  $\alpha=-\frac{1}{4}$ et $\beta=\frac{3}{4}$. De cette façon, en remplaçant ces valeurs  de $\alpha$ et $\beta$ dans (\ref{l_D}), Kolmogorov  a obtenu que  l'échelle de dissipation $\ell_D$ est donnée par 
\begin{equation}\label{Def_l_D}
\ell_D=\left( \frac{\nu^3}{\varepsilon_D}\right)^{\frac{1}{4}}.  
\end{equation}
Nous observons que cette échelle n'est pas un paramètre du modèle (contrairement à l'échelle d'injection d'énergie $\ell_0$) car $\ell_D$ dépend de la  quantité $\varepsilon_D$ qui elle même dépend finalement de la solution $\vu$ des équations (\ref{N-S}) (voir toujours la Définition \ref{Def_U_varepsilon_per}). \\
\\
Une fois que nous avons introduit  le modèle de cascade d'énergie qui est valable dans l'intervalle $]\ell_D,\ell_0[$ nous pouvons maintenant présenter la loi de dissipation d'énergie de Kolmogorov dont le domaine de validité est donné par ce même intervalle. 
Ainsi, conformément au modèle de cascade d'énergie, nous savons que l'énergie cinétique introduite dans le fluide   à l'échelle $\ell_0$  est transférée à un taux $\varepsilon_T>0$ par les forces inertielles du fluide jusqu'à l'échelle $\ell_D$. \\
\\
Le but de la loi de dissipation d'énergie est de donner une expression quantitative du taux  $\varepsilon_D>0$ auquel cette énergie cinétique est alors dissipée par les force visqueuses.\\
\\
Pour cela, tout d'abord, nous avons besoin de caractériser l'état turbulent de ce fluide et  nous allons introduire rapidement les nombres de Reynolds (voir  la  page \pageref{Def_U_varepsilon_per}, pour une exposition  plus complète sur ces  nombres). Pour  $U>0$ la vitesse caractéristique du fluide, $\nu>0$ sa constante de viscosité et $\ell_0>0$ l'échelle d'injection d'énergie,  le nombre de  Reynolds $Re$ peut être défini par la relation 
\begin{equation}\label{Re_lo_intro}
Re=\frac{U\ell_0}{\nu},
\end{equation} 
et ce  nombre  mesure l'importance relative des forces inertielles  qui transfèrent l'énergie cinétique (représentées par le numérateur $U\ell_0$)   sur les forces  visqueuses qui dissipent cette énergie (qui sont représentées par le dénominateur $\nu$).  Ainsi, l'état turbulent du fluide est caractérisé lorsque $$Re>>1,$$ ce qui représente le fait que les forces inertielles sont beaucoup plus fortes que les effets visqueux et alors la cascade d'énergie décrite ci-dessus a lieu. \\
\\
Dans ce régime turbulent des grandes valeurs du nombre de Reynolds,  on peut observer de façon expérimentale (voir par exemple les expériences physiques \cite{Houel} et \cite{Tennekes})  que le taux d'injection d'énergie  $\varepsilon_0>0$ à l'échelle $\ell_0$  est du même ordre de grandeur que le taux  de transfert $\varepsilon_T$ et est aussi du même ordre de grandeur que le taux de dissipation d'énergie $\varepsilon_D>0$ à l'échelle $\ell_D$, c'est à dire, nous avons $$ \varepsilon_0 \approx  \varepsilon_T\approx \varepsilon_D,$$ et alors, comme on veut donner une expression quantitative de $\varepsilon_D$, grâce à cette estimation\footnote{Il convient de préciser la notation $`` \approx"$. Pour $a,b\in \R$ nous dirons que le nombre $a$ est du même ordre de grandeur que le nombre $b$ s'il existent deux constates $c_1,c_2>0$, qui ne dépendent d'aucun paramètre physique  telles que $c_1 a \leq b \leq c_2 b$. Lorsque les nombres $a$ et $b$ sont du même ordre de grandeur nous écrirons $a\approx b$.} il est équivalent de donner une expression quantitative de $\varepsilon_0$. \\
\\
Ainsi, pour estimer $\varepsilon_0$ Kolmogorov considère les trois paramètres physiques suivants: l'échelle d'injection d'énergie $\ell_0>0$,  la vitesse caractéristique du fluide $U>0$ et la constante de viscosité du fluide $\nu>0$. \\
\\
\'Etant donné qu'on se trouve dans le régime turbulent caractérisé par $Re>>1$, alors les forces visqueuses du fluide sont négligeables par rapport aux  forces inertielles, ce qui amène  Kolmogorov à faire l'hypothèse  suivante: \emph{le taux d'injection d'énergie $\varepsilon_0$ est indépendant de la constante de viscosité $\nu$}. \\ 
\\
En supposant cette hypothèse, Kolmogorov introduit l'idée que le taux $\varepsilon_0$ doit alors  s'exprimer seulement en fonction des paramètres $U$ et $\ell_0$, ce qui l'amène à écrire l'estimation $\ds{\varepsilon_0 \approx U^{a}\ell^{b}_{0}}$. Pour déterminer les exposants $a, b \in \R$, en procédant par une analyse des dimensions physiques de la même façon que l'on a fait pour obtenir l'échelle $\ell_D$ donnée dans (\ref{Def_l_D})  on obtient  $\ds{\varepsilon_0 \approx \frac{U^3}{\ell_0}}$. \\
\\ 
De cette façon,  comme on a (expérimentalement) $\varepsilon_0\approx \varepsilon_D$, on obtient alors l'estimation du taux de dissipation d'énergie  $\varepsilon_D\approx \frac{U^3}{\ell_0}$. \\
\\
Finalement, pour simplifier la notation,  nous allons écrire dorénavant le taux de dissipation d'énergie $\varepsilon_D$ comme $\varepsilon$ et ainsi nous énonçons la \emph{loi de dissipation de Kolmogorov} pour des fluides en état turbulent: \vspace{5mm}
\begin{center}
\fbox{
\begin{minipage}[l]{90mm}
\begin{equation}\label{loi_dissipation_physique}
\text{si} \quad Re>>1\quad \text{alors}\quad \varepsilon \approx \frac{U^3}{\ell_0}. 
\end{equation}  
\end{minipage}
} 
\end{center}  
\vspace{5mm}
Nous rappelons que notre but est l'étude déterministe de cette loi de dissipation de Kolmogorov (\ref{loi_dissipation_physique})  pour un fluide visqueux et incompressible modélisé mathématiquement par les équations de Navier-Stokes (\ref{N-S})  dans l'espace $\Rt$ \emph{tout entier.} \\
\\
 Néanmoins, pour mieux comprendre les enjeux de cette théorie, dans la  section qui suit nous analyserons le cadre d'un fluide périodique en variable d'espace. Le passage par un cadre périodique nous permettra, tout d'abord, de faire une courte introduction de l'état actuel des connaissances sur l'étude de la loi (\ref{loi_dissipation_physique}) et ensuite de mettre en perspective  quelques difficultés lorsqu'on considère un fluide non périodique dans tout l'espace $\Rt$ comme nous le verrons à la fin de cette introduction.    
\subsection{Le cadre d'un fluide périodique}\label{Sec:cadre-periodique}
 Nous considérons donc une longueur $L>0$ (qui sera la période) et  les équations de Navier-Stokes  périodiques sur le cube $[0,L]^3\subset \Rt$:
\begin{equation}\label{N-S-per}
\left\lbrace \begin{array}{lc}\vspace{3mm}
\partial_t\vu=\nu\Delta \vu-(\vu\cdot \vec{\nabla}) \vu -\vec{\nabla}p+\fe, \quad div(\vu)=0,\quad \nu>0, &\\
\vu(0,\cdot)=\vu_0,&\\
\end{array}\right.
\end{equation}
où le champ de vitesse $\vu:[0,+\infty[\times [0,L]^3\longrightarrow \Rt$, la pression $p:[0,+\infty[\times [0,L]^3\longrightarrow \mathbb{R}$ sont toujours les inconnues, $\vu_0\in L^2([0,L]^3)$ est la donnée initiale, périodique et à divergence nulle et finalement  nous considérons une force extérieure $\fe \in L^2([0,L]^3)\cap \dot{H}^{-1}([0,L]^3)$ périodique, à divergence nulle et \emph{stationnaire}, c'est à dire, nous avons  $\fe(t,x)=\fe(x)$.  \\
\\
Ainsi, pour une force $\fe \in L^{2}([0,L]^3)\cap \dot{H}^{-1}([0,L]^3)$ et une donnée initiale $\vu_0\in L^2([0,L]^3)$  telles que   
\begin{equation}\label{moyenne_nulle_donnes}
\int_{[0,L]^3}\fe(x)dx=\int_{[0,L]^3}\vu_0(x)dx=0,
\end{equation}   nous considérons  les solutions de Leray (globales en temps) des équations périodiques (\ref{N-S-per}): 
$$\vu\in L^{\infty}([0,+\infty[,L^2([0,L]^3)\cap L^{2}_{loc}([0,+\infty[,\dot{H}^1([0,L]^3)),$$
où la condition (\ref{moyenne_nulle_donnes})   nous permet de construire des solutions $\vu$  qui satisfont, pour tout temps $t>0$, 
\begin{equation}\label{moyenne_nulle_vitesse}
\int_{[0,L]^3}\vu(t,x)dx=0,
\end{equation}
et de cette propriété de moment nul  nous pouvons en tirer l'inégalité de Poincaré: pour $c_L>0$ une  constante qui dépend de la période $L$ on a l'estimation
\begin{equation}\label{Poincare}
\Vert \vu(t,\cdot)\Vert^{2}_{L^2}\leq \frac{L}{2\pi}\Vert \vec{\nabla}\otimes \vu(t,\cdot)\Vert^{2}_{L^2}, 
\end{equation}  qui nous sera fondamentale par la suite (voir les livres \cite{BoyerFabrice} et \cite{FMRTbook} pour une preuve de cette inégalité). \\
\\ 
Ces solutions $\vu$ vérifient en plus les propriétés suivantes:     
  \begin{enumerate}
\item[(i)] pour tout temps $t>0$, $\vu(t,x)$ est périodique sur le cube $[0,L]^3$, 
\item[(ii)]  pour tout temps $t>0,$ on a l'inégalité d'énergie
\begin{equation}\label{ineg_ener_Leray_periodique}
\frac{d}{dt}\Vert \vu(t,\cdot)\Vert^{2}_{L^2}+2\nu \Vert \vec{\nabla} \otimes \vu(t,\cdot)\Vert^{2}_{L^2} \leq 2 \int_{[0,L]^3}\fe(x)\cdot \vu(t,x)\,dx, %
\end{equation} 
\end{enumerate}  Pour l'existence de ces solutions   voir les livres \cite{ConstDoer} et \cite{FMRTbook} de C. Foias \emph{et al.}, le mémoire \cite{Leray}  de J. Leray ainsi que le livre \cite{Temam} de R. Temam \emph{et al.}\\
\\
Finalement, observons que ces solutions appartiennent à l'espace  $L^{\infty}_{t}L^{2}_{x}\cap (L^{2}_{t})_{loc}\dot{H}^{1}_{x}$ où nous remarquons que 
 le fait que ces solutions soient  \emph{globalement} bornées en temps provient de  l'inégalité d'énergie ci-dessus et  de l'inégalité de Poincaré (\ref{Poincare}), tandis que,   le fait que ces solutions soient \emph{localement} carré intégrables en temps ($\vu \in (L^{2}_{t})_{loc}\dot{H}^{1}_{x}$) provient de l'hypothèse d'une force extérieure stationnaire. 

\subsubsection{1) Des quantit\'es physiques associ\'ees au fluide} 
Maintenant,  en suivant les articles \cite{DoerFoias}, \cite{FMRT2} et le livre \cite{FMRTbook} de  F. Foias,  R. Temam \emph{et al.},  (voir aussi les articles \cite{Layton}, \cite{OttoRamos})  nous allons introduire certaines quantités physiques  qui sont nécessaires pour l'étude déterministe de la loi de dissipation de Kolmogorov (\ref{loi_dissipation_physique}) dans ce cadre périodique.\\
\begin{enumerate}
\item[$A)$] \textbf{La longueur caractéristique du fluide et l'échelle d'injection d'énergie.} \\
\\
En mécanique des fluides, la longueur caractéristique du fluide $L>0$  est la taille de l'écoulement considéré où nous allons étudier son comportement turbulent. Par exemple, si nous considérons un fluide dans un tuyau de diamètre $L>0$, alors la longueur caractéristique de ce fluide correspond naturellement au diamètre $L$. \\
\\
Lorsqu'on considère un modèle plus artificiel d'un cadre périodique sur le cube $[0,L]^3\subset\Rt$  nous avons  aussi une définition naturelle de la longueur caractéristique du fluide:  
\begin{Definition}[Longueur caractéristique]\label{Def_L_periodique} Dans le cadre périodique nous définissons la longueur caractéristique  comme la période $L>0$.
\end{Definition}
Dans ce cadre périodique  nous voulons étudier le comportement turbulent d'un fluide dans un cube $[0,L]^3$ et en revenant  au modèle de cascade d'énergie  nous savons que pour obtenir un comportement turbulent  il faut qu'une force extérieure $\fe$ agite ce fluide.  L'action de la force sur le fluide est  effectuée à une échelle de longueur $\ell_0>0$ et nous avons 

\begin{Definition}[Échelle d'injection d'énergie]\label{Def_ell_0_periodique} Nous définissons l'échelle d'injection d'énergie $\ell_0>0$ comme l'échelle de longueur à laquelle la force $\fe$ agit sur le fluide.
\end{Definition} Comme nous avons déjà mentionné dans l'introduction nous pouvons penser, par exemple, à un agitateur de taille $\ell_0$ qui communique constamment de l'énergie cinétique  au fluide. \\
\\
Une fois que nous avons défini l'échelle d'injection d'énergie nous voulons donc  savoir quelle est la relation entre cette échelle et  la longueur caractéristique $L$ donnée dans la Définition \ref{Def_L_periodique}. Ainsi, étant donné que le fluide se trouve dans un cube $[0,L]^3$ alors  il est naturel de supposer que la force $\fe$ agit sur le fluide à l'intérieur de ce cube, ce qui nous amène à supposer que 
\begin{equation}\label{relation_el_L_intro}
\ell_0\leq L.
\end{equation}
\item[$B)$] \textbf{La vitesse caractéristique $U$ et le taux de dissipation d'énergie $\varepsilon$.} \\
\\
L'objectif ici est de donner une définition rigoureuse des quantités $U$, qui correspond à la vitesse caractéristique du fluide, et $\varepsilon$ qui est le taux de dissipation d'énergie.  Ainsi,  pour $\vu \in L^{\infty}_{t}L^{2}_{x}\cap (L^{2}_{t})_{loc}\dot{H}^{1}_{x}$ une solution des équations (\ref{N-S-per}) et $L>0$ la longueur caractéristique du fluide donnée par la Définition \ref{Def_L_periodique}, nous définissons  les quantités  
\begin{equation}\label{U(T)}
\frac{1}{L^3}\Vert \vu(t,\cdot)\Vert^{2}_{L^2},
\end{equation} et
\begin{equation}\label{varepsilon(T)}
\frac{\nu}{L^3}\Vert \vec{\nabla}\otimes \vu(t,\cdot)\Vert^{2}_{L^2},
\end{equation}
qui correspondent à \emph{l'énergie cinétique}  du fluide et au \emph{taux de dissipation d'énergie} à l'instant $t>0$ respectivement. En effet, si nous supposons pour l'instant que la solution $\vu$ ci-dessus est assez régulière, de sorte que  l'inégalité d'énergie (\ref{ineg_ener_Leray_periodique}) devient une égalité (voir le livre \cite{BoyerFabrice}), et si nous considérons pour simplifier une force nulle, alors, nous pouvons écrire l'identité $$\frac{d}{dt}\left( \frac{1}{L^3}\Vert \vu(t,\cdot)\Vert^{2}_{L^2}\right)=-\frac{\nu}{L^3}\Vert \vec{\nabla}\otimes \vu(t,\cdot)\Vert^{2}_{L^2},$$ 
ce qui justifie le nom de \emph{``taux de dissipation d'énergie"} pour la quantité $\frac{\nu}{L^3}\Vert \vec{\nabla}\otimes \vu(t,\cdot)\Vert^{2}_{L^2}$. \\
\\
Nous observons maintenant que les quantités associées au champ de vitesse, données par les formules (\ref{U(T)}) et (\ref{varepsilon(T)}) ci-dessus  dépendent de chaque instant du temps $t>0$ et dans des expériences physiques  ces quantités fluctuent fortement (voir les articles \cite{Houel,Tennekes,Wilcox}). Néanmoins, d'après la théorie K41 les moyennes en temps de ces quantités sont censées présenter un comportement universel et  afin de capturer ce comportement  nous considérons les moyennes en temps de la façon suivante: tout d'abord pour un temps $T>0$, dans les expressions  (\ref{U-per}) et (\ref{varepsilon-per }) ci-dessus nous considérons ses moyennes sur l'intervalle de temps $[0,T]$: 
$$ \frac{1}{T}\int_{0}^{T}\Vert \vu(t,\cdot)\Vert^{2}_{L^2}\frac{dt}{L^3}, $$ et 
$$ \frac{1}{T}\int_{0}^{T}\Vert \vec{\nabla}\otimes\vu(t,\cdot)\Vert^{2}_{L^2}\frac{dt}{L^3}.$$
Ensuite,  rappelons que  nous considérons  $\fe$  une force extérieure \emph{stationnaire}, alors $\fe$ agit sur le fluide en introduisant l'énergie cinétique indépendamment du temps et donc, une fois que le fluide est en régime turbulent nous allons supposer  qu'il  restera dans cet état au cours du temps. Pour cette raison, nous allons étudier le comportement des quantités  ci-dessus dans le régime asymptotique lorsque $T\longrightarrow+\infty$ et nous considérons ainsi les moyennes en temps  suivantes: 

\begin{equation}\label{U-per}
U^2=\limsup_{T\longrightarrow +\infty}\frac{1}{T}\int_{0}^{T}\Vert \vu(t,\cdot)\Vert^{2}_{L^2}\frac{dt}{L^3},
\end{equation} et
\begin{equation}\label{varepsilon-per }
\varepsilon=\nu \limsup_{T\longrightarrow +\infty}\frac{1}{T}\int_{0}^{T}\Vert \vec{\nabla}\otimes\vu(t,\cdot)\Vert^{2}_{L^2}\frac{dt}{L^3},
\end{equation}
où $U^2$ sera l'énergie cinétique moyenne du fluide et $\varepsilon$ est le taux moyen de dissipation ( ou plus  simplement le \emph{ taux de dissipation}.\\
\\
Cette moyenne en temps, $\ds{\limsup_{T\longrightarrow +\infty}\frac{1}{T}\int_{0}^{T}(\cdot)dt}$, également appelée \emph{moyenne en temps long}, est très utilisée dans la littérature pour l'étude déterministe de la turbulence (voir les articles \cite{DoerFoias}, \cite{FMRT2} \cite{Layton}, \cite{OttoRamos} et le livre \cite{FMRTbook}). \\
\begin{Remarque}\label{Remarque:sol_Leray} Observons  que comme nous allons considérer cette moyenne en temps long alors nous avons besoin de considérer des solutions de Leray des équations (\ref{N-S-per}) qui sont globales en temps. 
\end{Remarque}Il est nécessaire maintenant de justifier rigoureusement  que les quantités (\ref{U-per}) et (\ref{varepsilon-per }) ont bien un sens mathématique, c'est à dire,  pour $\vu$ une solution de Leray quelconque des équations (\ref{N-S-per}), nous voulons nous assurer que l'on a bien $U^2<+\infty$ et $\varepsilon<+\infty$. Nous avons donc la proposition suivante.
\begin{Proposition}\label{Remarque_Poincare} 
Soit la période $L>0$, soient  $\fe \in L^2([0,L]^3)\cap \dot{H}^{-1}([0,L]^3)$ une force donnée et  $\vu_0\in L^2([0,L]^3)$ une donnée initiale. Soit  $\vu\in L^{\infty}([0,+\infty[,L^2([0,L]^3)\cap L^{2}_{loc}([0,+\infty[,\dot{H}^1([0,L]^3))$ une solution de Leray des équations  (\ref{N-S-per}) associée aux données $(\vu_0, \fe)$.  Alors:
\begin{enumerate} 
\item[1)] $\ds{U^2=\limsup_{T\longrightarrow +\infty}\frac{1}{T}\int_{0}^{T}\Vert \vu(t,\cdot)\Vert^{2}_{L^2}\frac{dt}{L^3} <+\infty}$ et 
\item[2)] $\ds{\varepsilon=\nu \limsup_{T\longrightarrow +\infty}\frac{1}{T}\int_{0}^{T}\Vert \vec{\nabla}\otimes\vu(t,\cdot)\Vert^{2}_{L^2}\frac{dt}{L^3}<+\infty.}$
\end{enumerate}
\end{Proposition}
\pv \begin{enumerate}
\item[1)]  Soit $T>0$.  Dans l'inégalité  l'inégalité d'énergie (\ref{ineg_ener_Leray_periodique}) nous prenons l'intégrale sur l'intervalle de temps $[0,T]$ et nous obtenons  $$\Vert \vu(t,\cdot)\Vert^{2}_{L^2}+ 2\nu \int_{0}^{T}\Vert \vec{\nabla} \otimes \vu(t,\cdot)\Vert^{2}_{L^2}dt\leq \Vert \vu_0\Vert^{2}_{L^2}+2\int_{0}^{T}\int_{[0,L]^3}\fe(x)\cdot \vu(t,x)dxdt,$$ d'où, comme $\Vert \vu(t,\cdot)\Vert^{2}_{L^2}$ est une quantité positive nous pouvons écrire 
\begin{equation}\label{ineq_ener_aux_per}
2\nu \int_{0}^{T}\Vert \vec{\nabla} \otimes \vu(t,\cdot)\Vert^{2}_{L^2}dt\leq \Vert \vu_0\Vert^{2}_{L^2}+2\int_{0}^{T}\int_{[0,L]^3}\fe(x)\cdot \vu(t,x)dxdt,
\end{equation} et nous allons maintenant étudier le terme $2\nu \int_{0}^{T}\Vert \vec{\nabla} \otimes \vu(t,\cdot)\Vert^{2}_{L^2}dt$.  Par  l'inégalité de Poincaré (\ref{Poincare})  nous écrivons  $\ds{\frac{1}{L}\Vert \vu(t,\cdot)\Vert^{2}_{L^2}\leq \Vert \vec{\nabla}\otimes \vu(t,\cdot)\Vert^{2}_{L^2}}$, d'où en intégrant en temps et en multipliant par $2\nu$ nous avons  la minoration$$\ds{\frac{2\nu}{L}\int_{0}^{T}\Vert \vu(t,\cdot)\Vert^{2}_{L^2}dt\leq 2\nu \int_{0}^{T}\Vert \vec{\nabla} \otimes \vu(t,\cdot)\Vert^{2}_{L^2}dt}.$$ Ainsi, en remplaçant cette minoration dans le terme à gauche de l'estimation   (\ref{ineq_ener_aux_per}) nous pouvons écrire 
$$  \frac{2\nu}{L}\int_{0}^{T}\Vert \vu(t,\cdot)\Vert^{2}_{L^2}dt\leq  \Vert \vu_0\Vert^{2}_{L^2}+2\int_{0}^{T}\int_{[0,L]^3}\fe(x)\cdot \vu(t,x)dxdt.$$
Ensuite, nous divisons chaque terme de cette estimation par $T$, puis nous appliquons  l'inégalité de Cauchy-Schwarz et les inégalités de Young  dans le deuxième terme à droite ci-dessus et nous obtenons alors 
\begin{eqnarray*}
\frac{2\nu}{L}\frac{1}{T}\int_{0}^{T} \Vert \vu(t,\cdot)\Vert^{2}_{L^2} dt &\leq & \frac{1}{T }\Vert \vu_0\Vert^{2}_{L^2}+\frac{2}{T}\int_{0}^{T}\int_{[0,L]^3}\fe(x)\cdot \vu(t,x)\,dxdt\\
&\leq & \frac{1}{T }\Vert \vu_0\Vert^{2}_{L^2} + \frac{L}{\nu} \Vert \fe \Vert^{2}_{L^2} + \frac{\nu}{L}\frac{1}{T}\int_{0}^{T} \Vert \vu(t,\cdot)\Vert^{2}_{L^2} dt. 
\end{eqnarray*} Nous prenons  maintenant la limite $\ds{\limsup_{T\longrightarrow +\infty}}$ et nous avons l'estimation
\begin{equation}\label{estim-per-u}
u^2=\limsup_{T\longrightarrow +\infty}\frac{1}{T}\int_{0}^{T}\Vert \vu(t,\cdot)\Vert^{2}_{L^2} dt \leq \frac{L^2}{\nu^2}  \Vert \fe \Vert^{2}_{L^2}.
\end{equation} 
Finalement, nous divisons chaque terme de cette estimation par $L^3$ pour récupérer
ainsi la vitesse caractéristique $U^2$ et nous avons $U^2 \leq \frac{1}{L}  \frac{\Vert \fe \Vert^{2}_{L^2}}{\nu^2} <+\infty$. 
\item[2)] Comme $\fe \in  L^2([0,L]^3)\cap \dot{H}^{-1}([0,L]^3)$ et $\vu(t,\cdot)\in L^2([0,L]^3)\cap \dot{H}^{1}([0,L]^3)$,  par l'inégalité de Cauchy-Schwarz (en variable d'espace) nous écrivons  
$$ 2\int_{0}^{T}\int_{[0,L]^3}\fe(x)\cdot \vu(t,x)dxdt \leq  \int_{0}^{T} \Vert \fe \Vert_{\dot{H}^{-1}} \Vert \vu(t,\cdot) \Vert_{\dot{H}^{1}}dt,$$ d'où, par les inégalités de Young et comme $\fe$ est une fonction stationnaire nous obtenons 
$$ 2\int_{0}^{T}\int_{[0,L]^3}\fe(x)\cdot \vu(t,x)dxdt \leq  T\frac{\Vert \fe \Vert^{2}_{\dot{H}^{-1}}}{\nu}+\nu \int_{0}^{T}\Vert \vec{\nabla} \otimes \vu(t,\cdot)\Vert^{2}_{L^2}dt.$$ Maintenant, en remplaçant cette inégalité dans le deuxième terme à droite de (\ref{ineq_ener_aux_per})  nous pouvons écrire 
$$ 2\nu \int_{0}^{T}\Vert \vec{\nabla} \otimes \vu(t,\cdot)\Vert^{2}_{L^2}dt\leq \Vert \vu_0\Vert^{2}_{L^2}+T\frac{\Vert \fe \Vert^{2}_{\dot{H}^{-1}}}{\nu}+\nu \int_{0}^{T}\Vert \vec{\nabla} \otimes \vu(t,\cdot)\Vert^{2}_{L^2}dt,$$ d'où nous obtenons 
$$ \nu \int_{0}^{T}\Vert \vec{\nabla} \otimes \vu(t,\cdot)\Vert^{2}_{L^2}dt\leq \Vert \vu_0\Vert^{2}_{L^2}+ T\frac{\Vert \fe \Vert^{2}_{\dot{H}^{-1}}}{\nu},$$ et alors en divisant par $TL^3$  et en prenant la limite $\underset{T\longrightarrow +\infty}{\limsup}$ nous pouvons écrire 
$$\varepsilon = \nu \limsup_{T\longrightarrow +\infty}\frac{1}{T}\int_{0}^{T}\Vert \vec{\nabla}\otimes\vu(t,\cdot)\Vert^{2}_{L^2}\frac{dt}{L^3}\leq  \frac{ \Vert \fe \Vert^{2}_{\dot{H}^{-1}}}{\nu  L^3}<+\infty,$$ ce qui termine la preuve de la proposition. \finpv
\end{enumerate} 
\begin{Remarque} Observons que la preuve du point $1)$ de la Proposition \ref{Remarque_Poincare} est basée sur l'inégalité d'énergie (\ref{ineg_ener_Leray_periodique}) vérifiée par le champ de vitesse $\vu$ et sur l'inégalité de Poincaré (\ref{Poincare}) qui est valable dans ce cadre périodique. D'autre part, la preuve du point $2)$ de la Proposition \ref{Remarque_Poincare} repose uniquement sur l'inégalité d'énergie (\ref{ineg_ener_Leray_periodique}). 
\end{Remarque}
Ainsi nous avons 
\begin{Definition}\label{Def_U_varepsilon_per} Dans le cadre de la Proposition \ref{Remarque_Poincare} ci-dessus, où nous avons montré que les quantités $U$ et $\varepsilon$ ci-dessous sont bien définies, nous définissons les quantités moyennes:
\begin{enumerate}
\item[1)] La vitesse caractéristique du fluide: $\ds{U=\left(\limsup_{T\longrightarrow +\infty}\frac{1}{T}\int_{0}^{T}\Vert \vu(t,\cdot)\Vert^{2}_{L^2}\frac{dt}{L^3}\right)^{\frac{1}{2}}}$.
\item[2)] Le taux de dissipation d'énergie: $\ds{\varepsilon=\nu \limsup_{T\longrightarrow +\infty}\frac{1}{T}\int_{0}^{T}\Vert \vec{\nabla}\otimes\vu(t,\cdot)\Vert^{2}_{L^2}\frac{dt}{L^3}}$. \\
\end{enumerate}
\end{Definition}
\item[$C)$] \textbf{Les nombres de Reynolds.} \\  
\\
Plusieurs nombres sans dimensions permettent de caractériser le comportement de l'écoulement des fluides (voir le livre \cite{BoyerFabrice}) et ici nous allons nous intéresser au nombre de Reynolds qui caractérise le régime  laminaire ou turbulent d'un fluide comme nous l'expliquerons  dans cette section.  \\
\\
Ainsi, pour introduire le nombre de Reynolds nous avons besoin de considérer les équations de Navier-Stokes (\ref{N-S-physique})  que nous avons introduit dans l'introduction de ce chapitre: 
$$ \rho \left( \partial_t\vu+(\vu\cdot \vec{\nabla}) \vu \right)=  \mu\Delta \vu  - \vec{\nabla}p_r+\vec{f}_{ext}, \quad \rho \, div(\vu)=0, $$ où $\rho>0$ est la constante de densité du fluide et $\mu=\rho \nu$ est la constante de viscosité dynamique du fluide. \\
\\ 
Dans le cadre de ces équations, le nombre de Reynolds a été mis en évidence  par Osborne Reynolds dans l'année 1883 \cite{McDonough} et mesure  le rapport  entre l'ordre de grandeur du terme de transport  $\rho\,(\vu\cdot \vec{\nabla}) \vu$ et l'ordre de grandeur du terme de viscosité $\mu \Delta \vu$.  Nous allons voir que le nombre de Reynolds apparaît naturellement dans les équations de Navier-Stokes  ci-dessus.  \\%
\\
En effet, pour $L>0$ la longueur caractéristique du fluide dans le sens de la Définition \ref{Def_L_periodique} et $U>0$ la vitesse caractéristique du fluide donnée par la Définition    
\ref{Def_U_varepsilon_per}, nous définissons   les variables et opérateurs adimensionnels suivants
\begin{equation}\label{changement_variables} 
\vec{u'}=\frac{\vu}{U}, \,\,\, {p'}_{r}=\frac{1}{\rho U^2}p_r,\,\,\,\vec{f'}_{ext}=\frac{L}{\rho U^2}\vec{f}_{ext},\,\,\,\partial_{t'}=\frac{L}{U}\partial_t,\,\,\, \vec{\nabla}_{x'}=L\vec{\nabla},
\end{equation} et alors les équations de Navier-Stokes ci-dessus  se réécrivent, après simplification par le facteur $\frac{\rho U^2}{L}$, comme
$$ \partial_{t'}\vec{u'} +(\vec{u'} \cdot \vec{\nabla}_{x'}) \vec{u'} =\frac{\mu}{\rho LU}\Delta_{x'}\vec{u'}-\vec{\nabla}_{x'}{p'}_{r}+\rho \vec{f'}_{ext}.$$ 
 Le nombre de Reynolds $Re$ est alors défini comme l'inverse de la constante qui se trouve devant le terme de viscosité, c'est à dire: $\ds{\frac{\rho U L}{\mu}}$, mais comme la constante  de viscosité dynamique est définie par $\mu=\rho \nu$ (où $\nu>0$ est la constante de viscosité cinétique)  nous avons l'identité $\frac{\rho}{\mu}=\frac{1}{\nu}$ et donc nous avons
 \begin{Definition}\label{Re} Le nombre de Reynolds $Re$ est défini par
 \begin{equation}\label{formule_Re}
 Re=\frac{U L }{\nu}.
 \end{equation}
 \end{Definition} 
Ce nombre sert à caractériser la nature du régime du mouvement du fluide: laminaire ou turbulent. En effet, le régime laminaire est caractérisé par les faibles valeurs du nombre $Re$, où les forces visqueuses sont dominantes: deux particules du fluide qui étaient voisines à un instant donné resteront voisines à l'instant suivant et les couches du fluide maintiennent leur cohésion au cours du temps.  Par contre, si le fluide est en régime turbulent alors on s'attend avoir des valeurs élevées du nombre $Re$   (voir les résultats expérimentaux  dans les articles \cite{Houel}, \cite{Wilcox} et le livre \cite{McDonough}) ce qui exprime le fait que le fluide est dominé par les force d'inertie, qui tendent à produire des tourbillons chaotiques et autres instabilités: les particules qui étaient voisines à un instant donné ne  seront plus voisines à l'instant suivant. Voir la Figure  \ref{Fig:Reynolds} ci-dessous pour une image d'une expérience physique des régimes laminaire et turbulent correspondant à différentes valeurs du nombre de Reynolds. \\
\\
\begin{figure}[!h]\label{Fig:Reynolds}
\begin{center}
\includegraphics[scale=0.30]{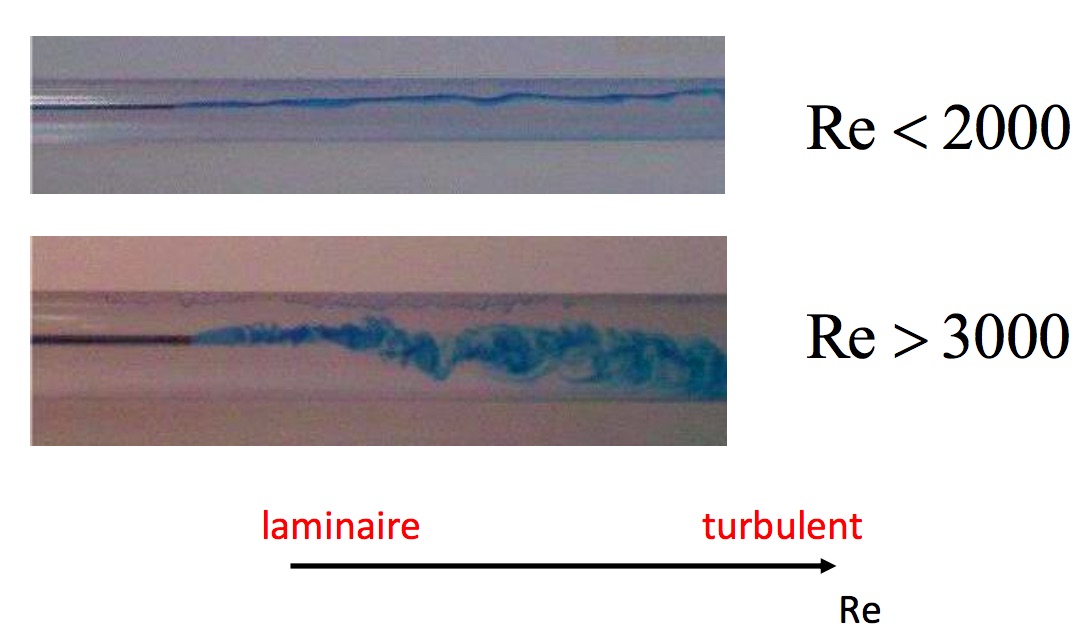}\caption{\footnotesize{ Expérience réalisée par N.H. Johannesen et C. Lowen avec de l'eau colorée introduite dans un tube, \cite{VanDyke}. Dans cette expérience, pour faire varier le nombre de Reynolds $Re=\frac{UL}{\nu}$, où la longueur caractéristique $L$ correspond au diamètre du tube et la constate de viscosité de l'eau est fixé par $\nu \approx 0,884\times 10^{-6} \,\frac{m^2}{s}$, la vitesse caractéristique du fluide $U$ est incrémentée par un agitateur extérieur.}}
\end{center}
\end{figure} \\
Une fois que nous avons définit le nombre de Reynolds $Re$ ci-dessus, il convient de faire les remarques suivantes.
\begin{Remarque}
\begin{enumerate}
\item[]
\item[$i)$] La définition du nombre de Reynolds n'est pas universelle dans le sens qu'elle dépend directement des  définitions de la vitesse caractéristique du fluide $U$ et de la longueur  $L$. En effet,  si dans la formule (\ref{changement_variables}) nous considérons  une autre définition de vitesse caractéristique  $U'$ et une autre définition de longueur  $L'$ alors par la formule (\ref{formule_Re}) nous obtenons le nombre de Reynolds $Re'$  donné par $Re'=\frac{U'L'}{\nu}$.
\item[$ii)$] Dans le cadre du point $i)$ ci-dessus, observons que pour définir le nombre de Reynolds $Re$ dans (\ref{formule_Re})  on a utilisé la longueur caractéristique $L$ tandis que pour définir  le nombre de Reynolds  (\ref{Re_lo_intro}) dans l'introduction de ce chapitre,  l'échelle d'injection d'énergie $\ell_0$ a été utilisée. Ce choix entre la longueur $L$  ou l'échelle $\ell_0$  n'empêche pas que le régime turbulent soit caractérisé par des grandes valeurs du nombre de Reynolds car,  étant donné que les quantités $L,\ell_0$ et $\nu$ sont fixes alors pour obtenir des grandes valeurs du nombre $Re$ c'est la vitesse caractéristique $U$ qui doit incrémenter.    
\item[$iii)$]  Comme  la vitesse caractéristique du fluide $U$ est définie à partir de la solution $\vu$ des équations  (\ref{N-S-per})  alors le nombre de Reynolds $Re$ dépend de cette solution $\vu$ et pour cette raison  ce nombre nous fourni une caractérisation \emph{a posteriori}  du régime du fluide, qu'il soit laminaire ou turbulent.  
\end{enumerate}
\end{Remarque}    
Revenons maintenant à la loi de dissipation d'énergie (\ref{loi_dissipation_physique}) proposée par la théorie K41. \\
\\
Nous savons que la loi de dissipation d'énergie (\ref{loi_dissipation_physique}) que nous souhaitons établir rigoureusement propose un encadrement du taux de dissipation d'énergie $\varepsilon$ lorsque le fluide est en régime turbulent. \\
\\
Nous venons de voir également que ce régime turbulent est caractérisé par des grands nombres de Reynolds; ce qui nous donne alors un cadre de travail assez naturel et nous supposerons souvent que $Re>>1$.  
\end{enumerate} 

\subsubsection{2) La loi de dissipation d'énergie dans le cadre périodique} 
Maintenant que nous avons les ingrédients nécessaires nous pouvons énoncer la loi de dissipation d'énergie. \\
\\  
Pour $\vu\in L^{\infty}([0,+\infty[,L^2([0,L]^3))\cap L^{2}_{loc}([0,+\infty[,\dot{H}^1([0,L]^3))$ une solution de Leray des équations de Navier-Stokes périodiques (\ref{N-S-per})  sur le cube $[0,L]^3$, on considère la vitesse caractéristique $U>0$ et le taux de dissipation $\varepsilon>0$ données par la Définition \ref{Def_U_varepsilon_per}. La longueur caractéristique du fluide $L$ est donnée par la Définition \ref{Def_L_periodique} et ainsi avec ces objets nous pouvons considérer  le nombre de Reynolds $Re=\frac{UL}{\nu}$. \\
\\
Donc, pour une échelle d'injection d'énergie $\ell_0>0$ donnée par la Définition \ref{Def_ell_0_periodique},  l'étude déterministe de la loi de dissipation de Kolmogorov (\ref{loi_dissipation_physique}) donnée page \pageref{loi_dissipation_physique} nous ramène à établir l'encadrement du taux de dissipation $\varepsilon$ suivant: 
\begin{equation}\label{Loi_Kolmogorov}
\text{si}\quad Re>>1\quad \text{alors}\quad  c_1\frac{U^3}{\ell_0} \leq \varepsilon \leq c_2 \frac{U^3}{\ell_0},
\end{equation} où $c_1>0$, $c_2>0$ sont des constantes indépendantes du nombre de Reynolds $Re$.\\ 
\\
Pour ce cadre périodique, lorsqu'on considère une échelle d'injection d'énergie $\ell_0$  égale à   la longueur caractéristique du fluide $L$ (c'est à dire $\ell_0 = L$)  il est alors possible de démontrer  la majoration du taux de dissipation  
\begin{equation}\label{major_varepsilon_periodique_1}
\varepsilon \leq c_2\frac{U^3}{L},
\end{equation} qui est obtenue dans différents contextes techniques où $c_2>0$ est une constante convenable qui reste bornée même dans le régime asymptotique lorsque le nombre de Reynolds est assez grand (voir l'article \cite{Childress} de S. Childress,  les articles  \cite{FMRT1}, \cite{FMRT2}, \cite{FMRT3} et le livre \cite{FMRTbook} de C. Foias, R. Temam \emph{et. al.}  ainsi que les articles \cite{Layton} de W. Layton  et \cite{Vigneron} de F. Vigneron).\\ 
\\
D'autre part, dans l'article \cite{DoerFoias} de C. Foias, on considère un cadre plus général où l'échelle d'injection d'énergie $\ell_0$ n'est pas forcément égale à la longueur caractéristique du fluide $L$. En effet, dans cet article la quantité $\ell_0>0$ est définie  par  $\ell_0=\frac{L}{n}$, avec $n$ un paramètre entier positif. De plus, en considérant une force extérieure particulière $\fe$, il est encore possible d'obtenir la majoration 
\begin{equation}\label{major_varepsilon_periodique_2} 
\varepsilon \leq c_2\frac{U^3}{\ell_0}, 
\end{equation}  avec $c_2>0$ une constante qui ne dépend pas des paramètres $\nu,n$ et $L$.\\
\\
Si nous comparons cette majoration ci-dessus du taux de dissipation d'énergie $\varepsilon$ avec la majoration donnée dans (\ref{major_varepsilon_periodique_1}) nous pouvons observer que cette majoration est plus proche de la loi de dissipation d'énergie (\ref{Loi_Kolmogorov}) car elle fait intervenir une échelle d'injection d'énergie $\ell_0$ qui n'est pas forcément du même ordre  que la longueur caractéristique du fluide $L$.  \\
\\
Maintenant, quant à la minoration du taux de dissipation $$ c_1\frac{U^3}{\ell_0}\leq \varepsilon,$$ cette estimation  reste encore un problème ouvert que l'on sait pas résoudre, même dans le cadre périodique. En effet, les majorations du taux de dissipation d'énergie (\ref{major_varepsilon_periodique_1}) et (\ref{major_varepsilon_periodique_2}) reposent essentiellement  sur l'inégalité d'énergie vérifiée par les solutions de Leray des équations de Navier-Stokes, mais, dans l'état actuel de nos connaissances, nous ne savons pas comment utiliser cette inégalité d'énergie pour étudier la minoration du taux de dissipation d'énergie ci-dessus.    
\subsection{Le  cadre non périodique}\label{sec:cadre_non_per}
Maintenant nous retournons au cadre d'un fluide non périodique posé dans l'espace $\Rt$ tout entier. Notre modèle déterministe pour étudier la loi de dissipation d'énergie  est alors donné par les équations de Navier-Stokes posées sur l'espace $\Rt$:
\begin{equation}\label{N-S-Rt}
\left\lbrace \begin{array}{lc}\vspace{3mm}
\partial_t\vu=\nu\Delta \vu-(\vu\cdot \vec{\nabla}) \vu -\vec{\nabla}p +\fe, \quad div(\vu)=0,\quad \nu>0, &\\
\vu(0,\cdot)=\vu_0\in L^2(\Rt),\quad div(\vu_0)=0,&\\
\end{array}\right.
\end{equation}
où, en suivant les idées exposées dans la section précédente, dans les équations ci-dessus  nous considérons une force extérieure, stationnaire et à divergence nulle telle que $\fe \in L^2(\Rt)\cap \dot{H}^{-1}(\Rt)$. Ainsi, si on considère les équations ci-dessus  nous avons  que les solutions de Leray $\vu$  vérifient $$\vu\in L^{\infty}_{loc}([0,+\infty[,L^2(\Rt))\cap L^{2}_{loc}([0,+\infty[,\dot{H}^{1}(\Rt)),$$ c'est à dire, ces solutions sont localement bornées temps et localement de carré intégrables en temps (voir le livre \cite{PGLR1} pour les détails). \\
\\
Une fois que nous avons rapidement  introduit les équations de Navier-Stokes (\ref{N-S-Rt}) sur tout l'espace $\Rt$,  nous voulons maintenant établir de façon rigoureuse  l'encadrement (\ref{Loi_Kolmogorov}) et il est important de souligner qu'il  y a très peu de références à ce sujet.  Une idée suggérée dans les notes de cours  \cite{Const} de P. Constantin est la suivante.\\
\\
Tout d'abord, on considère $\ell_0>0$ une échelle d'injection d'énergie donnée par la Définition \ref{Def_ell_0_periodique} page \pageref{Def_ell_0_periodique} et en s'inspirant de la Définition \ref{Def_U_varepsilon_per} page \pageref{Def_U_varepsilon_per}, pour $\vu$ une solution de Leray des équations (\ref{N-S-Rt}) on considère la vitesse caractéristique $U$ et le taux de dissipation $\varepsilon$  par les expressions: 
\begin{equation}\label{U,varepsilon_non_per} 
U=\left(\limsup_{T\longrightarrow +\infty}\frac{1}{T}\int_{0}^{T}\Vert \vu(t,\cdot)\Vert^{2}_{L^2}\frac{dt}{\ell^{3}_{0}}\right)^{\frac{1}{2}}\quad\text{et}\quad \varepsilon=\nu \limsup_{T\longrightarrow +\infty}\frac{1}{T}\int_{0}^{T}\Vert \vec{\nabla}\otimes\vu(t,\cdot)\Vert^{2}_{L^2}\frac{dt}{\ell^{3}_{0}}.
\end{equation}  
Ainsi, pour étudier l'encadrement du taux de dissipation $\varepsilon$ donné dans (\ref{Loi_Kolmogorov}) il nous manque un ingrédient et nous avons besoin de définir une longueur (que ce soit $\ell_0$ comme dans (\ref{Loi_Kolmogorov}) où $L$ comme dans  (\ref{major_varepsilon_periodique_1})) et ceci pose un certain nombre de problèmes. \\
\\
En effet,  rappelons tout d'abord que dans le cadre d'un fluide périodique sur le cube $[0,L]^3$ nous avons vu dans la Définition \ref{Def_L_periodique} que  la 
 \emph{longueur caractéristique}  du fluide est donnée de façon naturelle par la période $L>0$, mais, comme nous considérons maintenant un fluide posé sur l'espace $\Rt$ tout entier, alors on perd toute notion physique et mathématique de cette longueur caractéristique et la définition adéquate de cette longueur est une question qui n'est pas évidente à répondre.  \\
\\
Dans ce contexte, comme la force extérieure est  une donnée du modèle d'un fluide non périodique, toujours  dans \cite{Const} il est suggéré de considérer une longueur caractéristique $L_c$ en fonction de la force extérieure $\fe$ de la façon suivante: tout d'abord, pour l'échelle d'injection d'énergie $\ell_0>0$  nous supposons que la transformée de Fourier de la force extérieure $\fe$ est localisé aux fréquences $\vert \xi \vert \leq \frac{1}{\ell_0}$. Cette hypothèse sur la force extérieure représente le fait que, selon le modèle de cascade d'énergie, l'énergie cinétique est introduite dans le fluide  par la force extérieure $\fe$  aux échelles de longueur de l'ordre de $\ell_0$ et donc aux fréquences de l'ordre de $\frac{1}{\ell_0}$. \\
\\
Ensuite,  nous définissons la moyenne en norme $L^2$ de la force $\fe$ par la quantité   
\begin{equation}\label{F}
F=\frac{\Vert \fe \Vert_{L^2}}{\ell^{\frac{3}{2}}_{0}},  
\end{equation} et, toujours en suivant \cite{Const}, la longueur $L_c$ est définie comme
\begin{equation}\label{L_constantin}
L_c=\frac{F}{\Vert \vec{\nabla}\otimes \fe\Vert_{L^{\infty}}}.
\end{equation}
La signification physique de la longueur caractéristique $L_c$ n'est pas totalement claire  mais cette longueur apparaît  dans les calculs faits dans \cite{Const} et de cette façon, le but de P. Constantin dans ses notes  est de montrer que si l'on considère un fluide posé dans tout l'espace $\Rt$ alors, pour le taux de dissipation $\varepsilon$, pour la vitesse caractéristique $U$ (donnés dans (\ref{U,varepsilon_non_per})) et pour cette  longueur caractéristique $L_c$ ci-dessus, on peut obtenir l'estimation  du taux de dissipation suivante:
\begin{equation}\label{Ineq_taux_Constantin}
\varepsilon \leq c \frac{U}{L_c}\left( 1+\frac{1}{\sqrt{Re}} +\frac{3}{4 Re} \right),
\end{equation} où nous observons que, dans le régime asymptotique des nombres de Reynolds $\ds{Re}$ suffisamment grands, on peut alors obtenir l'estimation  du taux de dissipation $\ds{\varepsilon \leq c \frac{U^3}{L_c}}$  ce qui généralise au cadre $\Rt$ l'estimation (\ref{major_varepsilon_periodique_1}), page \pageref{major_varepsilon_periodique_1}, obtenue dans le cadre périodique.  \\ 
\\
Néanmoins,  l'estimation  (\ref{Ineq_taux_Constantin}) ci-dessus présente quelques lacunes d'un point de vue mathématique et met en évidence quelques contraintes techniques lorsqu'on considère un fluide posé sur tout l'espace $\Rt$. Dans la section qui suit nous expliquons plus en détail ces contraintes techniques.     
\subsection{Problèmes dans le cadre non périodique}\label{Sec:prob_non_periodique}
\subsubsection*{1) Une vitesse caractéristique potentiellement mal posée.}
La première contrainte technique relative à l'estimation (\ref{Ineq_taux_Constantin}) porte sur les solutions de Leray des équations (\ref{N-S-Rt}).  En effet, pour une solution de Leray $\vu$ quelconque nous devons considérer  les quantités moyennes suivantes
$$ U=\left(\limsup_{T\longrightarrow +\infty}\frac{1}{T}\int_{0}^{T}\Vert \vu(t,\cdot)\Vert^{2}_{L^2}\frac{dt}{\ell^{3}_{0}}\right)^{\frac{1}{2}}\quad\text{et}\quad \varepsilon=\nu \limsup_{T\longrightarrow +\infty}\frac{1}{T}\int_{0}^{T}\Vert \vec{\nabla}\otimes\vu(t,\cdot)\Vert^{2}_{L^2}\frac{dt}{\ell^{3}_{0}},
 $$où  $\ell_0>0$ est l'échelle d'injection d'énergie. Néanmoins (en \emph{toute} généralité) nous ne pouvons pas assurer que la vitesse caractéristique $U$ ci-dessus est  une quantité bien définie. En effet,  nous allons nous concentrer sur la moyenne en temps long qui apparaît dans la vitesse caractéristique $U$: 
 \begin{equation}\label{def-u}
 u=\left(\limsup_{T\longrightarrow +\infty}\frac{1}{T}\int_{0}^{T}\Vert \vu(t,\cdot)\Vert^{2}_{L^2} dt \right)^{\frac{1}{2}},
 \end{equation} où nous observons que  l'on a l'identité  
 \begin{equation}\label{U}
 U=\frac{u}{\ell^{\frac{3}{2}}_{0}},
 \end{equation}  et nous allons maintenant observer que cette moyenne en temps long $u$ est potentiellement mal posée. En effet, si nous  écrivons l'inégalité d'énergie vérifiée par la solution $\vu$:
\begin{equation}\label{ineq_energ_probleme_cadre_non_per}
\Vert \vu(t,\cdot)\Vert^{2}_{L^2}+2\nu \int_{0}^{t}\Vert \vec{\nabla} \otimes \vu(s,\cdot)\Vert^{2}_{L^2}ds\leq \Vert \vu_0\Vert^{2}_{L^2}+2\int_{0}^{t}\int_{\Rt}\fe(x)\cdot \vu(s,x)\,dxds, 
\end{equation} où, pour le dernier terme  à droite de cette inégalité,  comme 
  $\fe\in  \dot{H}^{-1}(\Rt)$ et $\vu(s,\cdot)\in \dot{H}^{1}(\Rt)$,  par l'inégalité de Cauchy-Schwarz nous pouvons écrire
$$  \int_{\Rt}\fe(x)\cdot \vu(s,x)\,dx \leq \Vert \fe \Vert_{\dot{H}^{-1}}\Vert \vu(s,\cdot)\Vert_{\dot{H}^1}\leq c\Vert \fe \Vert_{\dot{H}^{-1}}\Vert \vec{\nabla}\otimes \vu(s,\cdot)\Vert_{L^2},$$ et de plus, comme $\fe$ est une fonction stationnaire, par les inégalités de Young nous avons  
$$ 2\int_{0}^{t}\int_{\Rt}\fe(x)\cdot \vu(s,x)\,dxds \leq \frac{c t \Vert \fe \Vert^{2}_{\dot{H}^{-1}}}{2\nu}+2\nu \int_{0}^{t}\Vert \vec{\nabla}\otimes \vu(s,\cdot)\Vert^{2}_{L^2}ds.$$ De cette façon, en remplaçant cette estimation dans l'inégalité d'énergie (\ref{ineq_energ_probleme_cadre_non_per}) nous obtenons 
$$ \Vert \vu(t,\cdot)\Vert^{2}_{L^2}+2\nu \int_{0}^{t}\Vert \vec{\nabla} \otimes \vu(s,\cdot)\Vert^{2}_{L^2}ds\leq \Vert \vu_0\Vert^{2}_{L^2}+ \frac{c t \Vert \fe \Vert^{2}_{\dot{H}^{-1}}}{2\nu}+2\nu \int_{0}^{t}\Vert \vec{\nabla}\otimes \vu(s,\cdot)\Vert^{2}_{L^2}ds,$$ d'où nous pouvons en tirer le contrôle en temps suivant  pour la quantité $\Vert\vu(t,\cdot)\Vert^{2}_{L^2}$:
\begin{equation}\label{controle_temps_Leray}
\Vert \vu(t,\cdot)\Vert^{2}_{L^2}\leq \Vert \vu_0\Vert^{2}_{L^2}+\frac{c t}{2\nu}\Vert \fe \Vert^{2}_{\dot{H}^{-1}},
\end{equation}
 mais, lorsqu'on applique la moyenne en temps long à cette inégalité nous avons 
 \begin{eqnarray}\label{remarque_U} \nonumber
 u^2&=&\limsup_{T\longrightarrow +\infty}\frac{1}{T}\int_{0}^{T}\Vert \vu(t,\cdot)\Vert^{2}_{L^2} dt  \leq \limsup_{T \longrightarrow +\infty}\left[ \frac{1}{T}\int_{0}^{T}\Vert \vu_0 \Vert^{2}_{L^2} + \frac{1}{T}\int_{0}^{T} \frac{c t}{2\nu}\Vert \fe \Vert^{2}_{\dot{H}^{-1}} dt\right]\\
 & \leq &  \limsup_{T\longrightarrow +\infty} \left[ \Vert \vu_0\Vert^{2}_{L^2}+T \frac{\Vert \fe \Vert^{2}_{\dot{H}^{-1}}}{2\nu }\right] =+\infty,
 \end{eqnarray}  et nous ne connaissons pas un meilleur contrôle en temps du type (\ref{controle_temps_Leray}) de la quantité $\Vert \vu(t,\cdot)\Vert^{2}_{L^2}$ pour assurer que l'on a bien  $u^2<+\infty$.\\
\\
D'autre part, nous ne savons pas non plus si la quantité $u^2$ donnée dans (\ref{remarque_U}) diverge: dans l'état actuel de nos connaissances  nous ne savons pas  construire des solutions de Leray particuli\`eres telles qu'elles aient, par exemple,  le comportement suivant 
\begin{equation}\label{discussion_divergence_U_non_per}
t^{\alpha}\leq \Vert u(t)\Vert^{2}_{L^{2}},\quad \text{avec  }\quad 0<\alpha<1,   
\end{equation} où nous pouvons observer que si l'on prend la moyenne en temps long alors nous obtenons 
$$u^2= \limsup_{T\longrightarrow +\infty}\frac{1}{T}\int_{0}^{T} \Vert u(t)\Vert^{2}_{L^{2}} dt \geq  \limsup_{T\longrightarrow +\infty}\frac{1}{T}\int_{0}^{T}t^{\alpha} dt=+\infty,$$ et donc on pourrait  en conclure que la moyenne en temps long $u^2$ est effectivement mal posée lorsqu'on considère un fluide dans tout $\Rt$. \\ 
\\
Un des premiers résultats de cette thèse sera de donner un sens mathématique rigoureux à cette quantité $u$  et cet objectif sera atteint en utilisant un modèle particulier des équations de Navier-Stokes dans la Section  \ref{sec:N-S-amorties} ci-dessous.  
\subsubsection{2) Une longueur caractéristique mal posée}\label{longeur-mal-pose} 
Passons maintenant à la deuxième contrainte technique dans l'estimation (\ref{Ineq_taux_Constantin}) qui porte sur la définition de la longueur caractéristique $L_c$ donnée par l'expression (\ref{L_constantin}). Afin d'expliquer de façon plus claire cette contrainte, nous reprenons les grandes lignes des calculs faits dans \cite{Const}. En effet,  la preuve de l'estimation (\ref{Ineq_taux_Constantin})  repose sur l'inégalité 
\begin{equation}\label{ineq_F_Constantin}
F\leq \frac{U^2}{L_c}+\frac{1}{L_c}\sqrt{\nu}\sqrt{\varepsilon},
\end{equation} où la force moyenne $F$ est donnée dans l'expression (\ref{F}),  la vitesse caractéristique $U$  est donnée dans l'expression (\ref{U}), et le taux de dissipation $\varepsilon$  est  donné  par la formule (\ref{U,varepsilon_non_per}).  Pour prouver l'inégalité (\ref{ineq_F_Constantin}) ci-dessus, dans \cite{Const},  on considère les équations de Navier-Stokes 
$$ \partial_t\vu=\nu\Delta \vu-\P((\vu\cdot \vec{\nabla}) \vu) +\fe,$$ où $\vu\in L^{\infty}_{loc}([0,+\infty[,L^2(\Rt))\cap L^{2}_{loc}([0,+\infty[,\dot{H}^{1}(\Rt))$ est une solution de Leray et  $\fe\in \dot{H}^{-1}(\Rt) $ est une force extérieure régulière, stationnaire, à divergence nulle et telle que  $supp\left( \widehat{\fe}\right) \subset \lbrace \xi \in \Rt: \vert \xi \vert \leq \frac{1}{\ell_0}\rbrace$ pour $\ell_0>0$ une échelle d'injection d'énergie. On multiplie alors  ces équations  par $\fe$ puis on intègre en variable d'espace:   
\begin{eqnarray*}
\int_{\Rt} \partial_t\vu(t,x)\cdot \fe(x)dx&= &\int_{\Rt}\nu\Delta \vu(t,x)\cdot\fe(x)dx-\int_{\Rt}(\P(\vu\cdot\vec{\nabla}\vu(t,x)))\cdot \fe(x)dx+\Vert \fe \Vert^{2}_{L^2},  
 \end{eqnarray*} et on cherche à faire apparaître  les termes $F, U, \varepsilon$ et $L_c$ qui interviennent dans l'inégalité   (\ref{ineq_F_Constantin}). Dans l'identité ci-dessus, en utilisant des intégrations par parties, l'inégalité de H\"older et la moyenne $\ds{\limsup_{T\longrightarrow +\infty}\frac{1}{T}\int_{0}^{T}(\cdot)\frac{dt}{\ell^{3}_{0}}}$ on arrive à l'estimation suivante (voir les notes de cours \cite{Const} pour les détails)
\begin{eqnarray*}
\frac{\Vert \fe \Vert^{2}_{L^2}}{\ell^{3}_{0}} & \leq & \left(\limsup_{T\longrightarrow +\infty} \frac{1}{T}\int_{0}^{T}\Vert  \vu(t,\cdot)\Vert^{2}_{L^2}\frac{dt}{\ell^{3}_{0}} dt\right)\Vert \vec{\nabla} \otimes \fe \Vert_{L^{\infty}} \\
 & & +\sqrt{\nu}\left(\nu \limsup_{T\longrightarrow+\infty}\frac{1}{T}\int_{0}^{T}\Vert \vec{\nabla}\otimes\vu(t,\cdot)\Vert^{2}_{L^2}\frac{dt}{\ell^{3}_{0}} dt\right)^{\frac{1}{2}} \frac{\Vert \vec{\nabla}\otimes \fe \Vert_{L^2}}{\ell^{\frac{3}{2}}_{0}},
\end{eqnarray*}  d'où, par  les expressions des quantités $F,U$ et $\varepsilon$   on obtient   
$$ F^2 \leq \Vert \vec{\nabla}\otimes \fe \Vert_{L^{\infty}} + \sqrt{\nu}\sqrt{\varepsilon}\frac{\Vert \vec{\nabla}\otimes \fe \Vert_{L^2}}{ \ell^{\frac{3}{2}}_{0}}.$$
Dans  cette estimation nous écrivons  
$$ F\leq  U^2 \frac{\Vert \vec{\nabla}\otimes \fe \Vert_{L^{\infty}}}{F}+\sqrt{\nu}\sqrt{\varepsilon}\frac{\Vert \vec{\nabla}\otimes \fe \Vert_{L^2}}{ F\ell^{\frac{3}{2}}_{0}},$$ 
et, en utilisant le fait que la longueur $L_c$ est définie par $\ds{L_c=\frac{F}{\Vert \vec{\nabla}\otimes \fe \Vert_{L^{\infty}}}}$ alors, dans le premier terme de l'estimation ci-dessus nous pouvons écrire
\begin{equation}\label{ineq_F_Constantin_1}
F\leq   \underbrace{\frac{U^2}{L_c}}_{(a)} +\sqrt{\nu}\sqrt{\varepsilon}\underbrace{\frac{\Vert \vec{\nabla}\otimes \fe \Vert_{L^2}}{ F\ell^{\frac{3}{2}}_{0}}}_{(b)}.
\end{equation}
Si nous comparons cette estimation avec l'estimation recherchée (\ref{ineq_F_Constantin}) nous pouvons observer que  dans l'expression $(b)$ ci-dessus on veut faire apparaître le terme $\frac{1}{L_c}$ et pour cela, dans \cite{Const}, l'inégalité suivante est utilisée 
\begin{equation}\label{erreur_2_Const}
\frac{\Vert \vec{\nabla}\otimes \fe \Vert_{L^{2}}}{\ell^{\frac{3}{2}}_{0}} \leq \Vert \vec{\nabla}\otimes \fe \Vert_{L^{\infty}},
\end{equation}  indiquons rapidement que cette inégalité pose problème et nous y reviendrons dans les lignes qui suivent. Mais, en supposant pour l'instant que cette estimation est vraie, on peut écrire  l'estimation $$\frac{\Vert \vec{\nabla}\otimes \fe \Vert_{L^2}}{F\ell^{\frac{3}{2}}_{0}}\leq \frac{\Vert \vec{\nabla}\otimes \fe \Vert_{L^{\infty}}}{F} =\frac{1}{L_c},$$ et de cette façon, dans (\ref{ineq_F_Constantin_1}) on obtient l'inégalité cherchée (\ref{ineq_F_Constantin}). \\
\\
Revenons donc à l'estimation (\ref{erreur_2_Const}) et nous allons voir que  les calculs  ci-dessus présentent une lacune dans cette inégalité. En effet, la force $\fe$ est supposée régulière et sa transformée de Fourier $\widehat{\fe}$  vérifie  $\ds{supp(\widehat{\fe})\subset \left\lbrace \xi \in \Rt: \vert \xi \vert \leq \frac{1}{\ell_0} \right\rbrace }$; donc par les inégalités de Bernstein nous avons la majoration de la quantité $\Vert \vec{\nabla}\otimes \fe \Vert_{L^{\infty}}$: $$\Vert \vec{\nabla}\otimes \fe \Vert_{L^{\infty}}\leq c \frac{\Vert \vec{\nabla}\otimes \fe \Vert_{L^{2}}}{\ell^{\frac{3}{2}}_{0}},$$ mais, dans (\ref{erreur_2_Const}) nous avons besoin d'une minoration de la quantité $\Vert \vec{\nabla}\otimes \fe \Vert_{L^{\infty}}$ et nous ne savons pas déduire une telle estimation à partir des hypothèses de la fonction $\fe$.\\
\\
Ainsi, en revenant à la dernière expression de  l'inégalité (\ref{ineq_F_Constantin_1}) ci-dessus, nous pouvons observer que  la longueur caractéristique $L_c$ apparaît  naturellement dans le terme  $(a)$  mais non dans le terme $(b)$ de cette expression. \\
\\
\\
Donc, en résumé, nous pouvons  observer que l'estimation du taux de dissipation (\ref{Ineq_taux_Constantin}) proposée dans les notes de cours \cite{Const} possède deux contraintes techniques: d'une part pour $\vu$ une  solution de Leray des équations de Navier-Stokes posées sur tout $\Rt$, nous ne savons pas si la moyenne en temps long  $u$ donnée par l'expression (\ref{def-u})  est une quantité bien définie (et donc par l'identité (\ref{U}) on ne sait pas si $U<+\infty$).\\
\\
D'autre part, la longueur caractéristique $L_c$ ne convient pas pour obtenir l'inégalité (\ref{ineq_F_Constantin}) à partir de laquelle P. Constantin déduit  l'estimation du taux de dissipation (\ref{Ineq_taux_Constantin}). \\
\\
Pour régler le problème de la définition de la  moyenne en temps long $u$, dans la section qui suit nous proposons une modification des équations de Navier-Stokes (\ref{N-S-Rt}). Ensuite, dans la Section \ref{Sec:Theoreme principale} nous ferons une discussion sur la notion de  longueur caractéristique lorsqu'on travaille sur tout l'espace $\Rt$ où nous remarquerons le fait qu'une définition adéquate de telle longueur pour l'étude de l'encadrement du taux de dissipation (\ref{Loi_Kolmogorov}) semble actuellement hors de portée.

\section{Les équations de Navier-Stokes amorties}\label{sec:N-S-amorties}
Dans cette section, nous allons considérer un modèle particulier des équations de Navier-Stokes en introduisant un terme d'amortissement. \\
\\
L'étude de ce modèle nous permettra tout d'abord de bien définir la vitesse caractéristique $U$ et d'étudier ensuite l'encadrement du taux de dissipation d'énergie (\ref{Loi_Kolmogorov}), ce qui sera fait dans la section suivante.   
\subsection{Motivation du modèle}\label{Sec:motivation_alpha_modele} 
Dans la Section \ref{Sec:prob_non_periodique} ci-dessus, nous avons observé que lorsqu'on considère une solution de Leray $\vu$ des équations de Navier-Stokes posées dans  l'espace $\Rt$ tout entier, alors nous n'avons pas un  contrôle convenable en temps de la quantité $\Vert \vu (t,\cdot)\Vert^{^2}_{L^2}$ de sorte que l'on puisse assurer que la moyenne en temps  $$ u=\left( \limsup_{T\longrightarrow +\infty}\frac{1}{T}\int_{0}^{T}\Vert \vu(t,\cdot)\Vert^{2}_{L^2}dt \right)^{\frac{1}{2}},$$ soit bien une quantité finie. De cette façon, afin d'entraîner un contrôle sur la quantité $\Vert \vu(t,\cdot)\Vert^{2}_{L^2}$,  nous proposons ici de modifier les équations de Navier-Stokes en introduisant un terme additionnel $-\alpha \vu $  où $\alpha>0$ est un paramètre d'amortissement.  \\
\\
Ainsi le modèle sur lequel nous allons travailler dans tout ce chapitre est donné par le système d'équations suivant:
 \vspace{5mm}
\begin{center}
\fbox{
\begin{minipage}[l]{155mm}
\begin{equation}\label{NS-amortie-alpha} 
\left\lbrace \begin{array}{lc}\vspace{3mm}
\partial_t\vu = \nu\Delta \vu-(\vu\cdot\vec{\nabla}) \vu -\vec{\nabla}p +\fe-\alpha \vu, \quad div(\vu)=0,\quad div(\fe)=0, \quad \nu>0,\,\, \alpha>0,&\\
\vu(0,\cdot)=\vu_0, \quad div(\vu_0)=0.&
\end{array}\right. 
\end{equation}
\end{minipage}  
}
\end{center} 
\vspace{5mm}
Ce terme d'amortissement permet d'obtenir comme nous allons le voir un contrôle en temps de la quantité $\Vert \vu(t,\cdot)\Vert_{L^2}$  pour $\vu \in L^{\infty}_{t}L^{2}_{x}\cap (L^{2}_{t})_{loc}\dot{H}^{1}_{x}$ une solution faible des équations de Navier-Stokes amorties ci-dessus. En effet, dans le Théorème \ref{Theo:controle-temps} ci-dessous nous démontrons que les solutions faibles  de ces équations (solutions qui seront construites dans le Théorème \ref{Theo:existence-sol-ns-amortie}) vérifient le contrôle en temps:
\begin{equation}\label{controle_temps_alpha_modele_motivation} 
\Vert \vu(t,\cdot)\Vert^{2}_{L^2}\leq e^{-2\alpha t}\Vert \vu_0 \Vert^{2}_{L^2}+c\frac{\Vert \fe \Vert^{2}_{\dot{H}^{-1}}}{2\nu \alpha} \left( 1-e^{-2\alpha t}\right),
\end{equation} et dans cette estimation le terme $e^{-2\alpha t}$ (où $\alpha>0$) entraîne   $u^2<+\infty$ comme nous le verrons dans le  Corollaire \ref{Coro:theo_controle_temps_alpha_modele}. Ainsi, la vitesse la moyenne en temps long $u$ définies dans (\ref{def-u}) sera bien définie dans le cadre des équations de Navier-Stokes amorties (\ref{NS-amortie-alpha}).\\
\\
Avant de donner une preuve rigoureuse  de ces faits, ce qui sera fait dans la sous-section   \ref{Sec:Existence-prop-N-S-amortie} ci-dessous, il convient de faire maintenant une très courte discussion sur le terme d'amortissement $-\alpha \vu$ que nous venons d'ajouter aux équations de Navier-Stokes posées sur $\Rt$. \\
\begin{Remarque}\label{Poincare-artificiel}
Dans la  Remarque \ref{Remarque_Poincare} page \pageref{Remarque_Poincare} nous avons observé que si l'on considère un fluide périodique en variable d'espace alors l'inégalité de Poincaré nous permet d'assurer que la moyenne en temps  $u$ est bien une quantité finie. Pour assurer que la quantité $u$ est bien définie dans le cas de l'espace tout entier où l'on ne dispose pas d'un équivalent de l'inégalité de Poincaré, on verra que  le terme d'amortissement $-\alpha \vu$ peut dans un certain sens remplacer cette inégalité de Poincaré: c'est précisément l'utilité et l'intérêt d'introduire cette modification dans les équations de Navier-Stokes.\\
\end{Remarque}
Observons aussi que dans l'estimation  (\ref{controle_temps_alpha_modele_motivation})  nous avons  le terme $\ds{\frac{\Vert \fe \Vert^{2}_{\dot{H}^{-1}}}{2\nu \alpha}}$ et comme  $\alpha>0$ nous avons que ce terme n'est pas contrôlable lorsque $\alpha$ tend vers zéro et donc l'estimation (\ref{controle_temps_alpha_modele_motivation}) est seulement  valable dans le cadre des équation de Navier-Stokes amorties (\ref{NS-amortie-alpha}). \\
\\
Pour finir cette motivation, il est intéressant de souligner  que pour entraîner le contrôle $u^2<+\infty$ il est possible de considérer d'autres termes d'amortissement. Nous pouvons considérer  par exemple le terme d'amortissement $-\alpha P_{\kappa}(\vu)$, où le terme  $P_\kappa(\vu)$ est défini au niveau de Fourier par 
\begin{equation}\label{terme_amort2}
\widehat{P_{\kappa}(\vu)}(t,\xi)=\mathds{1}_{\vert \xi \vert <\kappa}(\xi)\widehat{\vu}(t,\xi),
\end{equation} où  $\kappa>0$ est  une fréquence de troncature et  $\widehat{\vu}$ dénote la transformée de Fourier de $\vu$ par rapport à la variable spatiale. Ainsi, dans un premier temps nous avons considéré le terme $-\alpha P_{\kappa}(\vu)$ pour obtenir un premier modèle d'équations de Navier-Stokes amorties:
\begin{equation}\label{N-S-amotrtie-terme2}
\partial_t\vu = \nu\Delta \vu-(\vu\cdot\vec{\nabla}) \vu-\vec{\nabla}p+\fe-\alpha P_{\kappa}(\vu), \quad div(\vu)=0,\quad \vu(0,\cdot)=\vu_0,
\end{equation} car ce terme supplémentaire entraîne un contrôle en temps de la quantité $\Vert \vu(t,\cdot)\Vert_{L^2}$ du même type que le contrôle (\ref{controle_temps_alpha_modele_motivation}) et donc les solutions $\vu$  des équations ci-dessus (qui sont construites de la même  façon que les solutions des équations (\ref{NS-amortie-alpha}))  vérifient aussi  $u^2<+\infty$. En effet, en suivant les mêmes lignes de la démonstration du Théorème \ref{Theo:controle-temps}  on peut  démontrer que toute solution  des équations (\ref{N-S-amotrtie-terme2}) vérifie le contrôle en temps suivant: pour  $\beta=\min(2\alpha,\nu \kappa^2)>0$,
\begin{equation}
\Vert \vu(t,\cdot)\Vert^{2}_{L^2}\leq e^{-\beta t}\Vert \vu_0 \Vert^{2}_{L^2}+c\frac{\Vert \fe \Vert^{2}_{\dot{H}^{-1}}}{2\nu \beta} \left( 1-e^{-\beta t}\right),
\end{equation}  et nous observons alors que ce contrôle en temps qui peut être obtenu  grâce au terme d'amortissement $-\alpha P_{\kappa}(\vu)$ est équivalent au contrôle   (\ref{controle_temps_alpha_modele_motivation}) qui sera obtenue par le biais du terme $-\alpha \vu$ et donc tous les résultats que nous allons obtenir dans le cadre des équations (\ref{NS-amortie-alpha}) peuvent être aussi obtenus dans le cadre des équations (\ref{N-S-amotrtie-terme2}).\\
\\
 Néanmoins,  nous allons préférer le terme d'amortissement $-\alpha \vu$ au lieu du terme $-\alpha P_{\kappa}(\vu)$ car ce premier terme est plus naturel  d'un point de vue physique. En effet, le terme d'amortissement $-\alpha \vu$, également appelé \emph{terme de friction}, a été considéré dans des modèles océaniques  \cite{Pedlosky} et ce terme  permet de modéliser la friction de l'eau avec le fond marin. D'autre part, le terme  $-\alpha P_{\kappa}(\vu)$ est un terme de troncature des hautes fréquences comme l'on peut observer dans la formule  (\ref{terme_amort2}) mais ce terme n'a, à notre connaissance, aucune signification physique. \\
\\
Notre étude se décomposera de la façon suivante: dans la Section \ref{Sec:Existence-prop-N-S-amortie} nous établissons les résultats de base par rapport aux équations (\ref{NS-amortie-alpha}) (existence, inégalités) et ensuite dans la Section  \ref{sec:parametre-amortissement} nous verrons  comment calibrer convenablement le paramètre $\alpha>0$  pour obtenir des résultats intéressants sur l'étude déterministe de  la loi de dissipation de Kolmogorov  (\ref{Loi_Kolmogorov}) dans le cadre des équations (\ref{NS-amortie-alpha}). \\
\subsection{Existence et propriétés}\label{Sec:Existence-prop-N-S-amortie}  
Dans cette section on commence par donner une preuve de l'existence de solutions faibles des équations de Navier-Stokes amorties (\ref{NS-amortie-alpha}). Dans le Théorème \ref{Theo:existence-sol-ns-amortie} ci-dessous nous construisons dans un premier temps des solutions $\vu$ localement bornées en temps pour ensuite vérifier dans le Théorème \ref{Theo:controle-temps}  que ces solutions satisfont le contrôle en temps (\ref{controle_temps_alpha_modele_motivation}) et qu'alors ces solutions sont globalement bornées en temps. 
\begin{Theoreme}\label{Theo:existence-sol-ns-amortie}
Soit $\vu_0\in L^2(\Rt)$ une donnée initiale à divergence nulle, soit $\fe \in \dot{H}^{-1}(\Rt)$ une force extérieure stationnaire et à divergence nulle.  Alors, pour tout $\alpha>0$  il existe des fonctions $\vu=\vu_\alpha\in L^{\infty}_{loc}([0,+\infty[,L^{2}(\Rt))\cap L^{2}_{loc}([0,+\infty[,\dot{H}^{1}(\Rt))$ et  $\ds{p=p_{\alpha}\in  L^{2}_{loc}([0,+\infty[,\dot{H}^{-\frac{1}{2}}(\Rt))}$ qui sont solution faible du système (\ref{NS-amortie-alpha}).
\end{Theoreme}  
La preuve de l'existence de ces solutions  suit  essentiellement les mêmes lignes que celle de l'existence des solutions de Leray des équations de Navier-Stokes classiques (voir le livre \cite{PGLR1}, Section $12.2$ pour tous les détails) et par conséquent nous détaillerons seulement les estimations réalisées sur le terme d'amortissement $-\alpha  \vu $.  \\
\\
\dm Nous appliquons  le projecteur de Leray aux équations (\ref{NS-amortie-alpha}) et comme $div(\vu)=0$ et $div(\fe)=0$ alors nous obtenons
\begin{equation}\label{N-S-amortie-alpha-Projec-Leray}   
\partial_t\vu = \nu\Delta \vu- \P \left( (\vu\cdot\vec{\nabla}) \vu \right) +\fe-\alpha \vu, \quad div(\vu)=0,  \quad \vu(0,\cdot)=\vu_0. 
\end{equation} 
Maintenant, soit $\theta \in \mathcal{C}^{\infty}_{0}(\Rt)$ une fonction positive telle que $\ds{\int_{\Rt}\theta(x)dx=1}$, pour $\delta>0$ on considère  la fonction $\theta_\delta$ donnée par $\theta_\delta(x)=\frac{1}{\delta^3}\theta \left( \frac{x}{\delta}\right)$ et on étudie alors  l'équation intégrale régularisée suivante 
\begin{eqnarray}\label{NS_regularised_point_fix}\nonumber
\vu(t,x)&=&h_{\nu t}\ast \vu_0(x) +\int_{0}^{t}h_{\nu(t-s)}\ast\fe(x)ds-\int_{0}^{t}h_{\nu(t-s)}\ast(\P(([\theta_\delta\ast\vu]\cdot \vec{\nabla})\vu)(s,x)ds-\\
& &-\alpha\int_{0}^{t}h_{\nu(t-s)}\ast \vu(s,x)ds. 
\end{eqnarray} Nous allons dans un premier temps appliquer un argument de point fixe dans l'espace\\ $L^{\infty}([0,T],L^2(\Rt))\cap L^2([0,T],\dot{H}^1(\Rt))$ muni de la norme $\Vert \cdot \Vert_T=\Vert \cdot \Vert_{L^{\infty}_{t}L^{2}_{x}}+\sqrt{\nu}\Vert \cdot \Vert_{L^{2}_{t}\dot{H}^{1}_{x}}$. Pour cela nous étudions la quantité 
\begin{eqnarray*} 
\Vert \vu \Vert_T &=& \left\Vert h_{\nu t}\ast \vu_0 +\int_{0}^{t}h_{\nu(t-s)}\ast\fe(\cdot)ds -\int_{0}^{t}h_{\nu(t-s)}\ast(\P(([\theta_\delta\ast\vu]\cdot \vec{\nabla})\vu)(s,\cdot)ds\right.\\
& &- \left.\alpha\int_{0}^{t}h_{\nu(t-s)}\ast \vu(s,\cdot)ds\right\Vert_T\\
&\leq &\underbrace{\left\Vert h_{\nu t}\ast \vu_0 +\int_{0}^{t}h_{\nu(t-s)}\ast\fe(\cdot)ds\right\Vert_T}_{(1)}  + \underbrace{\left\Vert \int_{0}^{t}h_{\nu(t-s)}\ast(\P(([\theta_\delta\ast\vu]\cdot \vec{\nabla})\vu)(s,\cdot)ds\right\Vert_T}_{(2)} \\
& & +\underbrace{\alpha  \left\Vert \int_{0}^{t}h_{\nu(t-s)}\ast \vu(s,\cdot)ds \right\Vert_T}_{(3)}.
\end{eqnarray*}
Les termes  $(1)$ et $(2)$ sont classiques à estimer. En effet, en utilisant le Théorème $12.2$ du livre \cite{PGLR1}  nous avons les estimations 
\begin{equation}\label{estimate_force_th_existence}
\Vert h_{\nu t}\ast \vu_0 \Vert_{T}\leq c \Vert \vu_0 \Vert_{L^2} \quad \text{et}\quad \left\Vert \int_{0}^{t}h_{\nu(t-s)}\ast\fe(\cdot)ds \right\Vert_T \leq c\left(\frac{1}{\nu}+T\sqrt{\nu} \right) \Vert \fe \Vert_{L^{2}_{t}H^{-1}_{x}},
\end{equation}  et comme la force $\fe \in\dot{H}^{-1}(\Rt)$ est indépendante de la variable de temps nous avons $$\ds{ \Vert \fe \Vert_{L^{2}_{t}H^{-1}_{x}} \leq \sqrt{T} \Vert \fe \Vert_{L^{\infty}_{t}H^{-1}_{x}}\leq \sqrt{T}\Vert \fe \Vert_{\dot{H}^{-1}}},$$  ce qui nous permet d'écrire la deuxième estimation de  (\ref{estimate_force_th_existence})  comme suit $$ \left\Vert \int_{0}^{t}h_{\nu(t-s)}\ast\fe(\cdot)ds \right\Vert_T \leq c\leq \left(\frac{1}{\nu}+T\sqrt{\nu} \right) \sqrt{T}\Vert \fe \Vert_{\dot{H}^{-1}}.$$ De cette façon, pour le terme $(1)$ ci-dessus nous pouvons écrire l'estimation 
\begin{equation}\label{estimate_terme_independante_existence_alpha_modele}
\left\Vert h_{\nu t}\ast \vu_0 +\int_{0}^{t}h_{\nu(t-s)}\ast\fe(\cdot)ds\right\Vert_T  \leq  c\Vert \vu_0 \Vert_{L^2}+c \left( \frac{1}{\sqrt{\nu}}+T\sqrt{\nu}\right)\sqrt{T}\Vert \fe \Vert_{\dot{H}^{-1}}.
\end{equation} Pour le terme $(2)$ nous avons directement l'estimation
\begin{equation}\label{estimate_forme_bilineaire_existence_alpha_modele}
\left\Vert \int_{0}^{t}h_{\nu(t-s)}\ast(\P(([\theta_\delta\ast\vu]\cdot \vec{\nabla})\vu)(s,\cdot)ds\right\Vert_T \leq c\frac{\sqrt{T}\delta^{-\frac{3}{2}}}{\sqrt{\nu}}\Vert \vu \Vert_T\,\Vert \vu \Vert_T,
\end{equation}
voir le livre \cite{PGLR1}  page $352$ pour les détails.  \\
\\
Donc, nous avons besoin d'étudier uniquement le terme $(3)$. En remplaçant $\fe$ par $\alpha \vu$ dans la deuxième estimation de (\ref{estimate_force_th_existence}) nous avons 
\begin{eqnarray}\label{estimate_terme_amortissement_existence_alpha_modele} \nonumber
\alpha  \left\Vert \int_{0}^{t}h_{\nu(t-s)}\ast \vu(s,\cdot)ds \right\Vert_T &\leq &\alpha c \left( \frac{1}{\sqrt{\nu}}+T\sqrt{\nu}\right)\Vert \vu \Vert_{L^{2}_{t}H^{-1}_{x}}\leq \alpha c \left( \frac{1}{\sqrt{\nu}}+T\sqrt{\nu}\right)\Vert \vu \Vert_{L^{2}_{t}H^{1}_{x}}\\ \nonumber
&\leq & \alpha c\left( \frac{1}{\sqrt{\nu}}+T\sqrt{\nu}\right)\left[\sqrt{T} \Vert \vu \Vert_{L^{\infty}_{t}L^{2}_{x}}+\frac{\sqrt{\nu}}{\sqrt{\nu}}\Vert \vu \Vert_{L^{2}_{t}\dot{H}^{1}_{x}}\right]\\
&\leq & \alpha c\left( \frac{1}{\sqrt{\nu}}+T\sqrt{\nu}\right)\max\left(\sqrt{T} ,\frac{1}{\sqrt{\nu}}\right) \Vert \vu \Vert_T.
\end{eqnarray} Une fois que nous avons les estimations (\ref{estimate_terme_independante_existence_alpha_modele}), (\ref{estimate_forme_bilineaire_existence_alpha_modele}) et (\ref{estimate_terme_amortissement_existence_alpha_modele}), pour un temps $T>0$ suffisamment petit et pour $\delta>0$, nous pouvons appliquer l'argument de point fixe de Picard pour construire une fonction $\vu_\delta$ telle que  $\vu_\delta \in L^{\infty}([0,T],L^2(\Rt))\cap L^2([0,T],\dot{H}^1(\Rt))$ et $\vu_\delta$ est  solution des équations approchées (\ref{NS_regularised_point_fix}).\\
\\
Une fois que l'on a construit une solution $\vu_\delta$ (locale en temps) nous prouvons  l'existence globale de cette solution: pour $\vu_\delta \in L^{\infty}([0,T],L^2(\Rt))\cap L^2([0,T],\dot{H}^1(\Rt))$ une solution des équations régularisées 
\begin{equation}\label{NS-regul-alpha}
\partial_t\vu_\delta= \nu\Delta \vu_\delta-\P(([\theta_\delta\ast\vu_\delta]\cdot\vec{\nabla}) \vu_\delta)+\fe-\alpha \vu_\delta,
\end{equation}
nous pouvons écrire 
\begin{eqnarray}\label{ineg-energ-delta}\nonumber 
\frac{d}{dt} \Vert \vu_\delta(t,\cdot)\Vert^{2}_{L^2} &= & 2\langle \partial_t \vu_\delta(t,\cdot),\vu_\delta(t,\cdot)\rangle_{\dot{H}^{-1}\times \dot{H}^1} = -2\nu \Vert \vu_\delta(t,\cdot)\Vert^{2}_{\dot{H}^1}+2\langle \fe,\vu_\delta(t,\cdot)\rangle_{\dot{H}^{-1}\times \dot{H}^1}\\
& &-2\alpha \Vert \vu_\delta(t,\cdot)\Vert^{2}_{L^2},
\end{eqnarray} d'où, en appliquant l'inégalité de Cauchy-Schwarz et puis les inégalités de Young  sur le terme $2\langle \fe,\vu_\delta(t,\cdot)\rangle_{\dot{H}^{-1}\times \dot{H}^1}$,  nous obtenons   
\begin{eqnarray}\label{controle-alpha} \nonumber 
\frac{d}{dt} \Vert \vu_\delta(t,\cdot)\Vert^{2}_{L^2} & \leq  & -2\nu \Vert \vu_\delta (t,\cdot)\Vert^{2}_{\dot{H}^1} +\nu \Vert \vu_\delta(t,\cdot)\Vert^{2}_{\dot{H}^1}+\frac{1}{\nu}\Vert f \Vert^{2}_{\dot{H}^{-1}} -2\alpha \Vert \vu_\delta(t,\cdot)\Vert^{2}_{L^2}\\
& \leq &  -\nu \Vert \vu_\delta (t,\cdot)\Vert^{2}_{\dot{H}^1} +\frac{1}{\nu}\Vert f \Vert^{2}_{\dot{H}^{-1}} -2\alpha \Vert \vu_\delta(t,\cdot)\Vert^{2}_{L^2},
\end{eqnarray} 
et alors, comme $-2\alpha \Vert \vu_\delta(t,\cdot)\Vert^{2}_{L^2}$ est une quantité négative,  nous avons donc  l'estimation 
$$ \frac{d}{dt} \Vert \vu_\delta(t,\cdot)\Vert^{2}_{L^2} \leq -\nu \Vert \vu_\delta (t,\cdot)\Vert^{2}_{\dot{H}^1} +\frac{1}{\nu}\Vert f \Vert^{2}_{\dot{H}^{-1}}.$$
Finalement, nous intégrons en temps sur l'intervalle $t\in [0,T]$ pour obtenir l'estimation
\begin{equation}\label{estim-norme-delta}
\Vert \vu_\delta (t,\cdot)\Vert^{2}_{L^2}+\nu\int_{0}^{t} \Vert \vu_\delta(s,\cdot)\Vert^{2}_{\dot{H}^1}ds \leq e^{\nu t}\left(\Vert\vu_0\Vert^{2}_{L^2}+\frac{T}{\nu}\Vert f \Vert^{2}_{\dot{H}^{-1}}\right),
\end{equation} et alors la solution locale en temps $\vu_\delta$ peut être étendue à l'intervalle   $[0,+\infty[$.\\
\\ 
Nous passons maintenant à la convergence vers une solution faible des équations (\ref{N-S-amortie-alpha-Projec-Leray}).  En effet, par le lemme de Rellich-Lions  (voir le Théorème $12.1$ du livre \cite{PGLR1}) il existe une suite de nombres positifs  $(\delta_n)_{n\in\mathbb{N}}$ et une fonction $\vu \in L^{2}_{loc}([0,+\infty[\times \Rt)$ telles que $(\vu_{\delta_n})_{n\in \mathbb{N}}$  converge fortement vers $\vu$ dans $L^{2}_{loc}([0,+\infty[\times \Rt)$. De plus, pour tout $T>0$, cette suite converge  vers $\vu$ dans la topologie faible étoile des espaces $L^{\infty}([0,T],L^2(\Rt))$ et $L^{2}([0,T],\dot{H}^1(\Rt))$ .\\ 
\\
Ainsi,  par les convergences ci-dessus nous obtenons que le terme $\P(([\theta_{\delta_n}\ast\vu_{\delta_n}]\cdot\vec{\nabla}) \vu_{\delta_n})$ converge  vers $\P((\vu\cdot\vec{\nabla}) \vu)$ dans la topologie faible étoile de l'espace $(L^{2}_{t})_{loc}(H^{-\frac{3}{2}}_{x})$ et donc la limite $\vu$ est une solution des équations de Navier-Stokes amorties  (\ref{N-S-amortie-alpha-Projec-Leray}). \\
\\
Une fois que l'on a construit une solution $\vu$ nous récupérons maintenant la pression $p$ reliée a cette solution. En effet, comme $\vu$ vérifie les équations  (\ref{N-S-amortie-alpha-Projec-Leray}) et comme $div(\vu)=0$ et $div(\fe)=0$  nous pouvons  écrire 
$$ \P \left[ \partial_t\vu -\nu\Delta \vu + \left( (\vu\cdot\vec{\nabla}) \vu \right) -\fe+\alpha \vu  \right]=0,$$ et par les propriétés du projecteur de Leray (voir le livre \cite{PGLR1}) il existe  $p\in \mathcal{D}^{'}(\Rt)$ telle que 
$$ \partial_t\vu -\nu\Delta \vu +  (\vu\cdot\vec{\nabla}) \vu  -\fe+\alpha \vu= \vec{\nabla}p. $$
De plus, en appliquant l'opérateur de divergence à chaque côté de cette identité nous obtenons la relation 
$$  p=\frac{1}{\Delta} div ((\vu\cdot\vec{\nabla}) \vu),$$ et comme $\vu \in (L^{\infty}_{t})_{loc}L^{2}_{x}\cap (L^{2}_{t})_{loc}\dot{H}^{1}_{x}$ alors nous avons $(\vu\cdot\vec{\nabla}) \vu \in (L^{2}_{t})_{loc}\dot{H}^{-\frac{3}{2}}_{x}$ et donc par la relation ci-dessus nous obtenons $p\in (L^{2}_{t})_{loc}\dot{H}^{-\frac{1}{2}}_{x}$.  \finpv \\
Une fois que nous avons montré l'existence globale des solutions faibles des  équations de Navier-Stokes amorties (\ref{NS-amortie-alpha}), dans la proposition qui suit nous montrons que ces solutions vérifient une inégalité d'énergie qui sera  exploitée tout au long de la Section \ref{Sec:Theoreme principale} ci-après. \\
\begin{Proposition}\label{Proposition_ineq_energie_alpha_modele} Dans le cadre du Théorème \ref{Theo:existence-sol-ns-amortie}, les solutions faibles des équations de Navier-Stokes amorties (\ref{NS-amortie-alpha}), $\vu\in L^{\infty}_{loc}(]0,+\infty[,L^2(\Rt))\cap L^{2}_{loc}(]0,+\infty[,\dot{H}^1(\Rt))$, vérifient l'inégalité d'énergie suivante: pour tout $T\geq 0$,
\begin{eqnarray}\label{energ_ineq_u}\nonumber
\Vert \vu(T,\cdot)\Vert^{2}_{L^2}+2\nu\int_{0}^{T} \Vert \nabla \otimes \vu(t,\cdot)\Vert^{2}_{L^2}dt &\leq &\Vert \vu_0 \Vert^{2}_{L^2}+2\int_{0}^{T} \langle \fe, \vu(t,\cdot) \rangle_{\dot{H}^{-1}\times \dot{H}^{1}}dt\\
& & -2\alpha \int_{0}^{T}\Vert \vu(t,\cdot) \Vert^{2}_{L^2}dt.
\end{eqnarray}
\end{Proposition}
 Pour montrer l'inégalité  ci-dessus nous suivrons les mêmes lignes de la preuve de l'inégalité d'énergie des solutions de Leray des équations de Navier-Stokes classiques faite dans le  Théorème $12.2$ du livre \cite{PGLR1}. \\
 \\
\pv Notre point de départ est l'identité (\ref{ineg-energ-delta}), d'où, pour $T>0$ en intégrant sur l'intervalle $[0,T]$ nous obtenons
\begin{eqnarray*}\label{energ_eq_delta_2} \nonumber
\Vert \vu_{\delta_n}(T,\cdot)\Vert^{2}_{L^2}+2\nu \int_{0}^{T} \Vert \vu_{\delta_n}(t,\cdot)\Vert^{2}_{\dot{H}^1}dt & =& \Vert \vu_0\Vert^{2}_{L^2}-2\alpha \int_{0}^{T}\Vert \vu_{\delta_n}(t,\cdot)\Vert^{2}_{L^2}dt \\
& & +2 \int_{0}^{T}\langle \fe, \vu_{\delta_n}(t,\cdot)\rangle_{\dot{H}^{-1}\times \dot{H}^1}dt.
\end{eqnarray*}
Maintenant, on régularise cette égalité en variable du temps et pour cela nous considérons la fonction test positive $w\in \mathcal{C}^{\infty}_{0}([-\eta,\eta])$ telle que $\ds{\int_{\mathbb{R}}w(t)dt=1}$. De cette façon, dans l'identité  ci-dessus  nous obtenons   
\begin{eqnarray}\label{eq_aux_step_4}\nonumber
& &\Vert w\ast \vu_{\delta_n}(T,\cdot)\Vert^{2}_{L^2}+ 2 w\ast \left(\nu \int_{0}^{T} \Vert \vu_{\delta_n}(t,\cdot)\Vert^{2}_{\dot{H}^1}dt+\alpha \int_{0}^{T}\Vert \vu_{\delta_n}(t,\cdot)\Vert^{2}_{L^2}dt\right)\\ 
& \leq & \Vert \vu_0\Vert^{2}_{L^2} +2w\ast \left( \int_{0}^{T}\langle \fe, \vu_{\delta_n}(t,\cdot)\rangle_{\dot{H}^{-1}\times \dot{H}^1}dt\right).
\end{eqnarray} 
\`A ce stade, nous avons besoin de vérifier que la suite $(\vu_{\delta_n})_{n\in\mathbb{N }}$ converge faiblement vers $\vu$ dans $L^{2}([0,T],L^2(\Rt))$. En effet, nous avons  $\Vert \vu_{\delta_n}\Vert_{L^{2}([0,T],L^2(\Rt))}\leq \sqrt{T}\Vert \vu_{\delta_n}\Vert_{L^{\infty}([0,T],L^2(\Rt))}$ et  alors par l'inégalité  (\ref{estim-norme-delta}) nous obtenons cette convergence. De plus, nous savons en plus que la suite $(\vu_{\delta_n})_{n\in \mathbb{R}}$ converge vers $\vu$ dans la topologie faible étoile de $L^{\infty}([0,T],L^2(\Rt))$ et $L^{2}([0,T],\dot{H}^1(\Rt))$  (pour tout $T>0$) et  alors  par l'inégalité (\ref{eq_aux_step_4}) nous pouvons écrire 
\begin{eqnarray*}
& &\Vert w\ast \vu(T,\cdot)\Vert^{2}_{L^2} + 2\nu w\ast \left( \int_{0}^{T} \Vert \vu(t,\cdot)\Vert^{2}_{\dot{H}^1}dt\right) +  2\alpha w\ast \left( \int_{0}^{T}\Vert \vu(t,\cdot)\Vert^{2}_{L^2}dt\right) \\
& \leq & \liminf_{n\longrightarrow +\infty} \Vert w\ast \vu_{\delta_n}(T,\cdot)\Vert^{2}_{L^2} +\liminf_{n\longrightarrow +\infty}2\nu w\ast \left( \int_{0}^{T} \Vert \vu_{\delta_n}(t,\cdot)\Vert^{2}_{\dot{H}^1}dt\right)\\
& &+\liminf_{n\longrightarrow +\infty}2\alpha w\ast \left( \int_{0}^{T}\Vert \vu_{\delta_n}(t,\cdot)\Vert^{2}_{L^2}dt\right)\\
& \leq & \Vert \vu_0\Vert^{2}_{L^2} +2w\ast \left( \int_{0}^{T}\langle \fe, \vu(t,\cdot)\rangle_{\dot{H}^{-1}\times \dot{H}^1}dt\right).
\end{eqnarray*}
Donc, pour $T>0$ un point de Lebesgue de l'application $t\mapsto \Vert \vu(t,\cdot)\Vert^{2}_{L^2}$ nous obtenons l'inégalité d'énergie (\ref{energ_ineq_u}) et en plus, cette égalité est étendue à tout temps $T>0$ par la continuité faible de l'application $t\mapsto \Vert \vu(t,\cdot)\Vert^{2}_{L^2}$ (voir  le Théorème $12.2$ du livre \cite{PGLR1} pour les détails).  \finpv  
\\
Nous avons maintenant à notre disposition tout les outils pour montrer que les solutions faibles des équations de Navier-Stokes amorties (\ref{NS-amortie-alpha}) vérifient le  contrôle en temps (\ref{controle_temps_alpha_modele_motivation}) où nous allons pouvoir apprécier le rôle du terme d'amortissement $-\alpha \vu$ dans ces estimations.  \\ 
\begin{Theoreme}\label{Theo:controle-temps} 
Soit $\vu_0\in L^2(\Rt)$ une donnée initiale à divergence nulle, soit $\fe \in \dot{H}^{-1}(\Rt)$ une force extérieure stationnaire et à divergence nulle et soit $\alpha>0$ le paramètre d'amortissement. Alors toute solution faible $\vu\in (L^{\infty}_{t})_{loc}L^{2}_{x}\cap (L^{2}_{t})_{loc}\dot{H}^{1}_{x}$ du système d'équations (\ref{NS-amortie-alpha}) obtenue par le biais du Théorème \ref{Theo:existence-sol-ns-amortie} vérifie $\vu \in L^{\infty}_{t}L^{2}_{x}\cap (L^{2}_{t})_{loc}\dot{H}^{1}_{x}$ et l'on a l'inégalité suivante: pour tout $t\geq 0$,
\begin{equation}\label{control_in_time_Th_2}
\Vert \vu(t,\cdot)\Vert^{2}_{L^2}\leq e^{-2\alpha t}\Vert \vu_0 \Vert^{2}_{L^2}+c\frac{\Vert \fe \Vert^{2}_{\dot{H}^{-1}}}{2\alpha\nu} \left( 1-e^{-2\alpha t}\right). 
\end{equation} 
\end{Theoreme} 
\dm   Nous allons montrer ce contrôle en temps pour la fonction $\vu_{\delta_n}$ (solution de l'équation régularisée (\ref{NS-regul-alpha})) et ensuite,  par  convergence faible étoile de  $(\vu_{\delta_n})_{n\in\mathbb{N}}$ vers $\vu$ dans l'espace  $L^{\infty}_{t}(L^{2}_{x})$  nous récupérons ce contrôle en temps pour la solution  $\vu$.\\
\\
On commence donc par l'inégalité (\ref{controle-alpha}):
$$ \frac{d}{dt} \Vert \vu_\delta(t,\cdot)\Vert^{2}_{L^2} \leq  -\nu \Vert \vu_\delta (t,\cdot)\Vert^{2}_{\dot{H}^1} +\frac{1}{\nu}\Vert f \Vert^{2}_{\dot{H}^{-1}} -2\alpha \Vert \vu_\delta(t,\cdot)\Vert^{2}_{L^2},$$ d'où, comme $-\nu \Vert \vu_\delta (t,\cdot)\Vert^{2}_{\dot{H}^1}$ est une quantité négative nous pouvons écrire 
$$ \frac{d}{dt} \Vert \vu_\delta(t,\cdot)\Vert^{2}_{L^2} \leq \frac{1}{\nu}\Vert f \Vert^{2}_{\dot{H}^{-1}} -2\alpha \Vert \vu_\delta(t,\cdot)\Vert^{2}_{L^2}, $$
 et par une application de l'inégalité de Gr\"onwall nous obtenons 
\begin{equation}\label{control_time_delta-alpha} 
\Vert \vu_\delta(t,\cdot)\Vert^{2}_{L^2}\leq e^{-2\alpha t}\Vert \vu_0 \Vert^{2}_{L^2}+c\frac{\Vert \fe \Vert^{2}_{\dot{H}^{-1}}}{2\alpha\nu} \left( 1-e^{-2\alpha t}\right),
\end{equation} pour tout temps $t\in [0,+\infty[$.\\  
\\
Maintenant, nous récupérons ce contrôle en temps pour la solution $\vu$ et pour cela on suit encore  l'argument utilisé dans la démonstration de la Proposition \ref{Proposition_ineq_energie_alpha_modele}.   En effet,  on régularise la quantité $\Vert \vu_{\delta_n}(t,\cdot)\Vert^{2}_{L^2}$  en variable de temps en prenant le produit de convolution avec une fonction positive $w\in\mathcal{C}^{\infty}_{0}([-\eta,\eta])$  telle que $\ds{\int_{\mathbb{R}}w(t)dt=1}$.  De cette façon, dans l'inégalité (\ref{control_time_delta-alpha}) ci-dessus nous avons   $$\Vert w\ast \vu_{\delta_n}(t,\cdot)\Vert_{L^2}\leq  w\ast \Vert \vu_{\delta_n}(t,\cdot)\Vert^{2}_{L^2}\leq w\ast \left( e^{-2 \alpha t}\Vert \vu_0 \Vert^{2}_{L^2}+c\frac{\Vert \fe \Vert^{2}_{\dot{H}^{-1}}}{ 2 \alpha \nu } \left( 1-e^{-2 \alpha  t}\right)\right).$$ 
Ensuite, comme $(\vu_{\delta_n})_{n\in\mathbb{N}}$ converge vers $\vu$ dans la topologie faible étoile de l'espace $L^{\infty}_{t}(L^{2}_{x})$ alors $w\ast \vu_{\delta_n}(t,\cdot)$ converge faiblement vers $w\ast\vu(t,\cdot)$ dans $L^2(\Rt)$ et de cette façon nous obtenons   
\begin{equation*}  
\Vert w\ast \vu(t,\cdot)\Vert^{2}_{L^2} \leq  \liminf_{n\longrightarrow +\infty} \Vert w\ast \vu_{\delta_n}(t,\cdot)\Vert^{2}_{L^2}\leq  w\ast \left( e^{-2 \alpha t}\Vert \vu_0 \Vert^{2}_{L^2}+c\frac{\Vert \fe \Vert^{2}_{\dot{H}^{-1}}}{ 2 \alpha \nu} \left( 1-e^{-2 \alpha  t}\right)\right).
\end{equation*} Ainsi, par la continuité faible de l'application $t\mapsto \Vert \vu(t,\cdot)\Vert^{2}_{L^2}$ (voir le Théorème $12.2$ du livre \cite{PGLR1},) nous avons le contrôle en temps cherché. \finpv
\\    
\begin{Corollaire}\label{Coro:theo_controle_temps_alpha_modele} Toute solution faible des équations (\ref{NS-amortie-alpha}), $\vu\in L^{\infty}_{t}L^{2}_{x}\cap (L^{2}_{t})_{loc}\dot{H}^{1}_{x}$, vérifie 
$$ u^2=\limsup_{T\longrightarrow +\infty}\frac{1}{T}\int_{0}^{T}\Vert \vu(t,\cdot)\Vert^{2}_{L^2}dt \leq \frac{\Vert \fe \Vert^{2}_{\dot{H}^{-1}}}{ \nu \alpha} <+\infty.$$ 
\end{Corollaire}
\pv Cette estimation est une conséquence directe du contrôle en temps  \begin{equation*}
\Vert \vu(t,\cdot)\Vert^{2}_{L^2}\leq e^{-2 \alpha t}\Vert \vu_0 \Vert^{2}_{L^2}+c\frac{\Vert \fe \Vert^{2}_{\dot{H}^{-1}}}{\nu \alpha}(1-e^{-2 \alpha t}),
\end{equation*} donné par le Théorème \ref{Theo:controle-temps}. En effet,   nous écrivons  $$\ds{\frac{1}{T}\int_{0}^{T} \Vert \vu(t,\cdot)\Vert^{2}_{L^2} dt\leq \frac{1}{T}\int_{0}^{T} e^{-2 \alpha  t}\Vert \vu_0 \Vert^{2}_{L^2}dt +\frac{1}{T}\int_{0}^{T}c\frac{\Vert \fe \Vert^{2}_{\dot{H}^{-1}}}{\nu \alpha} \left( 1-e^{-2\alpha t}\right)dt}$$ et comme la force $\fe$ est une fonction stationnaire alors en prenant la limite $\ds{\limsup_{T\longrightarrow +\infty}}$ nous obtenons l'estimation recherchée. \finpv  \\
Nous observons de cette façon que le terme d'amortissement $-\alpha \vu$  dans les équations (\ref{NS-amortie-alpha}) nous a permis d'obtenir un contrôle  en temps de la quantité $\Vert \vu(t,\cdot)\Vert^{2}_{L^2}$  donné par le biais du Théorème  \ref{Theo:controle-temps} et ce résultat est valable pour tout paramètre d'amortissement  $\alpha>0$.  Dans la sous-section qui suit nous allons fixer le paramètre $\alpha>0$ d'une façon convenable et cette valeur particulière de $\alpha$ va nous permettre de faire ensuite (dans la Section  \ref{Sec:Theoreme principale}) une discussion sur la loi de dissipation d'énergie dans le cadre des équations  (\ref{NS-amortie-alpha}).  
\subsection{Le paramètre d'amortissement}\label{sec:parametre-amortissement}  
Nous fixons ici le paramètre $\alpha>0$   dans le terme d'amortissement $-\alpha \vu$ des équations (\ref{NS-amortie-alpha}) et pour fixer les idées nous avons besoin de considérer pour l'instant le cadre d'un fluide périodique en variable d'espace qui a été introduit dans la Section \ref{Sec:cadre-periodique}. \\
\\
Rappelons que dans le cadre périodique nous considérons une période  $L\geq \ell_0$,   où $\ell_0>0$ est toujours une échelle d'injection fixe (voir l'estimation (\ref{relation_el_L_intro}) pour tous les détails ); et  nous considérons les équations  de Navier-Stokes périodiques sur le cube $[0,L]^3 \subset \Rt$. Rappelons maintenant que pour  toute  solution de Leray $\vu \in L^{\infty}_{t}L^{2}_{x}\cap (L^{2}_{t})_{loc}\dot{H}^{1}_{x}$,  par  l'inégalité d'énergie (\ref{ineg_ener_Leray_periodique}) et l'inégalité de Poincaré (\ref{Poincare}) nous  avons l'estimation  (\ref{estim-per-u}) qui nous rappelons ci-dessous:
\begin{equation*}
\frac{\nu}{L^2}\left( \limsup_{T\longrightarrow +\infty}\frac{1}{T} \int_{0}^{T}\Vert \vu(t,\cdot)\Vert^{2}_{L^2}dt\right) \leq \frac{\Vert \fe \Vert^{2}_{\dot{H}^{-1}}}{\nu}, 
\end{equation*} où $\nu>0$ est la constante de viscosité du fluide (voir toujours l'estimation (\ref{estim-per-u}) page \pageref{estim-per-u} pour tous les détails de la preuve de cette estimation). Nous observons ainsi que dans le cadre périodique nous avons toujours un contrôle sur la moyenne en temps long de la quantité $\Vert \vu(t,)\Vert^{2}_{L^2}$ et que ce contrôle est dû à l'inégalité de Poincaré. \\
\\
De cette estimation nous nous intéressons au terme $ \ds{\frac{\nu}{L^2}}$ devant la moyenne en temps long  car ce terme nous permettra de fixer une valeur assez naturelle du paramètre $\alpha>0$ dans le terme d'amortissement $-\alpha \vu$ des équations (\ref{NS-amortie-alpha}) comme nous l'expliquons tout de suite. \\
\\
Revenons à présent au cadre d'un fluide non périodique posé dans l'espace $\Rt$ tout entier et à notre modèle des équations de Navier-Stokes amorties (\ref{NS-amortie-alpha}). Dans ce cadre, par  le Corollaire \ref{Coro:theo_controle_temps_alpha_modele} nous savons que toute  solution faible $\vu \in L^{\infty}_{t}L^{2}_{x}\cap (L^{2}_{t})_{loc}\dot{H}^{1}_{x}$ des équations (\ref{NS-amortie-alpha})  vérifie l'estimation: 
\begin{equation}
\alpha \left( \limsup_{T\longrightarrow +\infty}\frac{1}{T} \int_{0}^{T}\Vert \vu(t,\cdot)\Vert^{2}_{L^2}dt\right) \leq \frac{\Vert \fe \Vert^{2}_{\dot{H}^{-1}}}{\nu}, 
\end{equation}et si nous comparons maintenant cette estimation avec l'estimation ci-dessus nous pouvons observer que le paramètre $\alpha$ joue le rôle du terme $\ds{\frac{\nu}{L^2}} $ qui apparaît dans le cadre d'un fluide périodique où $L\geq \ell_0$ est la période. De cette façon, par analogie au cadre périodique, pour  $\ell_0>0$ une échelle d'injection d'énergie fixe, nous intégrons à notre modèle une longueur $L\geq \ell_0$  qui est un paramètre du modèle tout comme la constante de viscosité du fluide $\nu$ et qui 
représente un analogue à la période dans le cadre d'un fluide périodique; et  nous allons fixer le paramètre d'amortissement $\alpha>0$ par l'expression
\begin{equation} 
\alpha=\frac{\nu}{L^2}.
\end{equation} Ainsi dans ce chapitre nous allons dorénavant travailler  avec le système d'équations de Navier-Stokes amorties suivant:

\begin{equation}\label{NS-amortie-alpha-fix} 
\left\lbrace \begin{array}{lc}\vspace{3mm}
\partial_t\vu = \nu\Delta \vu-(\vu\cdot\vec{\nabla}) \vu -\vec{\nabla}p +\fe-\ds{\frac{\nu}{L^{2}} \vu}, \quad div(\vu)=0,\quad div(\fe)=0, \quad \nu>0,\,\, L\geq \ell_0,&\\
\vu(0,\cdot)=\vu_0, \quad div(\vu_0)=0.&
\end{array}\right. 
\end{equation}
Dans la section qui suit nous allons faire une discussion plus précise sur la longueur $L\geq \ell_0$ ci-dessus  et nous allons  voir le rôle de cette  longueur  dans l'étude déterministe  de la loi de dissipation de d'énergie (\ref{loi_dissipation_physique}) dans le cadre des équations ci-dessus.        
 
\section[Discussion sur la loi de dissipation d'énergie]{Discussion sur  la loi de dissipation d'énergie dans les équations de Navier-Stokes amorties}\label{Sec:Theoreme principale}

Le but de cette section est de faire  une discussion rigoureuse sur l'étude de  l'encadrement du taux de dissipation d'énergie $\varepsilon$:
\begin{equation}\label{Kolmogorov_deter_2}
c_1\frac{U^3}{\ell_0}\leq \varepsilon \leq c_2 \frac{U^3}{\ell_0},
\end{equation} qui est censé être observé dans le régime turbulent des grandes valeurs  du nombre de Reynolds
\begin{equation}\label{Reynolds-discussion}
Re=\frac{U\ell_0}{\nu}
\end{equation} selon la loi de dissipation d'énergie (\ref{loi_dissipation_physique}) énoncée page \pageref{loi_dissipation_physique} et  proposée par la théorie K41.\\ 
\\ 
Dans cet encadrement  la vitesse caractéristique $U$ et le taux de dissipation $\varepsilon$ sont définies à partir des solutions faibles  $\vu$  des équations de Navier-Stokes amorties  données dans (\ref{NS-amortie-alpha-fix})  où $\fe$ est une force extérieure stationnaire et à divergence nulle sur laquelle nous ferons quelques hypothèses supplémentaires,  tandis que  $\ell_0$ est l'échelle d'injection d'énergie qui sera fixée par la force $\fe$ ci-après. \\
\\
Nous commençons donc  par définir quelques quantités dont nous aurons besoin pour faire notre étude  et la première chose à faire est de fixer la force extérieure $\fe$. Nous définissons  cette force en considérant $\fe\in L^2(\Rt)$ un champ de vecteurs stationnaire, c'est à dire $\fe(t,x)=\fe(x)$, à divergence nulle et tel que sa transformée de Fourier satisfait     
\begin{equation}\label{loc_fe}
supp\left(\widehat{\fe}\right) \subset \left\lbrace \xi\in \Rt: \frac{\rho_1}{ \ell_0}\leq \vert\xi \vert\leq \frac{\rho_2}{\ell_0}\right\rbrace,
\end{equation} pour une échelle d'injection d'énergie $\ell_0>0$ donnée et fixée une fois pour toute; et où $0<\rho_1<\rho_2$ sont deux constantes qui ne dépendent  d'aucun paramètre physique. Cette localisation de la transformée de Fourier de $\fe$ représente le fait que, selon le modèle de cascade d'énergie, l'énergie cinétique est introduite dans le fluide (par la force extérieure) uniquement aux échelles de longueur de l'ordre de $\ell_0$ et donc aux fréquences de l'ordre $\frac{1}{\ell_0}$. \\
\\
Nous voulons maintenant définir la vitesse caractéristique $U>0$ et le taux de dissipation d'énergie $\varepsilon>0$. Pour $\vu$ une solution faible des équations de Navier-Stokes amorties (\ref{NS-amortie-alpha-fix}) associée à la force $\fe$ ci-dessus il s'agit de moyenner cette fonction $\vu$ et toutes ses dérivées $\vec{\nabla}\otimes \vu$  tout d'abord en variable d'espace et ensuite en variable du temps  pour obtenir de cette façon une vitesse caractéristique $U$ et un taux moyen de dissipation  d'énergie $\varepsilon$ respectivement. Néanmoins, pour moyenner ces deux fonctions en variable d'espace nous allons voir que l'on rencontre  quelques difficultés. 

\subsection{La notion de longueur caractéristique du fluide}\label{Sec:problematique}

Dans cette section nous expliquons  que l'essentiel des difficultés que l'on rencontre lorsqu'on considère un  fluide posé sur tout l'espace $\Rt$  repose  sur la notion  de longueur caractéristique du fluide. \\
\\
Rappelons rapidement que 
dans le cadre d'un fluide périodique cette longueur caractéristique est définie de façon naturelle par la période $L$ (voir toujours la Définition \ref{Def_L_periodique}) mais si le fluide est posé dans tout l'espace nous pouvons alors observer  que l'on perd toute notion physique et mathématique de cette   longueur caractéristique. 
Ainsi,  nous nous posons alors  la question de comment  choisir  une longueur $L\geq \ell_0$ adéquate pour  définir la moyenne en espace des fonctions $\vu$ et $\vec{\nabla}\otimes \vu$  en termes de la norme $L^2:$ $\ds{\frac{\Vert \cdot \Vert_{L^2}}{L^{\frac{3}{2}}}}$, mais la réponse à cette question n'est pas évidente. \\
\\
Dans ce cadre,  les notes de cours \cite{Const} de P. Constantin suggèrent de considérer  $L=\ell_0$ (où $\ell_0>0$ est toujours l'échelle d'injection d'énergie définie das (\ref{loc_fe})) pour définir la moyenne en espace ci-dessus mais nous allons maintenant observer que le choix de cette longueur présente quelques lacunes.  En effet,  rappelons tout d'abord que la force $\fe$ est une fonction localisée aux fréquences de l'ordre de $\frac{1}{\ell_0}$ (voir l'expression (\ref{loc_fe})) et par  cette localisation fréquentielle nous observons que la quantité 
\begin{equation}\label{moyenne-fourier-f}
\left(\int_{\Rt} \vert \widehat{\fe}(\xi)\vert^2 \frac{dx}{\ell^{3}_{0}}\right)^{\frac{1}{2}}=\frac{\Vert \widehat{\fe} \Vert_{L^2}}{\ell^{\frac{3}{2}}_{0}},
\end{equation}
est une moyenne naturelle  en termes de la norme $L^2$ de la transformée de Fourier de $\fe$. Après,  par l'identité de Plancherel nous savons que $\Vert \fe \Vert_{L^2}=\Vert \widehat{\fe} \Vert_{L^2}$ et alors dans \cite{Const} on considère la moyenne en espace de la force $\fe$ comme  la quantité 
\begin{equation}\label{def:F}
F=\frac{\Vert \fe \Vert_{L^2}}{\ell^{\frac{3}{2}}_{0}}.
\end{equation} 
Mais en variable d'espace nous n'avons  aucune information supplémentaire sur la localisation de la fonction $\fe$ et alors la moyenne  $F>0$ ci-dessus n'est pas bien comprise. Dans la Section \ref{Sec:force ext: construction et proprietes} du chapitre suivant nous revisitons ce problème et nous construisons un exemple concret de force $\fe$ qui est une fonction bien localisée en variable de fréquence mais aussi en variable d'espace.\\
\\ 
Ensuite, toujours dans  \cite{Const}, il est suggéré  de considérer cette même moyenne en espace (\ref{def:F}) pour les fonctions $\vu$ et $\vec{\nabla}\otimes \vu$ et l'on a ainsi les quantités   $$\ds{\frac{\Vert \vu(t,\cdot)\Vert_{L^2}}{\ell^{\frac{3}{2}}_{0}} \quad \text{et}\quad \nu \frac{\Vert \vec{\nabla}\otimes\vu(t,\cdot)\Vert_{L^2}}{\ell^{\frac{3}{2}}_{0}}},$$ 
néanmoins, le choix de la longueur $\ell_0$ pour définir  cette moyenne en espace pour les fonctions $\vu$ et $\vec{\nabla}\otimes \vu$ n'est pas du tout  clair car il n'a  aucune explication rigoureuse que ce soit du point de vue  physique ou mathématique; et  comme nous voulons faire  une discussion aussi rigoureuse que possible de la loi de dissipation d'énergie (\ref{Kolmogorov_deter_2}) alors nous n'allons pas considérer ici cette moyenne en espace.\\
\\
Nous observons de cette façon que dans le cadre d'un fluide dans tout l'espace, une définition convenable de moyenne en espace pour le champ de vitesse $\vu$ et ses dérivées $\vec{\nabla}\otimes \vu$  est une question que l'on ne sait pas répondre de façon tout à fait satisfaisante car cette question est directement reliée au problème de trouver une définition convenable de longueur  caractéristique $L\geq \ell_0$. Ainsi, pour pouvoir faire une étude  rigoureuse  de l'encadrement du taux de dissipation (\ref{Kolmogorov_deter_2}), dans un premier temps nous allons seulement  considérer la moyenne en temps long des quantités $\Vert \vu(t,\cdot)\Vert_{L^2}$ et $\nu\Vert \vec{\nabla}\otimes \vu(t,\cdot)\Vert_{L^2}$ qui sont définies de la façon suivante. \\
\\
La moyenne en temps de la quantité $\Vert \vu(t,\cdot)\Vert_{L^2}$  a été introduite dans l'expression  (\ref{def-u}) page \pageref{def-u} mais pour la commodité du lecteur nous allons récrire cette moyenne donnée par l'expression:
\begin{equation}\label{u}
u=\left( \limsup_{T\longrightarrow +\infty} \frac{1}{T}\int_{0}^{T} \Vert \vu(t,\cdot)\Vert^{2}_{L^2} dt \right)^{\frac{1}{2}},
\end{equation}
où rappelons que, dans le cadre des ces équations amorties, le Corollaire \ref{Coro:theo_controle_temps_alpha_modele} nous assure que la moyenne en temps long $u$  est une quantité bien définie. \\
\\
Rappelons aussi que cette  moyenne en temps long a bien une signification mathématique: la force $\fe$ étant toujours une fonction stationnaire alors nous nous intéressons à étudier le comportement turbulent du fluide dans le régime asymptotique lorsque le temps $T$ tend vers l'infini (voir toujours  la Section \ref{Sec:cadre-periodique} page \pageref{U-per} pour plus de détails à ce sujet).\\
\\
Ensuite,  en suivant les idées  ci-dessous nous allons maintenant définir la moyenne en temps long de la quantité $\nu\Vert \vec{\nabla}\otimes \vu(t,\cdot)\Vert_{L^2}$ par l'expression: 

\begin{equation}\label{e}
 e=\nu \limsup_{T\longrightarrow +\infty} \frac{1}{T}\int_{0}^{T} \Vert  \vec{\nabla}\otimes \vu(t,\cdot)\Vert^{2}_{L^2} dt.\\
\end{equation} 
Une fois que nous avons introduit les quantités $u$ et $e$ ci-dessus expliquons de façon plus précise pourquoi nous allons considérer ici ces quantités pour l'étude de l'encadrement (\ref{Kolmogorov_deter_2}).  Nous pouvons observer que cet encadrement fait intervenir la vitesse caractéristique $U$ et le taux moyen de dissipation $\varepsilon$, qui dans le cadre d'un fluide périodique sur le cube $[0,L]^3$  s'écrivent de façon rigoureuse par les expressions:

\begin{equation}\label{U,var,per-motiv}
U=\left( \limsup_{T\longrightarrow +\infty} \frac{1}{T}\int_{0}^{T} \Vert \vu(t,\cdot)\Vert^{2}_{L^2}\frac{dt}{L^3}\right)^{\frac{1}{2}}\quad \text{et}\quad \varepsilon = \nu \limsup_{T\longrightarrow +\infty} \frac{1}{T}\int_{0}^{T} \Vert  \vec{\nabla}\otimes \vu(t,\cdot)\Vert^{2}_{L^2}\frac{dt}{L^3}.
\end{equation}  
Si nous comparons maintenant  les quantités $u$ et $e$ avec les quantités $U$ et $\varepsilon$ respectivement nous observons que si nous considérons n'importe quelle longueur $L\geq \ell_0$, qui pourrait  représenter un analogue à la période dans le cadre d'un fluide périodique, alors nous avons les relations
\begin{equation}\label{relation:U-u}
U=\frac{u}{L^{\frac{3}{2}}} \quad \text{et} \quad \varepsilon =\frac{e}{L^3}, 
\end{equation}  et alors  il s'agit donc d'étudier l'encadrement (\ref{Kolmogorov_deter_2}) tout d'abord dans le cadre rigoureux des quantités $u$ et $e$ et ceci sera fait dans le Théorème \ref{Theo:loi-kolmogorov} dans la Section \ref{sec:estim-rigoureuses}  ci-dessous. Ensuite dans la Section \ref{sec:conclusions} nous allons observer qu'à partir de ce résultat  et en considérant les quantités $U$ et $\varepsilon$, avec n'importe quelle longueur  $L> \ell_0$, nous obtenons alors une estimation du type  $\varepsilon \lesssim \frac{U^3}{\ell_0}$ bien qu'il s'agisse d'une estimation partielle par rapport à l'encadrement (\ref{Kolmogorov_deter_2}).
\subsection{Quelques estimations rigoureuses }\label{sec:estim-rigoureuses} 

Comme annoncé, dans cette section nous allons étudier l'encadrement du taux de dissipation (\ref{Kolmogorov_deter_2}): $ \varepsilon \approx \frac{U^3}{\ell_0}$, dans le cadre plus rigoureux des moyennes en temps long $u$ et $e$ introduites précédemment. Plus précisément, dans le Théorème \ref{Theo:loi-kolmogorov} ci-dessous nous allons montrer l'estimation:
\begin{equation}\label{estim-theo-motiv}
 e \lesssim   \frac{u_{\ell_0} u^2}{\ell_0} \left(\frac{\Vert \fe \Vert_{L^{\infty}}}{\Vert \fe \Vert_{L^{2}}} \right), 
\end{equation} où la quantité $u_{\ell_0}>0$ est définie de la façon suivante: rappelons tout d'abord que la force $\fe$ est une fonction localisée aux fréquences $\frac{\rho_1}{\ell_0}\leq \vert \xi \vert \leq \frac{\rho_2}{\ell_0}$ et alors nous définissons la fonction $\U_{\ell_0}$ comme  la localisation du champ de vitesse $\vu$ aux mêmes fréquences: 
\begin{equation}\label{def:U0}
\U_{\ell_0}(t,x)=\mathcal{F}^{-1} \left[ \mathds{1}_{\frac{\rho_1}{\ell_0}\leq \vert \xi \vert \leq \frac{\rho_2}{\ell_0}} (\xi) \widehat{\vu}(t,\cdot) \right](x),
\end{equation} où $\mathcal{F}^{-1}$ dénote la transformée de Fourier inverse;  et ensuite la quantité $u_{\ell_0}>0$ est définie comme la moyenne en temps long de cette fonction:
\begin{equation}\label{u_0}
u_{\ell_0}= \left( \limsup_{T \longrightarrow +\infty}\frac{1}{T} \int_{0}^{T} \left\Vert \U_{\ell_0}(t,\cdot) \right\Vert^{2}_{L^2}dt \right)^{\frac{1}{2}}. 
\end{equation} 
Avant d'entrer dans les détails techniques de la preuve de l'estimation (\ref{estim-theo-motiv}) nous allons tout d'abord  expliquer  cette estimation  et 
pour cela nous avons besoin de rappeler rapidement  l'étude de la loi de dissipation d'énergie (\ref{Kolmogorov_deter_2}) dans le cadre d'un fluide périodique qui a été exposée dans la Section \ref{Sec:cadre-periodique}.\\
\\
Rappelons que dans le cadre périodique on a une estimation du taux de dissipation $\varepsilon$: 
\begin{equation}\label{estim-var-motiv-per}
\varepsilon \lesssim \frac{U^3}{\ell_0},
\end{equation} et nous  maintenant expliquer les grandes lignes de la preuve de cette estimation (voir aussi l'article \cite{DoerFoias} de C. Doering et C. Foias pour tous les détails des calculs). Dans ce cadre périodique nous définissons la force moyenne $F>0$ comme $\ds{ F=\frac{\Vert \fe \Vert_{L^2}}{L^{\frac{3}{2}}}}$, le nombre de Reynolds $\ds{ Re=\frac{U \ell_0}{\nu}}$; et l'estimation (\ref{estim-var-motiv-per})  repose essentiellement sur les deux inégalités  techniques suivantes: 
\begin{equation}\label{estim:motiv-per1}
F \lesssim \frac{U^2}{\ell_0}\left(1 +\frac{1}{Re}\right),
\end{equation} et 
\begin{equation}\label{estim:motiv-per2}
\varepsilon  \lesssim U F. 
\end{equation} En effet,  dans le régime turbulent caractérisé lorsque $Re>>1$  dans l'inégalité (\ref{estim:motiv-per1})  nous obtenons alors  $F \lesssim \frac{U^2}{\ell_0}$,  %
et en multipliant par $U$ à chaque côté de cette inégalité nous avons $UF\lesssim \frac{U^3}{\ell_0}$, d'où,  
par l'inégalité (\ref{estim:motiv-per2}) nous pouvons  finalement écrire $\varepsilon \lesssim FU \lesssim \frac{U^3}{\ell_0}$ ce qui nous donne l'estimation (\ref{estim-var-motiv-per}). \\
\\
\`A ce stade  il est important de souligner que, dans l'état actuel de nos connaissances, l'étude déterministe de l'estimation (\ref{estim-var-motiv-per}) que l'on peut trouver dans la littérature suit essentiellement les grandes lignes expliquées ci-dessus et que les résultats que  l'on obtient (toujours dans le cadre périodique) dans  différentes contextes techniques sont toujours  similaires à cette estimation.  \\ 
\\ 
Avec ces idées en tête,  revenons à présent à notre cadre d'étude d'un fluide dans tout l'espace $\Rt$. Soit donc la force $\fe \in L^2(\Rt)$ définie par (\ref{loc_fe}) et soit $\vu$ une solution faible des équations de Navier-Stokes amorties (\ref{NS-amortie-alpha-fix}). 
Dans le Théorème \ref{Theo:loi-kolmogorov} nous allons vérifier l'estimation (\ref{estim-theo-motiv}) qui s'agit d'une estimation  analogue à l'estimation (\ref{estim-var-motiv-per}) (obtenue dans le cadre périodique)  mais en considérant seulement les quantités $\Vert \fe \Vert_{L^{\infty}}$, $\Vert \fe \Vert_{L^2}$, $u$, $u_0$ et $e$ qui ont bien un sens mathématique. \\
\\
 Ainsi, pour vérifier l'estimation (\ref{estim-theo-motiv}) nous allons suivre  les grandes lignes de la preuve de l'estimation (\ref{estim-var-motiv-per})  et alors il s'agit  de vérifier des inégalités analogues à celles données dans  
(\ref{estim:motiv-per1}) et (\ref{estim:motiv-per2})
(dans le cadre périodique).  Observons  tout d'abord que l'estimation (\ref{estim:motiv-per1}) fait intervenir le nombre de Reynolds $Re$ qui est défini à partir de  la vitesse caractéristique $U$ comme $$Re=\frac{U\ell_0}{\nu},$$ mais  comme  nous ne disposons pas ici d'une définition rigoureuse de cette vitesse caractéristique, en suivant les idées précédentes, nous allons alors remplacer cette quantité $U$ par la moyenne en temps long  $u$ et nous allons donc considérer  le nombre de Reynolds $ \mathcal{R}e$ défini comme 
\begin{equation}\label{Re-rigoureux}
\mathcal{R}e=\frac{u \ell_0}{\nu},
\end{equation} où $\ell_0>0$ est toujours l'échelle d'injection d'énergie fixée par la force $\fe$ dans (\ref{loc_fe}) et $\nu>0$ est toujours la constante de viscosité du fluide. \\
\\
\'Etudions  maintenant  la relation entre ce nouveau nombre de Reynolds $\mathcal{R}e$ et le nombre de Reynolds classique $Re$, où nous pouvons observer que pour toute longueur   $L\geq \ell_0$ on a 
\begin{equation}\label{Prop:relation-reynolds} 
Re=\frac{1}{L^{\frac{3}{2}}} \mathcal{R}e.
\end{equation}
%
En effet, il suffit de remarquer que pour une longueur $L\geq \ell_0$, par la relation  (\ref{relation:U-u}) on a $U= \frac{1}{L^{\frac{3}{2}}}u$ d'où nous avons directement $Re=\frac{U \ell_0}{\nu}= \frac{1}{L^{\frac{3}{2}}} \left(\frac{u \ell_0}{\nu}\right)= \frac{1}{L^{\frac{3}{2}}} \mathcal{R}e$. \\
\\
Maintenant que l'on dispose de cette identité, nous allons expliquer  comment le nombre $\mathcal{R}e$ nous permet aussi de caractériser le régime turbulent du fluide. Il s'agit de fixer le nombre $\mathcal{R}e$ suffisamment grand de sorte que ceci entraîne $Re>>1$ (qui caractérise le régime turbulent du fluide) et pour cela nous allons suivre le raisonnement suivant: pour $\ell_0>0$ l'échelle d'injection d'énergie rappelons que nous fixons le paramètre $L\geq \ell_0$ qui représente  la longueur caractéristique du fluide. Ensuite, si nous fixons le nombre $\mathcal{R}e$ tel que $\mathcal{R}e>>L$  alors on a $\frac{\mathcal{R}e}{L^{\frac{3}{2}}}>>1$ et donc, par la Proposition \ref{Prop:relation-reynolds} on obtient $Re>>1$. \\
\\
Nous observons ainsi  que le régime asymptotique des grandes valeurs du nombre $\mathcal{R}e$  ($\mathcal{R}e>>L$) entraîne le régime turbulent caractérisé par de grandes valeurs du nombre $Re$  et comme le nombre $\mathcal{R}e$ fait intervenir la moyenne en temps long $u$ au lieu de la vitesse caractéristique $U$ nous allons préférer ici ce nombre de Reynolds $\mathcal{R}e$ pour caractériser le régime turbulent. \\ 
\\  
Une fois que l'on a introduit le nombre de Reynolds $\mathcal{R}e$ ci-dessus, dans l'estimation suivante nous obtenons une inégalité analogue à l'inégalité (\ref{estim:motiv-per1}) obtenue dans le cadre périodique.

\begin{Proposition}[Première estimation dans le cadre non périodique]\label{Prop1} Soit $\ell_0>0$  une échelle d'injection d'énergie définie par la force $\fe\in L^2(\Rt)$ dans (\ref{loc_fe}). Soit  $\vu \in L^{\infty}_{t}L^{2}_{x}\cap (L^{2}_{t})_{loc}\dot{H}^{1}_{x}$ une solution faible des équations de Navier-Stokes amorties (\ref{NS-amortie-alpha-fix}) associée à cette force. Soit la moyenne en temps long $u$ définie à partir de  la solution $\vu$ dans l'expression (\ref{u}). Soit enfin le nombre de Reynolds $\mathcal{R}e$ défini dans l'expression (\ref{Re-rigoureux}). Alors on a l'estimation 
\begin{equation}\label{estim-rigoureux-motiv2}
\Vert \fe\Vert_{L^2} \leq \frac{u^2}{\ell_0} \left( \frac{\Vert \fe \Vert_{L^{\infty}}}{\Vert \fe \Vert_{L^2}} +\frac{1}{\mathcal{R}e} \right).
\end{equation}	
\end{Proposition}	
\pv  L'inégalité (\ref{estim-rigoureux-motiv2}) repose essentiellement sur l'estimation technique suivante.
\begin{Lemme}\label{lemme-tech1} Dans le cadre la Proposition \ref{Prop1} on a l'estimation: 
	\begin{equation}\label{ineq-aux-3}
	\Vert \fe \Vert^{2}_{L^2}\leq u^2 \Vert \vec{\nabla} \otimes \fe \Vert_{L^{\infty}} +\nu\, u  \Vert \Delta \fe \Vert_{L^2} +\frac{\nu}{\ell^{2}_{0}} u \Vert  \fe \Vert_{L^2}. 
	\end{equation}
\end{Lemme}	
La preuve de cette estimation suit les grandes lignes des notes de cours \cite{Const} de P. Constantin et pour la commodité du lecteur nous ferons tous les calculs en détail à la fin du chapitre. \\
\\
Nous allons maintenant étudier le premier et deuxième  terme à droite de l'estimation ci-dessus et alors, étant donnée que la force $\fe$ est localisée aux fréquences $\frac{\rho_1}{\ell_0}\leq \vert \xi \vert \leq \frac{\rho_2}{\ell_0}$, par les inégalités de Bernstein nous avons qu'il existe une constante $c>0$, qui ne dépend  d'aucun paramètre physique, telle que l'on a 
$$ \Vert \vec{\nabla} \otimes \fe \Vert_{L^{\infty}} \approx \frac{\Vert \fe \vert_{L^{\infty}}}{\ell_0}  \quad \text{et}\quad \Vert \Delta \fe \Vert_{L^2}\leq c \frac{1}{\ell^{2}_{0}} \Vert \fe \Vert_{L^2}.$$ Nous remplaçons maintenant les estimations ci-dessus dans (\ref{ineq-aux-3}) et nous obtenons l'estimation 
$$ \Vert \fe \Vert^{2}_{L^2}\leq c \frac{u^2}{\ell_0}   \Vert \fe \Vert_{L^{\infty}} + c\frac{\nu}{\ell^{2}_{0}} u \Vert \fe \Vert_{L^2} +  \frac{\nu}{\ell^{2}_{0}} u \Vert \fe \Vert_{L^2},$$ d'où nous pouvons écrire 
$$ \Vert \fe \Vert_{L^2} \leq c \frac{u^2}{\ell_0} \left( \frac{\Vert \fe \Vert_{L^{\infty}}}{\Vert \fe \Vert_{L^2}} + \frac{\nu}{u \ell_0}  \right), $$ et comme  l'on a définit le nombre de Reynolds $\mathcal{R}e$ (donné dans (\ref{Re-rigoureux}))  par l'expression $\mathcal{R}e=\frac{u \ell_0}{\nu}$ nous obtenons  l'estimation cherchée $\ds{\Vert \fe \Vert_{L^2}\leq c \frac{u^2}{\ell_0}\left( \frac{\Vert \fe \Vert_{L^{\infty}}}{\Vert \fe \Vert_{L^2}} +\frac{1}{\mathcal{R}e}\right)}$. \finpv
\\
Nous allons maintenant étudier une inégalité analogue à l'inégalité (\ref{estim:motiv-per2}): $\varepsilon\lesssim UF$, obtenue dans le cadre périodique et pour cela on commence par faire la remarque suivante: observons que dans l'inégalité (\ref{estim:motiv-per2}) interviennent les termes  $\varepsilon$, $U$ et $F$; et étant donné que ces  termes sont définis à partir des moyennes en temps long  $e$, $u$ (données dans les expression (\ref{e}) et (\ref{u}) respectivement),  la quantité $\Vert \fe \Vert_{L^2}$ et la période $L$  comme:  $\varepsilon=\frac{e}{L^{3}}$, $U=\frac{u}{L^{\frac{3}{2}}}$ et $F=\frac{\Vert \fe \Vert_{L^2}}{L^{\frac{3}{2}}}$; nous pouvons alors écrire  $\ds{\frac{e}{L^{3}} \leq \frac{u}{L^{\frac{3}{2}}} \frac{\Vert \fe \Vert_{L^2}}{L^{\frac{3}{2}}}}$, d'où nous obtenons l'estimation 
\begin{equation}\label{estim-per-motiv}
e \leq u \Vert \fe \Vert_{L^2}.
\end{equation}  Nous observons que la période $L$ ne joue aucun rôle dans cette estimation et ainsi,  par analogie au cadre périodique, nous voulons alors étudier une estimation du même type que celle  ci-dessus. 
\begin{Proposition}[Deuxième estimation dans le cadre non périodique]\label{Prop2} Soit $\ell_0>0$ une échelle d'injection d'énergie définie par la force $\fe\in L^2(\Rt)$ dans (\ref{loc_fe}). Soit  $\vu \in L^{\infty}_{t}L^{2}_{x}\cap (L^{2}_{t})_{loc}\dot{H}^{1}_{x}$ une solution faible des équations de Navier-Stokes amorties (\ref{NS-amortie-alpha-fix}) associée à cette force. Soient les moyennes en temps $e>0$ et $u_{\ell_0}>0$ définies à partir de la solution $\vu$ par les expressions (\ref{e}) et (\ref{u_0}) respectivement. Alors on a l'estimation:
	\begin{equation}\label{estim-rigoureux-motiv1}
	e \leq u_{\ell_0} \Vert \fe \Vert_{L^2}. \\
	\end{equation}
\end{Proposition}
Si nous comparons cette inégalité (\ref{estim-rigoureux-motiv1}) avec l'inégalité (\ref{estim-per-motiv}) nous pouvons observer que cette première est plus précise car cette inégalité fait intervenir la quantité $u_{\ell_0}$ au lieu de la quantité $u$ dans (\ref{estim-per-motiv}). En effet, rappelons que la quantité $u_{\ell_0}$ est définie dans (\ref{u_0}) est il correspond à  la moyenne en temps long de la fonction $\U_{\ell_0}$ qui est la localisation fréquentielle du champ de vitesse $\vu$ (voir l'expression (\ref{def:U0}) pour une définition de la fonction $\U_{\ell_0}$) et ainsi, comme l'on a $\Vert \U_{\ell_0}(t,)\Vert_{L^2} \leq \Vert \vu(t,)\Vert_{L^2}$ nous avons alors  
\begin{equation}\label{estim:u0-u}
u_{\ell_0} \leq u. 
\end{equation} 
\textbf{Preuve de la Proposition \ref{Prop2}.}  La preuve de cette estimation repose sur l'inégalité d'énergie vérifiée par la solution $\vu$ et qui a été obtenue dans la Proposition \ref{Proposition_ineq_energie_alpha_modele}: 
\begin{eqnarray*}
	\Vert \vu(T,\cdot)\Vert^{2}_{L^2}+2\nu\int_{0}^{T} \Vert \nabla \otimes \vu(t,\cdot)\Vert^{2}_{L^2}dt &\leq &\Vert \vu_0 \Vert^{2}_{L^2}+2\int_{0}^{T} \int_{\Rt} \fe(x)\cdot  \vu(t,\cdot) dx\,dt\\
	& & -2 \frac{\nu}{L^2} \int_{0}^{T}\Vert \vu(t,\cdot) \Vert^{2}_{L^2}dt,
\end{eqnarray*} d'où, étant donné que $\Vert \vu(T,\cdot)\Vert^{2}_{L^2}$ est une quantité positive et de plus, comme $2 \frac{\nu}{L^2} \int_{0}^{T}\Vert \vu(t,\cdot) \Vert^{2}_{L^2}dt$ est une quantité négative, nous pouvons écrire 
\begin{equation}\label{ineq-aux-1}
\nu\int_{0}^{T} \Vert \nabla \otimes \vu(t,\cdot)\Vert^{2}_{L^2}dt \leq \frac{1}{2} \Vert \vu_0 \Vert^{2}_{L^2}+\int_{0}^{T} \int_{\Rt} \fe(x)\cdot  \vu(t,\cdot) dx\,dt.
\end{equation} Nous allons maintenant étudier  le deuxième terme à droite de cette estimation et pour la fonction $\U_{\ell_0}$ définie dans l'expression (\ref{def:U0}) comme 
$$  \U_{\ell_0}=\mathcal{F}^{-1} \left[ \mathds{1}_{\frac{\rho_1}{\ell_0}\leq \vert \xi \vert \leq \frac{\rho_2}{\ell_0}} (\xi)\mathcal{F}[\vu](t,\cdot) \right],$$ nous allons montrer que l'estimation 
\begin{equation}\label{ineq-aux-2}
\int_{0}^{T} \int_{\Rt} \fe(x)\cdot  \vu(t,\cdot) dx\,dt \leq T^{\frac{1}{2}}\left( \int_{0}^{T}\Vert \U_{\ell_0}(t,\cdot)\Vert^{2}_{L^2} dt \right)^{\frac{1}{2}}\Vert \fe \Vert_{L^2}. 
\end{equation} 
En effet, comme la force $\fe$ est localisée aux fréquences: $supp\left(\widehat{\fe}\right) \subset \left\lbrace \xi\in \Rt: \frac{\rho_1}{ \ell_0}\leq \vert\xi \vert\leq \frac{\rho_2}{\ell_0}\right\rbrace,$ nous appliquons l'identité  de Parseval pour écrire
\begin{eqnarray*}
	\int_{\Rt} \fe(x)\cdot  \vu(t,\cdot) dx\,dt&=&  \int_{\Rt} \widehat{\fe} (\xi)\cdot  \widehat{\vu}(t,\xi) d\xi = \int_{\Rt} \left[ \mathds{1}_{\{  \frac{\rho_1}{ \ell_0}\leq \vert\xi \vert\leq \frac{\rho_2}{\ell_0}\}}(\xi) \widehat{\fe} (\xi)\right]\cdot  \widehat{\vu}(t,\xi) d\xi\\
	&=& \int_{\Rt}  \widehat{\fe} (\xi)\cdot \left[ \mathds{1}_{\{  \frac{\rho_1}{ \ell_0}\leq \vert\xi \vert\leq \frac{\rho_2}{\ell_0}\}}(\xi)\widehat{\vu}(t,\xi) \right] d\xi,
\end{eqnarray*} et ensuite,  par l'inégalité de Cauchy-Schwarz, l'identité de Plancherel  et par la définition de la fonctions $\U_{\ell_0}$ ci-dessus nous avons
$$ \int_{\Rt}  \widehat{\fe} (\xi)\cdot  \mathds{1}_{\{  \frac{\rho_1}{ \ell_0}\leq \vert\xi \vert\leq \frac{\rho_2}{\ell_0}\}}(\xi)\widehat{\vu}(t,\xi) d\xi \leq \Vert \U_{\ell_0}(t,\cdot)\Vert_{L^2}\Vert \fe \Vert_{L^2}.$$ 
Finalement, nous prenons l'intégrale sur l'intervalle de temps $[0,T]$ à chaque côté de cette inégalité où, étant  donné que $\fe$ est une fonction stationnaire et en appliquant l'inégalité de Cauchy-Schwarz (en variable de temps) nous pouvons alors écrire l'estimation (\ref{ineq-aux-2}). \\
\\
Une fois que l'on a cette estimation, nous la  remplaçons dans l'estimation (\ref{ineq-aux-1}) et nous obtenons 
$$ \nu\int_{0}^{T} \Vert \nabla \otimes \vu(t,\cdot)\Vert^{2}_{L^2}dt \leq \frac{1}{2} \Vert \vu_0 \Vert^{2}_{L^2}+T^{\frac{1}{2}}\left( \int_{0}^{T}\Vert \U_{\ell_0}(t,\cdot)\Vert^{2}_{L^2} dt \right)^{\frac{1}{2}}\Vert \fe \Vert_{L^2}.$$
Ainsi, nous divisons cette estimation par $T$, puis nous prenons la limite supérieure lorsque  $T \longrightarrow +\infty$ et par la définition des moyennes en temps long $e$ et $u_{\ell_0}$ nous obtenons finalement l'estimation cherchée 
$ e \leq u_{\ell_0} \Vert \fe \Vert_{L^2}.$ \finpv 
\\
Nous avons maintenant tous les ingrédients dont on a besoin pour étudier l'estimation (\ref{estim-theo-motiv}) et nous avons ainsi le résultat suivant.

\begin{Theoreme}[Loi de dissipation d'énergie dans le cadre non périodique]\label{Theo:loi-kolmogorov} Soit $\nu>0$ la constante de viscosité du fluide. Soit $\fe \in L^2(\Rt)$ une force à divergence nulle et qui vérifie la localisation fréquentielle 
$$supp\left(\widehat{\fe}\right) \subset \left\lbrace \xi\in \Rt: \frac{\rho_1}{ \ell_0}\leq \vert\xi \vert\leq \frac{\rho_2}{\ell_0}\right\rbrace,$$ pour une échelle d'injection d'énergie $\ell_0>0$ donnée et fixée; et où $0<\rho_1<\rho_2$ sont deux constantes qui ne dépendent pas d'aucun paramètre physique. Soit $\vu \in L^{\infty}_{t}L^{2}_{x}\cap (L^{2}_{t})\dot{H}^{1}_{x}$ une solution faible des équations de Navier-Stokes amorties
$$\partial_t\vu = \nu\Delta \vu-(\vu\cdot\vec{\nabla}) \vu -\vec{\nabla}p +\fe-\ds{\frac{\nu}{L^{2}} \vu}, \quad div(\vu)=0,\quad div(\fe)=0, \quad \nu>0,\,\, L\geq \ell_0,$$ obtenue par le biais du Théorème \ref{Theo:existence-sol-ns-amortie}, et à partir de laquelle on considère  $u$, $u_0$ et $e$ les moyennes en temps long  données par les expressions (\ref{u}), (\ref{e}) et (\ref{u_0}) respectivement. \\  
\\
Si $\mathcal{R}e>>L$, où $\ds{\mathcal{R}e=\frac{u \ell_0}{\nu}}$, alors on a l'estimation donnée dans  (\ref{estim-theo-motiv}):
\begin{equation*}
e \leq c \frac{u_{\ell_0} u^2}{\ell_0}\left(\frac{\Vert \fe \Vert_{L^{\infty}}}{\Vert \fe \Vert_{L^2}}\right),
\end{equation*} où $c>0$ est une constante numérique qui ne dépend  d'aucun paramètre physique ci-dessus. \\
\end{Theoreme}
\textbf{Démonstration.} Ce résultat repose sur les estimations (\ref{estim-rigoureux-motiv2}) et (\ref{estim-rigoureux-motiv1}) obtenues dans les Propositions \ref{Prop1} et \ref{Prop2} respectivement. En effet, observons tout d'abord que si l'on suppose que le nombre $\mathcal{R}e$ est suffisamment grand: $\mathcal{R}e>>L$,  ce qui caractérise le régime turbulent grâce à l'identité  \ref{Prop:relation-reynolds};  alors le terme $\frac{1}{\mathcal{R}e}$  devient négligeable et par l'estimation  (\ref{estim-rigoureux-motiv2}) nous pouvons écrire  $$\ds{\Vert \fe \Vert_{L^2}\leq c \frac{u^2}{\ell_0}\left(  \frac{\Vert \fe \Vert_{L^{\infty}} }{\Vert \fe \Vert_{L^2}}
	\right)}.$$
Ensuite, nous multiplions à chaque côté de cette estimation par la quantité $u_{\ell_0}>0$ et nous obtenons l'estimation suivante:
$$ u_{\ell_0} \Vert \fe \Vert_{L^2}\leq c \frac{u^2}{\ell_0}\left( \frac{\Vert \fe \Vert_{L^{\infty}} }{\Vert \fe \Vert_{L^2}} \right).$$
Finalement, par l'estimation (\ref{estim-rigoureux-motiv1}) nous savons que l'on a $\ds{e\leq u_{\ell_0} \Vert \fe \Vert_{L^2}}$ et ainsi nous pouvons écrire l'estimation cherché:  $ \ds{e \leq c \frac{u_{\ell_0} u^2}{\ell_0}\left(\frac{\Vert \fe \Vert_{L^{\infty}}}{\Vert \fe \Vert_{L^2}}\right),}$
où nous observons que la constante $c>0$ est indépendante de tout paramètre physique de notre modèle et alors il s'agit d'une constante universelle. \finpv
\\
Maintenant que l'on a vérifié l'estimation du terme de dissipation d'énergie  (\ref{estim-theo-motiv}), dans la section qui suit nous allons faire quelques remarques sur cette estimation par rapport à l'étude déterministe de la loi de dissipation de Kolmogorov $\varepsilon \approx \frac{U^3}{\ell_0}$ .   

 \subsection{Conclusions}\label{sec:conclusions}    
 Dans cette section nous allons faire une discussion sur l'estimation (\ref{estim-theo-motiv}) obtenue dans le Théorème \ref{Theo:loi-kolmogorov}.   Insistons tout d'abord sur le fait que cette estimation du terme de dissipation d'énergie $e$ est une estimation rigoureuse car tous les termes qui interviennent ont bien un sens mathématique. Nous souhaitons  maintenant d'expliquer comment, à partir de cette inégalité, on peut récupérer des estimations du taux moyen de dissipation d'énergie $\varepsilon$ selon la loi de Kolmogorov: 
\begin{equation}\label{estim-kolm-diss}
\varepsilon \lesssim \frac{U^3}{\ell_0}. 
\end{equation}
Rappelons rapidement  l'essentiel du problème:  étant donné qu'on travail sur tout l'espace $\Rt$ alors  la longueur caractéristique du fluide $L\geq \ell_0$ n'est pas rigoureusement définie et donc  la vitesse caractéristique $ \ds{U=\frac{u}{L^{\frac{3}{2}}}}$ et le taux moyen de dissipation d'énergie $\ds{\varepsilon=\frac{e}{L^3}}$ ne le sont pas non plus. Ainsi,  avant d'étudier l'estimation (\ref{estim-kolm-diss}) nous avons tout d'abord étudié l'estimation (\ref{estim-theo-motiv}) où  cette longueur $L$ n'intervient pas.\\
\\
Dans ce cadre, en suivant les idées de la Section \ref{sec:parametre-amortissement}, nous allons considérer la longueur $L\geq \ell_0$ tout simplement comme un paramètre pour définir (toujours formellement) les quantités moyennes $U$ et $\varepsilon$, et  à partir de l'estimation rigoureuse (\ref{estim-theo-motiv})  nous allons étudier l'estimation (\ref{estim-kolm-diss}). Toutes les estimations que nous allons obtenir sont des corollaires du Théorème \ref{Theo:loi-kolmogorov} et  étant donné que $L\geq \ell_0$ nous allons diviser notre étude en regardant deux cas: nous allons tout d'abord considérer $L=\ell_0$ pour ensuite étudier le cas lorsque $L>\ell_0$. 
\subsubsection{A) Le cas $L=\ell_0$} 
Rappelons rapidement que dans la formule (\ref{def:F}) page \pageref{def:F} nous avons expliqué que les notes de cours \cite{Const} de P. Constantin suggèrent de considérer la longueur $L=\ell_0$ (où $\ell_0>0$ est toujours une échelle d'injection d'énergie fixe)  pour définir les quantités moyennes $U$ et $\varepsilon$.  Nous suivons donc ici cette idée (même si elle n'a pas aucune explication rigoureuse)  et nous avons  les estimations suivantes  par rapport à la loi de Kolmogorov.
\begin{Proposition}\label{Prop:ell0-1} Soit $\ell_0>0$ une échelle d'injection d'énergie fixée par la force $\fe$ dans (\ref{loc_fe}). On fixe la longueur caractéristique  $L=\ell_0$. 
Soit   $\vu \in L^{\infty}_{t}L^{2}_{x}\cap (L^{2}_{t})_{loc}\dot{H}^{1}_{x}$ une solution faible des équations de Navier-Stokes amorties (\ref{NS-amortie-alpha-fix}) et soient les moyennes en temps $e>0$ et $u>0$ définies à partir de la solution $\vu$ par les expressions (\ref{e}) et (\ref{u_0}). Soient enfin la vitesse caractéristique $ \ds{U=\frac{u}{\ell^{\frac{3}{2}}_{0}}}$
et le taux de dissipation d'énergie $\ds{\varepsilon=\frac{e}{\ell^{3}_{0}}}$. Alors on a l'estimation: $$\ds{\varepsilon \leq c \frac{U^3}{\ell_0}},$$
où la constante $c>0$ ne dépend  d'aucun paramètre physique ci-dessus.
\end{Proposition}  
Nous observons ainsi que le choix $L=\ell_0$ nous donne une estimation du taux de dissipation $\varepsilon$ qui est bien en accord avec la loi de dissipation de Kolmogorov, néanmoins insistons sur le fait que cette estimation n'est pas tout à fait rigoureuse 
dans le sens que nous ne disposons  d'aucun argument supplémentaire (ni physique ni mathématique) pour justifier le choix de l'échelle $\ell_0$ pour définir les  quantités moyennes $U$ et $\varepsilon$ ci-dessus.  \\
\\
\pv Par l'estimation (\ref{estim-theo-motiv}) on commence par écrire   $\ds{ e \leq c \frac{u_{\ell_0} u^2}{\ell_0}\left(\frac{\Vert \fe \Vert_{L^{\infty}}}{\Vert \fe \Vert_{L^2}}\right)}$,  et nous allons étudier en plus le terme $\frac{\Vert \fe \Vert_{L^{\infty}}}{\Vert \fe \Vert_{L^2}}$. En effet, la force $\fe$ étant localisée aux fréquences de l'ordre de $\frac{1}{\ell_0}$ (voir toujours la formule (\ref{loc_fe})) alors par les inégalités de Bernstein nous avons l'estimation $ \frac{\Vert \fe \Vert_{L^{\infty}}}{\Vert \fe \Vert_{L^2}}   \lesssim  \frac{1}{\ell^{\frac{3}{2}}_{0}}$; et en remplaçant cette estimation dans l'estimation précédente nous avons $\ds{ e \leq c \frac{u_{\ell_0} u^2}{\ell_0} \frac{1}{\ell^{\frac{3}{2}}_{0}}}$. De plus,  par l'estimation (\ref{estim:u0-u}) nous avons  $u_{\ell_0}\leq u$ et dans l'estimation précédente nous pouvons écrire $\ds{e\leq \frac{u^3}{\ell_0} \frac{1}{\ell^{\frac{3}{2}}_{0}}}$. On divise chaque côté de cette estimation par $\ell^{3}_{0}$ et comme l'on a défini $\varepsilon=\frac{e}{\ell^{3}_{0}}$ et $U=\frac{u}{\ell^{\frac{3}{2}}_{0}}$ on obtient alors $\varepsilon \leq c \frac{U^3}{\ell_0}$. \finpv 
\\
Dans ce cas lorsque $L=\ell_0$ nous pouvons aussi en déduire l'estimation du taux de dissipation $\varepsilon$ faite dans \cite{Const} et qui a été exposée dans la Section \ref{sec:cadre_non_per}: rappelons rapidement que dans $\cite{Const}$ on considère les équation de Navier-Stokes classiques (sans terme d'amortissement) et l'on définit une longueur $L_c$ par le biais de la force $\fe$ comme:   
\begin{equation}\label{Lc}
L_c=\frac{F}{\Vert \vec{\nabla}\otimes \fe \Vert_{L^{\infty}}},
\end{equation}
avec $F=\frac{\Vert \fe \Vert_{L^2}}{\ell^{\frac{3}{2}}_{0}}$ (voir la formule (\ref{L_constantin}) pour tous les détails); et avec cette longueur, toujours dans \cite{Const},  on obtient l'estimation suivante $\ds{\varepsilon \lesssim \frac{U^3}{L_c}}$. \\
\\
Néanmoins, dans la Section \ref{Sec:prob_non_periodique} nous avons aussi expliqué que cette estimation présente quelques lacunes techniques et l'une de ces lacunes étant que la vitesse caractéristique $U=\frac{u}{\ell^{\frac{3}{2}}_{0}}$ est potentiellement mal posée dans le cadre des équations de Navier-Stokes classiques (voir l'estimation (\ref{remarque_U})  page \pageref{remarque_U} pour tous les détails à ce sujet). 
\begin{Proposition}\label{Prop:const-revisite} Soit   $\vu \in L^{\infty}_{t}L^{2}_{x}\cap (L^{2}_{t})_{loc}\dot{H}^{1}_{x}$ une solution faible des équations de Navier-Stokes amorties (\ref{NS-amortie-alpha-fix}) et soient soient les moyennes en temps $e>0$ et $u>0$ définies à partir de la solution $\vu$ par les expressions (\ref{e}) et (\ref{u_0}). Soit $\ell_0>0$ une échelle d'injection d'énergie fixée par la force $\fe$ dans (\ref{loc_fe}), soit la longueur $L=\ell_0$ et soient la vitesse caractéristique $ \ds{U=\frac{u}{\ell^{\frac{3}{2}}_{0}}}$ et le taux de dissipation d'énergie $\ds{\varepsilon=\frac{e}{\ell^{3}_{0}}}$. Soit enfin longueur  $L_c$ est définie dans  (\ref{Lc}). Alors on a l'estimation $$ \varepsilon \leq c \frac{U^3}{L_c}.$$
 \end{Proposition} 
\pv  Comme nous avons  $u_{\ell_0}\leq u$ (voir toujours  l'estimation (\ref{estim:u0-u})) alors par  l'estimation  (\ref{estim-theo-motiv})  nous avons $e\leq \frac{u^3}{\ell_0} \left( \frac{\Vert \fe \Vert_{L^{\infty}}}{\Vert \fe \Vert_{L^2}} \right)$ d'où nous écrivons 
$\ds{  e \leq c \frac{u^2}{\ell^{\frac{3}{2}}_{0}} \frac{\ell^{\frac{3}{2}}_{0}}{\Vert f \Vert_{L^2}}\frac{\Vert \fe \Vert_{L^{\infty}}}{\ell_0}}$, et nous allons maintenant vérifier  l'encadrement   $\ds{ \frac{\ell^{\frac{3}{2}}_{0}}{\Vert f \Vert_{L^2}}\frac{\Vert \fe \Vert_{L^{\infty}}}{\ell_0} \approx \frac{1}{L_c}}$.  Par la définition de la longueur $L_c$ ci-dessus nous savons que $\frac{1}{L_c}=\frac{\Vert \vec{\nabla}\otimes \fe \Vert_{L^{\infty}}}{F}$, et comme  $F=\frac{\Vert \fe \Vert_{L^{2}}}{\ell^{\frac{3}{2}}_{0}}$ alors nous écrivons 
$\frac{1}{L_c}=\Vert \vec{\nabla}\otimes \fe \Vert_{L^{\infty}} \frac{\ell^{\frac{3}{2}}}{\Vert \fe \Vert_{L^2}} $. De plus, par la localisation fréquentielle de la force $\fe$ (voir toujours la formule \ref{loc_fe}) et par les inégalités de Bernstein nous avons  
 de plus  $\Vert \fe \Vert_{L^{\infty}}\approx \frac{\Vert \fe \Vert_{L^{\infty}}}{\ell_0}$; et ainsi l'on a l'encadrement ci-dessus. Nous écrivons donc $ \ds{ e \leq c \frac{u^3}{\ell^{\frac{3}{2}}_{0}} \frac{1}{L_c}}$ puis nous divisons chaque terme par $\ell^{3}_{0}$ pour écrire $\varepsilon \leq c \frac{U^3}{L_c}$. \finpv\\
Soulignons maintenant le fait que cette estimation est encore moins rigoureuse que l'estimation obtenue dans la Proposition \ref{Prop:ell0-1} car ici on ne comprend pas tout à fait   la signification de la longueur $L_c$. \\ 
\\
En conclusion, nous observons que si l'on considère la longueur caractéristique $L=\ell_0$ alors l'estimation du taux de dissipation $\varepsilon$ donnée par la Proposition \ref{Prop:ell0-1} est préférable à l'estimation donnée dans la Proposition \ref{Prop:const-revisite} car cette première estimation est plus  en accord à ce qu'on s'attend selon la théorie K41.
\subsubsection{B) Le cas $L>\ell_0$}
Nous considérons ici un cas plus général où la longueur caractéristique $L$ n'est pas forcément égale à l'échelle d'injection d'énergie $\ell_0>0$. Dans ce cas nous avons une estimation du taux de dissipation $\varepsilon$ suivante:

\begin{Proposition}\label{Prop:estim-varep-L} Soit $\ell_0>0$ une échelle d'injection d'énergie fixée par la force $\fe$ dans (\ref{loc_fe}). On fixe la longueur caractéristique $L>\ell_0$. Soit   $\vu \in L^{\infty}_{t}L^{2}_{x}\cap (L^{2}_{t})_{loc}\dot{H}^{1}_{x}$ une solution faible des équations de Navier-Stokes amorties (\ref{NS-amortie-alpha-fix}) et soient les moyennes en temps $e>0$ et $u>0$ définies à partir de la solution $\vu$ par les expressions (\ref{e}) et (\ref{u_0}). Soient  la vitesse caractéristique $ \ds{U=\frac{u}{L^{\frac{3}{2}}}}$, le taux de dissipation d'énergie $\ds{\varepsilon=\frac{e}{L^{3}}}$ et soit la force moyenne $F=\frac{\Vert \fe \Vert_{L^{2}}}{L^{\frac{3}{2}}}$. Soit enfin $c>0$ la constante donnée dans l'estimation (\ref{estim-theo-motiv}) qui ne dépend  d'aucun paramètre physique ci-dessus. Alors on a l'estimation $$\ds{\varepsilon \leq c \frac{U^3}{\ell_0} \left(\frac{\Vert \fe \Vert_{L^{\infty}}}{F}\right)}.$$    
	\end{Proposition} 
Avant de donner une preuve de cette estimation il convient tout d'abord d'expliquer  ce résultat et nous allons maintenant observer que cette estimation donne, dans un certain sens, une généralisation des estimations obtenues dans les Propositions \ref{Prop:ell0-1} et \ref{Prop:const-revisite} dans le cas  $L=\ell_0$.\\
\\ 
En effet, observons tout d'abord que  si nous supposons  que la force $\fe$ vérifie en plus la propriété 
\begin{equation}\label{hyp-force}
\Vert \fe \Vert_{L^{\infty}}\approx F,
\end{equation}
alors par l'estimation du taux de dissipation d'énergie $\varepsilon$ ci-dessus nous pouvons écrire 
\begin{equation}\label{estim-ger}
\varepsilon \lesssim \frac{U^3}{\ell_0},
\end{equation} ce qui nous donne une  estimation  de $\varepsilon$ analogue à celle obtenue dans la Proposition \ref{Prop:ell0-1} et qui est en avec  la loi de Kolmogorov. \\
\\
Remarquons maintenant que l'hypothèse supplémentaire sur la force $\fe$ donnée dans (\ref{hyp-force})  peut être vérifié  dans le cadre de certains forces particulières. En effet, la Définition \ref{Definition_force_ext}  du chapitre suivant nous donnons un exemple concret de force $\fe$ qui vérifie cette propriété (voir  la Remarque \ref{Remarque_normes_fe} page \pageref{Remarque_normes_fe} pour plus de détails à ce sujet). \\
\\ 
D'autre part, quant à l'estimation du taux de dissipation d'énergie $\varepsilon$ obtenue dans la Proposition \ref{Prop:const-revisite}, nous allons maintenant observer que dans le cas lorsque $L>\ell_0$ nous pouvons  obtenir une estimation similaire. En effet, il suffit d'observer le fait que la force $\fe$ étant localisée aux fréquences de l'ordre de $\frac{1}{\ell_0}$ alors par les inégalités de Bernstein nous pouvons écrire  $\Vert \vec{\nabla}\otimes \fe \Vert_{L^{\infty}} \approx \frac{\Vert \fe \Vert_{L^{\infty}}}{\ell_0}$ et donc, par l'estimation de $\varepsilon$  obtenue dans la proposition ci-dessus nous avons $\varepsilon\lesssim U^3\left(\frac{\Vert \fe \Vert_{L^{\infty}}}{F}\right)$. Ensuite, rappelons que dans \cite{Const} on considère la longueur $L_c$ donnée dans (\ref{Lc}): $L_c=\frac{F}{\Vert \vec{\nabla}\otimes \fe \Vert_{L^{\infty}}}$; et nous obtenons ainsi l'estimation $$ \varepsilon \lesssim \frac{U^3}{L_c}.$$
Néanmoins, comme l'on a déjà expliqué dans le cas $L=\ell_0$, cette estimation est moins intéressante que l'estimation (\ref{estim-ger}) car l'on ne sait pas donner une interprétation rigoureuse à cette longueur $L_c$ mais elle  apparaît de façon relativement naturelle dans les calculs faits dans \cite{Const}.  \\
\\
\textbf{Preuve de la Proposition \ref{Prop:estim-varep-L}.} Toujours par l'estimation (\ref{estim-theo-motiv}) on commence par écrire $\ds{ e \leq c \frac{u_{\ell_0} u^2}{\ell_0}\left(\frac{\Vert \fe \Vert_{L^{\infty}}}{\Vert \fe \Vert_{L^2}}\right)}$, et comme nous avons  $u_{\ell_0} \leq u$ (voir l'estimation (\ref{estim:u0-u})) alors nous écrivons $\ds{ e\leq c \frac{u^3}{\ell_0}\left(\frac{\Vert \fe \Vert_{L^{\infty}}}{\Vert \fe \Vert_{L^2}}\right)= \frac{u^3}{\ell_0 L^{\frac{3}{2}}}\left(\frac{L^{\frac{3}{2}} \Vert \fe \Vert_{L^{\infty}}}{\Vert \fe \Vert_{L^2}}\right)}$, et comme l'on a défini la force moyenne $F=\frac{\Vert \fe \Vert_{L^2}}{L^{\frac{3}{2}}}$ nous avons $\ds{ e\leq \frac{u^3}{L^{\frac{3}{2}}}\left(\frac{\Vert \fe \Vert_{L^{\infty}}}{F}\right)}$. Finalement, nous divisons chaque terme de cette estimations par $L^3$ pour écrire $\varepsilon \leq c \frac{U^3}{\ell_0} \left(\frac{\Vert \fe \Vert_{L^{\infty}}}{F}\right)$. \finpv
\section[Lemme technique]{Lemme technique: preuve du Lemme \ref{lemme-tech1} page \pageref{lemme-tech1}}  
 Pour $\vu\in L^{\infty}_{t}(L^{2}_{x})\cap (L^{2}_{t})_{loc}(\dot{H}^{1}_{x})$ nous avons $\partial_t \vu \in (L^{2}_{t})_{loc}(H^{-\frac{3}{2}}_{x}) $,  $ \P((\vu\cdot \vec{\nabla}) \vu)\in (L^{2}_{t})_{loc}(H^{-\frac{3}{2}}_{x})$,  $\Delta \vu \in (L^{2}_{t})_{loc}(H^{-1}_{x})$ et de plus, la force $\fe$  étant  localisée aux fréquences $\frac{\rho_1}{\ell_0}\leq \vert \xi \vert \leq \frac{\rho_2}{\ell_0}$ alors $\fe$ appartient à tous les espaces de Sobolev $H^{s}(\Rt)$ ($s\in \mathbb{R}$) et donc, en multipliant les équations de Navier-Stokes amorties  par $\fe$ et en intégrant en variables d'espace  nous pouvons écrire 
\begin{eqnarray*}
	\int_{\Rt} \partial_t\vu(t,x)\cdot \fe(x)dx&= &\int_{\Rt}\nu\Delta \vu(t,x)\cdot\fe(x)dx-\int_{\Rt}(\P(\vu\cdot\vec{\nabla}\vu(t,x)))\cdot \fe(x)dx+\Vert \fe \Vert^{2}_{L^2}\\
	& & -\frac{\nu}{L^2} \int_{\Rt} \vu(t,x)\cdot \fe(x)dx.
\end{eqnarray*} Dans cette identité nous cherchons à faire apparaître les termes $\Vert \fe\Vert_{L^2}$, $u$, $\ell_0$ et $\mathcal{R}e$  et pour cela on commence par écrire 
\begin{eqnarray}\label{eq_aux_main_theorem}\nonumber
\Vert \fe \Vert^{2}_{L^2}&=&\int_{\Rt} \partial_t\vu(t,x)\cdot \fe(x)dx -\int_{\Rt}\nu\Delta \vu(t,x)\cdot\fe(x)dx+\int_{\Rt}(\P(\vu\cdot\vec{\nabla}\vu(t,x)))\cdot \fe(x)dx\\
& &+\frac{\nu}{L^2} \int_{\Rt} \vu(t,x)\cdot \fe(x)dx,
\end{eqnarray} 
et nous avons les remarques suivantes:  pour le premier terme à droite  ci-dessus nous pouvons écrire
$$ \int_{\Rt} \partial_t\vu(t,x)\cdot \fe(x)dx= \partial_t\int_{\Rt}\vu(t,x)\cdot \fe(x)dx,$$ car $\fe$ est stationnaire. 
Pour le deuxième terme de (\ref{eq_aux_main_theorem}), par une intégration par parties et en appliquant l'inégalité de Cauchy-Schwarz nous obtenons       
$$-\int_{\Rt}\nu\Delta \vu(t,x)\cdot\fe(x)dx =-\nu \int_{\Rt}\vu(t,x)\cdot\Delta\fe(x)dx\leq \nu  \Vert\vu(t,\cdot)\Vert_{L^2}\Vert \Delta \fe\Vert_{L^2}.$$ 
Pour le troisième terme de (\ref{eq_aux_main_theorem}), $\fe$  étant une fonction à divergence nulle et en utilisant les propriétés du projecteur de Leray, par une intégration par parties et  par l'inégalité de H\"older nous avons
\begin{eqnarray*}
	\int_{\Rt}(\P(\vu\cdot\vec{\nabla}\vu(t,x)))\cdot \fe(x)dx&=& \int_{\Rt}(\vu\cdot\vec{\nabla}\vu(t,x))\cdot \fe(x)dx=-\sum_{i,j=1}^{3}\int_{\Rt}u_{i}(t,x)u_{j}(t,x)\partial_{j}f_i(x)dx\\
	& \leq & \Vert \vu(t,\cdot)\Vert^{2}_{L^2}\Vert \vec{\nabla} \otimes\fe\Vert_{L^{\infty}}.
\end{eqnarray*}
Finalement, pour le quatrième terme de (\ref{eq_aux_main_theorem}), par l'inégalité de Cauchy-Schwarz et de plus  étant donné que l'on a $L\geq \ell_0$  alors  nous avons  $$ \frac{\nu}{\ell^{2}_{0}} \int_{\Rt} \vu(t,x)\cdot \fe(x)dx \leq \frac{\nu}{L^2} \Vert \fe \Vert_{L^2}\Vert \vu(t,\cdot)\Vert_{L^2}.$$ 

De cette façon, par les remarques ci-dessus, dans l'identité (\ref{eq_aux_main_theorem})  nous obtenons 
\begin{eqnarray*}
	\Vert \fe \Vert^{2}_{L^2} &\leq&  \partial_t \int_{\Rt} \vu(t,x)\cdot \fe(x)dx+ \Vert \vu(t,\cdot)\Vert^{2}_{L^2}\Vert \vec{\nabla} \otimes\fe\Vert_{L^{\infty}}+\nu\Vert\vu(t,\cdot)\Vert_{L^2}\Vert \Delta \fe\Vert_{L^2}\\
	& & +\frac{\nu}{\ell^{2}_{0}} \Vert \fe \Vert_{L^2}\Vert \vu(t,\cdot)\Vert_{L^2},
\end{eqnarray*}
et maintenant, pour  $T>0$ nous prenons la moyenne en temps $\ds{\frac{1}{T}\int_{0}^{T}(\cdot)dt}$ et comme $\fe$ est stationnaire nous avons 
\begin{eqnarray*}
	\Vert \fe \Vert^{2}_{L^2} &  \leq & \frac{1}{T}\left( \int_{\Rt}\vu(T,x)\cdot\fe(x)-\vu_0(x)\cdot\fe(x)dx\right)+ \left(\frac{1}{T}\int_{0}^{T}\Vert \vu (t,\cdot) \Vert^{2}_{L^2}dt\right) \Vert \vec{\nabla} \otimes \fe \Vert_{L^{\infty}}\\
	& & +\nu \left(\frac{1}{T}\int_{0}^{T} \Vert \vu(t,\cdot) \Vert_{L^2}dt\right) \Vert \Delta \fe \Vert_{L^2}+ \frac{\nu}{\ell^{2}_{0}} \nu \left(\frac{1}{T}\int_{0}^{T} \Vert \vu(t,\cdot) \Vert_{L^2}dt\right) \Vert  \fe \Vert_{L^2}
\end{eqnarray*}
ensuite, nous prenons la limite $\ds{\limsup_{T\longrightarrow+\infty}}$ et nous obtenons
\begin{eqnarray}\label{eq_2_aux_main_theorem} \nonumber
\Vert \fe \Vert^{2}_{L^2} & \leq & \limsup_{T\longrightarrow +\infty}\frac{1}{T}\left( \int_{\Rt}\vu(T,x)\cdot\fe(x)-\vu_0(x)\cdot\fe(x)dx\right)+\left(\limsup_{T\longrightarrow +\infty} \frac{1}{T}\int_{0}^{T}\Vert \vu (t,\cdot) \Vert^{2}_{L^2}dt\right) \times \\
& & \times \Vert \vec{\nabla} \otimes \fe \Vert_{L^{\infty}}+\nu \left( \limsup_{T\longrightarrow +\infty}\frac{1}{T}\int_{0}^{T} \Vert \vu(t,\cdot) \Vert_{L^2}dt\right) \Vert \Delta \fe \Vert_{L^2}\\
& & +\frac{\nu}{\ell^{2}_{0}} \left( \limsup_{T\longrightarrow +\infty}\frac{1}{T}\int_{0}^{T} \Vert \vu(t,\cdot) \Vert_{L^2}dt\right) \Vert  \fe \Vert_{L^2}.
\end{eqnarray}
Dans cette inégalité, pour le premier terme à droite, par le Théorème \ref{Theo:controle-temps} nous savons que la vitesse  $\vu$ vérifi l'estimation
$$ \Vert \vu(T,\cdot)\Vert^{2}_{L^2}\leq e^{-\beta T}\Vert \vu_0 \Vert^{2}_{L^2}+\frac{\Vert \fe \Vert^{2}_{\dot{H}^{-1}}}{\nu \beta} \left( 1-e^{-\beta T}\right),$$ et alors en utilisant l'inégalité de  Cauchy-Schwarz nous avons
\begin{eqnarray*}\nonumber
	& & \left\vert \frac{1}{T}\left( \int_{\Rt}\left[\vu(T,x)\cdot\fe(x)-\vu_0(x)\cdot\fe(x)\right] dx\right)\right\vert  \leq \frac{1}{T}\left( \Vert \vu(T,\cdot)\Vert_{L^2}+ \Vert \vu_0\Vert_{L^2}\right)\Vert \fe \Vert_{L^2}\\
	&\leq & \frac{1}{T} \left[ e^{-\beta T}\Vert \vu_0 \Vert^{2}_{L^2}+\frac{\Vert \fe \Vert^{2}_{\dot{H}^{-1}}}{\nu \beta} ( 1-e^{-\beta T}) + \Vert \vu_0\Vert_{L^2}\right]\Vert \fe \Vert_{L^2}, 
\end{eqnarray*}d'où  nous avons 
\begin{equation}\label{eq_3_aux_main_theorem}
\limsup_{T\longrightarrow+\infty}\left\vert \frac{1}{T}\left( \int_{\Rt}\left[\vu(T,x)\cdot\fe(x)-\vu_0(x)\cdot\fe(x)\right] dx\right)\right\vert \leq 0. 
\end{equation} D'autre part, pour le troisième et quatrième terme à droite de l'estimation (\ref{eq_2_aux_main_theorem}), en appliquant l'inégalité de  Cauchy-Schwarz en variable de temps nous avons
\begin{equation}\label{eq_4_aux_main_theorem}
\frac{1}{T}\int_{0}^{T}\Vert \vu (t,\cdot) \Vert_{L^2}dt\leq \left(\frac{1}{T}\int_{0}^{T}\Vert \vu (t,\cdot) \Vert^{2}_{L^2}dt \right)^{\frac{1}{2}}.
\end{equation}
De cette façon, en remplaçant les inégalités (\ref{eq_3_aux_main_theorem}) et (\ref{eq_4_aux_main_theorem}) dans (\ref{eq_2_aux_main_theorem}) nous obtenons
\begin{eqnarray*}
	\Vert \fe \Vert^{2}_{L^2}&\leq& \left(\limsup_{T\longrightarrow +\infty} \frac{1}{T}\int_{0}^{T}\Vert \vu(t,\cdot) \Vert^{2}_{L^2}dt\right) \Vert \vec{\nabla} \otimes \fe \Vert_{L^{\infty}}+\nu \left(\limsup_{T\longrightarrow+\infty}\frac{1}{T}\int_{0}^{T}\Vert \vu(t,\cdot)\Vert^{2}_{L^2}dt\right)^{\frac{1}{2}} \Vert \Delta \fe \Vert_{L^2}\\
	& &+ \frac{\nu}{\ell^{2}_{0}} \left( \limsup_{T\longrightarrow +\infty}\frac{1}{T}\int_{0}^{T} \Vert \vu(t,\cdot) \Vert^{2}_{L^2}dt\right)^{\frac{1}{2}} \Vert  \fe \Vert_{L^2},
\end{eqnarray*}	d'où, par la définition de la moyenne en temps long $u$ nous pouvons écrire l'estimation cherchée
\begin{equation*}
\Vert \fe \Vert^{2}_{L^2}\leq u^2 \Vert \vec{\nabla} \otimes \fe \Vert_{L^{\infty}} +\nu\, u  \Vert \Delta \fe \Vert_{L^2} +\frac{\nu}{\ell^{2}_{0}} u \Vert  \fe \Vert_{L^2}. 
\end{equation*} \finpv

	 \chapter{Les solutions stationnaires amorties} \label{Chap2} 
Dans le chapitre précédent nous avons introduit les équations de Navier-Stokes amorties 
\begin{equation}\label{n-s-amortie-motiv}
\partial_t\vu = \nu\Delta \vu-(\vu\cdot\vec{\nabla}) \vu -\vec{\nabla}p +\fe-\alpha \vu, \quad div(\vu)=0,\quad \alpha>0,
\end{equation}  où la force $\fe$ est une fonction stationnaire;  et  nous avons étudié la loi de dissipation  d'énergie de Kolmogorov dans le cadre des ces   équations où le terme d'amortissement nous a permis de donner un sens mathématique rigoureux aux quantités considérées.\\
\\
Le fait que la force $\fe$ ne dépende pas de la variable du temps suggère d'étudier les équations de Navier-Stokes amorties et \emph{stationnaires}: 
 \begin{equation}\label{NS-stationnaire-motiv}
-\nu \Delta \U +(\U\cdot \vec{\nabla})\U +\vec{\nabla}P  =\fe-\alpha \U, \qquad div(\U)=0, \quad \alpha>0,
\end{equation} où le champ de vitesse $\U=\U(x)$ et la pression $P=P(x)$ ne dépendent que de la variable spatiale et dans ce chapitre nous allons étudier un tout autre problème relié à la turbulence  dans le cadre de ces équations stationnaires. Plus précisément nous allons étudier ici la stabilité  et la décroissance en variable d'espace des solutions $\U$.\\
\\
En effet, si le fluide est en régime laminaire on s'attend à ce que la solution $\U$ soit stable au sens suivant: si nous considérons $\vu_0\in L^2(\Rt)$ n'importe quelle donnée initiale à divergence nulle et $\vu(t,x)$ une solution faible du problème de Cauchy des équations (\ref{n-s-amortie-motiv})  (obtenue par le biais du Théorème \ref{Theo:existence-sol-ns-amortie} page \pageref{Theo:existence-sol-ns-amortie}) alors on a: 
\begin{equation}\label{estb-intro}
\lim_{t\longrightarrow +\infty}\Vert \vu(t,\cdot)-\U \Vert_{L^2}=0,
\end{equation} ce qui exprime le fait qu'à partir de n'importe quelle donnée initiale $\vu_0$ l'évolution au cours du temps du champs de vitesse $\vu(t,\cdot)$ converge toujours vers la solution stationnaire $\U$ dans le régime asymptotique du temps long; et nous allons observer que cette propriété de stabilité de la solution $\U$, qui sera étudiée  plus tard dans le Théorème \ref{Theo:stabilite_sol_stationnaire_laminaire} (page \pageref{Theo:stabilite_sol_stationnaire_laminaire}),  est seulement valable dans le cadre d'un fluide en régime laminaire tandis que, dans le cadre plus général d'un fluide en régime turbulent, nous allons montrer une toute autre propriété des solutions stationnaires qui porte sur leur décroissance à l'infini en variable d'espace.  Plus précisément, dans le Théorème \ref{Theo:dec-U-turb} (page \pageref{Theo:dec-U-turb}) nous allons montrer  si la force $\fe$ est une fonction bien localisée en variable d'espace alors toute solution  $\U$  des équations (\ref{NS-stationnaire-motiv}) vérifie une décroissance  
\begin{equation}\label{dec-espace-intro}
\vert \U(x)\vert \lesssim \frac{1}{\vert x \vert^4}.
\end{equation}  Cette décroissance est intéressante car rappelons rapidement que dans le cadre des équations de Navier-Stokes classiques (sans terme d'amortissement) on s'attend à ce que les solutions n'aient pas une meilleure décroissance à l'infini que $\frac{1}{\vert x \vert^4}$ et nous allons montrer le terme d'amortissement $-\alpha \U$ introduit dans les équations (\ref{NS-stationnaire-motiv}) entraîne cette décroissance précise à l'infini. \\
\\
Dans la Section \ref{Sec:motivation} nous allons expliquer plus en détail l'intérêt d'étudier les équations stationnaires (\ref{NS-stationnaire-motiv}) et ensuite dans la Section \ref{Sec:existence} nous montrons un résultat général sur l'existence des solutions $(\U,P)$ des ces équations. \\
\\
Une fois que le système d'équations sera posé, nous allons étudier dans la Section  \ref{sec:regime-lam-turb} sous quelles conditions on se trouve dans un régime laminaire ou turbulent. Nous verrons ainsi  que pour faire une étude qui ne dépende pas explicitement des solutions des équations de Navier-Stokes, il sera préférable d'utiliser les nombres de Grashof au lieu des nombres de Reynolds introduits dans le chapitre précédent; et ces nombres nous permettront de déterminer le régime laminaire et le régime turbulent. En particulier, le régime laminaire sera caractérisé par un contrôle sur les nombres de Grashof et pour obtenir ce régime laminaire, nous avons besoin de choisir correctement une force extérieure et ceci sera fait dans la sous-section  \ref{Sec:force ext: construction et proprietes}. \\
\\
Une fois que nous avons tous les ingrédients de notre étude à disposition (les équations stationnaires, les nombres de Grashof et la force extérieure bien préparée) nous pourrons étudier dans la Section \ref{Sec:proprietes-sol-stat} la propriété de stabilité (\ref{estb-intro}) et la décroissance en variable d'espace (\ref{dec-espace-intro}) des solutions des équations stationnaires (\ref{NS-stationnaire-motiv}).

\section{Introduction}
Comme annoncé nous commençons par expliquer plus en détail notre intérêt pour étudier les équations de Navier-Stokes  amorties et stationnaires (\ref{NS-stationnaire-motiv}).
\subsection[Motivation]{Motivation: le problème en temps long}\label{Sec:motivation}
Le fait que la force $\fe$ soit une fonction stationnaire motive l'étude du  comportement des solutions $\vu(t,\cdot)$ lorsque le temps $t$ tend vers l'infini, ce qui est également appelé le problème en temps long pour les équations de Navier-Stokes; et pour expliquer comment  ce problème en temps long nous amène à l'étude des équations stationnaires (\ref{NS-stationnaire-motiv}) nous avons besoin de considérer pour l'instant les équations de Navier-Stokes avec le même terme d'amortissement mais en deux dimensions: 
\begin{equation}\label{n-s-amortie-motiv-2D}
\partial_t\vv = \nu\Delta \vv-(\vv\cdot\vec{\nabla}) \vv -\vec{\nabla}q +\vec{g}-\alpha \vv, \quad div(\vv)=0,\quad \alpha>0,\quad \text{sur}\quad [0,+\infty[\times \mathbb{R}^2,
\end{equation} où $\vec{g}=\vec{g}(x)\in \mathbb{R}^2$ est  toujours une force stationnaire. \\
\\
Le problème en temps long pour les équations (\ref{n-s-amortie-motiv-2D}) (en dimension 2) a été étudié par A. Ilyin \emph{et. al.} en $2015$ dans l'article \cite{Iliyin2} et nous allons expliquer très rapidement les grandes lignes de cette étude: tout d'abord, dans le Théorème $2.1$ de \cite{Iliyin2}, à partir d'une donnée initiale $\vv_0\in L^2(\mathbb{R}^2)$ à divergence nulle on montre l'existence  d'une \emph{unique} solution  $\vv \in L^{\infty}_{t}L^{2}_{x}\cap (L^{2}_{t})_{loc}\dot{H}^{1}_{x}$ des équations (\ref{n-s-amortie-motiv-2D}) qui vérifie 
\begin{equation}\label{cond-initiale-0}
\vv(0,\cdot)=\vv_0,
\end{equation}  et il s'agit d'étudier le comportement de cette solution $\vv(t,\cdot)$ lorsque $t\longrightarrow +\infty$. Cette étude faite dans \cite{Iliyin2} est plutôt technique mais nous expliquons tout de suite les idées générales.\\
\\
L'étude du  problème en temps long des équations de Navier-Stokes amorties et en deux dimensions 	(\ref{n-s-amortie-motiv-2D})  repose essentiellement sur deux ingrédients: la notion de solution éternelle   donnée dans la Définition  \ref{def-sol-eternelle} ci-dessous et  l'unicité de la solution $\vu$ du problème de Cauchy  pour ces équations. On commence donc par introduire la notion de solution éternelle:  
\begin{Definition}\label{def-sol-eternelle} Une fonction $\vv_e:]-\infty , +\infty[\times \mathbb{R}^2 \longrightarrow \mathbb{R}^2$ est une solution éternelle des équations de Navier-Stokes amorties (\ref{n-s-amortie-motiv-2D}) si $\vv_e  \in L^{\infty}_{t}L^{2}_{x}\cap (L^{2}_{t})_{loc}\dot{H}^{1}_{x}$ et si la fonction $\vv_e$ vérifie ces équations.
	\end{Definition}
Nous allons maintenant voir comment ces solutions éternelles nous permettent d'étudier le comportement en temps long de la solution $\vv$ (voir toujours l'article \cite{Iliyin2} pour tous les détails des calculs). On considère une suite $(t_n)_{n\in \mathbb{N}}$  telle que $t_0=0$ et telle que $t_n\longrightarrow +\infty$ lorsque $n\longrightarrow +\infty$;  et pour tout $n\in \mathbb{N}$ on considère le problème de Cauchy des équations (\ref{n-s-amortie-motiv-2D}) avec la condition initiale 
\begin{equation}\label{cond-initiale-n}
\vv_n(-t_n,\cdot)=\vv_0,
\end{equation} où nous pouvons observer que l'on prend ici la même donnée  initiale $\vv_0$ mais maintenant au temps $-t_n$ au lieu du temps $t_0=0$. Ainsi, toujours par le Théorème $2.1$ dans \cite{Iliyin2} nous savons qu'il existe une fonction $\vv_n  \in L^{\infty}_{t}L^{2}_{x}\cap (L^{2}_{t})\dot{H}^{1}_{x}$ définie sur $[-t_n,+\infty[\times \mathbb{R}^2$ qui est l'unique solution du problème de Cauchy avec la condition (\ref{cond-initiale-n}). Ensuite,  dans le Théorème $2.6$ dans \cite{Iliyin2} on montre que la suite de fonctions  $(\vv_n)_{n\in \mathbb{N}}$ converge dans la topologie forte de l'espace $(L^{2}_{t}L^{2}_{x})_{loc}$ vers une solution éternelle $\vv_e$ donnée dans la Définition \ref{def-sol-eternelle}. De plus (toujours dans le  Théorème $2.6$ de l'article \cite{Iliyin2})  on montre que l'on a aussi la convergence 
\begin{equation}\label{conv-1}
\lim_{n\longrightarrow +\infty}  \Vert \vv_n(0,\cdot)-\vv_e(0,\cdot)\Vert_{L^2}=0, 
\end{equation} et cette convergence et \emph{l'unicité} des solutions du problème de Cauchy des équations (\ref{n-s-amortie-motiv-2D}) vont nous permettre d'étudier le  comportement en temps long de la solution $\vv$.  \\
\\
En effet,  dans les expressions (\ref{cond-initiale-0}) et (\ref{cond-initiale-n}) nous pouvons observer que les solutions $\vv$ et $\vv_n$ sont construites à partir de la même donnée $\vv_0$ et ainsi, par l'unicité de la solution, nous observons  que la solution $\vv_n$ est en réalité un décalage de la solution $\vv$  au temps initial $-t_n$ et  nous pouvons ainsi écrire l'identité 
\begin{equation}\label{identite-2D}
\vv_n(0,\cdot)=\vv(t_n,\cdot),
\end{equation} pour tout $n\in \mathbb{N}$. Si nous remplaçons maintenant cette identité dans l'expression (\ref{conv-1}) nous pouvons alors écrire 
\begin{equation}\label{prob-temps-long-2D}
	\lim_{n\longrightarrow +\infty}  \Vert \vv(t_n,\cdot)-\vv_e(0,\cdot)\Vert_{L^2}=0,
	\end{equation} pour observer que la solution $\vv$ converge (via la suite $(t_n)_{n\in \mathbb{N}}$)  vers la solution éternelle $\vv_e$ des équations (\ref{n-s-amortie-motiv-2D}) lorsque le temps $t_n$ tend vers l'infini et ceci nous permet de comprendre le comportement en temps long de cette solution. \\
	\\
Revenons à présent à notre cas d'étude donné par les équations de Navier-Stokes amorties sur l'espace $\Rt$ (\ref{n-s-amortie-motiv}). Rappelons  que nous voulons étudier le comportement en temps long des solutions $\vu$ et pour cela il serait naturel de suivre les lignes exposées  ci-dessus (dans le cadre des équations en deux dimensions )  mais nous allons voir que l'on a ici une contrainte technique qui porte sur l'unicité des solutions et que nous expliquons tout suite. \\
\\
Soit donc $\vu_0 \in L^2(\Rt)$ une donnée initiale à divergence nulle et soit une suite  $(t_n)_{n\in \mathbb{N}}$ telle que $t_0=0$ et  $t_n\longrightarrow +\infty$ lorsque $n\longrightarrow +\infty$. En suivant les idées ci-dessus nous considérons le problème de Cauchy pour les équations (\ref{n-s-amortie-motiv}) avec une condition initiale 
	\begin{equation}
	\vu_n(-t_n,\cdot)=\vu_0,
	\end{equation} pour tout $n\in \mathbb{N}$; et alors par le Théorème \ref{Theo:existence-sol-ns-amortie} (page \pageref{Theo:existence-sol-ns-amortie}) nous savons qu'il existe une solution $\vu_n:[-t_n,+\infty[\times \Rt \longrightarrow \Rt$. De plus, en suivant exactement les mêmes lignes de la preuve du Théorème $2.6$ dans \cite{Iliyin2} (toutes les estimations s'adaptent sans aucun problème lorsqu'on considère les équations en trois dimensions) nous avons que cette suite des solutions $(\vu_n)_{n\in \mathbb{N}}$ converge vers une solution éternelle (au sens de la Définition \ref{def-sol-eternelle}) des équations (\ref{n-s-amortie-motiv}): $$\vu_e:]-\infty,+\infty[\times \Rt \longrightarrow \Rt,$$ et de plus (toujours par le Théorème $2.6$ dans \cite{Iliyin2}) nous avons  
	\begin{equation}\label{conv-2}
	\lim_{n\longrightarrow +\infty}  \Vert \vu_n(0,\cdot)-\vu_e(0,\cdot)\Vert_{L^2}=0.
	\end{equation} Mais, dans ce cas où l'on considère les équations de Navier-Stokes  posées sur $\Rt$ l'unicité des solutions du problème de Cauchy  est encore une question ouverte, que ce soit pour les équations de Navier-Stokes amorties ou pour les équations de Navier-Stokes classiques (sans terme d'amortissement); et cette fois-ci on ne peut pas écrire l'identité $\ds{\vu_n(0,\cdot)=\vu(t_n,\cdot)}$ pour ainsi obtenir un résultat similaire à celui donné dans l'expression (\ref{prob-temps-long-2D}) qui nous a permis de comprendre le comportement en temps long des équations amorties en deux dimensions. \\
	\\
Ainsi, l'étude du comportement en temps long des solutions des équations de Navier-Stokes amorties et en  dimension deux (\ref{n-s-amortie-motiv-2D}) fait dans l'article \cite{Iliyin2} ne peut pas être appliqué en toute généralité  à notre cadre des équations amorties en  dimension trois (\ref{n-s-amortie-motiv}) et l'essentiel du problème repose sur le manque d'information sur l'unicité des solutions du problème de Cauchy.  \\
\\
Néanmoins cette étude nous suggère de considérer les solutions éternelles $\vu_e$ des équations (\ref{n-s-amortie-motiv}) et un cas particulier des solutions éternelles sont les solutions stationnaires $\U$. En effet, comme la fonction $\U$ 
est constante en temps alors  cette fonction est bien évidemment définie pour tout temps $t\in \mathbb{R}$ et de plus, comme la fonction $\U$ vérifie les équations (\ref{NS-stationnaire-motiv}) et comme l'on a $\partial_t \U =0$, alors cette fonction vérifie aussi les équations (\ref{n-s-amortie-motiv}).\\
\\
Nous observons alors que l'étude du problème en temps long des équations  (\ref{n-s-amortie-motiv}) nous amène finalement à l'étude des équations stationnaires (\ref{NS-stationnaire-motiv}) et dans la Section \ref{Sec:Stabilite-laminaire} nous allons voir que  les solutions $\U$ de ces équations stationnaires nous permettent de comprendre le comportement en temps des solutions $\vu$ des équations non stationnaires  (\ref{n-s-amortie-motiv})  dans un cadre particulier lorsque le fluide est en régime laminaire. Mais, avant d'entrer dans le vif du sujet, nous allons tout d'abord étudier un peu plus les équations stationnaires (\ref{NS-stationnaire-motiv}) et nous allons montrer un résultat général sur l'existence des solutions de ces équations. 

\subsection{Existence}\label{Sec:existence}  
Comme annoncé  nous  considérons ici  le système de Navier-Stokes amorti et stationnaire (\ref{NS-stationnaire-motiv})  et nous allons construire  des solutions dans le cadre de l'espace de Sobolev $H^1(\Rt)$, mais,  avant de construire ces solutions il convient d'expliquer un peu plus en détails pourquoi nous cherchons à construire des solutions dans  cet espace fonctionnel. Le choix de l'espace $H^1(\Rt)$ est une conséquence des estimations  \emph{a priori} suivantes: si nous supposons que la vitesse $\U$ et la force $\fe$  sont suffisamment régulières et intégrables  alors, en multipliant les équations (\ref{NS-stationnaire-motiv}) par $\U$ et puis en intégrant en variable d'espace nous avons (formellement)
$$ -\nu\int_{\Rt}\Delta \U \cdot \U dx+ \int_{\Rt} \left[ (\U\cdot \vec{\nabla}) \U \right]\cdot \U dx +\int_{\Rt}\vec{\nabla}P\cdot \U dx=\int_{\Rt}\fe \cdot \U dx -\alpha \int_{\Rt}\U \cdot \U dx,$$ d'où, comme $div(\U)=0$, par une intégration par parties nous avons (toujours formellement) les identités 
$$ \int_{\Rt} \left[ (\U\cdot \vec{\nabla}) \U \right]\cdot \U dx=0 \quad \text{et}\quad  \int_{\Rt}\vec{\nabla}P\cdot\U dx=0,$$ et alors nous obtenons l'identité 
$$ \nu\int_{\Rt}\vert \vec{\nabla}\otimes \U \vert^2dx + \alpha \int_{\Rt}\vert\U \vert^2 dx =\int_{\Rt}\fe \cdot \U dx,$$ d'où nous pouvons écrire 
$$ \min(\nu,\alpha)\Vert \U \Vert^{2}_{H^1}\leq  \nu\int_{\Rt}\vert \vec{\nabla}\otimes U \vert^2dx + \alpha \int_{\Rt}\vert\U \vert^2 dx =\int_{\Rt}\fe \cdot \U dx.$$
Dans le dernier terme à droite ci-dessus nous observons que si l'on a $\fe \in H^{-1}(\Rt)$ et $\U\in H^1(\Rt)$ alors   en appliquant l'inégalité de Cauchy-Schwarz dans ce dernier terme   nous obtenons
$$ \min(\nu,\alpha)\Vert \U \Vert^{2}_{H^1}\leq \Vert \fe \Vert_{H^{-1}}\Vert \U \Vert_{H^{1}}, $$ d'où nous écrivons (au moins formellement) l'estimation \emph{a priori} suivante:
\begin{equation}\label{estimation-a-priori}
\Vert \U \Vert_{H^1} \leq \frac{1}{\min(\nu,\alpha)}\Vert \fe \Vert_{H^{-1}}. 
\end{equation}
Avec ces estimations formelles, nous pouvons observer que l'espace $H^{-1}(\Rt)$ est un espace naturel pour  la force $\fe$ et avec $\fe \in H^{-1}(\Rt)$  nous sommes censés obtenir des solutions dans l'espace $H^1(\Rt)$. \\
\\
Observons aussi que le terme d'amortissement $-\alpha \U$ entraîne un contrôle sur la norme $L^2$ des  solutions $\U$. Dans la Section \ref{Sec:proprietes-sol-stat} nous allons étudier quelques propriétés intéressantes des solutions $\U$ qui sont entraînées par ce terme $-\alpha \U$. \\
\\
Une fois que nous avons fait ces calculs préliminaires,  nous allons construire des solutions $(\U,P)$ des équations stationnaires (\ref{NS-stationnaire-motiv}) et pour cela nous procéderons de la façon suivante:  comme nous travaillons dans l'espace $\Rt$ tout entier, dans un premier temps  nous appliquons le projecteur de Leray aux équations  (\ref{NS-stationnaire-motiv}) et comme $div(\U)=0$ et $div(\fe)=0$ nous obtenons ainsi les équations
$$ -\nu \Delta \U +\P ((\U\cdot \vec{\nabla}) \U)=\fe -\alpha \U,$$ d'où nous pouvons écrire 
$$ (-\nu \Delta +\alpha I_d )\U + \P ((\U\cdot \vec{\nabla}) \U)=\fe,$$ où $I_d$ est l'opérateur identité.  Dans cette expression nous observons que ces équations s'écrivent formellement comme le problème de point fixe équivalent suivant:
\begin{equation}\label{point-fixe-stattionnaire}
\U = - \frac{1}{-\nu \Delta +\alpha I_d}\left[ \P ((\U\cdot \vec{\nabla}) \U)\right] + \frac{1}{-\nu \Delta +\alpha I_d}\left[\fe\right].
\end{equation}
Ainsi, nous utiliserons un principe de point fixe  pour résoudre ce problème et ensuite, à partir de la solution $\U$ nous pourrions récupérer le terme de pression $P$ qui est relié à cette fonction  par la relation $\ds{P=\frac{1}{-\Delta}div \left( (\U\cdot \vec{\nabla}) \U \right)}$.\\ 
\\
\`A ce stade il est important de souligner qu'il y a plusieurs façons de  résoudre le problème de point fixe (\ref{point-fixe-stattionnaire}).  En effet,  une manière classique de résoudre ce problème est d'utiliser le principe de contraction de Picard, où l'on obtient l'existence et l'unicité de la solution, mais avec cette méthode nous sommes obligés de contrôler la taille de la force $\fe$: $\Vert \fe \Vert_{H^{-1}}$ et alors nous avons l'existence et l'unicité de la solution seulement  pour des forces suffisamment petites. En effet, pour appliquer le principe de contraction de Picard nous devons contrôler le terme $\ds{\left\Vert  \frac{1}{-\nu \Delta +\alpha I_d}\left[\fe\right] \right\Vert_{H^{1}}}$ et pour cela nous devons finalement obtenir un contrôle sur la taille  de la force $\fe$:  l'opérateur $\ds{\frac{1}{-\nu \Delta +\alpha I_d}}$ s'exprime au niveau de Fourier par le symbole $\ds{\frac{1}{\nu \vert \xi \vert^2 +\alpha}}$ et comme nous avons l'estimation $\ds{\frac{1}{\nu \vert \xi \vert^2 +\alpha} \leq \frac{1}{\min(\alpha,\nu)}\frac{1}{\vert \xi \vert^1 +1}}$ nous obtenons alors l'estimation $\ds{\left\Vert  \frac{1}{-\nu \Delta +\alpha I_d}\left[\fe\right] \right\Vert_{H^{1}} \leq \frac{1}{\min(\alpha,\nu)}\Vert \fe \Vert_{H^{-1}}}$. \\
\\
\\ 
Une deuxième façon de résoudre le problème de point fixe (\ref{point-fixe-stattionnaire})  repose sur l'estimation \emph{a priori} (\ref{estimation-a-priori}) où en utilisant  le principe de point fixe de Schaefer (voir le Lemme \ref{Schaefer} ci-dessous)  nous pouvons construire des solutions $\U \in H^1(\Rt)$  pour \emph{n'importe quelle} force $\fe \in H^{-1}(\Rt)$ et donc il  s'agit d'une façon plus générale de construire ces solutions. Néanmoins comme nous allons le voir cette méthode nous ne fournit  aucune information supplémentaire sur l'unicité des solutions. \\
\\
Nous allons préférer ici cette deuxième façon de construire des solutions du problème (\ref{point-fixe-stattionnaire}) (et donc de construire des solutions des équations (\ref{NS-stationnaire-motiv})) car le fait de ne pas contrôler la taille force $\fe$ nous permettra après de caractériser  le régime turbulent du fluide, ce qui sera expliqué plus en détail dans la Section \ref{sec:regime-lam-turb}; et ensuite dans la Section \ref{Sec:localisation-espace} nous étudierons quelques propriétés des solutions $\U$ dans ce régime turbulent. \\
\\
Nous avons de cette façon le résultat suivant sur l'existence des solutions des équations de Navier-Stokes amorties et stationnaires. \\ 
\begin{Theoreme}\label{Theo:solutions_stationnaires_turbulent} Soit $\fe\in H^{-1}(\Rt)$ une force à divergence nulle.  Il existe $(\U,P)\in H^1(\Rt)\times H^{\frac{1}{2}}(\Rt)$ solution faible des équations de Navier-Stokes stationnaires et amorties (\ref{NS-stationnaire-motiv}). De plus, toute solution vérifie l'estimation (\ref{estimation-a-priori}).  \\
\end{Theoreme}
Avant de prouver ce théorème faisons les remarques suivantes. La preuve de ce théorème suit  les grandes lignes de la preuve du Théorème $16.2$ du livre \cite{PGLR1}; et comme l'on a déjà mentionné   cette preuve repose essentiellement sur l'estimation à priori (\ref{estimation-a-priori})  et le théorème de point fixe  de Schaefer  que nous énonçons en toute généralité comme suit  (voir  le Théorème $16.1$ du livre \cite{PGLR1} pour plus de références sur ce résultat). 
\begin{Lemme}\label{Schaefer} Soit $E$ un espace de Banach et $T:E\longrightarrow E$ un opérateur qui satisfait:
	\begin{enumerate}
		\item[1)] $T$ est un opérateur continu, 
		\item[2)] $T$ est un opérateur compact,
		\item[3)] Estimation \emph{a priori}: il existe une constante $C>0$ telle que, pour tout $\lambda \in [0,1]$, si $e\in E$ vérifie l'équation $e=\lambda T(e)$ alors on a $\Vert e \Vert_E\leq C$.
	\end{enumerate} Alors, il existe $e\in E$  une solution du problème de point fixe  $e=T(e)$.
\end{Lemme} 
Ainsi, dans le cadre de ce  lemme  nous définissons l'espace $E$ comme 
\begin{equation}\label{E_turbulent}
E=\lbrace \U \in H^1(\Rt): div(\U)=0\rbrace,
\end{equation} muni de la norme $\Vert \cdot \Vert_{H^1}$ et de plus, nous définissons l'opérateur $T$ par l'expression (\ref{point-fixe-stattionnaire}):
\begin{equation}\label{op_T_turbulent}
T(\U)= \frac{1}{-\nu \Delta +\alpha I_d}\left[  \P \left( \left((\U\cdot \vec{\nabla}) \U\right)\right)\right] - \frac{1}{-\nu \Delta +\alpha I_d}\left[ \fe\right],
\end{equation} et nous voulons maintenant vérifier que cet opérateur satisfait les hypothèses du Lemme \ref{Schaefer} mais nous trouvons ici une contrainte technique. En effet, le point $1)$ ci-dessus sera vérifié par le Lemme  \ref{Lemme_cont_forme_bilineaire_turbulent} ci-après, tandis que, le point $3)$  sera vérifié par l'estimation (\ref{estimation-a-priori}),  néanmoins, nous allons voir que le point $2)$  qui porte sur la compacité de l'opérateur $T$  pose des problèmes techniques.  En effet, si nous prenons $(\U_n)_{n\in\mathbb{N}}$ une suite bornée dans l'espace  $E$ défini dans (\ref{E_turbulent}), nous savons qu'il existe une sous suite $(\U_{n_k})_{k\in\mathbb{N}}$ qui converge faiblement dans cet espace et alors  par la continuité de l'opérateur $T$ nous avons seulement la convergence faible  la suite $(T(\U_{n_k}))_{k\in\mathbb{N}}$  et non pas sa convergence forte et ainsi  la compacité de  l'opérateur T  semble hors de portée. \\
\\
De ce cadre, nous allons contourner ce problème technique de la façon suivante: il s'agit d'approcher l'opérateur $T$ par une famille d'opérateurs compacts $(T_r)_{r>0}$ qui vérifient  les hypothèses du Lemme \ref{Schaefer} et alors nous obtiendrons une famille de solutions approchées $\ds{\U_r=T_r(\U_r)}$ où $U_r\in E$ pour tout $r>0$. Ensuite, par un lemme de Rellich-Lions nous prouvons que cette famille converge vers  une solution  des équations (\ref{NS-stationnaire-motiv}).\\ 
\\
\textbf{Démonstration  du Théorème \ref{Theo:solutions_stationnaires_turbulent}.} %
Soit  $\theta \in \mathcal{C}^{\infty}_{0}(\Rt)$  telle que $0\leq \theta(x)\leq 1$, $\theta(x)=1$ si $\vert x \vert  <1$ et $\theta(x)=0$ si $\vert x \vert >2$ et soit $r>0$. Nous définissons la fonction de troncature $\theta_r(x)$ par $\theta\left( \frac{x}{r}\right)$ et alors, dans le terme  bilinéaire de  (\ref{op_T_turbulent})  nous écrivons  $\ds{([\theta_r \U]\cdot \vec{\nabla})[\theta_r \U]}$,  et  nous définissons ainsi  l'opérateur approché $T_r$  par
$$ T_r(\U) =\frac{1}{-\nu \Delta +\alpha I_d}\left[ \P \left(\left(  ((\theta_r \U)\cdot \vec{\nabla})(\theta_r \U) \right)\right)\right] - \frac{-\nu \Delta}{-\nu \Delta +\alpha I_d}\left[( \fe)\right],$$ et pour mener à bien les estimations dont on a besoin nous récrivons cette opérateur comme

\begin{equation}\label{op_Tr_turbulent} 
T_r(\U) =\frac{-\nu \Delta}{-\nu \Delta +\alpha I_d}\left[\frac{1}{\nu} \P \left(\frac{1}{\Delta} \left(  ((\theta_r \U)\cdot \vec{\nabla})(\theta_r \U) \right)\right)\right] - \frac{-\nu \Delta}{-\nu \Delta +\alpha I_d}\left[\frac{1}{\nu \Delta}( \fe)\right].
\end{equation}
Vérifions maintenant que cet opérateur satisfait les hypothèses du Lemme \ref{Schaefer}.
\begin{enumerate}
	\item[$1)$] Continuité.   Soient  $\U,\V\in E$ et nous avons 
	$$ \left\Vert T_r(\U)-T_r(\V)\right\Vert_{H^1}= \left\Vert \frac{-\nu \Delta}{-\nu \Delta +\alpha I_d}\left[\frac{1}{\nu} \P \left(\frac{1}{\Delta} \left(  ((\theta_r \U)\cdot \vec{\nabla})[\theta_r \U] -(\theta_r \V)\cdot \vec{\nabla})(\theta_r \V)\right)\right)\right] \right\Vert_{H^1},$$ d'où, nous écrivons 
	\begin{equation}\label{identite_vec_1_turbulent}
	((\theta_r \U)\cdot \vec{\nabla})(\theta_r \U) -((\theta_r \V)\cdot \vec{\nabla})(\theta_r \V)=((\theta_r (\U-\V))\cdot \vec{\nabla})(\theta_r \U)+ ((\theta_r \V)\cdot \vec{\nabla})(\theta_r (\U-\V)),
	\end{equation} et alors  nous obtenons  
	\begin{eqnarray}\label{estim-aux} \nonumber \vspace{2mm}
		\left\Vert T_r(\U)-T_r(\V)\right\Vert_{H^1} &\leq & \left\Vert \frac{-\nu \Delta}{-\nu \Delta +\alpha I_d}\left[\frac{1}{\nu} \P \left(\frac{1}{\Delta} \left(  ([\theta_r (\U-\V)]\cdot \vec{\nabla})[\theta_r \U] \right)\right)\right] \right\Vert_{H^1}\\  \vspace{5mm}
		& & + \left\Vert \frac{-\nu \Delta}{-\nu \Delta +\alpha I_d}\left[\frac{1}{\nu} \P \left(\frac{1}{\Delta} \left(  ([\theta_r \V]\cdot \vec{\nabla})[\theta_r (\U-\V)] \right)\right)\right] \right\Vert_{H^1}. \vspace{5mm}
	\end{eqnarray}
	\`A ce stade nous avons besoin de l'estimation suivante:  
	\begin{Lemme}\label{Lemme_cont_forme_bilineaire_turbulent} Pour tout $\U,\V\in H^1(\Rt)$ tels que $div(\U)=0$ et $div(\V)=0$ on a: $$\ds{\left\Vert \frac{-\nu \Delta}{-\nu \Delta +\alpha I_d}\left[\frac{1}{\nu} \P \left(\frac{1}{\Delta} \left( (\U\cdot \vec{\nabla})\V \right)\right)\right] \right\Vert_{H^1}\leq c(\alpha,\nu)\Vert \U \Vert_{H^1} \Vert \U \Vert_{H^1}}.$$ 
	\end{Lemme} 
\pv  On commence par écrire
\begin{eqnarray*}
\left\Vert \frac{-\nu \Delta}{-\nu \Delta +\alpha I_d}\left[\frac{1}{\nu} \P \left(\frac{1}{\Delta} \left( (\U\cdot \vec{\nabla})\V \right)\right)\right] \right\Vert_{H^1}& \leq &c \left\Vert \frac{1}{-\nu \Delta +\alpha I_d}\left[ div (\U\otimes \V)\right] \right\Vert_{H^1} \\
& \leq & c(\alpha,\nu) \Vert div(\U \otimes \V)\Vert_{H^{-1}}\leq c(\alpha,\nu) \Vert \U\otimes \V \Vert_{L^2}\\
&\leq & c(\alpha,\nu) \Vert \U \otimes \V \Vert_{H^{\frac{1}{2}}}, 
\end{eqnarray*} où par les lois de produit (voir le livre \cite{PGLR1}) nous pouvons finalement écrire 
$$ c(\alpha,\nu) \Vert \U \otimes \V \Vert_{H^{\frac{1}{2}}} \leq c(\alpha,\nu)\Vert \U \Vert_{H^1}\Vert \V \Vert_{H^1}.$$ \finpv
Ainsi, nous appliquons cette  estimation à chaque terme à droite de l'inégalité (\ref{estim-aux}) et nous avons 
	\begin{eqnarray*}
		\left\Vert T_r(\U)-T_r(\V)\right\Vert_{H^1} & \leq & c(\alpha,\nu) \left[ \left\Vert \theta_r(\U-\V) \right\Vert_{H^1} \left\Vert \theta_r \U \right\Vert_{H^1}+\left\Vert \theta_r \V \right\Vert_{H^1} \left\Vert \theta_r (\U-\V) \right\Vert_{H^1}\right]\\ \vspace*{5mm}
		& \leq & c(\alpha,\nu) \left[ \left\Vert \theta_r \U \right\Vert_{H^1}+ \left\Vert \theta_r \V \right\Vert_{H^1} \right]\left\Vert \theta_r(\U-\V) \right\Vert_{H^1}.
	\end{eqnarray*} D'autre part, nous étudions un peu plus  le dernier terme à droite où l'on a l'estimation 
	\begin{eqnarray*}
		\left\Vert \theta_r(\U-\V) \right\Vert_{H^1} &=& \left\Vert \theta_r(\U-\V) \right\Vert_{L^2}+\left\Vert \vec{\nabla}\otimes ( \theta_r(\U-\V)) \right\Vert_{L^2}\\
		&\leq & \left[ \Vert \theta_r\Vert_{L^{\infty}}+\Vert \vec{\nabla}\theta_r\Vert_{L^{\infty}}\right] \Vert \U-\V \Vert_{H^1},  
	\end{eqnarray*} et alors  nous pouvons écrire
	$$\left\Vert T_r(\U)-T_r(\V)\right\Vert_{H^1} \leq c(\theta_r,\alpha,\nu)\left[ \left\Vert \theta_r \U \right\Vert_{H^1}+ \left\Vert \theta_r \V \right\Vert_{H^1} \right] \Vert \U-\V \Vert_{H^1},$$ d'où nous concluions la continuité de l'opérateur $T_r$ sur l'espace $E$.
	
	\item[2)] Compacité. Soit donc $(\V_n)_{n\in\mathbb{N}}$ une suite bornée dans $E$. Alors, la suite $(\theta_r \V_n)_{n\in \mathbb{N}}$ est aussi bornée dans $E$ et en plus, comme la fonction $\theta_r$ est définie par $\theta_r(x)=\theta\left( \frac{x}{r}\right)$ et la function de test $\theta$ satisfait $supp\,( \theta) \subset \lbrace x\in \Rt: \vert x \vert <2\rbrace$, alors, nous  avons $supp\, (\theta_r)\subset \lbrace x\in \Rt: \vert x \vert <2r\rbrace$ et donc , pour tout $n\in \mathbb{N}$ nous avons 
	$$ supp\, (\theta_r\V_n)\subset \lbrace x\in \Rt: \vert x \vert <2r\rbrace.$$
	De cette façon, par le lemme de Rellich-Lions il existe une sous-suite  $(\V_{n_k})_{k\in\mathbb{N}}$, qui  pour simplifier l'écriture sera notée  comme $(\V_n)_{n\in \mathbb{N}}$, et  il existe  $\V\in H^1(\Rt)$ telles que la suite $(\V_{n})_{n\in\mathbb{N}}$ converge fortement vers $\V$ dans $L^2(\Rt)$ mais aussi dans $L^p(\Rt)$ pour $2\leq p <6$.\\
	\\
	Ainsi, pour vérifier la compacité de l'opérateur $T_r$ nous allons montrer que la suite $(T_r(\V_{n}))_{n\in \mathbb{N}}$ converge fortement vers $T(\V)$ dans $E$. En effet, par l'identité (\ref{identite_vec_1_turbulent}) nous commençons par  écrire 
	\begin{eqnarray*}
		\Vert T(\V_n)-T(\V)\Vert_{H^1}&\leq & \left\Vert \frac{-\nu \Delta}{-\nu \Delta +\alpha I_d}\left[\frac{1}{\nu} \P \left(\frac{1}{\Delta} \left(  ([\theta_r (\V_n-\V)]\cdot \vec{\nabla})[\theta_r \V_n] \right)\right)\right] \right\Vert_{H^1}\\ \vspace{5mm}
		& & + \left\Vert \frac{-\nu \Delta}{-\nu \Delta +\alpha I_d}\left[\frac{1}{\nu} \P \left(\frac{1}{\Delta} \left(  ([\theta_r \V]\cdot \vec{\nabla})[\theta_r (\V_n-\V)] \right)\right)\right] \right\Vert_{H^1}, 
	\end{eqnarray*}
	et pour traiter les termes à droite ci-dessus nous avons l'estimation suivante:
	\begin{Lemme} On a $\ds{\left\Vert \frac{-\nu \Delta}{-\nu \Delta +\alpha I_d}\left[\frac{1}{\nu} \P \left(\vg\right)\right] \right\Vert_{H^1}\leq c(\alpha,\nu) \Vert \vg \Vert_{\dot{H}^{-1}}},$ pour tout $\vg \in \dot{H}^{-1}(\Rt)$.
	\end{Lemme} Nous prouverons ce lemme à la fin du chapitre et  maintenant, en appliquant ce lemme à chaque terme de l'identité précédente nous avons 
	\begin{equation*}
	\Vert T(\V_n)-T(\V)\Vert_{H^1}\leq  c \underbrace{\Vert ([\theta_r (\V_n-\V)]\cdot \vec{\nabla})[\theta_r \V_n] \Vert_{\dot{H}^{-1}}}_{(a)}+ c\underbrace{\Vert ([\theta_r \V]\cdot \vec{\nabla})[\theta_r (\V_n-\V)] \Vert_{\dot{H}^{-1}} }_{(b)}
	\end{equation*} et nous cherchons à montrer que les termes $(a)$ et $(b)$ ci-dessus convergent vers zéro lorsque $n$ tend à l'infini. 
	Pour traiter les termes $(a)$ et $(b)$  nous utiliserons l'identité suivante: soient  $\A=(A_1,A_2,A_3)$ et $\B=(B_1,B_2,B_3)$ deux fonctions vectorielles, alors on peut écrire  
	\begin{eqnarray}\label{identite_vec_2_turbulent}\nonumber
	((\theta_r \A)\cdot \vec{\nabla})(\theta_r \B)&=& \sum_{j=1}^{3}(\theta_r A_j)\partial_j (\theta_r \B)=\sum_{j=1}^{3}\left( \theta_r A_j(\partial_j \theta_r \B+\theta_{r} \partial_j \B)\right)\\
	&=& (\A \cdot \vec{\nabla}\theta_r)(\theta_r \B)+(\theta^{2}_{r}\A\cdot \vec{\nabla})\B.
	\end{eqnarray} Ainsi, pour étudier le  terme $(a)$, dans l'identité ci-dessus nous prenons $\A=\V_n-\V$  et $\B=\V_n$  et alors nous pouvons écrire 
	$$\Vert ([\theta_r (\V_n-\V)]\cdot \vec{\nabla})[\theta_r \V_n] \Vert_{\dot{H}^{-1}}\leq  \underbrace{\left\Vert  ((\V_n-\V)\cdot \vec{\nabla}\theta_r)[\theta_r \V_n]\right\Vert_{\dot{H}^{-1}}}_{(a.1)}+\underbrace{\left\Vert  (\theta^{2}_{r}(\V_n-\V)\cdot\vec{\nabla})\V_n  \right\Vert_{\dot{H}^{-1}}}_{(a.2)}.$$
	Nous avons encore besoin d'étudier les termes $(a.1)$ et $(a.2)$ ci-dessus. Pour le terme $(a.1)$,  par les inégalités de Hardy-Littlewood-Sobolev et l'inégalité de H\"older nous avons 
	\begin{eqnarray*}
		\left\Vert  ((\V_n-\V)\cdot \vec{\nabla}\theta_r)[\theta_r \V_n]\right\Vert_{\dot{H}^{-1}}&\leq & \left\Vert  ((\V_n-\V)\cdot \vec{\nabla}\theta_r)[\theta_r \V_n]\right\Vert_{L^{\frac{6}{5}}} \leq  \Vert  ((\V_n-\V)\cdot \vec{\nabla}\theta_r) \Vert_{L^2} \Vert \theta_r \V_n \Vert_{L^3}\\
		&\leq & \Vert \V_n-\V \Vert_{L^2}\Vert  \vec{\nabla}\theta_r \Vert_{L^{\infty}} \Vert \theta_r\Vert_{L^{6}}\Vert \V_n \Vert_{L^6},
	\end{eqnarray*} d'où, comme la suite $(\V_n)_{n\in \mathbb{N}}$ est bornée dans $E$ et comme $E\subset H^{1}(\Rt)$, nous avons donc que cette suite bornée dans $H^1(\Rt)$ et alors elle est aussi bornée dans $L^6(\Rt)$. De cette façon, par les estimations précédentes  nous pouvons écrire 
	$$  \left\Vert  ((\V_n-\V)\cdot \vec{\nabla}\theta_r)[\theta_r \V_n]\right\Vert_{\dot{H}^{-1}} \leq c(\theta_r)\Vert \V_n-\V \Vert_{L^2}, $$ d'où, comme la suite $(\V_n)_{n\in \mathbb{N}}$ converge fortement vers $\V$ dans $L^2$ nous avons  alors que le terme $(a.1)$ converge vers zéro lorsque $n$ tend vers l'infini. \\
	\\
	Pour le terme $(a.2)$, toujours par  les inégalités de Hardy-Littlewood-Sobolev et l'inégalité de H\"older nous pouvons écrire
	\begin{eqnarray*}
		\left\Vert  (\theta^{2}_{r}(\V_n-\V)\cdot\vec{\nabla})\V_n  \right\Vert_{\dot{H}^{-1}}&\leq & c \left\Vert  (\theta^{2}_{r}(\V_n-\V)\cdot\vec{\nabla})\V_n  \right\Vert_{L^{\frac{6}{5}}}\leq  \Vert \theta^{2}_{r}(\V_n-\V) \Vert_{L^3} \Vert \vec{\nabla}\otimes \V_n \Vert_{L^2}\\
		&\leq & \Vert \theta^{2}_{r}\Vert_{L^{12}}\Vert \V_n-\V \Vert_{L^4}\Vert \V_n \Vert_{\dot{H}^{1}} \leq c(\theta_r) \Vert \V_n-\V \Vert_{L^4},
	\end{eqnarray*} d'où, comme la suite comme la suite $(\V_n)_{n\in \mathbb{N}}$ converge fortement vers $\V$ dans $L^4(\Rt)$ nous avons  ainsi que le terme $(a.2)$ converge vers zéro lorsque $n$ tend vers l'infini et de cette façon  le terme $(a)$ n'annule si $n\longrightarrow +\infty$. \\
	\\
	La vérification  que le terme $(b)$ converge vers zéro lorsque $n\longrightarrow +\infty$ suit essentiellement les mêmes lignes que le terme $(a)$. \\
	\\
	Nous avons donc que  $(T_r(\V_{n}))_{n\in \mathbb{N}}$ converge fortement vers $T(\V)$ dans $E$ et alors $T_r$ est un opérateur compact sur cet espace.
	
	\item[3)] Estimation \emph{a priori}. Soit donc $\lambda \in [0,1]$, si $\U\in E$ satisfait l'équation de point fixe $\U=\lambda T_r(\U)$, par la définition de l'opérateur $T_r$ nous avons l'identité
	$$\U =\lambda \frac{-\nu \Delta}{-\nu \Delta +\alpha I_d}\left[\frac{1}{\nu} \P \left(\frac{1}{\Delta} \left(  ([\theta_r \U]\cdot \vec{\nabla})[\theta_r \U] \right)\right)\right] - \lambda\frac{-\nu \Delta}{-\nu \Delta +\alpha I_d}\left[\frac{1}{\nu \Delta}( \fe)\right],$$
	d'où nous pouvons écrire 
	$$ [-\nu \Delta +\alpha I_d]\U=-\lambda \P ([\theta_r \U]\cdot \vec{\nabla})[\theta_r \U]+\lambda \fe,$$ et alors nous observons que $\U$ satisfait l'équation 
	\begin{equation}\label{N-S-amortie-theta-turbulent}
	-\nu \Delta \U+\alpha \U= -\lambda \P ([\theta_r \U]\cdot \vec{\nabla})[\theta_r \U]+\lambda \fe.
	\end{equation}
	Maintenant, comme $\U \in E$ alors nous savons que $\U \in H^{1}(\Rt)$ et dans l'équation ci-dessus nous obtenons l'identité 
	$$-\nu \langle \Delta \U,\U\rangle_{H^{-1}\times H^1}+\alpha \langle\U, \U\rangle_{H^{-1}\times H^1}= -\lambda \langle\P ([\theta_r \U]\cdot \vec{\nabla})[\theta_r \U],\U \rangle_{H^{-1}\times H^1}+\lambda \langle\fe, \U\rangle_{H^{-1}\times H^1},$$ où, comme $div(\U)=0$ et par les propriétés du projecteur de Leray nous avons 
	$$\langle \P (([\theta_r \U]\cdot \vec{\nabla})[\theta_r \U]),\U \rangle_{H^{-1}\times H^{1}}= \langle ([\theta_r \U]\cdot \vec{\nabla})[\theta_r \U],\U \rangle_{H^{-1}\times H^{1}} =\langle (\U\cdot \vec{\nabla})[\theta_r \U],\theta_r \U \rangle_{ H^{-1}\times H^{1}} =0,
	$$ et de cette façon, par une intégration par parties nous écrivons 
	$$  \nu \Vert \U \Vert^{2}_{\dot{H}^{1}}+\alpha \Vert \U \Vert^{2}_{L^2}=\lambda \langle \fe, \U\rangle_{H^{-1}\times H^{1}},$$ d'où  nous avons
	$$ \min(\nu,\alpha)\Vert \U \Vert^{2}_{H^1} \leq \nu \Vert \U \Vert^{2}_{\dot{H}^{1}}+\alpha \Vert \U \Vert^{2}_{L^2}\leq \lambda \langle \fe, \U\rangle_{H^{-1}\times H^{1}}.$$
	De plus, comme $0\leq \lambda \leq 1$, par l'inégalité de Cauchy-Schwarz et les inégalités de Young nous pouvons écrire  
	$$ \lambda \langle \fe, \U\rangle_{H^{-1}\times H^{1}} \leq  \Vert \U \Vert_{H^1}\Vert \fe \Vert_{H^{-1}}\leq \frac{\min(\nu,\alpha)}{2}\Vert\U \Vert^{2}_{H^1} + \frac{2}{\min(\nu,\alpha)}\Vert \fe \Vert_{H^{-1}},$$ 
	%
	%
	et alors nous obtenons 
	\begin{equation}\label{control_Ur_turbulent}
	\Vert \U \Vert_{H^1}\leq c(\alpha,\nu)\Vert \fe \Vert_{H^{-1}},
	\end{equation} Finalement, nous posons la constante $C=c(\alpha,\nu)\Vert \fe \Vert_{H^{-1}}$ et nous vérifions de cette façon le point $3)$ du Lemme \ref{Schaefer}. \\
\end{enumerate}
Ainsi, par une application de ce lemme  il existe $\U_r\in E$ solution de l'équation de point fixe $\U_r =T_r(\U_r)$. \\
\\
Une fois que l'on a montré l'existe des solutions $\U_r$ (pour $r>0$ fixe) du problème de point fixe ci-dessus nous allons montrer que la famille de solutions $(\U_r)_{r>0}$ converge vers $\U$ une solution des équations de Navier-Stokes stationnaires et amorties. Tout d'abord, par l'estimation (\ref{control_Ur_turbulent}) nous savons que $(\U_r)_{r>0}$ est uniformément bornée dans l'espace $H^1(\Rt)$ et alors 
$$ \sup_{r>0}\Vert \varphi \U_r \Vert_{H^1}<+\infty,$$ pour tout $\varphi \in \mathcal{C}^{\infty}_{0}(\Rt).$ De cette façon, par le lemme de Rellich-Lions et en appliquant la méthode d'extraction diagonale de Cantor, il existe une suite de nombres positifs $(r_n)_{n\in \mathbb{N}}$, tels que $r_n \longrightarrow +\infty$, et il existe  $\U \in H^{1}_{loc}(\Rt)$ telle que la suite $(\U_{r_n})_{n\in\mathbb{N}}$ converge fortement vers $\U$ dans $L^{2}_{loc}(\Rt)$. Mais, toujours par l'estimation (\ref{control_Ur_turbulent}) nous avons que $(\U_{r_n})_{n\in\mathbb{N}}$  converge faiblement  dans $H^1(\Rt)$ ce qui nous donne  $\U\in H^1(\Rt)$. \\
\\
De plus, par les inégalités de Hardy-Littlewood-Sobolev nous savons que $(\U_{r_n})_{n\in\mathbb{N}}$ est uniformément bornée dans $L^6(\Rt)$ et donc  $(\U_{r_n})_{n\in \mathbb{N}}$ converge fortement vers $\U$ dans $L^{p}_{loc}(\Rt)$ pour $2\leq p <6$. \\
\\  
Maintenant, nous allons prouver que le terme non linéaire $((\theta_{r_n} \U_{r_n})\cdot \vec{\nabla})(\theta_{r_n} \U_{r_n})$ converge vers $(\U\cdot \vec{\nabla})\U$ dans $\mathcal{D}'$. En effet, pour $n\in \mathbb{N}$ fixé, comme $div(\U_{r_n})=0$, nous pouvons écrire 
$$ ((\theta_{r_n} \U_{r_n})\cdot \vec{\nabla})(\theta_{r_n} \U_{r_n})=\theta_r\left( (\U_{r_n}\cdot \vec{\nabla})(\theta_{r_n} \U_{r_n}) \right)=\theta_r div(\U_{r_n}\otimes (\theta_{r_n}\U_{r_n})),$$ et comme $\U_{r_n}$ converge fortement vers $\U$ dans $L^{4}_{loc}(\Rt)$ et de plus, étant donné que $\theta_{r_n}(x)=1$ lorsque $\vert x \vert < r_n$, alors nous avons que   $\theta_{r_n}\U_{r_n}$ converge fortement vers $\U$ dans $L^{4}_{loc}(\Rt)$ et de cette façon nous obtenons que $\U_{r_n}\otimes [\theta_{r_n}\U_{r_n}]$ converge fortement vers $\U\otimes \U$ dans $L^{2}_{loc}(\Rt)$. Ainsi, nous avons que $\theta_r div(\U_{r_n}\otimes [\theta_{r_n}\U_{r_n}])$ converge vers $div(\U\otimes \U)$ dans $\mathcal{D}'$ et comme, pour $\U \in H^1(\Rt)$ nous avons en plus que $div(\U)=0$, nous pouvons écrire $div(\U \otimes \U)=(\U\cdot \vec{\nabla})\U$. \\  
\\
Une fois que l'on a la convergence de  $((\theta_{r_n} \U_{r_n})\cdot \vec{\nabla})(\theta_{r_n} \U_{r_n})$  vers $(\U\cdot \vec{\nabla})\U$ dans $\mathcal{D}'$ alors nous avons l'identité suivante au sens des distributions
\begin{equation}\label{convergence_turbulent}
\lim_{r_n\longrightarrow +\infty} \left[ -\nu \Delta \U_{r_n} +((\theta_{r_n} \U_{r_n})\cdot \vec{\nabla})(\theta_{r_n} \U_{r_n}) -\fe +\alpha \U_{r_n}\right] = -\nu \Delta \U + ( \U \cdot \vec{\nabla}) \U-\fe +\alpha \U.
\end{equation}
D'autre part, comme $\U_{r_n}$ vérifie l'équation 
$$ -\nu \Delta \U_{r_n} + \P ((\theta_{r_n} \U_{r_n})\cdot \vec{\nabla})(\theta_{r_n} \U_{r_n}) =\fe-\alpha \U_{r_n}, $$ et comme $div(\U_{r_n})=0$ et $div(\fe)=0$ nous pouvons écrire 
$$ \P\left[  -\nu \Delta \U_{r_n} +((\theta_{r_n} \U_{r_n})\cdot \vec{\nabla})(\theta_{r_n} \U_{r_n}) -\fe +\alpha \U_{r_n}\right]=0,$$ et alors il existe $P_{r_n}\in \mathcal{D}'$  telle que 
$$ -\nu \Delta \U_{r_n} +((\theta_{r_n} \U_{r_n})\cdot \vec{\nabla})(\theta_{r_n} \U_{r_n}) -\fe +\alpha \U_{r_n}=\vec{\nabla}P_{r_n}.$$ \`A partir de cette identité nous pouvons observer que  le rotationnel du terme de gauche ci-dessus est égal à zéro et alors, par l'identité (\ref{convergence_turbulent})  nous en tirons que 
$$ \vec{\nabla}\wedge \left[  -\nu \Delta \U + ( \U \cdot \vec{\nabla}) \U-\fe +\alpha \U\right] =0,$$ et alors il existe $P\in \mathcal{D}'$
telle que $\ds{ -\nu \Delta \U + ( \U \cdot \vec{\nabla}) \U-\fe +\alpha \U=\vec{\nabla}P}$. De plus, toujours par la relation $P=\frac{1}{-\Delta}div (( \U \cdot \vec{\nabla}) \U)$ et comme $\U \in H^{1}(\Rt)$ nous avons que $P\in H^{\frac{1}{2}}(\Rt)$. \\
\\
Pour finir la preuve de ce théorème vérifions maintenant que toute solution $\U\in H^1(\Rt)$ vérifie l'estimation (\ref{estimation-a-priori}). En effet, il suffit de remarquer que si $\U\in H^1(\Rt)$ alors par les lois de produit nous avons $\U \otimes \U \in H^{\frac{1}{2}}(\Rt)$,  et  comme $div(\U)=0$ alors nous avons $(\U\cdot \vec{\nabla})\U =div(\U \otimes \U) \in H^{-\frac{1}{2}}(\Rt) \subset H^{-1}(\Rt)$. Ainsi, étant donné que $\U$ vérifie les équations (\ref{NS-stationnaire-motiv}) nous pouvons écrire $$ \langle -\nu \Delta \U + ( \U \cdot \vec{\nabla}) \U-\fe +\alpha \U+\vec{\nabla}P,U \rangle_{H^{-1}\times H^1}=0,$$ d'où nous pouvons en tirer l'estimation  (\ref{estimation-a-priori}).  \finpv
\\  
Nous observons ainsi que cette méthode nous permet de construire des solutions des équations  (\ref{NS-stationnaire-motiv}) sans faire aucune hypothèse supplémentaire sur la force extérieure. Mais, cette méthode ne fournit pas  d'information additionnelle sur l'unicité de ces solutions et  cette question est encore ouverte  dans le cadre d'une force extérieure $\fe$ non nulle. En revanche, si l'on considère le cas particulier  d'une force $\fe=0$ alors nous avons le corollaire suivant:

\begin{Corollaire} La solution triviale $\U=0$ est l'unique solution des équations 
$$	-\nu \Delta \U +(\U\cdot \vec{\nabla})\U +\vec{\nabla}P+ \alpha \U =0, \qquad div(\U)=0, \quad \alpha>0.$$
	\end{Corollaire}
\pv  Par le Théorème \ref{Theo:solutions_stationnaires_turbulent} nous savons que toute solution $\U\in H^1(\Rt)$ satisfait l'estimation (\ref{estimation-a-priori}): $\min(\alpha,\nu)\Vert \U \Vert_{H^1}\leq \Vert \fe \Vert_{H^{-1}}$, d'où nous observons que si $\fe=0$ alors $\U=0$ est l'unique solution des équations ci-dessus. \finpv 
\\
Néanmoins, dans ce chapitre nous allons  travailler avec  une force $\fe$ non nulle et cette force sera le fil conducteur dans notre étude  comme nous l'expliquerons  plus tard dans la Section \ref{Sec:force ext: construction et proprietes}. \\
\\
Une fois que nous avons construit des solutions des  équations (\ref{NS-stationnaire-motiv}) dans la Section \ref{Sec:proprietes-sol-stat} nous rentrons dans le vif du sujet et nous étudions quelques propriétés intéressantes de ces solutions: leur stabilité et la localisation en variable d'espace. \\
\\
Mais avant d'entrer dans les détails techniques, nous avons besoin de préciser encore un peu plus notre cadre de travail et c'est pour cette raison que dans la section suivante nous allons introduire deux objets qui seront à la base de notre étude.  

\section{Le régime laminaire et le régime turbulent}\label{sec:regime-lam-turb} 
 
Ces deux régimes du fluide seront caractérisés par le nombre de Grashof que nous introduisons tout de suite. 
\subsection{Les nombres de Grashof}\label{Sec:Grashof}
Dans la Section \ref{Sec:cadre-periodique} page \pageref{Sec:cadre-periodique} nous avons expliqué comment le nombre de Reynolds $Re$ sert à caractériser le régime du mouvement du fluide laminaire (ou turbulent) mais cette quantité  dépend des solutions des équations de Navier-Stokes et alors ce nombre nous donne une caractérisation \emph{a  posteriori} du régime laminaire ou turbulent du fluide.\\ 
\\
Dans ce cadre, afin d'obtenir une caractérisation \emph{a priori} du régime laminaire (ou turbulent) nous allons introduire une quantité \emph{équivalente} au nombre de Reynolds $Re$ et  en suivant une idée des articles \cite{DoerFoias}, \cite{FMRT1} et \cite{FMRT2}  de C. Foias, R. Temam \emph{et al.}  nous allons considérer ici le nombre  Grashof. Mais avant de définir  ce nombre nous avons besoin d'introduire tout d'abord quelques quantités physiques. \\
\\
En effet, pour définir le nombre de Grashof nous avons besoin de considérer une longueur caractéristique $L$ (comme dans la définition du nombre de Reynolds $Re$ ci-dessus) mais cette longueur est  très particulière car elle correspond à une signification physique bien précise: dans le cadre périodique elle représente la période et dans les autres cas elle est sensée représenter la dimension du domaine utilisé. \'Etant donné que cette notion physique est délicate à mettre en place lorsqu'on travaille sur l'espace tout entier, et en suivant les idées du chapitre précédent  nous allons intégrer cette longueur comme un paramètre du modèle: il s'agira alors d'une donnée du problème tout comme la constante de viscosité $\nu$.  Avec cette idée en tête nous commençons par fixer tous les paramètres physiques que nous considérerons dans ce chapitre. 
\begin{Definition}[Les paramètres physiques]\label{Param_physiques} 
\begin{enumerate}
\item[]
\item[1)] $\nu>0$ est toujours la constante de viscosité du fluide.
\item[2)] $L>0$ est la longueur caractéristique du fluide: la plus grande échelle de longueur sur laquelle on veut étudier le comportement laminaire ou turbulent du fluide.
\item[3)] $\ell_0>0$  (avec $L\geq \ell_0$) est l'échelle d'injection d'énergie: l'échelle de longueur à laquelle la force $\fe$ agira sur le fluide.  
\item[4)] $F>0$ est l'amplitude que l'on donnera à la force extérieure $\fe$. 
\end{enumerate}
\end{Definition} 
Une fois que l'on a fixé les paramètres physiques ci-dessus nous pouvons maintenant introduire les nombres de Grashof $G_{\theta}$  où $\theta$ sera un paramètre d'interpolation.  \\
\\
En effet, dans la littérature, il n'y a pas de façon standard pour définir le nombre de Grashof et  pour $(L,\ell_0,\nu, F)$ les paramètres physiques donnés dans la définition ci-dessus, dans les notes de cours \cite{Const} de P. Constantin, il est suggéré de définir le nombre de Grashof par  $$G_r=\frac{FL^3}{\nu^2},$$
où nous pouvons observer que seule la longueur caractéristique $L$ est considérée et   l'échelle d'injection d'énergie $\ell_0$ n'a pas été prise en compte.\\
\\
D'autre part, dans les articles \cite{DoerFoias}, \cite{FMRT1} et \cite{FMRT2} de C. Foias,  R. Temam \emph{et. al.}  les auteurs considèrent le nombre de Grashof $$G_r=\frac{F\ell^{3}_{0}}{\nu^2},$$ où cette fois-ci la longueur caractéristique $L$ a été négligée et seule l'échelle d'injection d'énergie $\ell_0$ est considérée. \\
\\
Remarquons  que comme  la longueur caractéristique $L$ vérifie toujours la relation $L\geq \ell_0$ (voir le point $3)$ de la Définition \ref{Param_physiques}) alors on obtient la relation  $$  \frac{F\ell^{3}_{0}}{\nu^2} \leq \frac{FL^3}{\nu^2},$$ et bien sûr ces deux quantités sont égales si l'on a $L=\ell_0$. \\
\\
Dans ce cadre, dans la définition du nombre de Grashof que nous allons mettre en œuvre  ici, nous allons  interpoler entre les échelles $\ell_0$ et $L$ et nous définissons une famille de nombres $G_\theta$ de la façon suivante.  
\begin{Definition}[Nombres de Grashof]\label{Grashof_theta} Soient $(\nu,L,\ell_0,F)$ les paramètres physiques donnés dans la Définition \ref{Param_physiques}. Pour $0\leq \theta \leq 3$  nous définissons le nombre de Grashof $G_\theta$ par $$G_\theta =\frac{F L^{\theta}\ell^{3-\theta}_{0}}{\nu^2}.$$
\end{Definition} 
Dans cette définition nous pouvons observer que si l'on prend le paramètre $\theta=0$ alors obtient le nombre $G_0=\frac{F \ell^{3}_{0}}{\nu^2}$ et si l'on fixe le paramètre $\theta=3$ on a bien  $G_3=\frac{F L^3}{\nu^2}$ ce qui correspond aux définitions du nombre de Grashof considérées dans la littérature.\\
\\
Une fois que l'on a définit la famille de nombres de Grashof $(G_\theta)_{0\leq \theta \leq 3}$ nous  étudions la relation entres ces nombres.
\begin{Proposition}\label{Proposition: ordre Grashof}Si $0\leq \theta_1\leq \theta_2\leq 3$ alors on a $G_{\theta_1}\leq G_{\theta_2}$.
\end{Proposition}
\pv En effet, nous pouvons écrire $\ds{ G_{\theta_1}=\frac{FL^{\theta_1}\ell^{3-\theta_1}_{0}}{\nu^2}=\frac{F}{\nu^2}  \left( \frac{L}{\ell^{1}_{0}}\right)^{\theta_1}\ell^{2}_{0}}$, d'où,
par le point $3)$ de la Définition \ref{Param_physiques} nous avons $\frac{L}{\ell_0}\geq 1$ et alors $\left( \frac{L}{\ell_0} \right)^{\theta_1}\leq \left( \frac{L}{\ell_0}\right)^{\theta_2}$. De cette façon, en revenant a l'identité précédente nous obtenons 
$$ G_{\theta_1}= \frac{F}{\nu^2}  \left( \frac{L}{\ell_0}\right)^{\theta_1} \ell^{2}_{0}\leq \frac{F}{\nu^2}  \left( \frac{L}{\ell_0}\right)^{\theta_2}\ell^{2}_{0}=G_{\theta_2}.$$ \finpv
Nous pouvons alors observer que si l'on contrôle un nombre de Grashof $G_{\theta_1}<\eta$,  alors on contrôle tous les nombres de Grashof $G_{\theta}$ avec $0\leq \theta \leq \theta_1$. \\ 
\\
Le régime laminaire que nous allons considérer tout au long de la Section \ref{Sec:Stabilite-laminaire}  correspond à un contrôle sur le nombre de Grashof $G_\theta<\eta$,  pour un certain paramètre $\theta \in [0,3]$ et une constante $\eta>0$ que nous fixerons plus tard.\\
\\
Observons maintenant qu'étant donné que la force extérieure $\fe$ des équations de Navier-Stokes amorties et stationnaires  (\ref{NS-stationnaire-motiv})  est une donnée de notre modèle, nous voulons à présent relier le contrôle sur le nombre de Grashof $G_{\theta}$ avec un contrôle sur cette force $\fe$ et donc, dans le section qui suit nous cherchons à construire une force extérieure $\fe$ de sorte que la taille de cette fonction puisse être directement reliée  avec le nombre de Grashof $G_\theta$.  L'idée sous-jacente étant de contrôler avec un seul paramètre la turbulence et la force.     
\subsection{Une force extérieure bien préparée}\label{Sec:force ext: construction et proprietes}   
Dans cette section nous définissons une force extérieure très particulière qui sera utilisée tout au long de ce chapitre et  pour cela nous avons besoin d'introduire l'ondelette $\vec{\phi}:\Rt\longrightarrow \Rt$ suivante.
\begin{Definition}\label{Ondelette de base} Soit $\vec{\phi}=(\phi_1, \phi_2, \phi_3)$ un champ de vecteurs dans la classe de Schwartz tel que: 
\begin{enumerate}
\item[1)] $\vec{\phi}$ est un vecteur  à divergence nulle. 
\item[2)] Pour $0<\rho_1<\rho_2$ deux constantes fixées, on a $supp\left(\widehat{\phi_i}\right)\subset \lbrace \xi \in \Rt: \rho_1 \leq \vert \xi \vert \leq  \rho_2 \rbrace$, pour tout $i=1,2,3$.
\item[3)] Pour tout $\vec{\psi}$ champ de vecteurs dans la classe de Schwartz qui vérifie la propriété $2)$ ci-dessus, on a 
$$\displaystyle{\int_{\Rt}\vec{\phi}(x-k)\cdot \vec{\psi}(x-m)dx}=\delta_{k,m},$$ pour tout $k,m\in \Zt,$ et où $\delta_{k,m}$ est la fonction delta de  Kronecker. 
\end{enumerate} 
\end{Definition} 
Dans le lemme qui suit nous énonçons quelques propriétés bien connues sur les ondelettes qui nous seront utiles par la suite, voir les livres \cite{KahLer} de J.P. Kahane \& P.G.  Lemarié-Rieusset et  \cite{Meyer} de Y. Meyer   pour une preuve de ces résultats classiques et pour plus de détails sur les ondelettes.
\begin{Lemme}\label{prop_ondelette_base} Soit $\vec{\phi}=(\phi_1,\phi_2,\phi_2)$ le champ de vecteurs donné par la Définition \ref{Ondelette de base}. Pour tout $s\in \mathbb{R}$ on définit $(-\Delta)^s \vec{\phi}= ((-\Delta)^s \phi_1, (-\Delta)^s \phi_2,(-\Delta)^s  \phi_3).$ Alors:
\begin{enumerate}
\item[1)]  Pour tout  $N\in \mathbb{N}$ il existe une constante $C_0(N,s,\vec{\phi})>0$ telle que, pour tout nombre réel $A\geq 1$,  on a la majoration: 
\begin{equation*}
\sum_{k\in \mathbb{Z}^3, \vert k \vert\leq A}\vert (-\Delta)^s \vec{\phi}(x-k)\vert \leq C_0 \left(\frac{A}{A+\vert x\vert}\right)^N.
\end{equation*}
\item[2)]  Pour  $1\leq p \leq +\infty$ il existe deux constantes $A_{s,p}, B_{s,p}>0$, qui ne dépendent que de $s,p$ et $\vec{\phi}$, telles que, pour toute suite  $(\lambda_k)_{k\in \mathbb{Z}^3}\in \ell^{p}(\mathbb{Z}^3)$ on a la \emph{propriété de presque-orthogonalité}:
\begin{equation*}
A_{s,p}\Vert (\lambda_k) \Vert_{\ell^{p}(\Zt)} \leq \left\Vert \sum_{k\in \mathbb{Z}^3}\lambda_k (-\Delta)^s \vec{\phi}(\cdot -k)\right\Vert_{L^p(\Rt)}\leq B_{s,p}\Vert (\lambda_k) \Vert_{\ell^{p}(\Zt)}.
\end{equation*}
 \end{enumerate}
\end{Lemme}
Une fois que nous avons introduit l'ondelette $\vec{\phi}$ ci-dessus nous pouvons construire la force extérieure $\fe$. Ainsi, pour $L>0$ la longueur caractéristique du fluide, donnée dans la Définition \ref{Param_physiques}, nous considérons le cube $[-L,L]^3\subset \Rt$. Comme l'échelle d'injection d'énergie $\ell_0>0$ est telle que $\frac{L}{\ell_0}\geq 1$  alors, dans le cube $[-L,L]^3$, nous allons considérer les points de la forme $\ell_0 k$ où $\vert \ell_0 k\vert\leq L$ avec $k\in \Zt$, voir la Figure $2.1$  ci-dessous.  \\
\\ 
\begin{figure}[!h]\label{dessin_force}
\begin{center}
\includegraphics[scale=0.35]{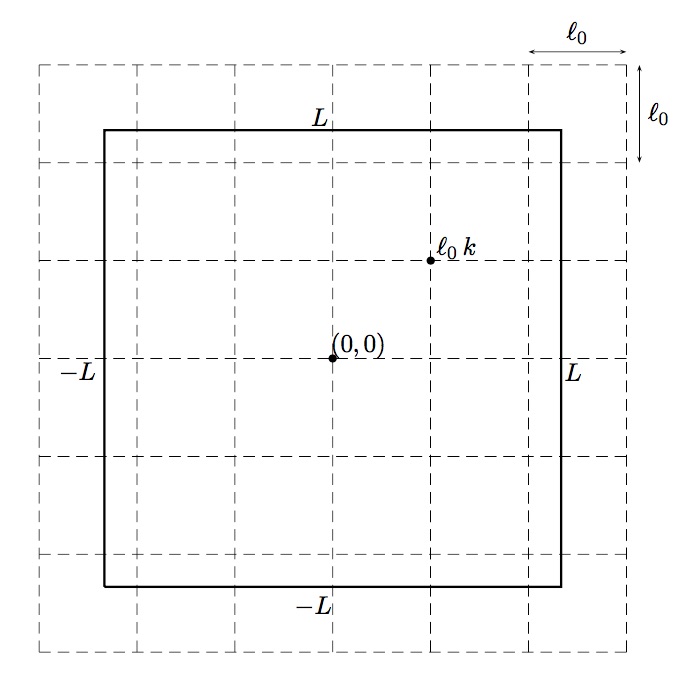}\caption{\footnotesize{Exemple de la construction de la force dans $\mathbb{R}^2$  où $2\ell_0<L$.}}
\end{center}
\end{figure} 
De cette façon, en suivant une idée de l'article \cite{DoerFoias} de C. Doering et C. Foias  (qui a été donnée dans le cadre périodique), nous construisons la force  $\fe$ par  translation de la fonction $\vec{\phi}$ en chaque point $\ell_0 k$ (avec $\vert \ell_0 k\vert\leq L$) et par  dilatation à l'échelle $\frac{1}{\ell_0}$. Nous avons donc:
\begin{Definition}[La force extérieure]\label{Definition_force_ext} Soient $L,\ell_0$ et $F$ les paramètres physiques donnés par la Définition \ref{Param_physiques} et soit $\vec{\phi}$ le champ de vecteurs donné par la Définition \ref{Ondelette de base}. Nous définissons le champ de vecteurs $\fe(x)$ par  l'expression:
$$ \fe(x)=F\sum_{\vert \ell_0 k \vert \leq L} \vec{\phi} \left(\frac{x}{\ell_0}-k\right), \quad \text{pour tout}\quad x\in \Rt.$$  
\end{Definition}
Grâce à cette définition nous observons directement que la force $\fe$ est un champ de vecteurs dans la classe de Schwartz et à divergence nulle. \\
\\
Dans les deux lemmes ci-dessous  nous étudions  maintenant quelques propriétés  de ce champ de vecteurs qui nous seront utiles par la suite.
\begin{Lemme}\label{lemma_loc_fe} Soient $L>0$ la longueur caractéristique du fluide donnée dans la Définition \ref{Param_physiques},  $\fe$ la force donnée par la Définition \ref{Definition_force_ext}  et  $s\in \mathbb{R}$.  Nous avons les points suivants:
\begin{enumerate}
\item[1)] \emph{Concentration en variable spatiale.}   Pour tout nombre réel $\mu\geq 1$ considérons le cube dilaté $[-\mu L,\mu L]^3$ et son complémentaire $([-\mu L,\mu L]^3)^c$. Alors la quantité 
\begin{equation*}
\ds{v(\mu)=\left(\int_{ ([-\mu L,\mu L]^3)^c}\vert (-\Delta)^s\fe(x)\vert^2dx\right)^{\frac{1}{2}}},
\end{equation*} est à décroissance rapide lorsque $\mu$ tend vers l'infini.
\item[2)] \emph{Support fréquentiel.}   On a  $supp\,\left((-\widehat{\Delta)^s\fe}\right)\subset \lbrace \xi \in \Rt: \frac{\rho_1}{\ell_0} \leq \vert \xi \vert\leq  \frac{\rho_2}{\ell_0}\rbrace$,  où   $0<\rho_1<\rho_2$ sont les constantes données dans le point  $2)$ de la Définition \ref{Ondelette de base}.
\end{enumerate} 
\end{Lemme}
\pv
\begin{enumerate}
\item[1)] Soit donc $\mu\geq 1$. Nous écrivons $\ds{v(\mu)=\left(\int_{([-\mu L,\mu L]^3)^c}\vert (-\Delta)^s \fe(x)\vert^2 dx\right)^{\frac{1}{2}}}$ et nous allons montrer  que pour tout $N\in \mathbb{N}$ tel que $N\geq 2$ il existe une constante $C(N)>0$ telle que 
$$ v(\mu)\leq \frac{C(N)}{\mu^{N-\frac{3}{2}}}.$$
En effet, par la Définition \ref{Definition_force_ext} et par le point $1)$ du  Lemme \ref{prop_ondelette_base}  (en prenant $s=0$ et $A=\frac{L}{\ell_0}\geq1$) nous obtenons directement la majoration
\begin{eqnarray*}
v^2(\mu)& \leq &F^2 \int_{\Rt \setminus [-\mu L,\mu L]^3}\left\vert \sum_{\vert k\vert\frac{L}{\ell_0}} \vec{\phi}\left( \frac{x}{\ell_0}-k \right)\right\vert^2dx \leq C_0(N) F^2 \int_{\Rt \setminus [-\mu L, \mu L]^3} \left(\frac{L}{L+\vert x \vert}\right)^N dx\\
& &  \leq C_0(N) F^2L^3\frac{1}{\mu^{2N-3}}=\frac{C(N)}{\mu^{2N-3}}. 
\end{eqnarray*}
\item[2)] Ce point est une conséquence directe de la Définition \ref{Definition_force_ext}. \finpv Il convient maintenant de faire la remarque suivante. 
\end{enumerate}
\begin{Remarque}\label{rmq:loc-fe}Le Lemme \ref{lemma_loc_fe} nous donne des informations intéressantes par rapport à la concentration  de la force extérieure $\fe$  en variable spatiale ou en fréquence.  En effet:
\begin{enumerate}
\item[i)] En variable d'espace. La quantité $v(\mu)$ mesure la concentration de la fonction $\fe$ hors du cube dilaté $[-\mu L, \mu L]^3$ et donc, comme la quantité $v(\mu)$ est à décroissance rapide lorsque $\mu$ tend vers l'infini nous observons donc que la fonction $\fe$  est \emph{essentiellement} concentrée sur le cube $[-L,L]^3$.
\item[ii)] En variable de fréquence. La transformée de Fourier de la force $\fe$ est localisée aux fréquences $ \frac{\rho_1}{\ell_0} \leq \vert \xi \vert \leq  \frac{\rho_2}{\ell_0}$ ce qui correspond au fait que, selon le modèle de cascade d'énergie expliqué dans l'introduction du chapitre précédent, cette force $\fe$ introduit l'énergie cinétique aux échelles de longueur de l'ordre de $\ell_0$ et donc aux fréquences de l'ordre de $\frac{1}{\ell_0}$. 
\end{enumerate} 
\end{Remarque}
Nous étudions maintenant une deuxième propriété utile de la force $\fe$. 
\begin{Lemme}\label{normes_fe} Soient les paramètres $L,\ell_0,F>0$ donnés dans la Définition \ref{Param_physiques} et la fonction $\fe$ donnée dans la Définition \ref{Definition_force_ext}. Pour $s\in \mathbb{R}$ et $1\leq p \leq +\infty$,  il existe deux constantes $C_{s,p}, D_{s,p}>0,$ indépendantes des paramètres physiques ci-dessus, telles que
\begin{equation*}
C_{s,p}\, FL^{\frac{3}{p}}\ell^{-2s}_{0}\leq \left\Vert (-\Delta)^s \fe \right\Vert_{L^p}\leq D_{s,p} \,FL^{\frac{3}{p}}\ell^{-2s}_{0}.
\end{equation*}
\end{Lemme}
\pv  Pour tout $s\in \mathbb{R}$ nous avons 
\begin{equation*}
(-\Delta)^s\fe(x)=F\ell^{-2s}_{0}\sum_{\vert \ell_0k\vert\leq L} (-\Delta)^s \vec{\phi}\left(\frac{x-\ell_0k}{\ell_0}\right)=F\ell^{-2s}_{0} \sum_{k\in \Zt}\lambda_k\left((-\Delta)^s \vec{\phi}\right)\left(\frac{x}{\ell_0}- k\right),
\end{equation*} 
 où la suite  $(\lambda_k)_{k\in \Zt}$ est définie par
 \begin{equation*}
 \lambda_k=\left\lbrace
\begin{array}{rl} 1&\text{si}\quad \vert \ell_0 k \vert\leq L, \\
0&\text{si non.}\end{array}\right. 
\end{equation*}   
Alors, en prenant la norme $L^p$ à chaque côté de l'identité ci-dessus et ensuite en utilisant le point $2)$ du Lemme \ref{prop_ondelette_base}, on a  qu'il existe deux constantes  $A_{s,p},B_{s,p}>0$ telles que
\begin{equation}\label{est_tech_1}
A_{s,p}\,F\ell^{-2s+\frac{3}{p}}_{0}\Vert (\lambda_k) \Vert_{\ell^{p}} \leq F\ell^{-2s+\frac{3}{p}}_{0}\left\Vert \sum_{k\in \mathbb{Z}^3}\lambda_k (-\Delta)^s \vec{\phi}(\cdot -k)\right\Vert_{L^p}\leq B_{s,p}\,F\ell^{-2s+\frac{3}{p}}_{0} \Vert (\lambda_k) \Vert_{\ell^{p}},
\end{equation} 
pour $1\leq p \leq+\infty$.  
D'autre part, par la définition de la suite $(\lambda_k)_{k\in \Zt}$ nous observons que $\Vert (\lambda_k)\Vert^{p}_{L^p}$ est le nombre de points de la forme $\ell_0k$ dans le cube $[-L,L]^3$ et alors, il existe deux constantes $c_1,c_2>0$ (indépendantes des paramètres $L,\ell_0$ et $F$) telles que  
\begin{equation*}
c_1\left(\frac{L}{\ell_0}\right)^{\frac{3}{p}} \leq \Vert (\lambda_k)\Vert_{\ell^{p}}\leq c_2\left(\frac{L}{\ell_0}\right)^{\frac{3}{p}},
\end{equation*}
 et donc, en posant les constantes $C_{s,p}=c_1A_{s,p}$ et  $D_{s,p}=c_2B_{s,p},$ par l'expression (\ref{est_tech_1}) nous obtenons l'estimation cherchée. \finpv
\begin{Remarque}\label{Remarque_normes_fe} Comme les constantes $C_{s,p}$ et $D_{s,p}$ sont indépendantes des paramètres  $(L,\ell_0,F,\nu,\alpha)$, nous écrirons par la suite
\begin{equation}\label{Formula_EstimateForce}
\left\Vert (-\Delta)^s\fe \right\Vert_{L^p}\approx FL^{\frac{3}{p}}\ell^{-2s}_{0},
\end{equation} 
d'où nous avons les estimations suivantes: si $s=0$, alors on a $\Vert \fe\Vert_{L^p} \approx FL^{\frac{3}{p}}$, pour tout $1\leq p\leq +\infty$. En particulier, pour $p=+\infty$ on obtient $\Vert \fe \Vert_{L^{\infty}}\approx F,$ et donc nous observons que l'amplitude de la force $\fe$ est bien donnée par le paramètre $F$ introduit dans la Définition \ref{Param_physiques}. De plus, comme $F\approx \frac{\Vert \fe\Vert_{L^p}}{L^{\frac{3}{p}}}$  et comme $\fe$ est essentiellement concentrée sur le cube  $[-L,L]^3$ (voir le point $1)$ du  Lemme \ref{lemma_loc_fe}) alors le paramètre  $F$ mesure aussi la moyenne de la force $\fe$ en norme $L^p$, avec $1\leq p <+\infty$.\\
\end{Remarque} L'intérêt principal de travailler avec cette force  $\fe$  repose sur le fait que grâce à l'estimation (\ref{Formula_EstimateForce}) ci-dessus nous allons pouvoir contrôler la taille de cette fonction par un contrôle direct sur  le nombre de Grashof $G_\theta$. En effet, nous avons le résultat suivant.

\begin{Proposition}\label{Prop:relation-force-grashof} Soit $\theta\in [0,3]$ et $G_{\theta}$ le nombre de Grashof donné dans la Définition \ref{Grashof_theta}. Soit $\fe$ la force extérieure donnée dans la Définition \ref{Definition_force_ext} et $\nu>0$ la constante de viscosité du fluide. Il existe  deux paramètres $s\in \mathbb{R}$ et $p\in [0,+\infty]$, qui dépendent du paramètre  $\theta$ ci-dessus, tels que l'on a l'estimation $\ds{\Vert (-\Delta)^s \fe \Vert_{L^p}\approx \nu^2 G_{\theta}}$. 
\end{Proposition}
\pv    Par le Lemme \ref{normes_fe}, pour $s\in \mathbb{R}$  et $1 \leq p \leq +\infty$ nous avons l'estimation 
$$ \Vert (-\Delta)^s \fe \Vert_{L^p}\approx F L^{\frac{3}{p}}\ell^{-2s}_{0}.$$
Maintenant, si  nous posons $\frac{3}{p}=\theta$ et $-2s =3-\theta$   alors nous obtenons la relation $$ \Vert (-\Delta)^s \fe \Vert_{L^p}\approx F L^{\frac{3}{p}}\ell^{-2s}_{0} = FL^{\theta}\ell^{3-\theta}_{0}=\nu^2 \frac{FL^{\theta}\ell^{3-\theta}_{0}}{\nu^2}=\nu^2 G_{\theta},$$ ce qui termine la  preuve de cette proposition.  \finpv \\ 
Ainsi,  nous pouvons écrire  $\ds{ \frac{\Vert (-\Delta)^s \fe \Vert_{L^p}}{\nu^2}\approx G_{\theta}},$ d'où nous observons que si nous  contrôlons le nombre de Grashof $G_{\theta}<\eta$ alors nous contrôlons la taille de la force $\fe$ dans l'espace de Sobolev $\dot{W}^{2s,p}(\Rt)$. \\
\\
Cette relation sera  utilisée dans la Section \ref{Sec:Stabilite-laminaire} où nous étudierons  quelques propriétés des solutions des équations de Navier-Stokes amorties et stationnaires (\ref{NS-stationnaire-motiv})  dans le régime laminaire.      


\section{Quelques propriétés des solutions stationnaires}\label{Sec:proprietes-sol-stat}
Maintenant que nous avons tous les ingrédients nécessaires nous allons rentrer dans le vif du sujet. Soit $(\U,P)\in H^1(\Rt)\times H^{\frac{1}{2}}(\Rt)$ une solution des équations de Navier-Stokes amorties et stationnaires 
 \begin{equation}\label{NS-stationnaire-prop}
-\nu \Delta \U +(\U\cdot \vec{\nabla})\U +\vec{\nabla}P  =\fe-\alpha \U, \qquad div(\U)=0, \quad \alpha>0,
\end{equation} obtenue par le biais du Théorème \ref{Theo:solutions_stationnaires_turbulent} et où $\fe$ est la force donnée dans la Définition \ref{Definition_force_ext}. \\
\\
Nous allons étudier maintenant  quelques  propriétés   de ces solutions. Dans  la Section \ref{Sec:Stabilite-laminaire} ci-dessous nous étudions une première propriété relative à la stabilité de la solution $\U$:  si nous considérons $\vu_0\in L^2(\Rt)$ une donnée initiale et $\vu(t,\cdot)$ une solution faible  du problème de Cauchy pour les équations de Navier-Stokes amorties (construite dans le Théorème \ref{Theo:existence-sol-ns-amortie} du chapitre précédent à partir de la donnée initiale $\vu_0\in L^2(\Rt)$)  alors nous voulons étudier la convergence 
\begin{equation}\label{stabilite} 
\lim_{t \longrightarrow +\infty} \Vert \vu(t,\cdot)-\U\Vert_{L^2}=0,
\end{equation} et cette convergence est également appelée la stabilité de la solution $\U$ car nous pouvons observer qu'à partir de n'importe quelle donnée initiale initiale $\vu_0$ alors toute solution faible $\vu(t,\cdot)$  des équations de Navier-Stokes amorties associée à cette donnée initiale converge toujours vers la solution stationnaire $\U$ lorsque le temps $t$ tend vers l'infini. \\
\\ 
Comme mentionné au début de ce chapitre, l'intérêt d'étudier la stabilité de la solution stationnaire $\U$ donnée dans (\ref{stabilite}) est que cette propriété nous permet de donner une réponse  au problème du comportement en temps long des solutions $\vu(t,\cdot)$  que l'on a introduit dans la Section \ref{Sec:motivation}. 
Plus précisément, dans le Théorème \ref{Theo:stabilite_sol_stationnaire_laminaire}  ci-après, nous montrons que si l'on a un contrôle sur le nombre de Grashof $G_{\theta}$ (avec un paramètre $\theta \in [0,3]$ choisi convenablement) alors on a la convergence (\ref{stabilite})  et nous pouvons ainsi observer que  toute solution faible $\vu(t,\cdot)$  des équations de Navier-Stokes amorties se comporte comme une solution stationnaire de ces équations $\U$ dans le régime asymptotique lorsque $t \longrightarrow +\infty$.\\ 
\\ 
Il est aussi important de remarquer que la propriété de stabilité  (\ref{stabilite})  est seulement valable lorsqu'on contrôle le nombre de Grashof $G_{\theta}$ (ce qui caractérise le régime laminaire du mouvement du fluide) et qu'une étude précise du comportement en temps long des solutions $\vu(t,\cdot)$ dans le cadre plus général d'un fluide en régime turbulent (où l'on ne contrôle pas le nombre de Grashof $G_{\theta}$) est une question ouverte (voir la Section \ref{Sec:motivation} pour tous les détails à ce sujet). \\
\\
\\
Dans la Section \ref{Sec:localisation-espace} nous étudions une toute autre  propriété des solutions stationnaires $\U$ qui porte sur la localisation spatiale de ces solutions: rappelons que la force $\fe$ donnée dans le Définition
\ref{Definition_force_ext} est une fonction qui appartient à la classe de Schwartz et donc cette force est bien localisée en variable d'espace; il s'agit alors d'étudier la localisation  spatiale en variable d'espace des solutions $\U$  associées à cette force et dans  ce cadre, dans le Théorème \ref{Theo:dec-U-turb},  nous allons prouver  que toute solution $\U$ des équations (\ref{NS-stationnaire-prop}) a une décroissance 
\begin{equation}\label{dec-espace-turb}
\vert \U(x)\vert \lesssim \frac{1}{\vert x \vert^4},
\end{equation}  lorsque $\vert x \vert$ est suffisamment grand. \\
\\
Soulignons maintenant que pour vérifier l'estimation (\ref{dec-espace-turb}) on n'a pas besoin de faire aucun contrôle sur le nombre de Grashof $G_{\theta}$ et ainsi ce résultat est valable que  dans le régime laminaire mais aussi dans le régime turbulent du fluide et donc il s'agit d'un résultat général sur la localisation spatiale des solutions des équations (\ref{NS-stationnaire-prop}). \\
\\
Expliquons maintenant l'intérêt d'étudier cette localisation spatiale. Dans la Section \ref{Sec:problematique} du chapitre précédent nous avons fait une discussion sur la difficulté de trouver une notion adéquate de longueur caractéristique $L$ (qui représente la plus grande échelle de longueur où l'on veut le comportement turbulent du fluide) lorsqu'on considère un fluide dans tout l'espace $\Rt$. Dans ce cadre, dans la Section \ref{Sec:Grashof} nous avons fixé cette longueur $L$ comme un paramètre du modèle (tout comme la constante de viscosité du fluide) et nous avons construit la force $\fe$ de sorte que cette fonction est essentiellement localisée en variable d'espace sur le cube $[-L,L]^3$ (voir le point $1)$ du Lemme \ref{lemma_loc_fe} et le point $1)$ de la Remarque \ref{rmq:loc-fe} pour tous les détails): l'idée sous-jacente étant que pour une longueur $L>0$ fixe nous voulons étudier le comportement turbulent du fluide dans le cube $[-L,L]^3$ et la force $\fe$ est bien localisée sur ce cube.\\
\\
Nous voulons maintenant savoir comment la solution stationnaire $\U$ associée à la force $\fe$ est localisée sur le cube ci-dessus et pour cela par la décroissance en variable d'espace de cette solution  donnée dans (\ref{dec-espace-turb}) nous pouvons en tirer l'estimation 

\begin{equation*}
 \int_{ ([- L, L]^3)^c}\vert \U(x)\vert dx \lesssim \frac{1}{ L^{3}},
\end{equation*}  où $([- L, L]^3)^c$ dénote l'ensemble complémentaire du cube $[-L,L]^3$, et cette estimation nous donne une mesure plus précise de la façon comment la solution stationnaire $\U$ est localisée sur le cube $[-L,L]^3$.  En effet, nous pouvons observer que le volume de la solution $\U$ en dehors de ce  cube, qui est mesuré par l'intégrale à gauche dans l'estimation ci-dessus, est de l'ordre de $\frac{1}{L^3}$. 

\subsection{Stabilité dans le cadre laminaire}\label{Sec:Stabilite-laminaire} 
Le but de cette section est de montrer la convergence (\ref{stabilite}) et pour cela nous aurons besoin de contrôler  le nombre de Grashof $G_{\theta}$. Ainsi,  dans cette section  nous supposons que nous avons le contrôle  $G_\theta<\eta$, avec un certain paramètre $\theta \in [0,3]$ qui apparaîtra dans les estimations dont on aura besoin;  et où $\eta>0$ est une constante qui ne dépend d'aucun paramètre physique donné dans la Définition \ref{Param_physiques} ni du paramètre d'amortissement $\alpha>0$ des équations (\ref{NS-stationnaire-prop}). \\
\\
Expliquons maintenant  pourquoi nous avons besoin de supposer ce contrôle sur le nombre de nombre de Grashof $G_{\theta}$. Dans le Théorème \ref{Theo:solutions_stationnaires_turbulent}  nous avons montré un résultat général d'existence des solutions $\U \in H^1(\Rt)$ des équations  stationnaires (\ref{NS-stationnaire-prop})  où ces solutions vérifient toujours l'estimation $\ds{\Vert \U \Vert_{H^1} \leq \frac{1}{\min(\nu,\alpha)}\Vert \fe \Vert_{H^{-1}}}$ (voir l'estimation (\ref{estimation-a-priori}) page \pageref{estimation-a-priori}),  néanmoins, pour étudier la convergence (\ref{stabilite}) nous avons besoin d'un contrôle plus précis sur la taille de la solution $\U$ et pour cela, dans le point $A)$ ci-dessous, nous allons tout d'abord montrer que si l'on  contrôle  le nombre de Grashof $G_{\theta}$  alors nous pouvons construire une solution $\U\in H^1(\Rt)$ des équations (\ref{NS-stationnaire-prop}) de sorte que la taille de cette solution peut être contrôlée convenablement par ce nombre $G_{\theta}$ et ceci sera fait dans le Théorème \ref{Prop:Existence_solutions_stationnaires_laminaire} ci-après. Ensuite,  dans le point $B)$, nous allons voir que ce  contrôle sur la taille de la solution $\U$ nous permettra finalement de vérifier la propriété de stabilité de cette solution donnée dans   (\ref{stabilite}) et ceci sera fait dans le Théorème \ref{Theo:stabilite_sol_stationnaire_laminaire} qui, comme l'on a déjà, mentionné est le résultat principal de cette section.

\subsubsection{A) Contrôle sur la solution stationnaire}
Si nous supposons un contrôle sur le nombre de Grashof $G_{\theta}$ nous allons voir que l'on peut utiliser le principe de contraction de Picard pour construire une solution $\U \in H^1(\Rt)$ des équations  (\ref{NS-stationnaire-prop}) telle que la taille de cette solution est contrôlée par ce nombre $G_{\theta}$. Dans ce cadre nous commençons donc par écrire les équations (\ref{NS-stationnaire-prop}) comme un problème de point fixe équivalent (voir l'expression (\ref{point-fixe-stattionnaire})  page \pageref{point-fixe-stattionnaire}): 
$$ \U = - \frac{1}{-\nu \Delta +\alpha I_d}\left[ \P ((\U\cdot \vec{\nabla}) \U)\right] + \frac{1}{-\nu \Delta +\alpha I_d}\left[\fe\right],$$ où, pour mener à bien les
estimations dont on aura besoin il convient de réécrire ce problème comme suit: 
  
\begin{equation}\label{N-S-stationnaire-alpha-pf}
\U =  \frac{-\nu \Delta}{-\nu \Delta +\alpha I_d}\left[\frac{1}{\nu} \P \left(\frac{1}{\Delta} \left((\U\cdot \vec{\nabla}) \U\right)\right)\right] - \frac{-\nu \Delta}{-\nu \Delta +\alpha I_d}\left[\frac{1}{\nu \Delta}( \fe)\right],
\end{equation}   
où l'opérateur  $\ds{\frac{-\nu \Delta}{-\nu \Delta +\alpha I_d}}$ est défini au niveau de Fourier  par le symbole   $\ds{\widehat{m}_{\alpha,\nu}(\xi)=\frac{\nu \vert \xi \vert^2}{\nu \vert \xi \vert^2 +\alpha}}$ et comme $\vert \widehat{m}_{\alpha,\nu}(\xi) \vert \leq 1$  pour tout $\xi \in \Rt$, nous avons que $\ds{\frac{-\nu \Delta}{-\nu \Delta +\alpha I_d}}$ est un opérateur linéaire et  borné dans tous les espaces de Sobolev $\dot{H}^s$ ($s\in \mathbb{R}$):  pour tout $g\in \dot{H}^s(\Rt)$ on a la majoration
\begin{equation}\label{borne_T_alpha_nu}
\left\Vert  \frac{-\nu \Delta}{-\nu \Delta +\alpha I_d}( g )\right\Vert_{\dot{H}^s}\leq \Vert g \Vert_{\dot{H}^s}.
\end{equation}
Nous cherchons à  trouver un espace fonctionnel $E\subset H^1(\Rt)$ dans lequel on puisse vérifier la continuité de la forme bilinéaire  
\begin{equation}\label{cont_forme_bilinaire}
\left\Vert  \frac{-\nu \Delta}{-\nu \Delta +\alpha I_d}\left[\frac{1}{\nu} \P \left(\frac{1}{\Delta} \left((\U\cdot \vec{\nabla}) \U\right)\right)\right] \right\Vert_E \leq \frac{C}{\nu}\Vert \U \Vert_E \Vert \U \Vert_E,
\end{equation}
 pour tout $\U \in  E$, avec $C>0$ une constante indépendante de $\nu$ et $\alpha$  et en plus, dans lequel on puisse relier la condition de petitesse sur le terme $\ds{\frac{-\nu \Delta}{-\nu \Delta +\alpha I_d}\left(\frac{1}{\nu\Delta}\fe\right)}$ par le contrôle sur le nombre de Grashof $G_\theta$:  
\begin{equation}\label{condition_Picard_grashof}
4\frac{C}{\nu} \left\Vert \frac{-\nu \Delta}{-\nu \Delta +\alpha I_d}\left(\frac{1}{\nu\Delta}\fe\right)\right\Vert_E\leq G_\theta,
\end{equation} avec une valeur convenable du paramètre $\theta \in [0,3]$. 
Donc, pour vérifier les points (\ref{cont_forme_bilinaire}) et (\ref{condition_Picard_grashof}) nous  introduisons  l'espace fonctionnel $E$ défini par 
\begin{equation}\label{E}
E=\{ \vec{g} \in H^1(\Rt): div(\vec{g})=0  \}
\end{equation} et  muni de la norme 
\begin{equation}\label{Norme_E}
\Vert \cdot \Vert_E= \frac{1}{L^{\frac{1}{2}}}\Vert \cdot \Vert_{L^2}+\frac{\ell_0}{L^{\frac{1}{2}}}\Vert \cdot \Vert_{\dot{H}^1}+\Vert \cdot \Vert_{L^3}, 
\end{equation} où $L>0$ est la longueur caractéristique du fluide et $\ell_0>0$ est l'échelle d'injection d'énergie. Nous observons que l'on a l'inclusion $E \subset H^1(\Rt)$ et alors:  
\begin{Theoreme}\label{Prop:Existence_solutions_stationnaires_laminaire} Soient $\nu,L,\ell_0$ et $F$ les paramètres physiques donnés dans la Définition \ref{Param_physiques}, soit $\fe$ la force donnée dans la Définition \ref{Definition_force_ext} et soit le nombre de Grashof $\ds{G_{\frac{3}{2}}=\frac{F L^{\frac{3}{2}} \ell^{\frac{3}{2}}_{0}}{\nu^2}}$. Il existe une constante $\eta_1>0$, qui ne dépend pas des paramètres physiques ci-dessus, telle que si $G_{\frac{3}{2}}<\eta_1$ alors il existe $\U\in E$ solution des équations de Navier-Stokes stationnaires et amorties (\ref{NS-stationnaire-prop})  qui est l'unique solution qui vérifie l'estimation 
	\begin{equation}\label{estim:U-Grashof}
	\Vert \U \Vert_E \leq 2C_1\nu G_{\frac{3}{2}},
	\end{equation} où $C_1>0$ est une constante indépendante des paramètres physiques.  \\
\end{Theoreme} 
Si nous comparons ce résultat avec le Théorème \ref{Theo:solutions_stationnaires_turbulent} page \pageref{Theo:solutions_stationnaires_turbulent} nous pouvons observer que si nous supposons un contrôle sur le nombre de Grashof $G_{\frac{3}{2}}$ alors nous pouvons construire une solution $\U$ des équations  (\ref{NS-stationnaire-prop}) qui vérifie en plus l'estimation (\ref{estim:U-Grashof}) et cette estimation va nous permettre d'étudier la stabilité cette solution dans le point $B)$ ci-après. \\
\\
\\
\dm \'Etant donné que $div(\U)=0$  nous écrivons le terme non linéaire $\ds{(\U\cdot \vec{\nabla}) \U}$ comme $\ds{div(\U\otimes \U)}$ et nous  étudions la quantité 
\begin{eqnarray*}
\|\U\|_{E} &\leq & \left\| \frac{-\nu \Delta}{-\nu \Delta +\alpha I_d}\left[ \frac{1}{\nu}\P\left(\frac{1}{\Delta} div(\U\otimes \U)\right) \right]\right\|_{E}+\left\|\frac{-\nu \Delta}{-\nu \Delta +\alpha I_d}\left[ \frac{1}{\nu \Delta}\fe\right]\right\|_{E}.
\end{eqnarray*} 
Nous commençons par  vérifier la continuité du terme bilinéaire sur l'espace $E$ et pour cela  nous avons besoin d'estimer chaque terme de la norme $\Vert \cdot \Vert_E$ donnée par l'expression (\ref{Norme_E}). Pour le premier terme qui compose la norme $\Vert \cdot \Vert_{E}$  
par la continuité de l'opérateur $\ds{\frac{-\nu \Delta}{-\nu \Delta +\alpha I_d}}$  (voir l'inégalité (\ref{borne_T_alpha_nu}))  et par la continuité du projecteur de Leray, nous avons 
$$\frac{1}{L^{\frac{1}{2}}}\left\Vert \frac{-\nu \Delta}{-\nu \Delta +\alpha I_d}\left[ \frac{1}{\nu}\P\left(\frac{1}{\Delta}div(\U\otimes \U)\right) \right]\right\Vert_{L^2}\leq  c_1\frac{1}{L^{\frac{1}{2}}}\left\Vert \frac{1}{\nu} \frac{1}{\Delta}div(\U\otimes \U)\right\Vert_{L^2}\leq
\frac{c_1}{\nu}\frac{1}{L^{\frac{1}{2}}}\left\Vert \U\otimes \U\right\Vert_{\dot{H}^{-1}},$$ 
d'où, par les inégalités de Hardy-Littlewood-Sobolev et l'inégalité de H\"older nous obtenons 
$$\frac{c_1}{\nu}\frac{1}{L^{\frac{1}{2}}}\left\Vert \U\otimes \U\right\Vert_{\dot{H}^{-1}} \leq \frac{c_1}{\nu}\frac{1}{L^{\frac{1}{2}}}\left\Vert \U\otimes \U\right\Vert_{L^{\frac{6}{5}}} \leq \frac{c_1}{\nu}\frac{1}{L^{\frac{1}{2}}} \Vert \U \Vert_{L^2}\Vert \U \Vert_{L^3}.$$ 
Ensuite, pour le deuxième terme de la norme $\Vert \cdot \Vert_{E}$, toujours par la continuité de l'opérateur $\ds{\frac{-\nu \Delta}{-\nu \Delta +\alpha I_d}}$ et du projecteur de Leray, nous avons  
$$ \frac{\ell_0}{L^{\frac{1}{2}}}\left\Vert \frac{-\nu \Delta}{-\nu \Delta +\alpha I_d}\left[ \frac{1}{\nu}\P\left(\frac{1}{\Delta}div(\U\otimes \U)\right) \right]\right\Vert_{\dot{H}^1}\leq \frac{c_2}{\nu} \frac{\ell_0}{L^{\frac{1}{2}}} \left\Vert \frac{1}{\Delta}div(\U\otimes \U) \right\Vert_{\dot{H}^1}\leq \frac{c_2}{\nu} \frac{\ell_0}{L^{\frac{1}{2}}} \left\Vert \U\otimes \U \right\Vert_{L^2},$$ d'où, par l'inégalité de H\"older et les inégalités de Hardy-Littlewood-Sobolev nous pouvons écrire
$$ \frac{c_2}{\nu} \frac{\ell_0}{L^{\frac{1}{2}}}\left\Vert \U\otimes \U \right\Vert_{L^2} \leq \frac{c_2}{\nu} \frac{\ell_0}{L^{\frac{1}{2}}}\Vert \U \Vert_{L^6}\Vert \U \Vert_{L^3}\leq  \frac{c_2}{\nu} \frac{\ell_0}{L^{\frac{1}{2}}}\Vert \U \Vert_{\dot{H}^1}\Vert \U \Vert_{L^3}.$$ 
Finalement, pour le troisième terme de la norme $\Vert \cdot \Vert_E$  par les inégalités de Hardy-Littlewood-Sobolev nous avons 
\begin{eqnarray*}
\left\Vert  \frac{-\nu \Delta}{-\nu \Delta +\alpha I_d}\left[ \frac{1}{\nu}\P\left(\frac{1}{\Delta}div(\U\otimes \U)\right) \right] \right\Vert_{L^3} &\leq & c_3\left\Vert \frac{-\nu \Delta}{-\nu \Delta +\alpha I_d}\left[ \frac{1}{\nu}\P\left(\frac{1}{\Delta}div(\U\otimes \U)\right) \right]\right\Vert_{\dot{H}^{\frac{1}{2}}},
\end{eqnarray*} d'où,
$$ \left\Vert \frac{-\nu \Delta}{-\nu \Delta +\alpha I_d}\left[\frac{1}{\nu}\P\left(\frac{1}{\Delta}div(\U\otimes\U)\right) \right]\right\Vert_{\dot{H}^{\frac{1}{2}}} \leq \frac{c_3}{\nu} \left\Vert \frac{1}{\Delta}div(\U\otimes \U)\right\Vert_{\dot{H}^{\frac{1}{2}}}\leq \frac{c_3}{\nu} \left\Vert \U\otimes \U\right\Vert_{\dot{H}^{-\frac{1}{2}}}.$$ Dans la dernière expression ci-dessus, nous appliquons encore les inégalités de Hardy-Littlewood-Sobolev et l'inégalité de H\"older pour obtenir 
$$ \frac{c_3}{\nu} \left\Vert \U\otimes \U\right\Vert_{\dot{H}^{-\frac{1}{2}}} \leq  \frac{c_3}{\nu} \left\Vert \U\otimes \U\right\Vert_{L^{\frac{3}{2}}}\leq \frac{c_3}{\nu} \left\Vert \U\right\Vert_{L^3}\left\Vert \U\right\Vert_{L^3}.$$
De cette façon, la propriété de continuité  du terme bilinéaire sur l'espace $E$ énoncée dans (\ref{cont_forme_bilinaire}) est bien vérifiée par les estimations ci-dessus, où on pose la constante 
\begin{equation}\label{const_bilin_forme}
C=\frac{1}{\nu}\max(c_1,c_2,c_3),
\end{equation}  
qui, comme nous pouvons observer, ne dépend d'aucun paramètre physique.\\   
\\ \vspace{5mm}
Vérifions maintenant le contrôle sur le terme $\ds{\left\Vert \frac{-\nu \Delta}{-\nu \Delta +\alpha I_d}\left[   \frac{1}{\nu\Delta}\fe\right]\right\Vert_{E}}$ donné dans (\ref{condition_Picard_grashof}) et pour cela nous avons besoin tout d'abord de vérifier  l'estimation 
 \begin{equation}\label{estim-C1}
 \left\Vert \frac{-\nu \Delta}{-\nu \Delta +\alpha I_d}\left[   \frac{1}{\nu\Delta}\fe\right]\right\Vert_{E}\leq \frac{C_1}{\nu}\Vert (-\Delta)^{-\frac{3}{4}}\fe \Vert_{L^2}.
 \end{equation} Nous allons estimer chaque terme de la norme   $\Vert \cdot \Vert_{E}$ donnée par l'expression (\ref{Norme_E}). Pour le premier terme de la norme $\Vert \cdot \Vert_E$, par la continuité de l'opérateur  $\ds{\frac{-\nu \Delta}{-\nu \Delta +\alpha I_d}}$ nous pouvons écrire 
$$ \frac{1}{L^{\frac{1}{2}}}\left\Vert \frac{-\nu \Delta}{-\nu \Delta +\alpha I_d}\left[   \frac{1}{\nu\Delta}\fe\right]\right\Vert_{L^2} \leq C_1 \frac{1}{L^{\frac{1}{2}}}\left\Vert \frac{1}{\nu\Delta}\fe \right\Vert_{L^2},$$ mais, par le point $1)$ du Lemme \ref{lemma_loc_fe}, page \pageref{lemma_loc_fe}, nous savons que la fonction $\ds{-\frac{1}{\nu\Delta}\fe}$ est essentiellement concentrée sur le cube $[-L,L]^3$ et alors  nous introduisons la fonction de troncature suivante: soit $\varphi \in \mathcal{C}^{\infty}_{0}(\Rt)$ telle que $\varphi$  est égale à $1$ sur la boule $\vert x \vert \leq 4$ et qui s'annule  en dehors de la boule $\vert x \vert\leq 8$ et nous considérons sa dilatation à l'échelle  $\frac{1}{L}$, notée $\varphi_L$, de sorte que nous avons  $\varphi_L\left(x\right)=1$ pour $\vert x \vert\leq 4L$  et  $\varphi_L\left(x\right)=0$ pour $\vert x \vert\geq8L$. De cette façon,   nous pouvons écrire
\begin{equation}\label{Estimate1} 
\left\Vert -\frac{1}{\nu \Delta}\fe\right\Vert_{L^2}\leq C_1 \left\Vert \varphi_L \left(-\frac{1}{\nu\Delta}\fe\right)\right\Vert_{L^2},
\end{equation} avec $C_1>0$ une constante qui ne dépend d'aucun paramètre. Donc, par l'inégalité de H\"older nous avons 
\begin{eqnarray*}
\left\Vert \varphi_L \left(-\frac{1}{\nu\Delta}\fe\right)\right\Vert_{L^2}&\leq  &\left\Vert \varphi_L\right\Vert_{L^6}\left\Vert-\frac{1}{\nu \Delta}\fe \right\Vert_{L^3} =L^{\frac{1}{2}}\Vert \varphi \Vert_{L^6}\left\Vert -\frac{1}{\nu \Delta}\fe \right\Vert_{L^3},
\end{eqnarray*}
 et alors nous obtenons 
 $$ \frac{1}{L^{\frac{1}{2}}}\left\Vert \frac{-\nu \Delta}{-\nu \Delta +\alpha I_d}\left[   \frac{1}{\nu\Delta}\fe\right]\right\Vert_{L^2} \leq C_1 \left\Vert -\frac{1}{\nu \Delta}\fe \right\Vert_{L^3}.$$
D'autre part, par les   inégalités de Hardy-Littlewood-Sobolev nous avons $\ds{\left\Vert -\frac{1}{\nu \Delta}\fe \right\Vert_{L^3}\leq \frac{c}{\nu}\Vert (-\Delta)^{-\frac{3}{4}}\fe \Vert_{L^2}}$, d'où, par les estimations ci-dessus nous écrivons 
\begin{equation}\label{Estimate2}
\frac{1}{L^{\frac{1}{2}}}\left\Vert -\frac{1}{\nu\Delta}\fe\right\Vert_{L^2}\leq  \frac{C_1}{\nu} \Vert (-\Delta)^{-\frac{3}{4}}\fe \Vert_{L^2}. 
\end{equation}
Pour le deuxième terme de la norme $\Vert \cdot \Vert_E$, toujours par la continuité de l'opérateur  $\ds{\frac{-\nu \Delta}{-\nu \Delta +\alpha I_d}}$ nous avons
$$ \frac{\ell_0}{L^{\frac{1}{2}}}\left\Vert \frac{-\nu \Delta}{-\nu \Delta +\alpha I_d}\left[   \frac{1}{\nu\Delta}\fe\right]\right\Vert_{\dot{H}^{1}} \leq C_1 \frac{\ell_0}{L^{\frac{1}{2}}}\left\Vert \frac{1}{\nu\Delta}\fe \right\Vert_{\dot{H}^1},$$ et en utilisant l'identité de Plancherel nous écrivons $ \ds{\left\Vert -\frac{1}{\nu \Delta}\fe\right\Vert_{\dot{H}^1}=\left\||\xi|\left(\widehat{-\frac{1}{\nu\Delta}\fe}\right)\right\|_{L^{2}}}$, mais par le point $2)$ du Lemme  \ref{lemma_loc_fe} nous savons  que $supp\,\left(\widehat{-\frac{1}{\nu\Delta}\fe}\right)\subset \left\{ \xi \in \Rt: \frac{\rho_1}{\ell_0} \leq \vert \xi \vert \leq  \frac{\rho_2}{\ell_0}\right\},$ 
et alors nous pouvons écrire
$$\left\Vert -\frac{1}{\nu \Delta}\fe\right\Vert_{\dot{H}^1}=\left\||\xi|\left(\widehat{-\frac{1}{\nu\Delta}\fe}\right)\right\|_{L^{2}}\leq \frac{\rho_{2}}{\ell_{0}}\left\|\left(\widehat{-\frac{1}{\nu\Delta}\fe}\right)\right\|_{L^{2}}=\frac{\rho_{2}}{\ell_{0}}\left\|-\frac{1}{\nu\Delta}\fe\right\|_{L^{2}},$$
ce qui équivaut à l'inégalité 
$$ \frac{\ell_{0}}{\rho_{2}}\left\Vert -\frac{1}{\nu \Delta}\fe\right\Vert_{\dot{H}^1}\leq\left\|-\frac{1}{\nu\Delta}\fe\right\|_{L^{2}},$$
et en suivant les mêmes estimations réalisées dans  (\ref{Estimate1}) et (\ref{Estimate2}), nous obtenons 
\begin{equation*}
\frac{\ell_0}{L^{\frac{1}{2}}}\left\Vert -\frac{1}{\nu \Delta}\fe\right\Vert_{\dot{H}^1}\leq  \frac{C_1}{\nu} \Vert (-\Delta)^{-\frac{3}{4}}\fe \Vert_{L^2}.
\end{equation*}  \'Etudions maintenant le troisième terme de la norme $\Vert \cdot \Vert_E$. Par les inégalités de Hardy-Littlewood-Sobolev et par la continuité de l'opérateur $\ds{\frac{-\nu \Delta}{-\nu \Delta +\alpha I_d}}$ nous écrivons 
\begin{eqnarray*}
\left\Vert \frac{-\nu \Delta}{-\nu \Delta +\alpha I_d}\left[ -\frac{1}{\nu \Delta}\fe\right]\right\Vert_{L^3}& \leq & C_1 
\left\Vert \frac{-\nu \Delta}{-\nu \Delta +\alpha I_d}\left[ -\frac{1}{\nu \Delta}\fe\right]\right\Vert_{\dot{H}^{\frac{1}{2}}}\leq  C_1\left\Vert 
 \frac{1}{\nu \Delta}\fe \right\Vert_{\dot{H}^{\frac{1}{2}}} \leq \frac{C_1}{\nu} \Vert (-\Delta)^{\frac{3}{4}}\fe \Vert_{L^2}. \\
\end{eqnarray*} 
De cette façon, par l'inégalité ci-dessus nous obtenons l'estimation  (\ref{estim-C1}). \\
\\
Une fois que nous disposons de cette estimation nous allons maintenant contrôler la quantité  $\ds{ C_1\Vert (-\Delta)^{-\frac{3}{4}}\fe \Vert_{L^2}}$. En effet,  par  la Proposition \ref{Prop:relation-force-grashof}  (en prenant $s=-\frac{3}{4}$ et $p=2$)  nous avons  
$$ \Vert (-\Delta)^{\frac{3}{4}}\fe \Vert_{L^2} \approx F L^p \ell^{-2s}_{0}= F L^{\frac{3}{2}}\ell^{\frac{3}{2}}_{0},$$ et alors, dans l'inégalité précédente nous obtenons  
$$ \left\Vert \frac{-\nu \Delta}{-\nu \Delta +\alpha I_d}\left[ -\frac{1}{\nu \Delta}\fe\right]\right\Vert_{E} \leq \frac{c}{\nu} F L^{\frac{3}{2}}\ell^{\frac{3}{2}}_{0}= C_1 \,\nu \left( \frac{F L^{\frac{3}{2}}\ell^{\frac{3}{2}}_{0}}{\nu^2}\right).$$ Nous observons que le nombre de Grashof $\ds{G_{\frac{3}{2}}=\frac{F L^{\frac{3}{2}}\ell^{\frac{3}{2}}_{0}}{\nu^2}}$ apparaît naturellement dans cette estimation et alors nous pouvons écrire 
$$ \left\Vert \frac{-\nu \Delta}{-\nu \Delta +\alpha I_d}\left[ -\frac{1}{\nu \Delta}\fe\right]\right\Vert_{E} \leq C_1\,\nu G_{\frac{3}{2}},$$
où la constante $C_1$ ne dépend d'aucun paramètre. \\ 
\\
Une fois que nous avons  les estimations $\ds{\left\Vert  \frac{-\nu \Delta}{-\nu \Delta +\alpha I_d}\left[\frac{1}{\nu} \P \left(\frac{1}{\Delta} \left((\U\cdot \vec{\nabla}) \U\right)\right)\right] \right\Vert_E \leq \frac{C}{\nu}\Vert \U \Vert_E \Vert \U \Vert_E,
}$ (où la constante $C>0$ est donnée dans (\ref{const_bilin_forme})) et  $ \ds{\left\Vert \frac{-\nu \Delta}{-\nu \Delta +\alpha I_d}\left[ -\frac{1}{\nu \Delta}\fe\right]\right\Vert_{E}\leq C_1 \nu G_{\frac{3}{2}}}$,  afin d'obtenir une solution $\U \in E$ des équations  (\ref{N-S-stationnaire-alpha-pf}), nous devons vérifier la condition $4\left(\frac{C}{\nu} \right) C_{1}\nu G_{\frac{3}{2}}<1,$ qui équivaut à l'inégalité
\begin{equation}\label{eta_1}
G_{\frac{3}{2}}<\frac{1}{4C C_1}=\eta_1.
\end{equation}
Nous observons que la dépendance du paramètre $\nu$ a disparu dans cette inégalité et que la constante $\eta_1$ ci-dessus ne dépend d'aucun paramètre physique ni du paramètre d'amortissement $\alpha>0$. De cette façon,  on obtient une fonction  $\U\in E$ qui est solution de  (\ref{N-S-stationnaire-alpha-pf}) et qui est l'unique solution dans la boule  $\Vert \U \Vert_E\leq 2C_1\nu G_{\frac{3}{2}}$. \finpv

\subsubsection{B) Stabilité de la solution stationnaire} 

Une fois que nous avons construit une solution $\U \in E \subset H^1(\Rt)$  qui vérifie l'estimation (\ref{estim:U-Grashof}) nous pouvons étudier maintenant la stabilité de cette solution.

\begin{Theoreme}\label{Theo:stabilite_sol_stationnaire_laminaire}  Soient $\nu,L,\ell_0$ et $F$ les paramètres physiques donnés dans la Définition \ref{Param_physiques}, soit $\fe$ la force donnée dans la Définition \ref{Definition_force_ext} et soit le nombre de Grashof $\ds{G_{\frac{3}{2}}=\frac{F L^{\frac{3}{2}} \ell^{\frac{3}{2}}_{0}}{\nu^2}}$. Soit le sous-espace $E\subset H^1(\Rt)$ donné dans (\ref{E}) et soit  $\U \in E$ 
 l'unique solution  des équations de Navier-Stokes amorties et stationnaires,
 $$  -\nu \Delta \U +(\U\cdot \vec{\nabla})\U +\vec{\nabla}P  =\fe-\alpha \U, \qquad div(\U)=0, \quad \alpha>0,$$ obtenue par le biais du Théorème \ref{Prop:Existence_solutions_stationnaires_laminaire}. \\
 \\
 D'autre part,  soit $\vu_0\in L^2(\Rt)$ une donnée initiale à divergence nulle et soit $\vu\in L^{\infty}_{t}L^{2}_{x}\cap (L^{2}_{t})_{loc}\dot{H}^{1}_{x}$  une solution faible  du problème de Cauchy  
 $$ \partial_t\vu = \nu\Delta \vu-(\vu\cdot\vec{\nabla}) \vu -\vec{\nabla}p +\fe-\alpha \vu, \quad div(\vu)=0,\quad  \vu(0,\cdot)=\vu_0,$$ construite dans le Théorème \ref{Theo:existence-sol-ns-amortie}.\\
 \\
  Si  $G_{\frac{3}{2}}< \frac{1}{2C_1}$, où $C_1>0$ est la constante donnée dans (\ref{estim-C1})  et qui ne dépend d'aucun paramètre physique   alors  on a  
 \begin{equation}\label{estatibilite}
\lim_{t\longrightarrow +\infty}\Vert \vu(t,\cdot)-\U\Vert_{L^2}=0.
\end{equation}
\end{Theoreme} 
Nous observons ainsi que si l'on a le contrôle sur le nombre de Grashof: $G_{\frac{3}{2}}< \frac{1}{2C_1}$, alors à partir de \emph{n'importe quelle} donnée initiale $\vu_0\in L^2(\Rt)$ on a que l'évolution au cours du temps du champ de vitesse $\vu(t,\cdot)$ associé à cette donnée initiale  converge toujours vers la solution stationnaire $\U$ lorsque le temps tend vers l'infini; et l'on obtient ainsi une description précise du comportement en temps long des solutions $\vu(t,\cdot)$ dans le cadre d'un fluide en régime laminaire.  \\
\\
\dm Soit $\vu\in L^{\infty}_{t}L^{2}_{x}\cap (L^{2}_{t})_{loc}\dot{H}^{1}_{x}$ une solution de faible des equations de Navier-Stokes amorties construite à partir d'une donnée initiale $\vu_0\in L^2(\Rt)$ et $\U \in E$ la solution stationnaires de ces équations.   On pose $$\vv(t,x)=\vu(t,x)-\U(x),$$ où $\vv \in L^{\infty}_{t}L^{2}_{x}\cap (L^{2}_{t})_{loc}\dot{H}^{1}_{x}$, et l'on cherche à montrer que $\ds{\lim_{t\longrightarrow +\infty}\Vert \vv(t,\cdot)\Vert_{L^2}=0}$. Pour cela  nous avons besoin de vérifier que la fonction $\vv $ satisfait l'inégalité d'énergie suivante:
\begin{Lemme}\label{Lemme:ineq_ener_stabilite} Pour tout  $t>0$ on a:
\begin{eqnarray}\label{theo_stat_ineq_6}\nonumber
\Vert \vv(t,\cdot)\Vert^{2}_{L^2}& \leq & \Vert \vu_0 -\U \Vert^{2}_{L^2} -2\nu \int_{0}^{t}\Vert \vv(s,\cdot)\Vert^{2}_{\dot{H}^{1}}ds-2\int_{0}^{t}\langle (\vv \cdot \vec{\nabla})\U , \vv \rangle_{\dot{H}^{-1}\times \dot{H}^1}ds\\
& & -2\alpha \int_{0}^{t} \Vert \vv(s,\cdot)\Vert^{2}_{L^2}ds. 
\end{eqnarray}
\end{Lemme} 
La preuve de ce lemme technique repose essentiellement sur l'inégalité d'énergie vérifiée par la solution $\vu$ et qui a été étudiée dans la Proposition \ref{Proposition_ineq_energie_alpha_modele} page  \pageref{Proposition_ineq_energie_alpha_modele};   et cette preuve  sera faite  à la fin du chapitre page \pageref{sec:preuve-lemme-tech-1}. Ainsi, nous supposons pour l'instant que l'inégalité d'énergie ci-dessus est vraie et dans le  troisième terme à droite de cette  inégalité,  en utilisant le fait  $div(\vv)=0$ et $div(\U)=0$, nous pouvons écrire 
$$ -2\int_{0}^{t}\langle (\vv \cdot \vec{\nabla})\U , \vv \rangle_{\dot{H}^{-1}\times \dot{H}^1}ds=-2 \int_{0}^{t}\langle (div(\vv \otimes \U) , \vv \rangle_{\dot{H}^{-1}\times \dot{H}^1}ds,$$ et comme $\U\in L^3(\Rt)$, par l'inégalité de Cauchy-Schwarz, l'inégalité de H\"older et les inégalités de Hardy-Littlewood-Sobolev   nous écrivons  
\begin{eqnarray*}
& & -2\int_{0}^{t}\langle (div(\vv \otimes \U) , \vv \rangle_{\dot{H}^{-1}\times \dot{H}^1}ds \leq 2\int_{0}^{t} \Vert div( \vv \otimes \U)\Vert_{\dot{H}^{-1}} \Vert \vv(s,\cdot)\Vert_{\dot{H}^1}ds\\
& \leq & 2 \int_{0}^{t} \Vert  \vv \otimes \U\Vert_{L^{2}}\Vert \vv(s,\cdot)\Vert_{\dot{H}^1}ds \leq 2 \int_{0}^{t} \Vert  \vv (s,\cdot)\Vert_{L^6}\Vert \U \Vert_{L^3} \Vert \vv(s,\cdot)\Vert_{\dot{H}^1}ds \\
& \leq & 2 \Vert \U \Vert_{L^3}\int_{0}^{t}\Vert \vv (s,\cdot)\Vert^{2}_{\dot{H}^{1}}ds. 
\end{eqnarray*} Alors, en remplaçant cette  dernière estimation dans (\ref{theo_stat_ineq_6}) nous obtenons 
\begin{eqnarray}\label{theo_stat_ineq_7}
\Vert \vv(t,\cdot)\Vert^{2}_{L^2}& \leq & \Vert \vu_0 -\U \Vert^{2}_{L^2} -2(\nu -\Vert \U \Vert_{L^3})\int_{0}^{t}\Vert \vv(s,\cdot)\Vert^{2}_{\dot{H}^{1}}ds -2\alpha \int_{0}^{t} \Vert \vv(s,\cdot)\Vert^{2}_{L^2}ds. 
\end{eqnarray} 
Mais,  par la définition de la norme $\Vert \cdot \Vert_E$ donnée par l'expression (\ref{Norme_E}) nous savons que $\ds{\Vert \U \Vert_{L^3}\leq \Vert \U \Vert_E}$ et de plus,  par le Théorème  \ref{Prop:Existence_solutions_stationnaires_laminaire} nous avons l'estimation $\ds{\Vert \U \Vert_E\leq 2C_1\nu G_{\frac{3}{2}}}$ et donc  nous obtenons  $\ds{\Vert \U \Vert_{L^3}\leq 2C_1\nu G_{\frac{3}{2}}}$. \\
\\
De cette façon, si nous supposons  le contrôle sur le  nombre de Grashof: $\ds{G_{\frac{3}{2}}<\frac{1}{2 C_1}}$, alors dans l'estimation précédente de la quantité $\Vert \U \Vert_{L^3}$  nous obtenons  $\ds{ \Vert \U \Vert_{L^3} < \nu}$; et alors nous avons que la quantité $\ds{-2(\nu -\Vert \U \Vert_{L^3})\int_{0}^{t}\Vert \vv(s,\cdot)\Vert^{2}_{\dot{H}^{1}}ds}$ est une quantité négative de sorte que, en revenant à l'estimation (\ref{theo_stat_ineq_7}),  nous pouvons écrire 
$$ \Vert \vv(t,\cdot)\Vert^{2}_{L^2} \leq  \Vert \vu_0 -\U \Vert^{2}_{L^2} -2\alpha \int_{0}^{t} \Vert \vv(s,\cdot)\Vert^{2}_{L^2}ds. $$ 
Finalement, nous appliquons le lemme de Gr\"onwall   pour obtenir l'estimation 
$$ \Vert \vv(t,\cdot)\Vert^{2}_{L^2} \leq  \Vert \vu_0 -\U \Vert^{2}_{L^2} e^{-2\alpha t},$$ d'où nous pouvons en tirer que 
$ \ds{\lim_{t\longrightarrow +\infty}\Vert \vv(t,\cdot)\Vert_{L^2}=0}$. \finpv\\
Nous pouvons observer que l'essentiel de cette preuve du Théorème \ref{Theo:stabilite_sol_stationnaire_laminaire}  repose sur le contrôle de  la taille de la solution stationnaire $\U$: $\ds{\Vert \U \Vert_{E}\leq 2C_1\nu G_{\frac{3}{2}}}$; et ce contrôle repose finalement sur un contrôle sur le nombre de Grashof $G_{\frac{3}{2}}$. \\
\\
Dans ce cadre, nous observons que l'étude de la stabilité de la solution stationnaire $\U$ que nous venons de faire dans le Théorème   \ref{Theo:stabilite_sol_stationnaire_laminaire} est seulement valable dans le cadre d'un régime laminaire et une étude analogue de cette propriété de stabilité  dans un cadre général où l'on ne contrôle pas le nombre de Grashof $G_{\frac{3}{2}}$ semble actuellement hors de portée. \\
\\
Dans la section qui suit nous étudions une toute autre propriété des solutions stationnaires $\U$ et
nous allons voir que, contrairement à la propriété de stabilité, la  localisation spatiale de la solution $\U$, valable dans le cadre général où l'on ne fait  aucun contrôle sur le nombre de Grashof $G_{\frac{3}{2}}$. 

\subsection{Localisation en variable d'espace dans le cadre turbulent}\label{Sec:localisation-espace} 
Dans cette section nous considérons $\U \in H^1(\Rt)$ une solution des équations de Navier-Stokes amorties et stationnaires (\ref{NS-stationnaire-prop}) où $\fe$ est toujours la force donnée dans la Définition \ref{Definition_force_ext}; et nous allons étudier ici la localisation spatiale de la solution $\U$ obtenue dans le Théorème \ref{Theo:solutions_stationnaires_turbulent}. Plus  précisément, la force $\fe$ étant une fonction bien localisée en variable d'espace (car $\fe$ appartient à la classe de Schwartz) il s'agit d'étudier  une estimation du type 
\begin{equation}\label{dec-U-motiv}
\vert \U(x)\vert \lesssim \frac{1}{\vert x \vert^\beta},
\end{equation} pour un certain paramètre $\beta>0$  et pour $\vert x \vert$ suffisamment grand. \\
\\
Lorsqu'on veut étudier cette décroissance en variable d'espace,  il est naturel de se poser la question de quel est le plus grand paramètre $\beta>0$  que l'on peut espérer dans (\ref{dec-U-motiv}) et dans ce cadre, pour motiver le résultat obtenu dans le Théorème \ref{Theo:dec-U-turb} ci-dessous  nous allons maintenant expliquer pourquoi on s'attend à ce que la solution $\U$ n'ait pas  une meilleure  décroissance à l'infini que $\ds{\frac{1}{\vert x \vert^4}}$, même si la force $\fe$ appartient à la classe de Schwartz; et pour expliquer ce fait nous avons besoin de rappeler (rapidement) quelques résultats classiques sur l'étude de la décroissance (\ref{dec-U-motiv}) qui ont été obtenus dans le cadre des équations de Navier-Stokes stationnaires sans terme d'amortissement (lorsque $\alpha=0$).\\
\\
Considérons donc pour l'instant le paramètre d'amortissement $\alpha=0$ et les équations de Navier-Stokes stationnaires classiques
\begin{equation}\label{N-S-dec}
-\nu \Delta \U +(\U\cdot \vec{\nabla})\U +\vec{\nabla}P  =\fe, \qquad div(\U)=0.
\end{equation} Si l'on suppose en plus que la force $\fe$ vérifie la condition de petitesse suivante: pour $a \in \mathbb{N}^{3}$ un multi-indice et $\eta>0$ une constante (suffisamment petite),
\begin{equation}\label{cond-force}
\sup_{\vert a \vert \leq 2}\sup_{x\in \Rt} (1+\vert x \vert)^4\vert \partial^{a}\fe(x)\vert \leq \eta \nu^2,
\end{equation} alors par le Théorème $4.10$ du livre \cite{PGLR1} nous avons qu'il existe $(\U,P)\in \mathcal{C}^2(\Rt)$ une solution (classique) des équations (\ref{N-S-dec}) qui vérifie 
\begin{equation}\label{cond-sol-estat-classique}
 \sup_{\vert a \vert \leq 2}\sup_{x\in \Rt} (1+\vert x \vert)\vert \partial^{a}\U(x)\vert <+\infty,
\end{equation}
et alors nous avons le résultat classique suivant:
 \begin{Proposition} Soit $\U \in  \mathcal{C}^2(\Rt) $ la solution classique des équations de Navier-Stokes (\ref{N-S-dec})  où la force $\fe$ vérifie la condition de petitesse (\ref{cond-force}). Alors cette solution ne peut pas avoir une meilleure décroissance que   $\ds{\frac{1}{\vert x \vert^4}}$.
 	\end{Proposition}
 \pv   Supposons pour l'instant que la solution $\U$ a une décroissance à l'infini  $\vert \U(x)\vert \lesssim \frac{1}{\vert x \vert^{\beta} }$, avec $\beta>4$, et on va obtenir  une contradiction.  Si nous considérons la solution stationnaire $\U \in \mathcal{C}^{2}(\Rt)$ ci-dessus comme étant la donnée initiale du problème de Cauchy 
 \begin{equation}\label{cauchy-n-s-dec-motiv}
 \partial_t \vu + (\vu \cdot \vec{\nabla})\vu -\nu \Delta\vu +\vec{\nabla}p =\fe, \quad div(\vu)=0, \quad \vu(0,\cdot)=\U,
 \end{equation} et comme $\U$ vérifie la propriété (\ref{cond-sol-estat-classique}) alors par le Théorème $4.7$ du livre \cite{PGLR1} nous avons qu'il existe un temps $T_0>0$ (qui dépend de $\U$ et $\fe$) et une fonction $\vu \in \mathcal{C}([0,T_0[, \mathcal{C}^2(\Rt))$ qui est l'unique solution du problème de Cauchy (\ref{cauchy-n-s-dec-motiv}). De plus, comme l'on a supposé que $\vert \U(x)\vert \lesssim \frac{1}{\vert x \vert^{\beta} }$, avec $\beta>4$ nous avons  $\ds{\lim_{\vert x \vert \longrightarrow+\infty}\vert x \vert^4 \vert \U(x)\vert=0}$ et alors nous pouvons appliquer le Théorème $4.12$
 du livre \cite{PGLR1} pour obtenir que la solution $\vu$ ne peut décroître à l'infini plus rapidement que $\frac{1}{\vert x \vert^4}$. Mais,  observons que la solution stationnaire $\U$ est aussi une solution du problème de Cauchy (\ref{cauchy-n-s-dec-motiv}) et alors, par l'unicité de la solution $\vu$ nous avons $\vu=\U$ et donc la solution stationnaire $\U$ ne décroît pas à l'infini plus rapidement que $\ds{\frac{1}{\vert x \vert^4}}$. \finpv 
\\ 
Nous observons ainsi que si l'on suppose une condition de petitesse sur la force $\fe$ et qui est donnée dans (\ref{cond-force}), alors les solutions classiques des équations (\ref{N-S-dec}) ne peuvent pas avoir une meilleure décroissance à l'infini que $\frac{1}{\vert x \vert^4}$. Ce fait est un résultat classique (voir toujours le Chapitre $4$ du livre \cite{PGLR1} pour plus de détails) et dans le résultat suivant nous étudions cette décroissance précise à l'infini dans le cadre des équations stationnaires et amorties (\ref{NS-stationnaire-prop}).\\ 
\\
Ainsi, nous allons montrer que le terme d'amortissement $-\alpha \U$ introduit dans ces équations entraîne une décroissance à l'infini du type $\ds{\frac{1}{\vert x \vert^4}}$ pour toute solution faible des équations (\ref{NS-stationnaire-prop}) et sans aucune condition de petitesse sur la force $\fe$ et donc  sans aucun contrôle sur le nombre de Grashof $G_{\frac{3}{2}}$ et alors le résultat que nous allons montrer ci-dessous est valable que ce soit dans le régime laminaire et aussi dans le régime turbulent du fluide. \\
\begin{Theoreme}\label{Theo:dec-U-turb} Soit $\fe$ la force donnée dans la Définition \ref{Definition_force_ext}. Alors toute solution $\U \in H^{1}(\Rt)$  des équations de Navier-Stokes amorties et stationnaires  (\ref{NS-stationnaire-prop}) vérifie la décroissance suivante:  
	\begin{equation}\label{dec-theo}
	\left\vert \U(x)\right\vert \leq  \frac{c}{1+\vert x \vert^4}, 
	\end{equation}  avec une constante $c=c(\U,\fe,\nu,\alpha)>0$  qui dépend de la solution $\U$, de la force $\fe$, de la constante de viscosité $\nu>0$ et du paramètre d'amortissement $\alpha>0$. 
	\end{Theoreme} 
\textbf{Démonstration.}  La preuve de ce résultat repose essentiellement sur une étude préliminaire de la localisation spatiale de la solution $\U$ qui est faite dans la proposition suivante:

\begin{Proposition}\label{Prop:dec-espace-CKN} Dans 
	le cadre du Théorème \ref{Theo:dec-U-turb}, la solution $\U\in H^1(\Rt)$ des équations (\ref{NS-stationnaire-prop})  vérifie la décroissance en variable d'espace: pour tout $x\in \Rt$, 
	\begin{equation}\label{estim-CKN}
	\vert \U(x)\vert \leq c \frac{1}{1+\vert x\vert },
	\end{equation} où $c=c(\U,\nu,\alpha)>0$ est une constante qui dépend de la solution $\U$, de la constante de viscosité $\nu>0$ et du paramètre d'amortissement $\alpha>0$. 
\end{Proposition} 
Pour prouver cette proposition nous avons besoin de passer par le cadre de la théorie de la régularité partielle développée principalement dans l'article \cite{CKN} de Caffarelli-Kohn et Niremberg (voir aussi les Chapitres $13$ et $14$ du livre \cite{PGLR1} pour plus de détails sur cette théorie) et ceci sera fait en  détail dans l'appendice à la fin du chapitre. Plus précisément, dans le Théorème \ref{Theo:CKN} page \pageref{Theo:CKN}  nous allons remarquer que cette théorie de la régularité s'adapte facilement au cadre des équations de Navier-Stokes amorties  et à partir de ce résultat   nous pourrons  vérifier la décroissance de la solution $\U$ donnée dans  (\ref{estim-CKN}).  Ainsi, nous allons maintenant supposer la Proposition \ref{Prop:dec-espace-CKN} et nous allons voir comment à partir de ce résultat nous pouvons démontrer le Théorème \ref{Theo:dec-U-turb}. \\ 
\\
 Soit donc $\U \in H^1(\Rt)$ une solution des équations (\ref{NS-stationnaire-prop}).  Pour étudier la décroissance spatiale (\ref{dec-theo}) on commence par écrire cette solution comme la solution du problème de point fixe équivalent (voir l'expression (\ref{point-fixe-stattionnaire}) page \pageref{point-fixe-stattionnaire}): 

\begin{equation}\label{point-fixe-stattionnaire-aux}
\U = - \frac{1}{-\nu \Delta +\alpha I_d}\left[ \P ((\U\cdot \vec{\nabla}) \U)\right] + \frac{1}{-\nu \Delta +\alpha I_d}\left[\fe\right],
\end{equation} et nous allons étudier la décroissance spatiale des deux termes ci-dessus. Plus précisément, nous allons montrer qu'il existe deux constantes  $c_1=c_1(\U,\nu, \alpha)>0$  et $c_2=c_2(\fe,\nu,\alpha)>0$ telles que l'on a 
\begin{equation}\label{estimation1}
\left\vert- \frac{1}{-\nu \Delta +\alpha I_d}\left[ \P ((\U\cdot \vec{\nabla}) \U)\right](x)\right\vert \leq \frac{c_1}{\vert x \vert^4}
\end{equation} et 
\begin{equation}\label{estimation2}
\left\vert   \frac{1}{-\nu \Delta +\alpha I_d}\left[\fe\right](x) \right\vert \leq \frac{c_2}{\vert x \vert^4}
\end{equation} pour tout $\vert x \vert > 8$. \\
\\
Nous allons tout d'abord étudier l'estimation (\ref{estimation2})  car les estimations que nous allons faire avec ce terme nous seront aussi utiles pour étudier après l'estimation (\ref{estimation1}).  	\\
\\
L'estimation (\ref{estimation2}) sera une conséquence  du lemme suivant:
\begin{Lemme}\label{Lemme:bessel} Soit $g \in L^1(\Rt)$ qui vérifie une décroissance $\ds{\vert g(x)\vert \leq \frac{c}{\vert x \vert^n}}$, pour $n \in \mathbb{N}$  et pour tout $\vert x \vert >8$.  Alors il existe une constante $c_2=c_2(g,\nu,\alpha)>0$ telle que l'on a: pour tout $\vert x \vert >8$, 
	\begin{equation}\label{loc-espace-Bessel}
	\left\vert \frac{1}{-\nu \Delta +\alpha I_d}\left[g\right](x) \right\vert \leq \frac{c_2}{\vert x \vert^n}. 
	\end{equation}
	\end{Lemme}
La preuve de ce lemme sera faite à la fin du chapitre page \pageref{sec:preuve-lemme-bessel} et elle repose essentiellement sur le fait que l'on peut écrire $$ \frac{1}{-\nu \Delta +\alpha I_d}\left[g\right](x) = G_{\nu,\alpha}\ast g (x),$$ où $G_{\nu,\alpha}(x)>0$ est le noyau du  potentiel de Bessel $\ds{\frac{1}{-\nu \Delta +\alpha I_d}\left[\cdot\right]}$; et ce noyau vérifie de bonnes propriétés de  décroissance. \\
\\
Une fois que nous avons énoncé le Lemme \ref{Lemme:bessel} nous pouvons maintenant en tirer l'estimation (\ref{estimation2}) directement de ce lemme. En effet, rappelons que la force $\fe=(f_1,f_2,f_3)$ donnée dans la Définition \ref{Definition_force_ext} est un champ de vecteurs dans la classe de Schwartz et alors par le Lemme  \ref{Lemme:bessel} nous avons que la fonction $\ds{\frac{1}{-\nu \Delta +\alpha I_d}\left[ f_i\right]}$ appartient  aussi à la classe de Schwartz (avec $1\leq i \leq 3$ ) et donc en prenant $n=4$ dans l'estimation (\ref{loc-espace-Bessel}) nous avons l'estimation (\ref{estimation2}). \\
\\
\\ 
Une fois que nous disposons de l'estimation (\ref{estimation2}) nous pouvons maintenant vérifier l'estimation (\ref{estimation1}) et pour mener à bien les estimations dont on a besoin on commence par écrire le terme $$ \ds{ \frac{1}{-\nu \Delta +\alpha I_d}\left[ \P ((\U\cdot \vec{\nabla}) \U)\right]},$$ de la façon équivalente suivante: tout d'abord, étant donnée que $\U \in H^1(\Rt)$ et que $div(\U)=0$ alors on peut écrire le terme bilinéaire $(\U\cdot \vec{\nabla}) \U$ comme $div(\U \otimes \U)$ et ensuite nous écrivons 
$$ - \frac{1}{-\nu \Delta +\alpha I_d}\left[ \P (div(\U \otimes \U))\right]=  \frac{\nu \Delta}{-\nu \Delta +\alpha I_d}\left[ \P \frac{1}{-\nu \Delta} (div(\U \otimes \U))\right],$$ et nous allons maintenant étudier le  terme  à droite de l'identité ci-dessus.  Plus précisément, nous allons montrer que ce terme peut s'écrire comme le produit de convolution avec un noyau  et pour définir ce noyau (ce qui sera fait dans l'expression (\ref{def:K}) ci-après)  nous devons tout d'abord étudier 
 le terme $\ds{\P \frac{1}{-\nu \Delta} (div(\U \otimes \U))}$. \\
 \\
  Nous allons maintenant observer que le terme ci-dessus   s'écrit comme le produit de convolution  $\ds{m \ast (\U \otimes \U)} $, où $m=(m_{i,j,k})_{1\leq i,j,k \leq 3}$ est un tenseur avec $m_{i,j,k}\in \mathcal{C}^{\infty}(\Rt \setminus \{0\})$ une fonction homogène de degré $-2$. En effet,  étant donné que le projecteur de Leray $\P$ est défini en variable d'espace comme $$ \P(\vec{\varphi})= \vec{\varphi} + (\vec{R} \otimes \vec{R})\vec{\varphi},$$ où $\vec{R}=(R_i)_{1\leq i \leq 3}$ avec $R_i =\frac{\partial_i}{\sqrt{-\Delta}}$ la i-ème transformée de Riesz; et de plus, étant donné que $$ \frac{1}{-\nu \Delta} div(\U \otimes \U)= \left( \sum_{k=1}^{3} \frac{\partial_k}{-\nu \Delta}  (U_i U_k) \right)_{1\leq i \leq 3},$$ nous écrivons alors 
\begin{eqnarray}\label{estim:dec-2} \nonumber
 \P \frac{1}{-\nu \Delta} (div(\U \otimes \U)) &=& \left( \sum_{k=1}^{3} \frac{\partial_k}{-\nu \Delta}(U_i U_k)+ \sum_{j=1}^{3}\sum_{k=1}^{3} R_i R_j \frac{\partial_k}{-\nu \Delta}(U_j U_k)  \right)_{1\leq i \leq 3}	\\
&=& \left( \sum_{j=1}^{3}\sum_{k=1}^{3} \left[ \delta_{i,j} + R_i R_j \right]  \frac{\partial_k}{-\nu \Delta}(U_j U_k)  \right)_{1\leq i \leq 3},
	\end{eqnarray} et en prenant la transformée de Fourier dans chaque terme de la somme ci-dessus nous obtenons $\ds{\mathcal{F}\left[ \left[ \delta_{i,j} + R_i R_j \right] \frac{\partial_k}{-\nu \Delta}(U_j U_k) \right](\xi)  = \left[ \delta_{i,j} + \frac{\xi_i \xi_j}{\vert \xi \vert^2} \right]  \frac{ i \xi_k}{\nu \vert \xi \vert^2}(\mathcal{F}[U_j] \ast \mathcal{F}[U_k] )(\xi)}$.\\
	 \\ 
	 Nous définissons alors la fonction $m_{i,j,k}$ au niveau de Fourier par
	\begin{equation}\label{def:m-fourier}
	   \mathcal{F}\left[ m_{i,j,k}\right](\xi)= \left[ \delta_{i,j} + \frac{\xi_i \xi_j}{\vert \xi \vert^2} \right]  \frac{ i \xi_k}{\nu \vert \xi \vert^2},
	   \end{equation}
	    où nous pouvons observer que $\mathcal{F}\left[ m_{i,j,k}\right](\xi)$ est une fonction homogène de degré $-1$ et de classe $\mathcal{C}^{\infty}$ en dehors de l'origine et alors  $m_{i,j,k}\in \mathcal{C}^{\infty}(\Rt \setminus \{0\})$ est  une fonction homogène de degré $-2$.\\
	\\
	De cette façon, en revenant à l'identité (\ref{estim:dec-2}) nous écrivons 
	$$  \P \frac{1}{-\nu \Delta} (div(\U \otimes \U))  = \left( \sum_{j=1}^{3}\sum_{k=1}^{3} m_{i,j,k} \ast (U_j U_k)  \right)_{1\leq i \leq 3}, $$ où, pour simplifier l'écriture on pose le tenseur 
	\begin{equation}\label{def:m}
	m=(m_{i,j,k})_{1\leq i,j,k\leq 3},
	\end{equation} et par un abus de notation nous allons écrire dorénavant 
	\begin{equation}\label{estim:dec-3}
	 \P \frac{1}{-\nu \Delta} (div(\U \otimes \U))  =m\ast (\U \otimes \U).
	\end{equation}
	Une fois que l'on a définit le tenseur $m$ ci-dessus nous définissons maintenant  le noyau $K_{\nu,\alpha}$ comme le tenseur
	\begin{equation}\label{def:K}
	K_{\nu,\alpha}= \frac{\nu \Delta}{-\nu \Delta +\alpha I_d} \left[  m\right], 
	\end{equation} c'est à dire, $[K_{\nu,\alpha}]_{i,j,k}= \frac{\nu \Delta}{-\nu \Delta +\alpha I_d} \left[  m_{i,j,k}\right]$ pour $1\leq i,j,k\leq 3$; et alors par les identités (\ref{estim:dec-3})  et (\ref{def:K}) nous pouvons alors écrire 
	\begin{equation}\label{estim:dec-1}
		\frac{\nu \Delta}{-\nu \Delta +\alpha I_d} \left[ \P \frac{1}{-\nu \Delta} (div(\U \otimes \U)) \right]= \frac{\nu \Delta}{-\nu \Delta +\alpha I_d} \left[  m\ast (\U \otimes \U) \right] = K_{\nu,\alpha}\ast (\U\otimes \U).
		\end{equation} 
Une fois que nous disposons de cette identité   nous observons  que  l'estimation (\ref{estimation1}) et  alors équivalente a l'estimation
\begin{equation}\label{dec-espace-terme-bilin}
\vert K_{\nu,\alpha}\ast (\U\otimes \U)(x) \vert \leq \frac{c_2}{\vert x \vert^4},
\end{equation} pour tout $\vert x \vert >8$ et où $c_2=c_2(\U,\nu,\alpha)>0$ est une constante. Vérifions donc l'estimation (\ref{dec-espace-terme-bilin}).  La première chose  à faire c'est étudier la décroissance en variable d'espace du noyau $K_{\nu,\alpha}$ ci-dessus et  nous avons le lemme technique suivant:
\begin{Lemme}\label{Lemme:dec-K} Soit le noyau $K_{\nu,\alpha}$ défini par l'expression (\ref{def:K}). Ce noyau vérifie une décroissance suivante:
\begin{equation}\label{dec-noyau-K}
K_{\nu,\alpha}(x)\leq c_{\nu,\alpha} \left\lbrace  \begin{array}{cl} \vspace{2mm}
\frac{1}{\vert x \vert},& \text{si} \,\, \vert x \vert \leq 4, \\
 \frac{1}{\vert x \vert^4}, & \text{si} \,\, \vert x \vert > 4,
\end{array}	 \right.
\end{equation} où $c_{\nu,\alpha}>0$ est une constante qui dépend seulement de $\nu>0$ et $\alpha>0$.	
\end{Lemme}	Expliquons rapidement les grandes lignes de la preuve de ce lemme qui, pour la commodité du lecteur, sera faite  en détail à la fin du chapitre. Cette preuve suit essentiellement les idées de la preuve du Lemme \ref{Lemme:bessel} mais nous devons contourner quelques contraintes techniques. En effet, observons 
 par la définition du noyau $K_{\nu,\alpha}$ donnée dans l'expression (\ref{def:K}) nous savons que 
$\ds{K_{\nu,\alpha}= \frac{\nu \Delta}{-\nu \Delta +\alpha I_d} \left[  m\right]}$ (où le tenseur $m$ est défini par les expressions (\ref{def:m-fourier}) et (\ref{def:m})) et alors  nous pouvons  écrire $\ds{K_{\nu,\alpha}= \frac{1}{-\nu \Delta +\alpha I_d} \left[ \nu \Delta m\right]}$, néanmoins, on ne peut pas appliquer  directement ici le Lemme \ref{Lemme:bessel}:  pour appliquer ce lemme il faut que la fonction $\nu \Delta m$ appartienne à l'espace $L^1(\Rt)$ mais ceci n'est pas possible  car nous allons montrer que cette fonction a une décroissance $\ds{\vert \nu \Delta m(x)\vert \lesssim \frac{1}{\vert x \vert^4} }$ et alors elle  n'est pas intégrable à l'origine. Dans la page \pageref{sec:preuve-lemme-dec-K} nous allons voir comment contourner ce problème  et nous donnons une preuve du Lemme \ref{Lemme:dec-K}.\\
\\
%
Une fois que nous disposons  de l'information nécessaire sur la décroissance du noyau $K_{\nu,\alpha}$, nous revenons maintenant au terme $\ds{K_{\nu,\alpha}\ast (\U\otimes \U)}$ pour vérifier l'estimation (\ref{dec-espace-terme-bilin}). 	L'idée pour vérifier cette estimation est la suivante:   nous allons tout d'abord  montrer que la solution $\U$ vérifie la décroissance ci-dessous : 
\begin{Lemme}\label{Lemme:dec-U-aux} Pour tout $\vert x \vert >4$ on a $\ds{\vert \U (x)\vert \leq \frac{c_2}{\vert x \vert^2}}$. 
	\end{Lemme} La preuve de ce lemme repose essentiellement sur le fait que par les identités (\ref{point-fixe-stattionnaire-aux}) et (\ref{estim:dec-1}) nous pouvons alors  écrire   la solution $\U$ des équations (\ref{NS-stationnaire-prop})  comme  $$ \U= K_{\nu,\alpha}\ast (\U\otimes \U)+\frac{1}{-\nu \Delta +\alpha I_d}\left[\fe\right].$$ Ensuite,  par la Proposition \ref{Prop:dec-espace-CKN} nous savons que la solution stationnaire $\U$ vérifie une décroissance  
\begin{equation*}
\vert \U (x)\vert \leq c \frac{1}{1+\vert x \vert},
\end{equation*} et par cette décroissance,  l'estimation du noyau $K_{\nu,\alpha}$ donnée dans le Lemme \ref{Lemme:dec-K} et l'estimation (\ref{estimation2}), nous pouvons montrer  que la solution $\U$ vérifie  la décroissance donnée dans le Lemme \ref{Lemme:dec-U-aux}. Tous les calculs seront faits en détail à la fin du chapitre page \pageref{sec:preuve-lemme-dec-U-aux}; et maintenant  nous allons nous servir de ce lemme  pour vérifier l'estimation cherchée (\ref{dec-espace-terme-bilin}).\\
\\
En effet, pour $\vert x \vert >8$ fixe on commence par écrire 
\begin{eqnarray}\label{estim1:th} \nonumber
\vert K_{\nu,\alpha}\ast (\U \otimes \U)(x)\vert &\leq & \int_{\Rt} \vert K_{\nu,\alpha}(x-y)\vert \vert (\U \otimes \U)(y)\vert dy=\int_{\vert y \vert \leq \frac{\vert x \vert}{2}} \vert K_{\nu,\alpha}(x-y)\vert \vert (\U \otimes \U)(y)\vert dy \\
& & + \int_{\vert y \vert >2 \frac{\vert x \vert}{2}} \vert K_{\nu,\alpha}(x-y)\vert \vert (\U \otimes \U)(y)\vert dy= I_1+I_2,
\end{eqnarray} et l'on cherche à estimer les termes $I_1$ et $I_2$ ci-dessus. \\
\\
Pour le terme $I_1$, comme nous considérons ici $\vert y \vert \leq \frac{\vert x \vert}{2}$ alors nous avons les inégalités suivantes: $\vert x-y\vert \geq \vert x \vert -\vert y \vert \geq \frac{\vert x \vert}{2}$, d'où, étant donné que $\vert x \vert >8$  nous écrivons $\frac{\vert x \vert }{2}>4$ pour obtenir $\vert x-y \vert >4 $ et ainsi, par l'estimation (\ref{dec-noyau-K}) obtenue dans le Lemme \ref{Lemme:dec-K} nous avons $\ds{\vert K_{\nu,\alpha}(x-y) \vert \leq \frac{c_{\nu,\alpha}}{\vert x-y\vert^4}  \leq \frac{c_{\nu,\alpha}}{\vert x\vert^4}}$, et  nous pouvons écrire 
\begin{equation}\label{estim:I_1:th}
I_1  =\int_{\vert y \vert \leq \frac{\vert x \vert}{2}} \vert K_{\nu,\alpha}(x-y)\vert \vert (\U \otimes \U)(y)\vert dy\leq \frac{c_{\nu,\alpha}}{\vert x \vert^4} \int_{\vert y \vert \leq \frac{\vert x \vert}{2}} \vert  (\U \otimes \U)(y)\vert dy \leq \frac{c_{\nu,\alpha}}{\vert x \vert^4} \Vert \U \Vert^{2}_{L^2}. 
\end{equation} 
Pour le terme $I_2$, on commence par écrire 
\begin{eqnarray}\label{estim:I2:th}\nonumber
I_2&=&\int_{\vert y \vert > \frac{\vert x \vert}{2},\, \vert x-y\vert \leq 2} \vert K_{\nu,\alpha}(x-y)\vert \vert (\U \otimes \U)(y)\vert dy+\int_{\vert y \vert > \frac{\vert x \vert}{2},\, \vert x-y\vert > 2} \vert K_{\nu,\alpha}(x-y)\vert \vert (\U \otimes \U)(y)\vert dy\\
&=& (I_2)_a+(I_2)_b,
\end{eqnarray} et l'on doit estimer chaque terme de l'identité ci-dessus. Pour estimer le  terme $(I_2)_{a}$ nous allons utiliser la décroissance de la solution $\U$ donnée dans le Lemme  \ref{Lemme:dec-U-aux}. En effet, comme nous considérons ici $\vert y \vert >\frac{\vert x \vert }{2}$ (où comme $\vert x \vert >8$ alors on a  $\vert y \vert >4$) par ce lemme   nous  pouvons alors  écrire  $\ds{\vert (\U\otimes \U)(y)\vert \leq  \frac{c}{\vert y \vert^4}\leq \frac{c}{\vert x \vert^4}}$. D'autre part, comme nous considérons aussi $\vert x-y \vert \leq 2$, toujours par  l'estimation (\ref{dec-noyau-K}) (obtenue dans le Lemme \ref{Lemme:dec-K}) nous savons que $\ds{\vert K_{\nu,\alpha}(x-y) \vert \leq \frac{c_{\nu,\alpha}}{\vert x-y\vert}}$; et nous obtenons de cette façon 
\begin{equation}\label{estim:I2a:th}
(I_2)_a =\int_{\vert y \vert > \frac{\vert x \vert}{2},\, \vert x-y\vert \leq 2} \vert K_{\nu,\alpha}(x-y)\vert \vert (\U \otimes \U)(y)\vert dy \leq \frac{c_{\nu,\alpha}}{\vert x \vert^4} \int_{\vert x-y\vert \leq 2} \frac{dy}{\vert x-y\vert} \lesssim  \frac{c_{\nu,\alpha}}{\vert x \vert^4}. 
\end{equation}	  Pour estimer  le terme  $(I_2)_b$, rappelons que nous considérons ici $\vert x-y\vert >2$ et ainsi, toujours par l'estimation (\ref{dec-noyau-K}), nous savons que $\ds{K_{\nu,\alpha}(x-y)\leq \frac{c_{\nu,\alpha}}{\vert x-y\vert^4}}$ et alors, en suivant les mêmes lignes que l'estimation (\ref{estim:I_1:th}) nous obtenons 
\begin{equation}\label{estim:I2b:th}
(I_2)_b \leq  \frac{c_{\nu,\alpha}}{\vert x \vert^4}\Vert \U \Vert^{2}_{L^2}.
\end{equation} Une fois que nous disposons des estimations (\ref{estim:I2a:th}) et (\ref{estim:I2b:th}) nous posons maintenant la constante $c_{2}=\max(c_{\nu,\alpha},c_{\nu,\alpha}\Vert \U \Vert^{2}_{L^2})>0$ et en revenant  à l'identité (\ref{estim:I2:th}) nous pouvons écrire $\ds{I_2\leq \frac{c_{2}}{\vert x \vert^4}}$. Ainsi, par cette estimation et par l'estimation (\ref{estim:I_1:th}) nous revenons maintenant à l'identité (\ref{estim:I_1:th}) où nous obtenons l'estimation cherchée  (\ref{dec-espace-terme-bilin}). \finpv

\section{Appendice}  
Dans la Section \ref{sec:CKN} ci-dessous nous revisitons la théorie de la régularité partielle de Caffarelli, Kohn et Nirenberg dans le cadre des équations  de Navier-Stokes stationnaires et amorties; et  nous donnons une preuve de la Proposition \ref{Prop:dec-espace-CKN}. Ensuite, dans la Section \ref{sec:lemme-tech} nous prouvons les lemmes techniques. 
\subsection{La théorie de la régularité de Caffarelli, Kohn et Nirenberg}\label{sec:CKN}  

Le but de cette section est de donner une preuve de la Proposition \ref{Prop:dec-espace-CKN}  (page \pageref{Prop:dec-espace-CKN}) et comme annoncé  nous avons besoin de passer par le cadre de la théorie de la régularité partielle  de Caffarelli, Kohn et Niremberg \cite{CKN} (CKN). Pour  exposer les idées de la preuve de cette proposition d'une façon plus claire  nous avons divisé cet appendice comme suit: dans le point $A)$ ci-dessous  on commence par expliquer les grandes lignes de la preuve de la Proposition \ref{Prop:dec-espace-CKN} et nous expliquons comment nous allons utiliser ici la théorie CKN. Ensuite, dans le point $B)$  nous faisons un rappel sur un résultat classique de la  théorie CKN dont on a besoin et après nous adaptons ce résultat au cadre des équations de Navier-Stokes stationnaires et amorties; et ceci sera fait dans le Théorème \ref{Theo:CKN}. Finalement, dans le point $C)$, à l'aide du Théorème \ref{Theo:CKN}   nous donnons une preuve de la  Proposition \ref{Prop:dec-espace-CKN}. 

\subsubsection{A) Les idées de la preuve de la Proposition \ref{Prop:dec-espace-CKN}}
Rappelons rapidement que dans cette proposition nous considérons les équations  
\begin{equation}\label{N-S-alpha-aux}
-\nu \Delta \U +(\U\cdot \vec{\nabla})\U +\vec{\nabla}P  =\fe-\alpha \U, \qquad div(\U)=0, \quad \alpha>0,
\end{equation} avec la  force $\fe$ donnée dans la Définition \ref{Definition_force_ext};  et nous voulons montrer que toute solution $\U \in H^{1}(\Rt)$ (obtenue par le Théorème \ref{Theo:solutions_stationnaires_turbulent}) vérifie une décroissance en variable d'espace:
\begin{equation}\label{dec-u}
\vert \U(x)\vert \leq C \frac{1}{1 +\vert x \vert},
\end{equation}  où $C=C(\U,\nu,\alpha)>0$ est une constante qui dépend de la solution $\U$, la constante de viscosité $\nu>0$ et le paramètre d'amortissement $\alpha>0$.\\
\\
Dans l'expression (\ref{dec-u}) nous pouvons observer  que cette estimation équivaut  aux estimations suivantes:
\begin{equation}\label{Estm1}
\vert \U(x)\vert \leq c,
\end{equation} et 
\begin{equation}\label{Estim2}
\vert \U(x)\vert \leq \frac{c}{\vert x \vert};
\end{equation} et il s'agit alors de vérifier ces deux estimations. L'estimation (\ref{Estm1}) sera  vérifié dans le Lemme \ref{lemme:U-L-infty}  du point $C)$ ci-après (où nous allons montrer que la solution $\U$ appartient à l'espace $L^{\infty}(\Rt)$) et cette estimation ne  présente  aucun problème particulier. Mais, l'estimation (\ref{Estim2}) est plus délicate à vérifier car \emph{a priori} nous ne disposons  d'aucune information supplémentaire pour que la solution ait une décroissance comme celle donnée dans  (\ref{Estim2}); et pour vérifier cette décroissance nous allons adapter un résultat de la théorie CKN au cadre des équations (\ref{N-S-alpha-aux}). 

\subsubsection{B) La théorie de la régularité CKN  pour les équations amorties}\label{sec:theorie-CKN} 
On commence donc par faire un très court  rappel sur un résultat de cette théorie. Pour un exposé complet sur l'état de l'art de la théorie de la régularité CKN  voir les Chapitres $13$ et $14$ du livre \cite{PGLR1}. \\ 
\\
Le résultat dont on a besoin est obtenu dans le cadre des  équations de Navier-Stokes classiques (sans terme d'amortissement): 
\begin{equation}\label{n-s-ckn}
\partial_t\vu = \nu\Delta \vu-(\vu\cdot\vec{\nabla}) \vu -\vec{\nabla}p +\vec{g}, \quad div(\vu)=0,
\end{equation} où $\vec{g}=\vec{g}(t,x)$ est une force qui vérifie $g \in L^2([0,T[, H^{-1}(\Rt))$ pour un temps $T>0$; et nous allons maintenant introduire quelques définitions et notations. Pour $x_0 \in \Rt$, $r_0>0$ et $t_0>0$ nous considérons le domaine $Q_{0}=]t_0-r^{2}_{0},t_0[\times B(x_0,r_0)\subset ]0,+\infty [\times \Rt$ et nous avons les définitions suivantes:
\begin{Definition}[Solution faible et solution adaptée]\label{def:CKN} 
	\begin{enumerate}
		\item[] 
		\item[1)] Le couple $(\vu,p)$ est une solution faible des équations (\ref{n-s-ckn}) sur le domaine $Q_{0}$ ci-dessus si l'on a $\vu \in L^{\infty}_{t}L^{2}_{x}\cap L^{2}_{t}\dot{H}^{1}_{x}(Q_{0})$ et $p \in \mathcal{D}^{'}(Q_{0})$.
		\item[2)] Une solution faible $\vu$ des équations (\ref{n-s-ckn})  est une solution adaptée si cette solution vérifie l'inégalité d'énergie locale dans $\mathcal{D}^{'}(Q_0)$: 
		\begin{equation}\label{ineq-ener-local}
		\partial_t \vert \vu \vert^2 \leq \nu \Delta (\vert \vu \vert^2)-\nu \vert \vec{\nabla}\otimes \vu \vert^2 -div((2p-\vert \vu \vert^2)\vu)+\vu\cdot \vec{g}. 
		\end{equation}
 		\end{enumerate}
\end{Definition}  Maintenant que l'on a introduit ces définitions  nous pouvons énoncer le résultat suivant qui s'agit d'un critère de régularité de la théorie CKN  (pour une preuve de ce  résultat voir le Théorème $14.4$ page $472$ du livre \cite{PGLR1})
\begin{Theoreme}\label{Theo:CKN-classique} Soit le domaine $Q_{0}=]t_0-r^{2}_{0},t_0 [\times B(x_0,r_0)$. Soit $(\vu,p) \in  L^{\infty}_{t}L^{2}_{x}\cap L^{2}_{t}\dot{H}^{1}_{x}(Q_{0}) \times \mathcal{D}^{'}(Q_{0})$ une solution faible des équations de Navier-Stokes (\ref{n-s-ckn}). Si:
	\begin{enumerate}
		\item[1)] $p \in L^{\frac{3}{2}}_{t}L^{\frac{3}{2}}_{x}(Q_{0})$,
		\item[2)] $\vec{g}\in L^{q}_{t}L^{q}_{x}(Q_{0})$, avec $q>\frac{5}{2}$,
		\item[3)] $\vu$	est une solution adapté au sens du point $2)$ de la Définition \ref{def:CKN},
	\end{enumerate}	 alors il existe deux constantes $\eta>0$ et $c>0$, qui ne dépendent que de la constante de viscosité  $\nu$ et de $q$, telles que: si  pour $0\leq \lambda \leq \eta$ on a: 
	\begin{enumerate}
		\item[4)] $\ds{\int_{Q_{0}}(\vert \vu(t,x)\vert^3+\vert p(t,x)\vert^{\frac{3}{2}})dx\,dt \leq \lambda^3 r^{2}_{0}} $ et 
		\item[5)]  $\ds{\int_{Q_{0}} \vert \vec{g}(t,x)\vert^{q}dx\,dt \leq \lambda^{2 q} r^{5-3q}_{0}}$;
	\end{enumerate}
	alors $\vu$ est bornée sur le sous-domaine $Q_{1}=]t_0-\frac{r^{2}_{0}}{4} ,t_0 [\times B(x_0,\frac{r_0}{2}) \subset Q_{0}$ et l'on a 
	\begin{equation}\label{estim-ckn}
	\sup_{(t,x)\in Q_1} \vert \vu(t,x)\vert \leq c \lambda \frac{1}{r_{0}}. 
	\end{equation}
\end{Theoreme}	 
Ce résultat de la théorie CKN  a été aussi étudié dans différentes  contextes  (voir par exemple le Chapitre $30$ du livre \cite{PGLR2}) et nous n'allons pas faire ici une discussion sur ce sujet.  Notre intérêt au Théorème \ref{Theo:CKN-classique}  repose essentiellement sur l'estimation (\ref{estim-ckn}) comme nous l'expliquons tout de suite: dans cette estimation nous observons que pour tout $(t,x)\in Q_1=]t_0-\frac{r^{2}_{0}}{4} ,t_0 [\times B(x_0,\frac{r_0}{2})$ nous avons 
$\ds{\vert \vu(t,x)\vert \leq  c \lambda \frac{1}{r_0}}$, et alors, si pour $x_0\in \Rt$, où $\vert x_0 \vert>0 $, on pose $r_0=\vert x_0 \vert$ alors par l'estimation (\ref{estim-ckn})  nous observons que  pour $(t,x_0)$ nous pouvons  écrire  
\begin{equation}\label{dec-ckn-n-s-classique}
\vert \vu(t,x_0) \leq c\lambda \frac{1}{\vert x_0 \vert}.
\end{equation}  Cette estimation est intéressante car nous pouvons observer qu'ici l'on obtient une décroissance similaire à la décroissance cherchée (\ref{Estim2}) et nous pouvons aussi observer que cette décroissance est valide pour tout $\vert x_0\vert>0$ car la constante $c\lambda$  ne dépend pas de $x_0$. \\
\\
En suivant ces idées, nous voulons donc adapter le Théorème \ref{Theo:CKN-classique} au cadre des équations de Navier-Stokes amorties et stationnaires (\ref{N-S-alpha-aux}) pour ensuite pouvoir vérifier l'estimation (\ref{Estim2}). Nous avons ainsi le résultat suivant:

\begin{Theoreme}\label{Theo:CKN}  Soit la force $\fe$ donnée dans la Définition \ref{Definition_force_ext} et soit $(\U,P) \in H^1(\Rt) \times H^{\frac{1}{2}}(\Rt)$ une solution des équations  (\ref{N-S-alpha-aux})   obtenue dans le Théorème \ref{Theo:solutions_stationnaires_turbulent}. Soit $x_0\in \Rt$, $r_0>0$ et  soit la boule $B_0=B(x_0,r_0)\subset \Rt$. Si: 
	\begin{enumerate}
		\item[1)] $P\in L^{\frac{3}{2}}(B_0)$ et 
		\item[2)] $\fe \in L^q(B_0)$, avec $q>\frac{5}{2}$,
	\end{enumerate} alors il existe deux constantes $\eta>0$ et $c>0$, qui ne dépendent que de la constante de viscosité  $\nu>0$, le paramètre $q>\frac{5}{2}$ et le paramètre d'amortissement $\alpha>0$, telles que: si pour $0\leq \lambda \leq \eta$ on a: 
	\begin{enumerate}
		\item[3) ] $\ds{\int_{B_0} \left(\vert \U(x)\vert^3 + \vert P(x)\vert^{\frac{3}{2}}\right)dx \leq \lambda^3 }$ et 
		\item[4)]  $\ds{\int_{B_0} \vert \fe (x)\vert^qdx \leq \lambda^{2q}r^{\frac{3-3q}{0}}}$,
	\end{enumerate}	alors la solution $\U$ est bornée sur la boule $B_1=B(x_0,\frac{r_0}{2}) \subset B_0$ et l'on a l'estimation: 
	\begin{equation}\label{estim-CKN-stat}
	\sup_{x\in B_1}\vert \U(x)\vert \leq c \lambda \frac{1}{r_0}. 
	\end{equation}
\end{Theoreme}	
\dm 
Comme le champ de vitesse $\U$ est une fonction stationnaire alors nous avons $\partial_t \U=0$; et de plus, comme le couple $(\U,P) \in H^1(\Rt) \times H^{\frac{1}{2}}(\Rt)$ est une solution des équations de Navier-Stokes stationnaires et amorties (\ref{N-S-alpha-aux}) alors $(\U,P)$ est aussi une solution des équations de Navier-Stokes amorties:
\begin{equation}\label{ns-evol-aux-theo}
\partial_t\U = \nu\Delta \U-(\U\cdot\vec{\nabla}) \U -\vec{\nabla}P +\fe-\alpha\U, \quad div(\U)=0;     
\end{equation} et alors  nous allons  maintenant montrer  que l'on peut appliquer le Théorème \ref{Theo:CKN-classique} aux équations ci-dessus. \\
\\
 Soit un temps  $t_0>0$ et soit donc le domaine $Q_{0}=]t_0-r^{2}_{0},t_0 [\times B(x_0,r_0)$.  Observons tout d'abord que $(\U,P)$ est une solution faible des équations (\ref{ns-evol-aux-theo}) sur le domaine $Q_0$ au sens du point $1)$ de la Définition \ref{def:CKN}.  En effet, comme $\U\in H^1(\Rt)$ ne dépend que de la variable spatiale et comme nous avons ici un intervalle du temps borné: $]t_0-r^{2}_{0},t_0 [$, alors nous avons $\U\in L^{\infty}_{t} L^{2}_{x} \cap L^{2}_{t}\dot{H}^{1}_{x}(Q_0)$. De même façon, comme $P\in H^{\frac{1}{2}}(\Rt)$ est une fonction stationnaire   nous avons $P\in \mathcal{D}^{'}(Q_0)$. \\
\\
Observons maintenant que  le champs de vitesse $\U$, la pression $P$ et la force $\fe$ vérifient les hypothèses $1)$, $2)$ et $3)$ du Théorème  \ref{Theo:CKN-classique}. En effet, Les points $1)$ et $2)$ du Théorème \ref{Theo:CKN-classique}  sont une conséquence directe des hypothèses  $1)$ et $2)$ données dans le  Théorème \ref{Theo:CKN} et du fait que les fonctions $P$ et $\fe$ ne dépendent pas de la variable du temps. Quant au  point $3)$  du Théorème \ref{Theo:CKN-classique},  nous allons maintenant montrer  que la solution $\U$ est une fonction adaptée au sens du point $2)$ de la Définition \ref{def:CKN}. \\
\\
 Il s'agit de montrer que la solution $\U$ vérifie l'inégalité locale (\ref{ineq-ener-local}) et pour cela nous remarquons que comme le couple $(\U,P)$ vérifie les équations (\ref{ns-evol-aux-theo}) et comme nous avons en plus $\U \in H^1(\Rt)$ et $P\in H^{\frac{1}{2}}(\Rt)$ alors les fonctions $\U$ et $P$ sont suffisamment régulières et nous pouvons écrire l'identité au sens des distributions suivante:
\begin{equation}\label{ident-ener-local-U}
\partial_t \vert \U \vert^2  = \nu \Delta (\vert \U \vert^2)-\nu \vert \vec{\nabla}\otimes \U \vert^2 -div((2P-\vert \U \vert^2)\vu) + \U\cdot \vec{\fe}-\alpha \vert \U \vert^2.
\end{equation} Observons maintenant que  le terme $-\alpha \vert \U \vert^2$ ci-dessus a un signe négatif et alors par cette identité nous pouvons aussi écrire l'inégalité (toujours au sens des distributions) suivante:
$$ \partial_t \vert \U \vert^2  \leq  \nu  \Delta (\vert \U \vert^2)-\nu \vert \vec{\nabla}\otimes \U \vert^2 -div((2P-\vert \U \vert^2)\vu) + \U\cdot \vec{\fe}, $$ qui est l'inégalité d'énergie locale (\ref{ineq-ener-local}); et donc la solution $\U$ est bien une solution adaptée au sens de la Définition \ref{def:CKN}. \\
\\
Une fois que l'on a vérifié les points $1)$, $2)$ et $3)$ du Théorème \ref{Theo:CKN-classique}, par ce théorème  avons alors la conclusion suivante: il existe deux constantes $\eta>0$ et $c>0$, qui ne dépendent que $\nu>0$, $q>\frac{5}{2}$ et de $\alpha>0$, telles que: si  pour $0\leq \lambda \leq \eta$ on a 
\begin{enumerate}
	\item[i)] $\ds{\int_{Q_{0}}(\vert \U(x)\vert^3+\vert P(x)\vert^{\frac{3}{2}})dx\,dt  \leq \lambda^3 r^{2}_{0}} $ et 
	\item[ii)]  $\ds{\int_{Q_{0}} \vert \fe (x)\vert^{q}dx\,dt \leq \lambda^{2 q} r^{5-3q}_{0}}$;
\end{enumerate}
alors $\U$ est bornée sur le sous-domaine $Q_{1}=]t_0-\frac{r^{2}_{0}}{4} ,t_0 [\times B(x_0,\frac{r_0}{2}) \subset Q_{0}$ et l'on a 
\begin{equation}\label{estim-ckn-amortie}
\sup_{(t,x)\in Q_1} \vert \U(x)\vert \leq c \lambda \frac{1}{r_{0}}. 
\end{equation}  Observons  finalement que cette conclusion est la même conclusion énoncée dans le Théorème \ref{Theo:CKN}: étant donné que  les fonctions $\U$, $P$ et $\fe$ sont toujours des fonctions stationnaires alors  les estimations données dans les points $i)$ et $ii)$ ci-dessus sont équivalentes aux estimations énoncées dans les points $3)$ et $4)$ du Théorème \ref{Theo:CKN}. En effet, comme $Q_0=]t_0-r^{2}_{0},t_0 [\times B(x_0,r_0)$ alors par  l'estimation donnée dans  $i)$  nous pouvons écrire 
\begin{eqnarray*}
r^{2}_{0} \int_{B(x_0,r_0)}  \vert \U(x)\vert^3+\vert P(x)\vert^{\frac{3}{2}})dx&=&\int_{t_0-r^{2}_{0}}^{t_0} \int_{B(x_0,r_0)}  \vert \U(x)\vert^3+\vert P(x)\vert^{\frac{3}{2}})dx\,dt \\
&=& \int_{Q_{0}}(\vert \U(x)\vert^3+\vert P(x)\vert^{\frac{3}{2}})dx\,dt \leq   \lambda^3 r^{2}_{0},
\end{eqnarray*}  ce qui équivaut à écrire   $\ds{\int_{B_0} \left(\vert \U(x)\vert^3 + \vert P(x)\vert^{\frac{3}{2}}\right)dx \leq \lambda^3}$, qui est l'estimation donnée dans le point $3)$ du Théorème \ref{Theo:CKN}. En suivant le même raisonnement nous avons que l'estimation donnée dans le point $ii)$ équivaut à l'estimation $4)$ de ce théorème. De plus, nous pouvons aussi observer que l'estimation (\ref{estim-ckn-amortie}) équivaut à l'estimation (\ref{estim-CKN-stat}) et alors le Théorème \ref{Theo:CKN} est maintenant démontré. \finpv
Maintenant que l'on dispose de l'estimation (\ref{estim-CKN-stat}) nous pouvons montrer que la solution $\U$ vérifie l'estimation (\ref{dec-u}).

\subsubsection{C) Preuve de la Proposition \ref{Prop:dec-espace-CKN}}\label{sec:fin-Prop-CKN}
Comme expliqué dans le point $A)$ on commence par vérifier l'estimation (\ref{Estm1}).
\begin{Lemme}\label{lemme:U-L-infty} Soit une $\U \in H^1(\Rt)$ une solution des équations (\ref{N-S-alpha-aux})  associée à la force $\fe$ donnée dans la Définition \ref{Definition_force_ext}. Alors $\U \in L^{\infty}(\Rt)$. 
	\end{Lemme}
\pv  Comme nous avons l'inclusion $H^2(\Rt)\subset L^{\infty}(\Rt)$ nous allons donc montrer que la solution $\U$ appartient à l'espace $H^2(\Rt)$. Nous écrivons cette solution comme la solution du problème de point fixe suivant:
\begin{equation}\label{eq-aux-pf}
\U = - \frac{1}{-\nu \Delta +\alpha I_d}\left[ \P ((\U\cdot \vec{\nabla}) \U)\right] + \frac{1}{-\nu \Delta +\alpha I_d}\left[\fe\right],
\end{equation}  voir l'identité (\ref{point-fixe-stattionnaire}) pour tous les détails des calculs; où nous observons tout d'abord que comme la force $\fe$ appartient à la classe de Schwartz alors  $\ds{\frac{1}{-\nu \Delta +\alpha I_d}\left[\fe\right] \in H^2(\Rt) }$ et il s'agit alors de vérifier que l'on a aussi $\ds{ \frac{1}{-\nu \Delta +\alpha I_d}\left[ \P ((\U\cdot \vec{\nabla}) \U)\right]\in H^2(\Rt)}$. En effet, comme $div(\U)=0$ nous écrivons  $\ds{(\U\cdot \vec{\nabla}) \U= div ( \U \otimes  \U)}$ et  nous avons alors l'identité 
\begin{equation}\label{eq-aux} 
\frac{1}{-\nu \Delta +\alpha I_d}\left[ \P ((\U\cdot \vec{\nabla}) \U)\right]= \frac{1}{-\nu \Delta +\alpha I_d}\left[ \P (div ( \U \otimes  \U))\right], 
\end{equation} et nous allons maintenant vérifier que le terme à droite de cette identité appartient à l'espace $H^2(\Rt)$.  \'Etant donné que $\U \in H^1(\Rt)$ alors par les lois de produit nous avons $\U\otimes \U \in H^{\frac{1}{2}}(\Rt)$ d'où nous obtenons  $\frac{1}{-\nu \Delta +\alpha I_d}\left[ \P (div ( \U \otimes  \U))\right] \in H^{\frac{3}{2}}(\Rt) $; et alors  par l'identité (\ref{eq-aux}) et l'identité  (\ref{eq-aux-pf}) nous avons $\U \in H^{\frac{3}{2}}(\Rt)$. Avec l'information $\U \in H^{\frac{3}{2}}(\Rt)$, toujours par le lois de produit nous avons $\U \otimes \U \in H^{\frac{3}{2}}(\Rt)\subset H^1(\Rt)$  et  nous avons $ \frac{1}{-\nu \Delta +\alpha I_d}\left[ \P (div ( \U \otimes  \U))\right]\in H^2(\Rt)$.\\
\\
Ainsi,  toujours par l'identité (\ref{eq-aux-pf}) nous obtenons $\U \in H^2(\Rt)$ et donc $\U \in L^{\infty}(\Rt)$. \finpv   
Dans le lemme suivant (qui s'agit d'un corollaire du Théorème \ref{Theo:CKN}) nous vérifions  maintenant l'estimation (\ref{Estim2}).
\begin{Lemme}\label{lemme-ckn} Soit une $\U \in H^1(\Rt)$ une solution des équations (\ref{N-S-alpha-aux})  associée à la force $\fe$ donnée dans la Définition \ref{Definition_force_ext}. Il existe deux constantes  $c>0$ et $M>0$, qui ne dépendent que de la constante de viscosité du fluide $\nu>0$ et le paramètre d'amortissement $\alpha>0$, telles que la solution $\U$ vérifie l'estimation $$\vert \U(x)\vert \leq \frac{c}{\vert x \vert},$$ pour tout $\vert x \vert >M$. 
\end{Lemme}
\pv   La première chose à faire est de vérifier les points $1)$ et $2)$ du Théorème \ref{Theo:CKN} pour obtenir l'estimation (\ref{estim-CKN-stat}). Soient $x_0\in \Rt$, $r_0>0$ et soit la boule $B_0=B(x_0,r_0)\subset \Rt$. Pour vérifier le point $1)$ du Théorème \ref{Theo:CKN} il suffit de montrer  que la pression $P$ appartient à l'espace $L^{\frac{3}{2}}(\Rt)$ et pour cela nous écrivons 
$\ds{P=-\frac{1}{\Delta} div (div(\U\otimes \U))= \sum_{i,j=1}^{3}R_{i}R_{j} (U_i,U_j),} $ où $R_i=\frac{\partial_i}{\sqrt{-\Delta}}$ est toujours la transformée de Riesz; et comme l'opérateur $R_iR_j$ est borné dans l'espace $L^{\frac{3}{2}}(\Rt)$ alors pour obtenir $P\in L^{\frac{3}{2}}(\Rt)$ nous allons montrer que $U_i,U_j \in L^{\frac{3}{2}}(\Rt)$ pour $1\leq i,j\leq 3$, ce qui équivaut à montrer $\U \otimes \U \in L^{\frac{3}{2}}(\Rt)$. En  effet, comme $\U \in L^2(\Rt)$ par les inégalités de H\"older nous avons $\U \otimes \U \in L^{1}(\Rt)$. D'autre par le Lemme   \ref{lemme:U-L-infty}  nous avons aussi $\U \in L^{\infty}(\Rt)$ et alors, comme $\U \in L^2(\Rt)$, nous avons  $\U \otimes \U \in L^{2}(\Rt)$; et par interpolation nous obtenons    $\U \otimes \U \in L^{\frac{3}{2}}(\Rt)$.\\
\\
Pour vérifier  le point $2)$ Théorème \ref{Theo:CKN} il suffit de remarquer que  la force $\fe$ appartient à la classe de Schwartz (voir toujours la Définition  \ref{Definition_force_ext}) et nous avons directement $\fe \in L^q(B_0)$ avec $q>\frac{5}{2}$. Ainsi, par le Théorème \ref{Theo:CKN} il existe deux constantes $\eta>0$ et $c>0$, qui ne dépendent que de $\nu>0$, $q>\frac{5}{2}$ et $\alpha>0$, telle que: si pour $0\leq \lambda \leq \eta$, on a :
\begin{equation}\label{H1}
\int_{B_0} \left(\vert \U(x)\vert^3 + \vert P(x)\vert^{\frac{3}{2}}\right)dx \leq \lambda^3
\end{equation} et 
\begin{equation}\label{H2}
\int_{B_0} \vert \fe (x)\vert^q dx \leq \lambda^{2q}r^{3-3q}_{0},
\end{equation} alors la solution $\U$ vérifie l'estimation (\ref{estim-CKN-stat}); et pour pouvoir écrire cette estimation nous devons maintenant vérifier les estimation (\ref{H1}) et (\ref{H2}). \\
\\
Fixons tout d'abord $\lambda=\eta$. Pour vérifier l'estimation (\ref{H1}) nous allons montrer les estimations suivantes
\begin{equation}\label{estim-aux-1}
\int_{B_0}\vert \U(x)\vert^3 \leq \frac{\eta^3}{2}\quad \text{et}\quad \int_{B_0}\vert P(x)\vert^{\frac{3}{2}}dx \leq \frac{\eta^3}{2}.
\end{equation} En effet, pour monter la première estimation de (\ref{estim-aux-1})  observons tout d'abord que comme $\U\in H^1(\Rt)=L^2\cap \dot{H}^{1}(\Rt)$ alors par les inégalités de Hardy-Littlewood-Sobolev nous avons $\U \in L^2\cap L^6(\Rt)$ et ensuite par les inégalités de interpolation nous obtenons $\U \in L^3(\Rt)$. Ainsi, pour $\frac{\eta^3}{2}>0$ il existe $M_1>0$ tel que si $\vert x_0 \vert>M_1$ alors nous avons $\ds{\int_{B_0}\vert \U(x)\vert^3 \leq \frac{\eta^3}{2}}$. \\
\\
Pour vérifier la deuxième estimation de (\ref{estim-aux-1}), rappelons que l'on a $P\in L^{\frac{3}{2}}(\Rt)$  et alors pour $\frac{\eta^3}{2}>0$ il existe $M_2>0$ tel que si $\vert x_0 \vert>M_2$ alors nous avons l'estimation  $\ds{\int_{B_0}\vert P\vert^{\frac{3}{2}}dx \leq \frac{\eta^3}{2}}$. \\
\\
Ainsi, on pose alors la constante 
\begin{equation}\label{M0}
M_0=\max(M_1,M_2)>0,
\end{equation} et pour $\vert x_0 \vert >M_0$ nous pouvons alors écrire les estimations données dans (\ref{estim-aux-1}) et donc nous avons vérifié l'estimation (\ref{H1}).\\
\\
Vérifions maintenant l'estimation  (\ref{H2}). Pour $\vert x_0 \vert >M_0$ nous prenons ici le rayon $r_0>0$ de la boule $B_0=B(x_0,r_0)$ comme 
\begin{equation}\label{r0}
r_0=\frac{\vert x_0 \vert}{2}, 
\end{equation} et alors l'estimation (\ref{H2}) s'écrit comme
\begin{equation}\label{H2-aux}
\left[\frac{\vert x_0 \vert}{2}\right]^{3q-3}\int_{B(x_0,\frac{\vert x_0\vert}{2})} \vert \fe \vert^q dx \leq \eta^{2q}, 
\end{equation} et pour vérifier cette estimation nous devons étudier un peu plus le terme à gauche ci-dessus. Toujours par le fait que  la force $\fe$ appartient à la classe de Schwartz  alors pour $\beta>3$ nous pouvons écrire  $\vert \fe (x)\vert \leq \frac{c}{\vert x \vert^{\beta}}$, pour tout $x\in \Rt$; et alors par cette estimation nous avons
$$  \left[\frac{\vert x_0 \vert}{2}\right]^{3q-3}\int_{B(x_0,\frac{\vert x_0\vert}{2})} \vert \fe \vert^q dx \leq c \left[\frac{\vert x_0 \vert}{2}\right]^{3q-3}\int_{B(x_0,\frac{\vert x_0\vert}{2})} \frac{dx}{\vert x \vert^{q\beta}},$$ et en passant à des cordonnées radiales nous obtenons 
$$ c \left[\frac{\vert x_0 \vert}{2}\right]^{3q-3}\int_{B(x_0,\frac{\vert x_0\vert}{2})} \frac{dx}{\vert x \vert^{q\beta}} \leq c \left[\frac{\vert x_0 \vert}{2}\right]^{3q-3}\int_{\frac{\vert x_0 \vert}{2}}^{ 3\frac{\vert x_0\vert}{2}} \rho^{2-q\beta}d \rho= c \vert x_0 \vert^{q(3-\beta)} =\frac{c}{\vert x_0 \vert^{q(\beta-3)}},$$ et alors nous pouvons finalement écrire 
$$ \left[\frac{\vert x_0 \vert}{2}\right]^{3q-3}\int_{B(x_0,\frac{\vert x_0\vert}{2})} \vert \fe \vert^q dx \leq \frac{c}{\vert x_0 \vert^{q(\beta-3)}}. $$ 
Par cette estimation  nous observons que pour vérifier l'estimation (\ref{H2-aux}) il suffit de vérifier l'estimation $\ds{\frac{c}{\vert x_0 \vert^{q(\beta-3)}} \leq \eta^{2q}}$ et ceci équivaut à écrire   $\ds{\frac{c}{\eta^{\frac{2q}{q(\beta-3)}}} \leq \vert x_0\vert}$. \\
\\
Ainsi,  on pose la constante $\ds{M=\max\left(M_0, \frac{c}{\eta^{\frac{2q}{q(\beta-3)}}} \right)>0}$, où la constante $M_0>0$ est donnée dans (\ref{M0}),  et pour  tout $\vert x_0 \vert >M$ et $r_0=\frac{\vert x_0 \vert}{2}$ nous avons l'estimations (\ref{H1}) et (\ref{H2}) et donc nous pouvons écrire (\ref{estim-CKN-stat}): $$ \sup_{x\in B(x_0,\frac{\vert x_0 \vert}{4})} \vert \U(x)\vert \leq  \frac{c\eta}{\vert x_0 \vert}.$$ Ainsi, comme la constante $c\eta>0$ ne dépend pas de $x_0$ et en écrivant maintenant $x_0=x$ alors  pour tout $\vert x  \vert >M$ par l'estimation ci-dessus nous pouvons écrire l'estimation cherchée $\ds{\vert \U(x) \vert \leq \frac{c\eta}{\vert x \vert}}$. \finpv 
\\
De cette façon par les Lemmes \ref{lemme:U-L-infty} et \ref{lemme-ckn} on pose maintenant la constante $C=\max(\Vert \U \Vert_{L^{\infty}}, c\eta)>0$ et pour tout $x\in \Rt$ nous pouvons écrire $\vert \U(x)\vert \leq C \frac{1}{1+\vert x \vert}$, qui est l'estimation cherchée (\ref{dec-u}) et alors la Proposition \ref{Prop:dec-espace-CKN} est maintenant montrée . \finpv

 \subsection{Les lemmes techniques}\label{sec:lemme-tech}  
\subsubsection*{Preuve du Lemme \ref{Lemme:ineq_ener_stabilite} page \pageref{Lemme:ineq_ener_stabilite}}\label{sec:preuve-lemme-tech-1} 
Nous commençons  par écrire 
\begin{equation}\label{theo_stat_ineq_0}
\Vert \vv(t,\cdot)\Vert^{2}_{L^2}= \Vert \vu(t,\cdot)\Vert^{2}_{L^2}-2\langle \vu(t,\cdot),\U\rangle_{L^2\times L^2} +\Vert \U \Vert^{2}_{L^2},
\end{equation} où le deuxième terme à droite peut s'écrire comme 
$$ -2\langle \vu(t,\cdot),\U\rangle_{L^2\times L^2} = -2\langle \vu_0,\U\rangle_{L^2\times L^2} -2 \int_{0}^{t} \langle \partial_t\vv(t,\cdot),\U\rangle_{H^{-2}\times H^2}ds.$$
En effet,  comme $\U$ est stationnaire nous avons $\partial_t(\vu\cdot \U)=(\partial_t \vu)\cdot \U$ et de plus,  comme  $\U \in H^2$  nous pouvons écrire $$ \partial_t \langle \vu(t,\cdot),\U\rangle_{L^2 \times L^2}=\langle \partial_t \vu(t,\cdot),\U\rangle_{H^{-2}\times H^2}.$$ Mais, comme $\vv=\vu-\U$ alors nous avons $\vu=\vv+\U$ et donc nous obtenons $\partial_t \vu=\partial_t \vv$. Ainsi,  en remplaçant cette identité dans l'identité précédente nous avons 
$$ \partial_t \langle \vu(t,\cdot),\U\rangle_{L^2 \times L^2}=\langle \partial_t \vv(t,\cdot),\U\rangle_{H^{-2}\times H^2}.$$ 
Finalement, on intègre par rapport à la variable de temps et nous obtenons de cette façon 
$$  \langle \vu(t,\cdot),\U\rangle_{L^2\times L^2} = \langle \vu_0,\U\rangle_{L^2\times L^2} + \int_{0}^{t} \langle \partial_t\vv(t,\cdot),\U\rangle_{H^{-2}\times H^2}ds.$$
De cette façon, en appliquant cette identité au deuxième terme de (\ref{theo_stat_ineq_0}) nous avons 
\begin{equation}\label{theo_stat_ineq_8}
\Vert \vv(t,\cdot)\Vert^{2}_{L^2}= \Vert \vu(t,\cdot)\Vert^{2}_{L^2} -2\langle \vu_0,\U\rangle_{L^2\times L^2} -2 \int_{0}^{t} \langle \partial_t\vv(t,\cdot),\U\rangle_{H^{-2}\times H^2}ds +\Vert \U \Vert^{2}_{L^2}.  
\end{equation} Maintenant, comme la solution  de Leray $\vu$ vérifie  l'inégalité d'énergie obtenue dans la Proposition \ref{Proposition_ineq_energie_alpha_modele}:    
\begin{eqnarray*}
\Vert \vu(t,\cdot)\Vert^{2}_{L^2} &\leq & \Vert \vu_0 \Vert^{2}_{L^2} +2\int_{0}^{t} \langle\fe, \vu\rangle_{H^{-1}\times H^1} ds -2\nu\int_{0}^{t} \Vert  \vu(s,\cdot)\Vert^{2}_{\dot{H}^1}ds-2\alpha \int_{0}^{T}\Vert \vu (t,\cdot) \Vert^{2}_{L^2}dt, 
\end{eqnarray*} alors en remplaçant cette estimation du terme $\Vert \vu(t,\cdot)\Vert^{2}_{L^2}$  dans (\ref{theo_stat_ineq_8}) nous pouvons écrire 
 
\begin{eqnarray}\label{theo_stat_ineq_1} \nonumber
\Vert \vv(t,\cdot)\Vert^{2}_{L^2}& \leq & \underbrace{\Vert \vu_0 \Vert^{2}_{L^2}  -2\langle \vu_0,\U\rangle_{L^2\times L^2} +\Vert \U \Vert^{2}_{L^2}}_{(a)}    -2\nu\int_{0}^{t} \Vert \vu(s,\cdot)\Vert^{2}_{\dot{H}^1}dt  -2\alpha \int_{0}^{t}\Vert \vu (s,\cdot) \Vert^{2}_{L^2}dt  \\
& &  \underbrace{2\int_{0}^{t} \langle \fe \cdot \vu\rangle_{H^{-1}\times H^1} ds}_{(b)}  -2 \underbrace{\int_{0}^{t} \langle \partial_t\vv(t,\cdot),\U\rangle_{H^{-2}\times H^2}ds}_{(c)},  
\end{eqnarray} où  nous avons  besoin d'étudier les termes $(a)$, $(b)$ et $(c)$ ci-dessus. Pour le terme $(a)$ nous écrivons directement 
\begin{equation}\label{theo_stat_ineq_2}
\Vert \vu_0 \Vert^{2}_{L^2}  -2\langle \vu_0,\U\rangle_{L^2\times L^2} +\Vert \U \Vert^{2}_{L^2}= \Vert \vu-\U \Vert^{2}_{L^2}.
\end{equation} Ensuite, pour le terme $(b)$ nous allons vérifier que l'on a l'identité suivante
\begin{eqnarray}\label{theo_stat_ineq_3}\nonumber
2\int_{0}^{t}\langle \fe,\vu\rangle_{H^{-1}\times H^1} &= &2\nu \int_{0}^{t}\langle \vec{\nabla}\otimes \U, \vec{\nabla}\otimes \vu\rangle_{L^2 \times L^2}ds+2\int_{0}^{t}\langle (\U\cdot \vec{\nabla})\U, \vu\rangle_{H^{-1}\times H^1}ds\\
& & +2\alpha \int_{0}^{t}\langle\U, \vu\rangle_{L^2 \times L^2}ds.
\end{eqnarray} En effet, comme la solution stationnaire $\U$ vérifie l'équation $\ds{-\nu \Delta \U +\P((\U\cdot \vec{\nabla})\U)+\alpha \U=\fe}$, nous avons 
$$ \langle \fe,\vu\rangle_{H^{-1}\times H^1}=\langle -\nu \Delta \U +\P((\U\cdot \vec{\nabla})\U+\alpha \U,\vu\rangle_{H^{-1}\times H^1}, $$ où, nous intégrons par parties le premier terme à droite et de plus, nous utilisons les propriétés du projecteur de Leray dans le deuxième terme à droite pour obtenir l'identité
$$ \langle \fe,\vu\rangle_{H^{-1}\times H^1} = \nu \langle \vec{\nabla}\otimes \U, \vec{\nabla}\otimes \vu\rangle_{H^{-1}\times H^1} + \langle (\U\cdot \vec{\nabla})\U,  \vu\rangle_{H^{-1}\times H^1}+\alpha \langle \U,\vu\rangle_{L^2\times L^2}.$$ Finalement nous intégrons par rapport au temps et nous avons de cette façon l'identité recherchée  (\ref{theo_stat_ineq_3}). \\
\\
\'Etudions maintenant le terme $(c)$ et pour cela, comme $\vv=\vu-\U$ nous utilisons la relation 
$$ \partial_t \vv =\nu \Delta \vv -\P\left( (\vu\cdot \vec{\nabla})\vu- (\U\cdot \vec{\nabla})\U \right)-\alpha \vv,$$ mais, afin de mener à bien les calculs plus bas, il convient d'écrire le terme non linéaire comme  
$$ (\vu\cdot \vec{\nabla})\vu- (\U\cdot \vec{\nabla})\U= (\vv \cdot \vec{\nabla})\vv +(\U\cdot \vec{\nabla})\vv + (\vv \cdot \vec{\nabla})\U,$$  de sorte que $\vv$ vérifie l'équation 
$$ \partial_t \vv =\nu \Delta \vv -\P\left(  (\vv \cdot \vec{\nabla})\vv +(\U\cdot \vec{\nabla})\vv+ (\vv \cdot \vec{\nabla})\U \right)-\alpha \vv. $$
De cette façon nous avons
\begin{eqnarray*}
\langle \partial_t \vv, \U\rangle_{H^{-2}\times H^2} &=& \langle \nu \Delta \vv, \U\rangle_{H^{-2}\times H^2}- \langle \P((\vv \cdot \vec{\nabla})\vv), \U\rangle_{H^{-2}\times H^2}-\langle \P(\U\cdot \vec{\nabla})\vv), \U\rangle_{H^{-2}\times H^2}\\
& &-\langle \P((\vv \cdot \vec{\nabla})\U), \U\rangle_{H^{-2}\times H^2}-\alpha \langle \vv,\U\rangle_{L^2\times L^2}.
\end{eqnarray*} Dans cette identité, dans les trois premiers termes à droite nous intégrons par parties et nous utilisons les propriétés du projecteur de Leray. De plus, pour le quatrième terme a droite, comme $div(\U)=0$ alors nous avons $\langle \P((\vv \cdot \vec{\nabla})\U), \U\rangle_{H^{-2}\times H^2}=0$ et donc nous obtenons  
\begin{eqnarray*}
\langle \partial_t \vv, \U\rangle_{H^{-2}\times H^2} &=&- \nu \langle \vec{\nabla}\otimes \vv,\vec{\nabla}\otimes \U\rangle_{L^{2}\times L^2} + \langle (\vv \cdot \vec{\nabla})\U, \vv\rangle_{H^{-1}\times H^1}+\langle \U\cdot \vec{\nabla})\U, \vv\rangle_{H^{-1}\times H^1}\\
& &-\alpha \langle \vv,\U\rangle_{L^2\times L^2}.
\end{eqnarray*} Finalement, nous intégrons par rapport à la variable de temps pour obtenir 
\begin{eqnarray}\label{theo_stat_ineq_4} \nonumber
-2\int_{0}^{t} \langle \partial_t \vv, \U\rangle_{H^{-2}\times H^2} ds &=& 2\nu\int_{0}^{t} \langle \vec{\nabla}\otimes \vv,\vec{\nabla}\otimes \U\rangle_{L^{2}\times L^2}ds  - 2\int_{0}^{t}\langle (\vv \cdot \vec{\nabla})\U, \vv\rangle_{H^{-1}\times H^1}ds\\
& &- 2\int_{0}^{t}\langle \U\cdot \vec{\nabla})\U, \vv\rangle_{H^{-1}\times H^1} +2\alpha \int_{0}^{t}\langle \vv,\U\rangle_{L^2\times L^2}.  
\end{eqnarray}
Une fois que nous avons les estimations (\ref{theo_stat_ineq_2}), (\ref{theo_stat_ineq_3}) et (\ref{theo_stat_ineq_4}), nous les remplaçons dans (\ref{theo_stat_ineq_1}) et nous obtenons 

\begin{eqnarray}\label{theo_stat_ineq_5} \nonumber
\Vert \vv(t,\cdot)\Vert^{2}_{L^2} 
& \leq &  \Vert \vu_0-\U \Vert^{2}_{L^2} -2\nu \int_{0}^{t}\Vert \vu(s,\cdot)\Vert^{2}_{\dot{H}^1}ds
+4\nu \int_{0}^{t}\langle \vec{\nabla}\otimes \vv, \vec{\nabla}\otimes \U\rangle_{L^2 \times L^2}ds\\ \nonumber
& & -2\int_{0}^{t}\langle (\vv\cdot \vec{\nabla})\U,\vv\rangle_{H^{-1}\times H^1}-2\alpha\int_{0}^{t}\Vert \vu(s,\cdot)\Vert^{2}_{L^2}ds +2\alpha \int_{0}^{t}\langle \U,\vu\rangle_{L^2\times L^2} ds\\
& & +2\alpha \int_{0}^{t}\langle \vv,\U\rangle_{L^2\times L^2} ds.
\end{eqnarray}
\`A ce stade nous avons les remarques suivantes: tout d'abord, étant donné que $\vv=\vu-\U$ alors, dans le troisième terme à droite de l'estimation ci-dessus  nous pouvons écrire 
$$ 4\nu \int_{0}^{t}\langle \vec{\nabla}\otimes \vv, \vec{\nabla}\otimes \U\rangle_{L^2 \times L^2}ds= 4\nu \int_{0}^{t}\langle \vec{\nabla}\otimes \vu, \vec{\nabla}\otimes \U\rangle_{L^2 \times L^2}ds-4\nu \int_{0}^{t} \Vert \U \Vert^{2}_{L^2}ds, $$ de sorte que 
\begin{eqnarray*}
& & -2\nu \int_{0}^{t}\Vert \vu(s,\cdot)\Vert^{2}_{\dot{H}^1}ds +4\nu \int_{0}^{t}\langle \vec{\nabla}\otimes \vv, \vec{\nabla}\otimes \U\rangle_{L^2 \times L^2}ds \\
& \leq & -2\nu \int_{0}^{t}\Vert \vu(s,\cdot)\Vert^{2}_{\dot{H}^1}ds + 4\nu \int_{0}^{t}\langle \vec{\nabla}\otimes \vu, \vec{\nabla}\otimes \U\rangle_{L^2 \times L^2}ds-4\nu \int_{0}^{t} \Vert \U \Vert^{2}_{L^2}ds\\
& \leq & -2\nu \int_{0}^{t}\Vert \vu(s,\cdot)-\U \Vert^{2}_{\dot{H}^1}ds \leq -2\nu \int_{0}^{t}\Vert \vv(s,\cdot)\Vert^{2}_{\dot{H}^1}ds.    
\end{eqnarray*}
D'autre part, toujours par la relation $\vv=\vu-\U$, dans le dernier terme de (\ref{theo_stat_ineq_5}) nous pouvons écrire
$$ +2\alpha \int_{0}^{t}\langle \vv,\U\rangle_{L^2\times L^2} ds=2\alpha \int_{0}^{t}\langle \vu,\U\rangle_{L^2\times L^2} ds-2\int_{0}^{t} \Vert
\U \Vert^{2}_{L^2} ds,$$ et alors 
\begin{eqnarray*} 
& & -2\alpha\int_{0}^{t}\Vert \vu(s,\cdot)\Vert^{2}_{L^2}ds +2\alpha \int_{0}^{t}\langle \U,\vu\rangle_{L^2\times L^2} ds+2\alpha\int_{0}^{t}\langle \vv,\U\rangle_{L^2\times L^2} ds\\
&\leq & -2\alpha\int_{0}^{t}\Vert \vu(s,\cdot)\Vert^{2}_{L^2}ds +4\alpha \int_{0}^{t}\langle \vu,\U\rangle_{L^2\times L^2} ds -2\alpha\int_{0}^{t} \Vert \U \Vert^{2}_{L^2} ds\\
&\leq & -2\alpha\int_{0}^{t}\Vert \vu(s,\cdot)-\U\Vert^{2}_{L^2}ds \leq -2\alpha \int_{0}^{t}\Vert \vv(s,\cdot)\Vert^{2}_{L^2}ds.
\end{eqnarray*}  
De cette façon, en remplaçant les estimations ci-dessus dans (\ref{theo_stat_ineq_5}) nous obtenons l'inégalité d'énergie recherchée. \finpv

\subsubsection{Preuve du Lemme \ref{Lemme:bessel} page \pageref{Lemme:bessel}}\label{sec:preuve-lemme-bessel}
Soit  $G_{\nu,\alpha}>0$ le noyau du potentiel de Bessel $\frac{1}{-\nu \Delta +\alpha I_d}[\cdot]$ qui vérifie l'estimation suivante:  pour $c_{\nu,\alpha}>0$  une constante qui dépend seulement de $\nu>0$ et $\alpha>0$ on a 
\begin{equation}\label{dec-noyau-bessel-aux}
G_{\nu,\alpha}(x)\leq c_{\nu,\alpha} \left\lbrace  \begin{array}{cl} \vspace{2mm}
\frac{1}{\vert x \vert},& \text{si} \,\, \vert x \vert \leq 2, \\
e^{-\frac{\vert x \vert}{2}}, & \text{si} \,\, \vert x \vert > 2.
\end{array}	 \right.
\end{equation}  Soit $g \in L^1(\Rt)$ qui vérifie une décroissance 
\begin{equation}\label{dec-g}
\vert g(x)\vert \leq \frac{c}{\vert x \vert^n},
\end{equation}  pour $n \in \mathbb{N}$ avec $n\geq 4$ et pour tout $\vert x \vert >8$; et soit enfin $x\in \Rt$ fixe tel que $\vert  x \vert >8$. Nous écrivons 
\begin{equation}\label{estim}
\vert G_{\nu,\alpha}\ast g (x)\vert  =  \int_{\vert y \vert \leq \frac{\vert x \vert}{2}} \vert G_{\nu,\alpha}(x-y)\vert \vert g(y)\vert dy + \int_{\vert y \vert > \frac{\vert x \vert}{2}} \vert G_{\nu,\alpha}(x-y)\vert \vert g(y)\vert dy=I_1+I_2,
\end{equation}
 et nous allons estimer chaque terme ci-dessus. Pour le terme $I_1$, comme l'on considère $\vert y \vert \leq \frac{\vert x \vert }{2}$ alors l'on a $\vert x-y \vert \geq \frac{\vert x \vert}{2}$, et de plus, étant donné que $\vert x \vert >8$ alors par l'inégalité précédente l'on a aussi $\vert x-y \vert >4 $ et ainsi, par l'estimation (\ref{dec-noyau-bessel-aux}) et comme $g\in L^1(\Rt)$ nous pouvons écrire 
 \begin{equation}\label{estim-I1}
 I_1 \leq c_{\nu,\alpha} \int_{\vert y \vert \leq \frac{\vert x \vert}{2}} e^{-\frac{\vert x-y \vert}{2}}\vert g(y)\vert dy \leq c_{\nu,\alpha} e^{-\frac{\vert x \vert}{2}}\Vert g \Vert_{L^1}. 
 \end{equation} Pour le terme $I_2$, nous écrivons:
 $$ I_2=  \int_{\vert y \vert > \frac{\vert x \vert}{2}, \vert x-y\vert\leq 2} \vert G_{\nu,\alpha}(x-y)\vert \vert g(y)\vert dy + \int_{\vert y \vert > \frac{\vert x \vert}{2},\vert x-y\vert >2} \vert G_{\nu,\alpha}(x-y)\vert \vert g(y)\vert dy=(I_2)_a+(I_2)_b,$$ et nous devons encore estimer chaque terme ci-dessus. Pour le terme $(I_2)_a$, comme l'on considère ici $\vert x-y \vert \leq 2 $ alors par l'estimation (\ref{dec-noyau-bessel-aux}) nous savons que l'on a $G_{\nu,\alpha}(x-y)\leq \frac{c_{\nu,\alpha}}{\vert x-y \vert}$. De plus, comme l'on considère aussi $ \vert y \vert > \frac{\vert x \vert}{2} $ alors pour l'estimation (\ref{dec-g}) nous savons que l'on a $\vert g (y )\vert \leq \frac{c}{\vert y \vert^n}\leq \frac{c}{\vert x \vert^n}$; et alors nous pouvons écrire 
 $$ (I_2)_a \leq  \int_{\vert y \vert > \frac{\vert x \vert}{2}, \vert x-y\vert\leq 2}\vert G_{\nu,\alpha}(x-y)\vert \vert g(y)\vert dy \leq \frac{c_{\nu,\alpha}}{\vert x \vert^n}\int_{\vert x-y\vert \leq 2} \frac{dy}{\vert x-y\vert}\lesssim  \frac{c_{\nu,\alpha}}{\vert x \vert^n}.$$  
 Ensuite, pour le terme  $(I_2)_b$, comme nous considérons ici $\vert x-y \vert >2$ alors nous écrivons
 $$ (I_2)_b \leq \int_{\vert y \vert > \frac{\vert x \vert}{2},\vert x-y\vert >2} \vert G_{\nu,\alpha}(x-y)\vert \vert g(y)\vert dy \leq \frac{c_{\nu,\alpha}}{\vert x \vert^n}\int_{\vert x-y\vert>2} e^{-\frac{\vert x-y\vert}{2}}dy \lesssim \frac{c_{\nu,\alpha}}{\vert x \vert^n};$$ et par les estimations ci-dessus nous avons 
 \begin{equation}\label{estim-I2}
 I_2 \leq \frac{c_{\nu,\alpha}}{\vert x \vert^n}.
 \end{equation} Une fois que l'on a vérifié les estimations
 (\ref{estim-I1}) et (\ref{estim-I2}), on pose la constante $c_2=\max(c_{\nu,\alpha}\Vert g \Vert_{L^1},c_{\nu,\alpha})>0$, et nous revenons  à l'identité (\ref{estim}) pour écrire $\ds{\vert G_{\nu,\alpha}\ast g (x)\vert \leq \frac{c_2}{\vert x \vert^n}}$.  \finpv

\subsubsection{Preuve du Lemme \ref{Lemme:dec-K} page \pageref{Lemme:dec-K}}\label{sec:preuve-lemme-dec-K}
Nous allons montrer ici l'estimation (\ref{dec-noyau-K}). Rappelons tout d'abord que par la définition du noyau $K_{\nu,\alpha}$ donnée dans l'expression (\ref{def:K}) nous savons que 
$\ds{K_{\nu,\alpha}= \frac{\nu \Delta}{-\nu \Delta +\alpha I_d} \left[  m\right]}$, où le tenseur $m=(m_{i,j,k})_{1\leq i,j,k\leq 3})$ est défini par les expressions (\ref{def:m-fourier}) et (\ref{def:m}); et alors  nous pouvons  écrire $\ds{K_{\nu,\alpha}= \frac{1}{-\nu \Delta +\alpha I_d} \left[ \nu \Delta m\right]}$. Ensuite,  
étant donné que l'action du potentiel de Bessel $\ds{\frac{1}{-\nu \Delta +\alpha I_d} \left[ \cdot\right]}$ s'écrite comme un produit de convolution avec le noyau $G_{\nu,\alpha}$ donné (\ref{dec-noyau-bessel-aux}) nous écrivons alors $\ds{K_{\nu,\alpha}=G_{\nu,\alpha}\ast (\nu \Delta m)}$. Ainsi, pour vérifier l'estimation (\ref{dec-noyau-K}) nous allons maintenant  montrer l'estimation 
\begin{equation}\label{estima1}
\vert G_{\nu,\alpha}\ast (\nu \Delta m)(x)\vert \leq \frac{c_{\nu,\alpha}}{\vert x \vert},\quad \text{pour}\quad \vert x \vert \leq 4;
\end{equation}
et ensuite nous allons vérifier l'estimation  
\begin{equation}\label{estima2}
\vert G_{\nu,\alpha}\ast (\nu \Delta m)(x)\vert \leq \frac{c_{\nu,\alpha}}{\vert x \vert^4},\quad \text{pour}\quad \vert x \vert >4.
\end{equation}
Pour montrer l'estimation (\ref{estima1}) nous écrivons $\ds{G_{\nu,\alpha}\ast (\nu \Delta m)(x)=\Delta G_{\nu,\alpha}\ast (\nu m)(x)}$. D'autre part, observons que l'on a les estimations 
\begin{equation}\label{estimG}
\vert \Delta G_{\nu,\alpha}(x)\vert \leq \frac{1c_{\nu,\alpha}}{\vert x \vert^2},\quad \text{pour}\quad \vert x \vert>0 \quad \text{et}\quad \vert x \vert \neq 2,
\end{equation} et 
\begin{equation}\label{estimM}
\vert m(x)\vert \leq \frac{c}{\vert x \vert^2},\quad \text{pour}\quad \vert x \vert >0.
\end{equation}
 En effet, l'estimation (\ref{estimG})  est une conséquence directe de la définition du noyau $G_{\nu,\alpha}$ donnée dans (\ref{dec-noyau-K}) tandis que l'estimation (\ref{estimM})  vient  du fait que le tenseur le tenseur $m$ est défini comme $m=(m_{i,j,k})_{1\leq i,j,k\leq 3}$ où pour chaque $1\leq i, j,k\leq 3$ la fonction $m_{i,j,k}\in \mathbb{C}^{\infty}(\Rt \setminus \{0\} ) $ définie dans (\ref{def:m-fourier}) est une fonction homogène de degré $-2$. De cette façon, pour tout $\vert x \vert >0$ nous pouvons écrire 
$$ \vert  \Delta G_{\nu,\alpha}\ast (\nu m)(x) \vert \leq c_{\nu,\alpha} \int_{\Rt} \frac{dy}{\vert x-y \vert^2 \vert y \vert^2} \leq \frac{c_{\nu,\alpha}}{\vert x \vert},$$ et en particulier cette estimation est valide pour $0<\vert x \vert \leq 4$, ce qui nous donne l'estimation (\ref{estima1}). \\
\\
Vérifions maintenant l'estimation (\ref{estima2}).  Soit $\vert x \vert >4$ fixe. On définit la fonction de troncature $\varphi_x$ comme $\varphi_x \in \mathbb{C}^{\infty}_{0}(\Rt)$ telle que: $0\leq \varphi_x \leq 1$, $\varphi_x(y)=1$ si $\vert y \vert \leq \frac{\vert x \vert}{4}$ et $\varphi_x(y)=0$ si $\vert y \vert \geq \frac{\vert x \vert}{2}$, $ \ds{\Vert \vec{\nabla}\varphi_x\Vert_{L^{\infty}} \leq \frac{c}{\vert x \vert}}$ et $\ds{\Vert \Delta \varphi_x\Vert_{L^{\infty}}\leq \frac{c}{\vert x \vert^2}} $; et nous écrivons
\begin{equation}\label{estimat1}
G_{\nu,\alpha}\ast (\nu \Delta m)=G_{\nu,\alpha}\ast [\varphi_x (\nu \Delta m)]+G_{\nu,\alpha}\ast [(1-\varphi_x) (\nu \Delta m)]=I_1+I_2,
\end{equation} où nous allons monter que $\ds{I_1\leq \frac{c_{\nu,\alpha}}{\vert x \vert^4}}$ et $\ds{I_2\leq \frac{c_{\nu,\alpha}}{\vert x \vert^4}}$  et pour cela nous allons suivre les idées  de la preuve du Lemme  \ref{Lemme:bessel}. En effet, nous écrivons chaque produit de convolution ci-dessus comme une intégrale sur $\Rt$ et nous divisons cette intégrale en deux morceaux: $\{\vert y \vert \leq \frac{\vert x \vert}{2}\}$ et $\{\vert y \vert >\frac{\vert x \vert}{2}\}$. Pour le terme $I_1$, par la définition de la fonction $\varphi_x$ nous savons que $\varphi_x(y)=0$ si $\vert y \vert >\frac{\vert x \vert}{2}$ et alors nous avons 
$$ I_1=\int_{\vert y \vert\leq \frac{\vert x \vert}{2}} G_{\nu,\alpha}(x-y)[\varphi_x(y)(\nu \Delta m (y))]dy, $$ et en intégrant par parties nous pouvons écrire 
$$ I_1= \int_{\vert y \vert\leq \frac{\vert x \vert}{2}}\Delta[G_{\nu,\alpha}(x-y)\varphi_x(y)]m(y)dy.$$
Ensuite, nous développons le terme $\Delta[G_{\nu,\alpha}(x-y)\varphi_x(y)]$ pour écrire 
\begin{eqnarray}\label{estimat2}\nonumber
I_1&=&\int_{\vert y \vert \leq \frac{\vert x \vert}{2}} \Delta G_{\nu,\alpha}(x-y)\varphi_x(y)m(y)dy+ 2 \int_{\vert y \vert \leq \frac{\vert x \vert}{2}} [\vec{\nabla} G_{\nu,\alpha}(x-y) \cdot  \vec{\nabla}\varphi_x(y)] m(y)dy \\
& & +\int_{\vert y \vert \leq \frac{\vert x \vert}{2}}  G_{\nu,\alpha}(x-y)\Delta \varphi_x(y) m(y)dy=(I_1)_a+(I_1)_b+(I_1)_c.
\end{eqnarray} Nous devons encore estimer ces trois termes ci-dessus mais avant de faire ça il convient de rappeler tout d'abord que l'on considère ici $\vert y \vert \geq \frac{\vert x \vert}{2}$ et comme l'on a fixé $\vert x \vert >4$ alors on a $\vert x-y \vert \leq \frac{\vert x \vert}{2} >2$ et alors par l'estimation (\ref{dec-noyau-bessel-aux}) on a
\begin{equation}\label{estimg}
G_{\nu,\alpha}(x-y)\leq c_{\nu,\alpha}e^{-\frac{\vert x-y\vert}{2}}.
\end{equation}  Avec cette estimation en tête nous étudions maintenant les trois termes ci-dessus. \\
\\
Pour le terme $(I_1)_a$: par l'estimation (\ref{estimg}) nous pouvons écrire  $\vert \Delta G_{\nu,\alpha}(x-y)\vert  \leq \frac{c_{\nu,\alpha}}{\vert x-y \vert^5}$ et comme $\vert x-y \vert \geq \frac{\vert x \vert}{2}$ nous avons $\vert \Delta G_{\nu,\alpha}(x-y)\vert \leq \frac{c_{\nu,\alpha}}{\vert x \vert^5}$. Ensuite, par l'estimation (\ref{estimM}) nous savons que $\vert m(y)\vert \leq \frac{c}{\vert y \vert^2}$. De plus, étant donné que $\Vert \varphi_x \Vert_{L^{\infty}}\leq 1$, alors nous écrivons 
\begin{equation}\label{I1a}
(I_1)_a \leq \frac{c_{\nu,\alpha}}{\vert x \vert^5}  \int_{\vert y \vert \leq \frac{\vert x \vert}{2}} \frac{dy}{\vert y \vert^2}\leq \frac{c_{\nu,\alpha}}{\vert x \vert^4}.
\end{equation}
Pour le terme $(I_1)_b$: en suit les mêmes idées pour estimer le terme $(I_1)_a$ ci-dessus: toujours par l'estimation (\ref{estimg}) nous pouvons écrire $\vert \vec{\nabla}G_{\nu,\alpha}(x-y)\vert \leq \frac{c_{\nu,\alpha}}{\vert x \vert^4}$ et comme nous avons aussi  $\Vert \vec{\nabla}\varphi_x\Vert_{L^{\infty}}\leq \frac{c}{\vert x \vert}$ par la définition de $\varphi_x$, alors nous écrivons  
\begin{equation}\label{I1b}
(I_1)_b \leq \frac{c_{\nu,\alpha}}{\vert x \vert^5} \int_{\vert y \vert \leq \frac{\vert x \vert}{2}} \frac{dy}{\vert y \vert^2}\leq \frac{c_{\nu,\alpha}}{\vert x \vert^4}. 
\end{equation} 
Pour le terme $(I_1)_c$: nous écrivons ici $\vert G_{\nu,\alpha}(x-y)\vert \leq \frac{c_{\nu,\alpha}}{\vert x-y \vert^3}\leq \frac{c_{\nu,\alpha}}{\vert x \vert^3}$ et comme $\Vert \Delta \varphi_x\Vert _{L^{\infty}}\leq \frac{c}{\vert x \vert^2}$ alors nous avons:
\begin{equation}\label{I1c}
(I_1)_b \leq \frac{c_{\nu,\alpha}}{\vert x \vert^5} \int_{\vert y \vert \leq \frac{\vert x \vert}{2}} \frac{dy}{\vert y \vert^2}\leq \frac{c_{\nu,\alpha}}{\vert x \vert^4}.
\end{equation} Ainsi, en revenant à l'estimation (\ref{estimat2}), par les estimations  (\ref{I1a}), (\ref{I1b}) et (\ref{I1c}) nous avons $\ds{I_1\leq \frac{c_{\nu,\alpha}}{\vert x \vert^4}} $.\\
\\
Vérifions maintenant que l'on a $\ds{I_2\leq \frac{c_{\nu,\alpha}}{\vert x \vert^4}}$, où rappelons que le terme $I_2$ est donné dans (\ref{estimat1}) par l'expression $\ds{I_2=}G_{\nu,\alpha}\ast [(1-\varphi_x) (\nu \Delta m)]$. Observons tout d'abord que dans ce terme nous avons la fonction $1-\varphi_x$, où étant donné que $\varphi_x(y)=1$ si $\vert y \vert \leq \frac{\vert x \vert}{4}$ alors on a que  $\vert y \vert \leq \frac{\vert x \vert}{4}(y)=0$ si  $\vert y \vert \leq \frac{\vert x \vert}{4}$ et donc nous écrivons 
$$ I_2= \int_{\vert y \vert > \frac{\vert x \vert}{4}} G_{\nu,\alpha}(x-y)(1-\varphi_x(y))\nu \Delta m(dy)dy.$$ Ensuite, nous avons besoin d'étudier un peu plus le terme $\Delta m$, où,  étant donné que pour chaque $1\leq i, j,k\leq 3$ la fonction $m_{i,j,k}\in \mathbb{C}^{\infty}(\Rt \setminus \{0\} ) $  est une fonction homogène de degré $-2$ alors  la fonction $\Delta m_{i,j,k}$ est une fonction homogène de degré $-4$ et donc, comme le tenseur $m$ est définie par $m=(m_{i,j,k})_{1\leq i,j,k \leq 3}$,  nous avons l'estimation $\vert \nu \Delta m(y)\vert \lesssim \frac{1}{\vert y \vert^4};
$ et comme nous considérons ici $\vert y \vert > \frac{\vert x \vert}{4}$ nous avons l'estimation $\vert \nu \Delta m(y)\vert \lesssim \frac{1}{\vert y \vert^4} \lesssim \frac{1}{\vert x \vert^4}.$ Avec cette estimation en tête nous revenons à l'identité ci-dessus pour écrire  $\ds{I_2 \leq \frac{1}{\vert x \vert^4} \Vert G_{\nu,\alpha}\Vert_{L^1}\leq \frac{c_{\nu,\alpha}}{\vert x \vert^4}} $. \finpv
\subsubsection{Preuve du Lemme \ref{Lemme:dec-U-aux}  page \pageref{Lemme:dec-U-aux}}\label{sec:preuve-lemme-dec-U-aux}
Expliquons tout d'abord l'idée de la preuve de ce lemme.  Nous allons tout d'abord montrer que le terme $\ds{K_{\nu,\alpha}\ast (\U\otimes \U)}$ vérifie une décroissance 
\begin{equation}\label{estim-K-UU-aux}
\vert K_{\nu,\alpha}\ast (\U\otimes \U)(x)\vert \leq \frac{c_2}{\vert x \vert^2}, 
\end{equation} pour $\vert x \vert >4$. D'autre part, par l'estimation (\ref{estimation2}) nous savons  pour $\vert x \vert >4$ l'on a 
\begin{equation}\label{estim:force}
\left\vert \frac{1}{-\nu \Delta +\alpha I_d}\left[\fe\right](x)\right\vert \leq \frac{c_1}{\vert x \vert^4},
\end{equation} d'où nous pouvons en tirer l'estimation 
\begin{equation}\label{estim:force:aux}
\left\vert \frac{1}{-\nu \Delta +\alpha I_d}\left[\fe\right](x)\right\vert \leq \frac{c_1}{\vert x \vert^2},
\end{equation} et étant donné que l'on a  l'identité 
\begin{equation}\label{ident:U}
\U= K_{\nu,\alpha}\ast (\U\otimes \U)+\frac{1}{-\nu \Delta +\alpha I_d}\left[\fe\right],
\end{equation} nous obtiendrons ainsi que la solution $\U$ vérifie une décroissance:  pour tout $\vert x \vert >4$,
\begin{equation*}\label{estim-aux-U}
\vert \U (x)\vert \leq \frac{c_2}{\vert x \vert^2},
\end{equation*} 
Vérifions donc l'estimation (\ref{estim-K-UU-aux}).  Pour $\vert x \vert >4$ fixe on commence par écrire 
\begin{eqnarray}\label{estim1} \nonumber
\vert K_{\nu,\alpha}\ast (\U \otimes \U)(x)\vert &\leq & \int_{\Rt} \vert K_{\nu,\alpha}(x-y)\vert \vert (\U \otimes \U)(y)\vert dy=\int_{\vert y \vert \leq \frac{\vert x \vert}{2}} \vert K_{\nu,\alpha}(x-y)\vert \vert (\U \otimes \U)(y)\vert dy \\
& & + \int_{\vert y \vert >2 \frac{\vert x \vert}{2}} \vert K_{\nu,\alpha}(x-y)\vert \vert (\U \otimes \U)(y)\vert dy= I_1+I_2,
\end{eqnarray} et l'on cherche à estimer les termes $I_1$ et $I_2$ ci-dessus. \\
\\
Pour le terme $I_1$, comme nous considérons ici $\vert y \vert \leq \frac{\vert x \vert}{2}$ alors nous avons les inégalités suivantes: $\vert x-y\vert \geq \vert x \vert -\vert y \vert \geq \frac{\vert x \vert}{2}$, d'où, étant donné que $\vert x \vert >4$  nous écrivons $\frac{\vert x \vert }{2}>2$ pour obtenir $\vert x-y \vert >2 $ et ainsi, par l'estimation (\ref{dec-noyau-K}) obtenue dans le Lemme \ref{Lemme:dec-K} nous avons $\ds{K_{\nu,\alpha}(x-y)\leq \frac{c_{\nu,\alpha}}{\vert x-y\vert^4} \leq \frac{c_{\nu,\alpha}}{\vert x\vert^4}}$, et  nous pouvons écrire 
\begin{equation}\label{estim:I_1}
I_1  =\int_{\vert y \vert \leq \frac{\vert x \vert}{2}} \vert K_{\nu,\alpha}(x-y)\vert \vert (\U \otimes \U)(y)\vert dy\leq \frac{c_{\nu,\alpha}}{\vert x \vert^4} \int_{\vert y \vert \leq \frac{\vert x \vert}{2}} \vert \vert (\U \otimes \U)(y)\vert dy \leq \frac{c_{\nu,\alpha}}{\vert x \vert^4} \Vert \U \Vert^{2}_{L^2}. 
\end{equation}  Ainsi, comme l'on a $\vert x \vert >4$ nous pouvons écrire 
\begin{equation}\label{estim:I1-aux}
I_1 \leq \frac{c_{\nu,\alpha}}{\vert x \vert^2} \Vert \U \Vert^{2}_{L^2}.
\end{equation}
Pour le terme $I_2$, on commence par écrire 
\begin{eqnarray}\label{estim:I2}\nonumber
I_2&=&\int_{\vert y \vert > \frac{\vert x \vert}{2},\, \vert x-y\vert \leq 2} \vert K_{\nu,\alpha}(x-y)\vert \vert (\U \otimes \U)(y)\vert dy+\int_{\vert y \vert > \frac{\vert x \vert}{2},\, \vert x-y\vert > 2} \vert K_{\nu,\alpha}(x-y)\vert \vert (\U \otimes \U)(y)\vert dy\\
&=& (I_2)_a+(I_2)_b,
\end{eqnarray} et on doit estimer chaque terme de l'identité ci-dessus. Pour estimer le premier terme $(I_2)_{a}$ nous allons utiliser la décroissance de la solution $\U$ donnée dans (\ref{dec-u}). En effet, par cette estimation nous pouvons écrire $\ds{ \vert (\U\otimes \U)(y)\vert \leq  c \left( \frac{1}{1+\vert y \vert}\right)^2}$, et comme nous considérons ici $\vert y \vert >\frac{\vert x \vert }{2}$ (et comme $\vert x \vert >4$ alors on a aussi $\vert y \vert >2$) alors nous avons l'estimation $\ds{\vert (\U\otimes \U)(y)\vert \leq  \frac{c}{\vert y \vert^2}\leq \frac{c}{\vert x \vert^2}}$. D'autre part, comme nous considérons aussi $\vert x-y \vert \leq 2$, toujours par  l'estimation (\ref{dec-noyau-K}) (obtenue dans le Lemme \ref{Lemme:dec-K}) nous savons que $\ds{\vert K_{\nu,\alpha}(x-y) \vert \leq \frac{c_{\nu,\alpha}}{\vert x-y\vert}}$; et nous obtenons de cette façon 
\begin{equation}\label{estim:I2a}
(I_2)_a =\int_{\vert y \vert >2 \frac{\vert x \vert}{2},\, \vert x-y\vert \leq 2} \vert K_{\nu,\alpha}(x-y)\vert \vert (\U \otimes \U)(y)\vert dy \leq \frac{c_{\nu,\alpha}}{\vert x \vert^2} \int_{\vert x-y\vert \leq 2} \frac{dy}{\vert x-y\vert} \lesssim  \frac{c_{\nu,\alpha}}{\vert x \vert^2}. 
\end{equation}	  Pour estimer maintenant le terme $(I_2)_b$, rappelons que nous considérons ici $\vert x-y\vert >2$ et alors toujours par l'estimation (\ref{dec-noyau-K}) nous savons que $\ds{K_{\nu,\alpha}(x-y)\leq \frac{c_{\nu,\alpha}}{\vert x-y\vert^4}}$ et ainsi, et suivant les mêmes lignes que l'estimation (\ref{estim:I_1}) nous obtenons 
\begin{equation}
(I_2)_b \leq  \frac{c_{\nu,\alpha}}{\vert x \vert^4}\Vert \U \Vert^{2}_{L^2},
\end{equation} et comme $\vert x \vert >4$ nous pouvons alors écrire
\begin{equation}\label{estim:I2b}
(I_2)_b \leq  \frac{c_{\nu,\alpha}}{\vert x \vert^2}\Vert \U \Vert^{2}_{L^2}.
\end{equation} Une fois que nous disposons des estimations (\ref{estim:I2a}) et (\ref{estim:I2b}) nous posons maintenant la constante $c_{2}=\max(c_{\nu,\alpha},c_{\nu,\alpha}\Vert \U \Vert^{2}_{L^2})$ et en revenant  à l'identité (\ref{estim:I2}) nous pouvons écrire $\ds{I_2\leq \frac{c_{2}}{\vert x \vert^2}}$. Ainsi, par cette estimation et par l'estimation (\ref{estim:I1-aux}) nous revenons maintenant à l'identité (\ref{estim1}) où nous obtenons l'estimation cherchée (\ref{estim-K-UU-aux}). \finpv

	\chapter{Décroissance fréquentielle}
Dans le chapitre précédent nous avons étudié les équations de Navier-Stokes stationnaires  avec un terme d'amortissement supplémentaire et ce terme nous  a permis d'étudier quelques propriétés intéressantes des solutions de ces équations.
Néanmoins, il y a d'autres propriétés des solutions, comme celle que nous étudierons dans ce chapitre,  que l'on peut étudier que ce soit dans le cadre des équations amorties ou dans le cadre des équations classiques (sans terme d'amortissement); et alors nous  allons préférer les équations classiques car ces équations correspondent à un modèle physique plus réaliste que les équations amorties.\\ 
\\ 
Dans ce chapitre nous revenons donc aux équations de Navier-Stokes stationnaires mais cette fois-ci sans aucun terme supplémentaire d'amortissement:  
\begin{equation}\label{N-S-stationnaire}
-\nu \Delta \U+(\U\cdot \vec{\nabla}) \U+\vec{\nabla} P =\fe, \qquad div(\U)=0,\\
\end{equation} et dans le cadre de ces équations nous allons étudier un autre problème relié à l'étude déterministe de la turbulence qui porte sur la décroissance de la transformée de Fourier de la solution $\U$. \\
\\
Comme nous l'expliquerons plus en détail ci-après, si le fluide est en régime laminaire  on s'attend à observer que la quantité $\ds{\left\vert \widehat{\U}(\xi) \right\vert}$ a une décroissance du type suivant
\begin{equation}\label{dec_fourier_U}
\left\vert \widehat{\U}(\xi) \right\vert \lesssim e^{-\vert \xi \vert},
\end{equation} aux fréquences $\vert \xi \vert >0$ et cette décroissance  correspond au fait que dans ce régime les forces de viscosité dissipent rapidement l'énergie cinétique introduite  par la force $\fe$, tandis que  si le fluide est en régime turbulent alors les forces de viscosité du fluide n'interviennent qu'aux hautes  fréquences et ainsi on s'attend à observer  cette même décroissance (\ref{dec_fourier_U}) mais seulement lorsque $\vert \xi \vert>>1$. En effet, selon la théorie K41, dans le régime turbulent on s'attend à observer un intervalle de fréquences  (nommé l'intervalle d'inertie) où la quantité $\left\vert \widehat{\U}(\xi) \right\vert$ est censée  avoir tout d'abord  une  décroissance polynomiale pour ensuite avoir une décroissance exponentielle aux très hautes fréquences.\\ 
\\
Ainsi, le but de ce chapitre est alors d'étudier la décroissance (\ref {dec_fourier_U}) dans le cadre des équations (\ref{N-S-stationnaire}) sans aucun  terme d'amortissement supplémentaire et pour cela dans la Section \ref{Sec:existence-sol-stattionnaires}  nous allons  tout d'abord rappeler  un résultat classique d'existence  des solutions $\U$. Ensuite, dans la Section \ref{Sec:motivation-K41} nous allons  revenir  à la théorie de la turbulence K41 pour expliquer  de façon plus précise pourquoi la solution $\U$ est censée avoir la décroissance fréquentielle (\ref{dec_fourier_U}).\\
\\
Dans la Section \ref{Sec:decroissance-freq} on s'attaquera finalement à l'étude déterministe de la décroissance  (\ref{dec_fourier_U}) et en suivant les idées du chapitre précédent nous allons étudier cette  décroissance  en regardant tout d'abord le régime laminaire et ensuite le régime turbulent. \\
\\
Rappelons que  le régime laminaire  a été caractérisé dans le chapitre précédent par un contrôle sur le nombre de Grashof $G_{\theta}$ qui induisait comme nous l'avons vu  un contrôle sur la force $\fe$, tandis que dans le régime turbulent  on ne suppose  aucun contrôle  sur le nombre $G_{\theta}$ ce qui laissera la force libre; et nous pouvons alors observer  que l'idée sous-jacente pour caractériser le régime laminaire ou turbulent du fluide repose essentiellement sur le fait  de contrôler ou non la force extérieure.\\
\\
Dans ce chapitre nous allons en revanche caractériser le régime laminaire du fluide en supposant un contrôle \emph{direct} sur la force $\fe$  tandis que dans le régime turbulent nous ne ferons aucun contrôle sur cette fonction et nous n'utiliserons plus les nombres de Grashof.  Cette façon de  caractériser le mouvement du fluide va nous permettre de considérer des forces $\fe$ dans des différents cadres fonctionnels et ainsi d'étudier la décroissance (\ref{dec_fourier_U}) des solutions $\U$  dans ces  cadres fonctionnels associés aux forces.  

\section{Introduction} 
Comme annoncé  nous allons tout d'abord rappeler un résultat général  d'existence  des solutions des équations de Navier-Stokes stationnaires (\ref{N-S-stationnaire}) pour ensuite expliquer le type de décroissance fréquentielle qu'il est possible d'obtenir dans ce cadre.

\subsection{Les équations de Navier-Stokes stationnaires}\label{Sec:existence-sol-stattionnaires}
Dans la Section \ref{Sec:existence} du chapitre précédent nous avons vérifie un résultat général sur l'existences des solutions des équations de Navier-Stokes stationnaire et amorties; et dans cette section nous suivons ces mêmes idées mais dans le cadre des équations classiques (\ref{N-S-stationnaire}). \\
\\
Pour trouver un espace fonctionnel où l'on puisse construire des solutions $(\U,P)$ des ces équations  nous allons faire les estimations \emph{a priori}  suivantes: si l'on suppose pour l'instant que la vitesse $\U$ est une fonction suffisamment  régulière et intégrable;  et si l'on multiplie l'équation ci-dessus par $\U$ et puis l'on intègre en variable d'espace on obtient (formellement) l'identité 
$$ -\nu \int_{\Rt} \Delta \U \cdot \U dx+ \int_{\Rt} \left[ (\U\cdot \vec{\nabla}) \U \right] \cdot \U dx+ \int_{\Rt} \vec{\nabla} P \cdot \U dx=\int_{\Rt} \fe \cdot \U dx, $$ d'où comme $div(\U)$ on a  (toujours formellement) les identités $$ \int_{\Rt} \left[ (\U\cdot \vec{\nabla}) \U \right] \cdot \U dx=0 \quad \text{et}\quad  \int_{\Rt} \vec{\nabla} P \cdot \U dx=0,$$ et ainsi, par une intégration par parties nous écrivons 
$$ \nu \int_{\Rt} \vert \vec{\nabla}\otimes \U \vert^2dx=\int_{\Rt} \fe \cdot \U dx.$$
Dans cette identité nous pouvons observer que si l'on considère une force $\fe \in \dot{H}^{-1}(\Rt)$ et si $\U \in\dot{H}^{1}(\Rt) $ alors en appliquant l'inégalité de Cauchy-Schwarz dans le terme à droite ci-dessus on a $$ \nu\Vert \U \Vert^{2}_{\dot{H}^{1}} \leq \Vert \fe \Vert_{\dot{H}^{-1}}\Vert \U \Vert_{\dot{H}^1},$$ d'où on obtient l'estimation \emph{a priori} 
\begin{equation}\label{estimation-a-priori-sol-estationnaires}
\nu \Vert \U \Vert_{\dot{H}^1}\leq \Vert \fe \Vert_{\dot{H}^{-1}}.
\end{equation}
Nous observons donc que l'espace $\dot{H}^{-1}(\Rt)$ est un espace fonctionnel naturel pour la force $\fe$ et dans ce cadre on veut construire  des solutions $\U$ dans l'espace $\dot{H}^{1}(\Rt)$. De plus, étant donné que la pression $P$ est toujours reliée à la vitesse $\U$ par la relation $P=\frac{1}{-\Delta}div((\U\cdot \vec{\nabla}) \U)$ alors on a $P \in \dot{H}^{\frac{1}{2}}(\Rt)$. En effet, comme $div(\U)=0$ nous écrivons le terme non linéaire $(\U\cdot \vec{\nabla}) \U=div(\U \otimes \U)$ et si $\U\in \dot{H}^{1}(\Rt)$ alors par le lois de produit on a $\U\otimes \U \in \dot{H}^{\frac{1}{2}}(\Rt)$ et donc 
\begin{equation}\label{terme-non-lin}
div (\U \otimes \U)\in \dot{H}^{-\frac{1}{2}}(\Rt); 
\end{equation}
et ainsi,  par la relation ci-dessus on obtient (formellement)  $P \in \dot{H}^{\frac{1}{2}}(\Rt)$.   \\
\\
L'existence des solutions $(\U,P) \in \dot{H}^{1}(\Rt)\times \dot{H}^{\frac{1}{2}}(\Rt)$ pour n'importe quelle force $\fe \in \dot{H}^{-1}(\Rt)$ est un résultat classique qui a été développé dans les travaux de R. Finn \cite{Finn} en 1961 et les travaux de O. Ladyzhenskaya \cite{Ladyzhenskaya1} en 1959 et \cite{Ladyzhenskaya2} en 1963. Ces résultats sont basés sur l'estimation (\ref{estimation-a-priori-sol-estationnaires})  et nous avons:
\begin{Theoreme}\label{Theo:exstence-sol-stationnaires} Soit $\fe \in \dot{H}^{-1}(\Rt)$ une force à divergence nulle. Il existe $(\U, P)\in \dot{H}^{1}(\Rt
)\times \dot{H}^{\frac{1}{2}}(\Rt)$ solution des équations de Navier-Stokes stationnaires  (\ref{N-S-stationnaire}) qui vérifie l'estimation (\ref{estimation-a-priori-sol-estationnaires}). 
\end{Theoreme}
Il s'agit donc d'un résultat général  d'existence de solutions où l'on ne fait aucune hypothèse supplémentaire sur la force $\fe$. Néanmoins, comme nous pouvons observer ce résultat  ne donne aucune information sur l'unicité de ces solutions. \\ 
\\
D'autre part, si nous comparons ce résultat avec le Théorème \ref{Theo:solutions_stationnaires_turbulent} où l'on considère les équations amorties, nous pouvons observer le terme d'amortissement ($-\alpha \U$ avec $\alpha>0$) entraîne que les solutions des équations amorties appartiennent  à l'espace $\dot{H}^{1}\cap L^2(\Rt)$, tandis que dans le résultat ci-dessus nous observons que les solutions des équations (\ref{N-S-stationnaire}) appartiennent seulement à l'espace $\dot{H}^1(\Rt)$. \\  
\\   
Pour une démonstration de ce résultat voir le livre \cite{PGLR1} (Théorème $16.2$ page 530).  L'idée de la preuve suit  les grandes lignes suivantes: on écrit les équations (\ref{N-S-stationnaire}) comme le problème de point fixe équivalent: 
\begin{equation}\label{n-s_etat_point-fixe}
\U=\frac{1}{\nu}\P\left(\frac{1}{\Delta} div(U\otimes U)\right)-\frac{1}{\nu \Delta}\fe
\end{equation} où en utilisant l'estimation (\ref{estimation-a-priori-sol-estationnaires}) et le principe de point fixe de Shaefer (voir le Lemme \ref{Schaefer}) on construit tout d'abord une solution $\U \in \dot{H}^{1}(\Rt)$ pour ensuite récupérer  la pression $P \in \dot{H}^\frac{1}{2}(\Rt)$. \\
\\
Une fois que nous avons énoncé ce résultat  d'existence des solutions des équations de Navier-Stokes stationnaires nous allons maintenant expliquer le type de décroissance fréquentielle que nous voulons étudier pour ces solutions.

\subsection{La décroissance fréquentielle selon la théorie K41}\label{Sec:motivation-K41}

Pour expliquer la décroissance fréquentielle  que nous voulons étudier, nous avons besoin de revenir à la théorie de la turbulence K41 dont on a déjà parlé tout au début du Chapitre  \ref{Chap.1}. Plus précisément, nous allons rappeler rapidement   l'idée phénoménologique de cette théorie qui se base sur le modèle de cascade d'énergie et qui explique le processus d'introduction, transfert et dissipation de l'énergie cinétique dans un fluide en état turbulent. \\
\\ 
En effet,  on considère un fluide visqueux et incompressible sur lequel agit une force extérieure $\fe$ et cette force introduit de l'énergie cinétique dans le fluide à une échelle de longueur $\ell_0>0$ qui a été appelée l'échelle d'injection d'énergie.  Nous allons supposer que les effets de cette force  sont suffisamment forts de sorte que le fluide se trouve dans un état turbulent et alors le modèle de cascade d'énergie explique que l'énergie cinétique introduite dans le fluide à l'échelle d'injection $\ell_0$ est  transférée aux échelles de longueur plus petites $\ell>0$ (avec $\ell<<\ell_0$) et ce transfert d'énergie est effectué  jusqu'à que l'on arrive à l'échelle de dissipation $\ell_D$ (voir (\ref{l_D}) page \pageref{l_D} pour la définition  de cette échelle) avec $\ell_D<<\ell << \ell_0$. Ainsi, à partir de l'échelle d'échelle $\ell_D$ les forces de viscosité dissipent sous forme de chaleur l'énergie cinétique provenant des  échelles plus grandes.\\
\\
Nous pouvons observer que le domaine de validité du modèle de cascade d'énergie est donné par l'intervalle $]\ell_D,\ell_0[$ et nous pouvons  aussi observer que pour expliquer  ce modèle nous avons pris en compte  des différents échelles de longueur. Néanmoins, pour étudier la décroissance fréquentielle et comme nous l'expliquerons ci-après, il est plus intéressant ici d'étudier ce modèle par une approche différente  et alors  au lieu de considérer les échelles de longueur $\ell$ nous allons maintenant considérer les fréquences  $\kappa=\frac{1}{\ell}$.  L'intérêt de regarder le modèle de cascade d'énergie d'un point de vue fréquentiel repose sur  le fait que la théorie de la turbulence K41 propose une  caractérisation intéressante de ce modèle  que nous expliquons tout de suite (voir également les articles \cite{Kolm1}, \cite{Kolm2} et \cite{Kolm3} de A.N. Kolmogorov ainsi que le livre \cite{JacTab} de L. Jacquin et P. Tabeling).\\
\\
Pour  une échelle d'injection d'énergie  $\ell_0>0$ nous définissons $\kappa_0=\frac{1}{\ell_0}$ la fréquence d'injection d'énergie,
pour $\ell_D$ l'échelle de dissipation nous définissons $\kappa_D=\frac{1}{\ell_D}$ la fréquence de dissipation et alors le modèle de cascade d'énergie se traduit au niveau fréquentiel de la façon suivante: dans un fluide qui se trouve en régime turbulent,  l'énergie cinétique qui est introduite dans ce fluide (toujours par l'action d'une force $\fe$)  à la fréquence $\kappa_0$ est alors transférée aux plus hautes fréquences $\kappa\in ]\kappa_0,\kappa_D[$  jusqu'à la fréquence de dissipation $\kappa_D$ à partir de laquelle cette énergie est dissipée comme chaleur.\\
\\
Dans ce cadre, il s'agit tout d'abord d'étudier la quantité d'énergie cinétique $E(\kappa)$   aux  fréquences $\kappa \in ]\kappa_0,\kappa_D[$ où  d'après le modèle ci-dessus  l'énergie cinétique est transférée; et ensuite il s'agit d'étudier la  quantité $E(\kappa)$ aux hautes  fréquences $\kappa >\kappa_D$ où  l'énergie cinétique est dissipée comme chaleur. Cette quantité  $E(\kappa)$ sera  appelée  le spectre d'énergie et elle est  définie par l'expression 
\begin{equation}\label{def_E(k)}  
E(\kappa)=\int_{\vert \xi \vert=\kappa}\left\vert\widehat{\U}(\xi)\right\vert^2d\sigma(\xi),
\end{equation}   où $\widehat{\U}$ est la transformée de Fourier du champ de vitesse $\U$ et $ d\sigma$ est la mesure de la sphère unité. \\
\\
On souhaite donc étudier la décroissance de cette quantité $E(\kappa)$ lorsque $\kappa \in ]\kappa_0,\kappa_D[$. Pour cela nous avons besoin de faire un court rappel de quelques quantités physiques qui interviennent:  pour caractériser le régime  turbulent dans la théorie K41 on considère  le nombre de Reynolds $Re$  où nous rappelons que  pour $U>0$ la vitesse caractéristique du fluide, $\ell_0>0$ l'échelle d'injection d'énergie et $\nu>0$ la constante de viscosité du fluide ce nombre $Re$ est défini par 
  \begin{equation}
Re=\frac{U \ell_0}{\nu}.  
\end{equation}
De cette façon le régime turbulent sera caractérisé lorsque $Re>>1$. \\
\\
Ainsi, pour estimer la quantité $E(\kappa)$, Kolmogorov considère les trois paramètres physiques suivants: la constante de viscosité $\nu$, la fréquence $\kappa \in ]\kappa_0,\kappa_D[$ et le taux de dissipation d'énergie $\varepsilon>0$; et étant donné qu'on se trouve dans le régime turbulent  alors les forces visqueuses du fluide sont négligeables par rapport aux  forces inertielles, ce qui amène a Kolmogorov à faire l'hypothèse  suivante: \emph{la quantité d'énergie cinétique  $E(\kappa)$ est indépendante de la constante de viscosité $\nu$}. \\ 
\\
En supposant cette hypothèse, Kolmogorov introduit l'idée que cette quantité  $E(\kappa)$ doit alors  s'exprimer seulement en fonction des paramètres $\varepsilon$ et $\kappa$, ce qui lui amène à écrire l'estimation $\ds{E(\kappa) \approx  \varepsilon^{\alpha}\kappa^{\beta}}$. Pour déterminer les exposants $\alpha, \beta \in \R$ on fait l'analyse  de dimensions physiques suivante: comme le taux de dissipation $\varepsilon$ a une dimension physique $\frac{longueur^{2}}{temps^{3}}$ et la fréquence $\kappa$ a une dimension physique $\frac{1}{longueur}$ alors la quantité $\varepsilon^{\alpha}\kappa^{\beta}$ a une dimension physique  $\frac{longueur^{2\alpha -\beta}}{temps^{3 \alpha}}$. D'autre part, comme le spectre d'énergie $E(\kappa)$ mesure la quantité d'énergie cinétique à chaque fréquence $\kappa$  alors cette quantité a  une dimension physique $\frac{energie}{frequence}$  et donc une dimension physique  $\frac{longueur^{3}}{temps^{2}}$.       Ainsi, dans l'estimation $\ds{E(\kappa) \approx  \varepsilon^{\alpha}\kappa^{\beta}}$ il faut que l'on ait les valeurs $\alpha=\frac{2}{3}$ et $\beta=-\frac{5}{3}$.  \\
\\
Nous obtenons de cette façon l'estimation du spectre d'énergie $E(\kappa)$ pour les fréquences $\kappa \in ]\kappa_0,\kappa_D[$: 
\vspace{5mm}
\begin{center}
	\fbox{
		\begin{minipage}[l]{90mm}
			\begin{equation}\label{loi-cinq-tiers}
			\text{si}\quad Re>>1\quad \text{alors}\quad E(\kappa)\approx \varepsilon^{\frac{2}{3}} \kappa^{-\frac{5}{3}}, 
			\end{equation}
			\end{minipage}
	} 
\end{center}  
\vspace{5mm}
qui est connue comme la loi $-\frac{5}{3}$ de Kolmogorov et qui est valable seulement dans l'intervalle d'inertie $]\kappa_0,\kappa_D[$. \\
\\
Maintenant, il est important de remarquer que l'estimation (\ref{loi-cinq-tiers}), partiellement observée dans des expériences physiques et numériques (voir les articles  \cite{Kida},\cite{Mart} ainsi que les livres \cite{McDonough}, \cite{Tennekes}), a une explication purement phénoménologique et que l'étude déterministe rigoureuse de cette estimation dans le cadre des équations de Navier-Stokes est encore hors de portée. En effet, il y a très peu de références à ce sujet et la difficulté principale de cette étude déterministe porte sur le fait que les outils qu'on connait pour étudier les équations de Navier-Stokes ne suffissent pas pour étudier  l'encadrement (\ref{loi-cinq-tiers}) et ce problème ne sera donc pas considéré ici. \\
\\
Dans ce chapitre nous allons plutôt nous concentrer à l'étude déterministe du spectre d'énergie $E(\kappa)$ dans \emph{l'intervalle de dissipation}, c'est à dire, pour les hautes fréquences $\kappa>\kappa_D$. En effet, dans le modèle de cascade d'énergie ci-dessus et une fois que l'on arrive à la fréquence de dissipation $\kappa_D$, on a que les forces de viscosité du fluide sont censées annuler tous les mouvements de l'écoulement aux fréquences supérieures à $\kappa_D$ et donc l'énergie cinétique est \emph{rapidement} dissipée sous forme de chaleur. \\
\\
Dans les expériences physiques et numériques   \cite{Houel},\cite{Wilcox} et \cite{Mart}  il est souvent observé  que dans cet intervalle de dissipation  $\kappa>\kappa_D$ alors le spectre d'énergie $E(\kappa)$ est  à décroissance rapide et cette quantité présente le comportement
\begin{equation}\label{exp_decay_freq}
E(\kappa)\approx e^{-\kappa}. 
\end{equation} 
De cette façon, étant donné que la quantité $E(\kappa)$ est définie la transformée de Fourier du champ de vitesse $\U$ dans l'expression (\ref{def_E(k)}), pour les fréquences $ \vert \xi \vert=\kappa>\kappa_D$  nous allons  alors  étudier la décroissance (\ref{dec_fourier_U}):
\begin{equation*}
\left\vert\widehat{\U}(\xi) \right\vert \lesssim e^{-\vert \xi \vert}.   
\end{equation*} 
Mais avant d'entrer dans le vif du sujet il convient tout d'abord  de préciser un peu plus notre cadre d'étude: le régime laminaire et le régime turbulent. 

\subsection{De retour  au régime laminaire et turbulent}\label{Sec:regime-lam-turb}
Dans la section précédente nous avons expliqué que la décroissance  (\ref{dec_fourier_U}) est censée être observée  aux  fréquences  $\vert \xi \vert \in ]\kappa_D,+\infty[$, où $\kappa_D$ est la fréquence de dissipation d'énergie; et  nous allons  maintenant étudier en peu plus en détail cet intervalle des fréquences $]\kappa_D,+\infty[$ dans le régime laminaire et dans le régime turbulent. \\ 
\\
Lorsqu'on considère un fluide en régime laminaire, à un nombre de Reynolds borné $Re\leq C$, on s'attend à ce la fréquence de dissipation d'énergie soit du même ordre de grandeur que la fréquence d'injection d'énergie $\kappa_0$.  En effet, les fluides en régime laminaire peuvent êtres caractérisés par un taux de dissipation d'énergie de l'ordre $\varepsilon \approx \frac{1}{Re}\frac{U^3}{\ell_0}$ (voir le livre \cite{FMRTbook}, Chapitre II de R. Temam \emph{et. al.} ainsi que le livre \cite{JacTab} de L. Jacquin and P. Tabeling) et en utilisant maintenant cette valeur de $\varepsilon$ dans la formule  (\ref{Def_KDF}) nous avons 
\begin{equation}\label{KDrelacionReynLam}
\kappa_D\approx Re^{\frac{1}{2}}\kappa_0  
\end{equation} 
et comme $Re\leq C$  nous obtenons les inégalités $\kappa_D\leq \sqrt{C}\kappa_0$ et donc la fréquence de dissipation 
$\kappa_D$ est du même  ordre de grandeur que la fréquence d'injection d'énergie $\kappa_0$.\\
\\
De cette façon nous observons que  pour les fluides en régime laminaire la décroissance (\ref{dec_fourier_U}  doit être observée  dans l'intervalle de fréquences $]\kappa_0\approx  \kappa_{D}, +\infty[$  qui commence à partir de la fréquence d'injection d'énergie $\kappa_0$ et non seulement à partir de la fréquence de dissipation $\kappa_D$.\\
\\
D'autre part, dans le cadre d'un fluide en régime turbulent  nous allons observer que l'on a $\kappa_D>>\kappa_0$ et donc la décroissance (\ref{dec_fourier_U}) a lieu aux hautes fréquences. En effet,  observons tout d'abord que la fréquence $\kappa_D=\frac{1}{\ell_D}$ est définie comme l'inverse de la longueur de dissipation  $\ell_D$ et cette longueur  a été définie dans (\ref{Def_l_D}) page \pageref{Def_l_D}  par l'expression  $\ds{\ell_D=\left( \frac{\nu^3}{\varepsilon}\right)^{\frac{1}{4}}}$,  nous avons de cette façon 
 \begin{equation}\label{Def_KDF}
 \kappa_D=\left( \frac{\varepsilon}{\nu^3}\right)^{\frac{1}{4}},
 \end{equation} où $\varepsilon$ est le taux de dissipation d'énergie et $\nu$ est toujours la constante de viscosité du fluide. Mais, dans l'introduction de Chapitre \ref{Chap.1} nous avons vu que la loi de dissipation de Kolmogorov explique que si le nombre de Reynolds $Re$ est suffisamment grand alors le taux de dissipation $\varepsilon$ est estimé par  $\varepsilon \approx \frac{U^3}{\ell_0}$ et 
 en remplaçant cette valeur de  $\varepsilon$ dans   (\ref{Def_KDF}), comme $Re=\frac{U\ell_0}{\nu}$ et $\kappa_0=\frac{1}{\ell_0}$ nous obtenons 
\begin{equation}\label{KDrelacionReynTurb}
\kappa_D=\left( \frac{\varepsilon}{\nu^3}\right)^{\frac{1}{4}}\approx \left( \frac{U^3}{\ell_0 \nu^3}\right)^{\frac{1}{4}}= \left( \frac{U^3\ell^{3}_{0}}{ \nu^3}\right)^{\frac{1}{4}}\frac{1}{\ell_0} =  Re^{\frac{3}{4}}\kappa_0.
\end{equation}
Ainsi, lorsque $1<<Re,$ on a $\kappa_0<<\kappa_D$.\\
\\
Dans la section qui suit nous allons  étudier ces deux cadres du mouvement du fluide où l'on obtiendra une décroissance fréquentielle (\ref{dec_fourier_U})  des solutions des équations  de Navier-Stokes stationnaires (\ref{N-S-stationnaire}) et cette décroissance sera obtenue à partir des bases fréquences $\vert \xi \vert>0$ dans le régime laminaire et seulement aux  hautes fréquences $\vert \xi \vert >>1$ dans le régime turbulent. \\
\\
Insistons sur le fait que dans ce chapitre nous allons caractériser le mouvement du fluide par une approche différente par rapport au chapitre précédent et ici le régime laminaire  sera représenté par un contrôle directe sur la taille de la force $\fe$ tandis que dans le régime turbulent nous ne ferons aucun contrôle la taille de  cette fonction.\\
%
%
\section{Décroissance fréquentielle des équations de Navier-Stokes stationnaires}\label{Sec:decroissance-freq}
On commence donc par étudier la décroissance fréquentielle dans le cadre d'un fluide en régime laminaire.
  \subsection{Quelques résultats dans le régime laminaire}\label{Sec:decroissance-freq-lam}
 Pour  étudier la décroissance fréquentielle  (\ref{dec_fourier_U}) nous proposons dans cette section deux approches différentes. Dans la première approche qui sera étudiée dans le point $A)$ ci-dessous  nous supposerons tout d'abord que la force $\fe \in \dot{H}^{-1}(\Rt)$ est une fonction localisée en  fréquence  et dans ce cadre nous utiliserons   la méthode développée par Oseen pour la construction des solutions classiques des équations de Navier-Stokes  (voir le livre \cite{PGLR1}, Chapitre $4$ pour plus de détails sur cette méthode) et cette méthode sera combinée avec  la notion de nombres de Catalan utilisés en combinatoire (voir le livre \cite{Comtet} pour plus de références) pour construire une solution $\U \in \dot{H}^{1}(\Rt)$ des équations de Navier-Stokes stationnaires  (\ref{N-S-stationnaire})  qui  vérifie la  décroissance fréquentielle ponctuelle $\ds{\vert \widehat{\U}(\xi)\vert\lesssim e^{-\vert\xi \vert}}$ pour $\vert \xi\vert$ supérieur à une certaine fréquence.\\
\\
Ensuite, dans la deuxième approche qui sera étudiée dans le point $B)$ ci-après, on considère une force  $\fe \in \dot{H}^{-1}(\Rt)$ qui n'est pas localisée en fréquence mais qui vérifie une décroissance fréquentielle exponentielle  et en utilisant la notion de  classe de Gevrey (voir la Définition \ref{Classe de Gevrey} page \pageref{Classe de Gevrey} pour une définition précise)  nous construirons une solution $\U \in \dot{H}^{1}(\Rt)$ des équations (\ref{N-S-stationnaire}) qui vérifie une décroissance fréquentielle $\left\vert \widehat{\U}(\xi)\right\vert \lesssim \frac{e^{-\vert\xi \vert}}{\vert \xi \vert}$, pour tout $\xi\neq 0.$ \\
\\ 
Dans ces deux approches nous aurons besoin de supposer des conditions de petitesse sur la taille de la force  $\fe$: $\Vert \fe \Vert_{\dot{H}^{-1}}<\eta$ et pour cette raison ces résultats sont valables seulement dans  le cadre d'un fluide en régime laminaire. 
\subsection*{A) Première approche: localisation fréquentielle de la force}
Soit $\fe \in \dot{H}^{-1}(\Rt)$ une force et soit $\ell_0>0$ une échelle d'injection d'énergie. Nous allons supposer que cette force vérifie 
\begin{equation}\label{loc-force}
supp \left( \widehat{\fe} \right)\subset \left\lbrace \xi \in \Rt: \frac{\rho_1}{\ell_0}\leq \vert \xi \vert \leq \frac{\rho_2}{\ell_0}\right\rbrace,
\end{equation} où $0<\rho_1<\rho_2$ sont deux constantes fixes
et nous allons voir maintenant comment on peut construire une solution $\U$ des équations  (\ref{N-S-stationnaire}) qui ait une décroissance fréquentielle exponentielle. \\
\\
Pour cela nous allons adapter à ce cadre stationnaire la méthode d'Oseen qui fût utilisée pour construire des solutions classiques  des équations de Navier-Stokes non-stationnaires et nous obtiendrons une solution $\U\in \dot{H}^{1}(\Rt)$  qui s'écrit comme une série  $\ds{\U=\sum_{n=1}^{+\infty}}\U_n$ où le terme $\U_n$ sera donné par l'expression (\ref{def_U_n}) ci-dessous. Ainsi, la force $\fe$ étant localisée aux fréquences nous observerons que cette localisation entraîne une localisation fréquentielle convenable  de chaque terme $\U_n$ de la série ci-dessus et cette propriété nous permettra obtenir une décroissance fréquentielle exponentielle pour cette solution $\U$,  qui sera étudiée dans le Théorème \ref{Theo:Ossen_decroissance}  ci-dessous.  \\
\\
Une fois que nous avons expliqué les grandes lignes de notre étude  nous allons maintenant entrer dans les détails et la première chose à faire est de construire une solution $\U$ des équations (\ref{N-S-stationnaire}) dans un sous-espace de $\dot{H}^{1}(\Rt)$. En effet, comme l'on travaille sur tout $\Rt$ on considère les équations  (\ref{N-S-stationnaire}) comme un problème de point fixe équivalent 
\begin{equation}\label{N-S-stationnaire-pf}
\U=\frac{1}{\nu} \P\left( \frac{1}{\Delta}div(\U \otimes \U) \right)-\frac{1}{\nu \Delta}\fe.
\end{equation} Pour résoudre ce problème nous allons utiliser le principe de contraction de Picard et alors nous cherchons un espace fonctionnel $E\subset \dot{H}^{1}(\Rt)$ où l'on puise vérifier la continuité de la forme bilinéaire 
\begin{equation}\label{cont_forme_bilin}
\left\Vert \frac{1}{\nu} \P\left( \frac{1}{\Delta}div(\U \otimes \U) \right) \right\Vert_E \leq \frac{C}{\nu} \Vert \U \Vert_E \Vert \U \Vert_E,
\end{equation}  pour tout $\U \in E$ et  où l'on puise relier la condition de petitesse sur le terme $\ds{\left\Vert \frac{1}{\nu \Delta}\fe \right\Vert_{E}}$ par un contrôle direct  sur la norme $\Vert \fe \Vert_{\dot{H}^{-1}}$:
\begin{equation}\label{controle-norme-force}
\left\Vert \frac{1}{\nu \Delta}\fe \right\Vert_{E}\leq C_{(\nu,\ell_0)}\Vert \fe \Vert_{\dot{H}^{-1}}, 
\end{equation} où $\ell_0>0$ est fixée par (\ref{loc-force}).  Donc, pour vérifier les points  (\ref{cont_forme_bilin}) et (\ref{controle-norme-force}) nous introduisons l'espace 
\begin{equation}\label{E}
E=L^2(\Rt)\cap \dot{H}^{1}(\Rt)
\end{equation} muni de la norme 
\begin{equation}\label{norme_E}
\Vert \cdot \Vert_{E}=  \Vert \cdot \Vert_{L^2}+ \Vert \cdot \Vert_{\dot{H}^1}+\Vert \cdot \Vert_{L^3}.
\end{equation}  \`A ce stade, il convient de souligner que par les inégalités de Hardy-Littlewood-Sobolev et par les inégalités d'interpolation dans l'espaces de Lebesgue  nous avons $E=L^2(\Rt)\cap \dot{H}^{1}(\Rt) \subset L^3(\Rt)$ et dans la norme $\Vert \cdot \Vert_{E}$ ci-dessus nous observons que nous aurons besoin de contrôler la quantité $\Vert \cdot \Vert_{L^3}$. \\
\\
Observons aussi que l'on a $E\subset \dot{H}^{1}(\Rt)$ et nous avons:

\begin{Proposition}\label{Prop:existence_sol_lam_force_loc} Soit $\fe\in \dot{H}^{-1}(\Rt)$ une force à divergence nulle  qui vérifie (\ref{loc-force}). Il existe une constante $\eta>0$ telle que si $\Vert \fe \Vert_{\dot{H}^{-1}}<\eta$ alors il existe $\U \in E$ solution des équations de Navier-Stokes stationnaires  (\ref{N-S-stationnaire}) qui est l'unique solution dans la boule $\Vert \U \Vert_{E}<\frac{\nu}{2C}$, où $C>0$ est la constante donnée dans (\ref{cont_forme_bilin}). 
\end{Proposition}
\pv Il s'agit de vérifier tout d'abord le point (\ref{cont_forme_bilin}) et pour cela  nous avons besoin d'estimer chaque terme de la norme $\Vert \cdot \Vert_E$ donnée par l'expression (\ref{norme_E}). Pour le premier qui compose la norme $\Vert \cdot \Vert_{E}$  
par la continuité  du projecteur de Leray et par les inégalités de Hardy-Littlewood-Sobolev et l'inégalité de H\"older nous obtenons 
\begin{equation}\label{estim1}
\left\Vert \frac{1}{\nu} \P\left( \frac{1}{\Delta}div(\U \otimes \U) \right)  \right\Vert_{L^2}\leq \frac{c_1}{\nu}\left\Vert \U\otimes \U\right\Vert_{\dot{H}^{-1}} \leq \frac{c_1}{\nu} \left\Vert \U\otimes \U\right\Vert_{L^{\frac{6}{5}}} \leq \frac{c_1}{\nu} \Vert \U \Vert_{L^2}\Vert \U \Vert_{L^3}.
\end{equation}
  
Ensuite, pour le deuxième terme de la norme $\Vert \cdot \Vert_{E}$, en suivant le même raisonnement ci-dessus nous avons  
\begin{equation}\label{estim2}
\left\Vert \frac{1}{\nu} \P\left( \frac{1}{\Delta}div(\U \otimes \U) \right)  \right\Vert_{\dot{H}^{1}}\leq \frac{c_2}{\nu} \left\Vert \U\otimes \U \right\Vert_{L^2} \leq \frac{c_2}{\nu} \Vert \U \Vert_{L^6}\Vert \U \Vert_{L^3}\leq  \frac{c_2}{\nu} \Vert \U \Vert_{\dot{H}^1}\Vert \U \Vert_{L^3}.
\end{equation}  
Finalement, pour le troisième terme de la norme $\Vert \cdot \Vert_E$, par les inégalités de Hardy-Littlewood-Sobolev nous avons  
\begin{equation*}
\left\Vert \frac{1}{\nu}\P\left(\frac{1}{\Delta}div(\U\otimes \U)\right)  \right\Vert_{L^3} \leq  \frac{c_3}{\nu}\left\Vert  \P\left(\frac{1}{\Delta}div(\U\otimes \U)\right) \right\Vert_{\dot{H}^{\frac{1}{2}}}
\end{equation*} d'où, toujours par la continuité du projecteur de Leray nous écrivons 
\begin{equation}\label{estim3}
\frac{c_3}{\nu} \left\Vert \P\left(\frac{1}{\Delta}div(\U\otimes\U)\right) \right\Vert_{\dot{H}^{\frac{1}{2}}} \leq \frac{c_3}{\nu} \left\Vert \U\otimes \U\right\Vert_{\dot{H}^{-\frac{1}{2}}}.
\end{equation} Dans la dernière expression ci-dessus, nous appliquons encore les inégalités de Hardy-Littlewood-Sobolev et l'inégalité de H\"older pour obtenir 
\begin{equation}\label{estim4}
\frac{c_3}{\nu} \left\Vert \U\otimes \U\right\Vert_{\dot{H}^{-\frac{1}{2}}} \leq  \frac{c_3}{\nu} \left\Vert \U\otimes \U\right\Vert_{L^{\frac{3}{2}}}\leq \frac{c_3}{\nu} \left\Vert \U\right\Vert_{L^3}\left\Vert \U\right\Vert_{L^3}.
\end{equation}
De cette façon, la propriété de continuité  du terme bilinéaire sur l'espace $E$ énoncé dans (\ref{cont_forme_bilin}) est bien vérifiée par les estimations ci-dessus où l'on pose la constante  $C=\max(c_1,c_2,c_3)>0$.\\
\\
Nous allons étudier  le point  (\ref{controle-norme-force}) et pour cela nous avons le lemme technique suivant.
\begin{Lemme}\label{lemme:controle-force}  Soit $\fe \in \dot{H}^{-1}(\Rt)$ une force qui vérifie la localisation fréquentielle (\ref{loc-force}) pour $\ell_0>0$ une échelle d'injection d'énergie. Soient $\nu>0$ la constante de viscosité du fluide et la norme $\Vert \cdot \Vert_E$ donnée dans (\ref{norme_E}). Alors il existe une constante  $C_{(\nu,\ell_0)}>0$ telle que l'on a l'estimation  $$ \ds{\left\Vert \frac{1}{\nu \Delta}\fe \right\Vert_{E}\leq C_{(\nu,\ell_0)}\Vert \fe \Vert_{\dot{H}^{-1}}.}$$
\end{Lemme}
\pv  Nous allons estimer chaque terme de la quantité $\ds{\left\Vert \frac{1}{\nu \Delta}\fe\right\Vert_{E}}$. En effet, pour le premier terme de cette quantité, comme $\fe$ est localisée aux fréquences $\frac{\rho_1}{\ell_0}\leq \vert \xi \vert \leq  \frac{\rho_2}{\ell_0}$ nous pouvons écrire 
\begin{eqnarray*}
	\left\Vert \frac{1}{\nu \Delta}\fe \right\Vert^{2}_{L^2}&=&\frac{1}{\nu^2}\int_{\Rt} \frac{1}{\vert \xi \vert^4}\vert \widehat{\fe}(\xi)\vert^{2}d\xi= \frac{1}{\nu^2}\int_{\Rt} \mathds{1}_{ \lbrace{ \xi\in \Rt :\frac{\rho_1}{\ell_0}\leq \vert \xi \vert  \leq \frac{\rho_2}{\ell_0} \rbrace}}\frac{1}{\vert \xi \vert^4}\vert \widehat{\fe}(\xi)\vert^{2}d\xi\\
	&=& \frac{1}{\nu^2}\int_{\Rt} \mathds{1}_{ \lbrace{ \xi\in \Rt :\frac{\rho_1}{\ell_0}\leq \vert \xi \vert  \leq \frac{\rho_2}{\ell_0} \rbrace}}\frac{1}{\vert \xi \vert^2}\left[ \frac{1}{\vert \xi \vert^2}\vert \widehat{\fe}(\xi)\vert^{2}\right] d\xi \leq \frac{\rho^{2}_{1}}{\nu^2 \ell^{2}_{0}}\int_{\Rt}  \frac{1}{\vert \xi \vert^2}\vert \widehat{\fe}(\xi)\vert^{2} d\xi \\
	&\leq & \frac{\rho^{2}_{1}}{\nu^2 \ell^{2}_{0}} \Vert \fe \Vert^{2}_{\dot{H}^{-1}},
\end{eqnarray*} d'où nous obtenons l'estimation  
\begin{equation}\label{estim_force_1}
\left\Vert \frac{1}{\nu \Delta}\fe \right\Vert_{L^2} \leq \frac{\rho_1}{\nu \ell_0} \Vert  \fe \Vert_{\dot{H}^{-1}}. 
\end{equation}
Ensuite, pour le deuxième terme de la quantité $\ds{\left\Vert \frac{1}{\nu \Delta}\fe\right\Vert_{E}}$ nous pouvons écrire directement
\begin{equation}\label{estim_force_2}
\left\Vert  \frac{1}{\nu \Delta}\fe\right\Vert_{\dot{H}^{1}} \leq \frac{1}{\nu} \Vert \fe \Vert_{\dot{H}^{-1}}.
\end{equation}
Finalement, pour le troisième terme de cette quantité, par les inégalités de Hardy-Littlewood-Sobolev nous avons 
$$ \left\Vert \frac{1}{\nu \Delta}\fe\right\Vert_{L^3} \leq c \left\Vert \frac{1}{\nu \Delta}\fe\right\Vert_{\dot{H}^{\frac{1}{2}}}, $$ et alors, toujours en utilisation la localisation fréquentielle de $\fe$ nous écrivons 
\begin{eqnarray*}
	\left\Vert \frac{1}{\nu \Delta}\fe \right\Vert^{2}_{\dot{H}^{\frac{1}{2}}}&=&\frac{1}{\nu^2}\int_{\Rt} \vert \xi \vert \left[\frac{1}{\vert \xi \vert^4} \vert \widehat{\fe}(\xi)\vert^{2}\right] d\xi= \frac{1}{\nu^2}\int_{\Rt} \mathds{1}_{ \lbrace{ \xi\in \Rt :\frac{\rho_1}{\ell_0}\leq \vert \xi \vert  \leq \frac{\rho_2}{\ell_0} \rbrace}}\frac{1}{\vert \xi \vert}\left[ \frac{1}{\vert \xi \vert^2}\vert \widehat{\fe}(\xi)\vert^{2}\right] d\xi\\
	&\leq &  \frac{\rho_{1}}{\nu^2 \ell_{0}}\int_{\Rt}  \frac{1}{\vert \xi \vert^2}\vert \widehat{\fe}(\xi)\vert^{2} d\xi  \leq \frac{\rho_{1}}{\nu^2 \ell_{0}} \Vert \fe \Vert^{2}_{\dot{H}^{-1}},
\end{eqnarray*} d'où nous pouvons en tirer l'estimation 
\begin{equation}\label{estim_force_3}
\left\Vert \frac{1}{\nu \Delta}\fe\right\Vert_{L^3} \leq c \frac{\rho^{\frac{1}{2}}_{0}}{\nu \ell^{\frac{1}{2}}_{0}} \Vert \fe \Vert_{\dot{H}^{-1}}. 
\end{equation}  Ainsi, par les estimations (\ref{estim_force_1}), (\ref{estim_force_2}) et (\ref{estim_force_3}) on fixe la constante $C_{(\nu,\ell_0)}=\max \left( \frac{\rho_1}{\nu \ell_0},  \frac{1}{\nu}, c \frac{\rho^{\frac{1}{2}}_{0}}{\nu \ell^{\frac{1}{2}}_{0}}\right)$ et nous écrivons alors l'estimation   $\ds{\left\Vert \frac{1}{\nu \Delta}\fe \right\Vert_{E}\leq C_{(\nu,\ell_0)}\Vert \fe \Vert_{\dot{H}^{-1}}.}$ \finpv
\\
Ainsi, si $\Vert \fe \Vert_{\dot{H}^{-1}}<\eta$ pour $\eta>0$ une constante suffisamment petite  alors par le lemme ci-dessus et par l'estimation (\ref{cont_forme_bilin}) nous pouvons appliquer le principe de contraction de Picard pour obtenir $\U\in E$ une solution des équations (\ref{N-S-stationnaire}) et qui est l'unique solution dans la boule $\Vert \U \Vert_{E}<\frac{\nu}{2C}$. La Proposition \ref{Prop:existence_sol_lam_force_loc} est maintenant vérifiée. \finpv 
\\
Une fois que nous avons construit la solution $\U\in E$ nous allons maintenant revisiter la méthode d'Oseen dans ce cadre des équations de Navier-Stokes stationnaires. Remarquons au passage que cette méthode a été aussi utilisée dans l'article \cite{Aucher} de P. Aucher et P. Tchamitchian  pour obtenir des solutions des équations de Navier-Stokes non-stationnaires dans le cadre de différents espaces fonctionnels comme les espaces de Sobolev, Besov et Morrey.\\

\begin{Proposition}\label{Prop:descomp_U} Dans le cadre de la Proposition \ref{Prop:existence_sol_lam_force_loc}, la solution $\U\in E$ des équations de Navier-Stokes stationnaires (\ref{N-S-stationnaire})  s'écrit comme une série convergente dans la topologie forte de l'espace $E$, $\ds{\U=\sum_{n=1}^{+\infty}\U_n}$, où le premier terme est donné par $\U_1=-\frac{1}{\nu\Delta}\fe$ et pour tout $n\geq 2$ on a
\begin{equation}\label{def_U_n}
 \U_n=\sum_{k=1}^{n-1} \frac{1}{\nu} \P\left( \frac{1}{\Delta}div(\U_k \otimes \U_{n-k}) \right).
\end{equation} 
\end{Proposition} 
\pv  Tout d'abord, pour simplifier l'écriture nous allons écrire le terme $\frac{1}{\nu} \P\left( \frac{1}{\Delta}div(\U \otimes \U) \right)$ comme  $\frac{1}{\nu}B(\U,\U)$ et en suivant l'idée d'Oseen on considère l'équation 
\begin{equation}\label{prob_approch_eps}
\U= \frac{1}{\nu}B(\U ,\U )-\frac{1}{\nu \Delta}\fe
\end{equation} 
et l'on développe la solution  $\U$ comme la série: 
\begin{equation}\label{descomp_serie_eps}
\U =\sum_{n=1}^{+\infty} \U_n.
\end{equation} 
Pour déterminer les termes $\U_n$ en remplace (\ref{descomp_serie_eps}) dans (\ref{prob_approch_eps}) et nous obtenons
\begin{equation}\label{Equation_Serie}
\sum_{n=1}^{+\infty}  \U_n= \frac{1}{\nu}B\left( \sum_{n=1}^{+\infty} \U_n,\sum_{n=1}^{+\infty}\U_n\right)-\frac{1}{\nu \Delta}\fe,
\end{equation}
et comme $ \U_1=-\frac{1}{\nu\Delta}\fe$ nous pouvons écrire
$$\sum_{n=2}^{+\infty} \U_n= \frac{1}{\nu}B\left( \sum_{n=1}^{+\infty} \U_n,\sum_{n=1}^{+\infty} \U_n\right).$$
En utilisant les propriétés du terme bilinéaire nous obtenons
$$\frac{1}{\nu}B\left( \sum_{n=1}^{+\infty} \U_n,\sum_{n=1}^{+\infty}\U_n\right)=\sum_{n=1}^{+\infty}\sum_{k=1}^{+\infty} \frac{1}{\nu}B\left( \U_n,  \U_k\right)
=\sum_{n=2}^{+\infty}\,\left(\sum_{k=1}^{n-1}\frac{1}{\nu}B(\U_k,\U_{n-k})\right),$$  
et nous trouvons que le terme $\U_n$ et bien donné par la formule (\ref{def_U_n}).\\
\\ 
Nous allons maintenant prouver  que la série (\ref{descomp_serie_eps}) converge fortement vers la solution $\U$ dans l'espace $E$ et pour cela,  nous prouvons tout d'abord que  $\sum_{n=1}^{+\infty}\U_n$ converge  fortement dans $E$ et ensuite nous allons observer que cette série est en fait une solution de l'équation (\ref{prob_approch_eps}).\\
\\
Pour montrer la convergence forte de série ci-dessus nous introduisons les nombres de Catalan d'une façon légèrement détournée et adaptée à notre problème: par l'expression  (\ref{def_U_n}) nous pouvons observer que  pour tout entier $n\geq 2,$ le terme $\U_n$ est la somme de toutes les façons possibles de multiplier  $n$ fois le terme $\U_1$ en utilisant la forme bilinéaire $\frac{1}{\nu} B(\cdot, \cdot)$. En effet, en prenant $n=2$ dans la formule (\ref{def_U_n}) nous avons $\U_2=\frac{1}{\nu}B(\U_1,\U_1)$, pour $n=3$ nous avons $\U_3=\frac{1}{\nu}B\left(\U_1,\U_{2}\right)+\frac{1}{\nu}B\left(\U_{2}, \U_1\right)$, et dans cette identité en remplaçant $\U_{2}$ par $\frac{1}{\nu}B(\U_1,\U_1)$ nous obtenons $\U_3=\frac{1}{\nu}B\left(\U_1, \frac{1}{\nu}B(\U_1,\U_1)\right)+\frac{1}{\nu}B\left(\frac{1}{\nu}B(\U_1,\U_1), \U_1\right)$ et ainsi de suite pour $n>3$. Une fois que l'on a remarqué ce fait, pour $n\geq 2$ nous définissons les ensemble $W_n$ comme l'ensemble de tous les produits de  $n$ fois le terme $\U_1$ en utilisant $\frac{1}{\nu} B(\cdot , \cdot)$, par exemple:
\begin{eqnarray*}
W_{2}&=&\left\{\frac{1}{\nu}B(\U_1,\U_1)\right\},\\
 W_{3}&=&\left\{\frac{1}{\nu}B\left(\U_1, \frac{1}{\nu}B(\U_1,\U_1)\right), \; \frac{1}{\nu}B\left(\frac{1}{\nu}B(\U_1,\U_1), \U_1\right)\right\},\\
W_{4}&=&\left\{\frac{1}{\nu}B\left(\U_1, \frac{1}{\nu}B(\U_1, \frac{1}{\nu}B(\U_1,\U_1))\right), \; \frac{1}{\nu}B\left(\U_1, \frac{1}{\nu}B(\frac{1}{\nu}B(\U_1,\U_1), \U_1)\right),\cdots\right\},\\
\end{eqnarray*}
et donc, avec cette définition des ensembles  $W_{n}$, nous pouvons écrire $\ds{\U_n=\sum_{w\in W_n}w}$. \\
\\
Maintenant, nous définissons les nombres $A_1=1$ et  pour tout $n\geq 2$ le nombre $A_n=card(W_n)$ est le  cardial de l'ensemble $W_n$. Cette suite de nombres  $(A_n)_{n\geq 1}$ est connue en combinatoire par les  \emph{nombres de Catalan}.  Pour une application des nombres de Catalan dans le problème de Navier-Stokes voir le livre \cite{PGLR1}, Théorème $5.1$ et voir aussi le livre \cite{Comtet} pour les propriétés générales de ces nombres.\\
\\
Dans notre étude nous aurons besoin uniquement de deux propriétés de ces nombres. Tout d'abord,  la suite $(A_n)_{n\geq 1}$ vérifie la relation  de récurrence:
\begin{equation}\label{def_A_n}
A_n=\left\lbrace
\begin{array}{cl}\vspace{3mm}
1,&\text{si}\quad n=1,\\ 
\ds{\sum_{k=1}^{n-1}A_k\, A_{n-k}},&\text{si}\quad n\geq 2.\\
\end{array}\right.
\end{equation} 
Ensuite, par cette formule de récurrence il est possible de montrer (voir le livre \cite{Comtet}) que sa série génératrice vérifie : 
\begin{equation}\label{suite_gener}
\sum_{n=1}^{+\infty}A_n z^n=\frac{1-\sqrt{1-4z}}{2},\quad \text{pour tout }\quad  \vert z \vert \leq \frac{1}{4}.
\end{equation}
Maintenant nous avons tous les outils pour prouver que $\ds{\left\Vert \sum_{n=1}^{+\infty}\U_n\right\Vert_E}<+\infty$. En effet, par la continuité de la forme bilinéaire  $\frac{1}{\nu}B(\cdot , \cdot)$ dans l'espace $E$ étudiée dans (\ref{cont_forme_bilin})   nous obtenons l'estimation suivante pour tout $w\in W_n$:
$$\Vert w \Vert_E\leq \left( \frac{C}{\nu}\right)^{n-1}\Vert \U_1 \Vert^{n}_{E},$$ et de plus, comme $A_n=card(W_n)$ nous écrivons
\begin{equation}\label{major_catalan}
\Vert \U_n \Vert_{E}\leq \sum_{w\in W_n}\Vert w \Vert_E\leq  A_n \left( \frac{C}{\nu}\right)^{n-1} \Vert \U_1 \Vert^{n}_{E}=\frac{\nu}{C} A_n\left(\frac{C}{\nu} \Vert \U_1 \Vert_{E}\right)^n,
\end{equation} d'où nous obtenons 
\begin{equation*}  
\left\Vert  \sum_{n=1}^{+\infty}\U_{n} \right\Vert_{E}\leq \sum_{n=1}^{+\infty}\Vert \U_n \Vert_{E} \leq \frac{\nu}{C}\sum_{n=1}^{+\infty}A_n\left( \frac{C}{\nu} \Vert \U_1 \Vert_{E}\right)^n.
\end{equation*} 
Ainsi, comme $\U_1=-\frac{1}{\nu\Delta}\fe$, par le Lemme \ref{lemme:controle-force} nous savons que $\Vert \U_1\Vert_{E}\leq C_{(\nu,\ell_0)}\Vert \fe \Vert_{\dot{H}^{-1}}$ et alors si $\Vert \fe \Vert_{\dot{H}^{-1}}<\eta$ pour une constante $\eta>0$ suffisamment petite nous obtenons $\frac{C}{\nu} \Vert \U_1 \Vert_{E}<\frac{1}{4}$ et de cette façon, par la formule (\ref{suite_gener}) nous obtenons  
$$\left\Vert  \sum_{n=1}^{+\infty}\U_{n} \right\Vert_{E}\leq \sum_{n=1}^{+\infty}\Vert \U_n \Vert_{E} \leq\frac{1-\sqrt{1-4 \frac{C}{\nu} \Vert \U_1 \Vert_{E}}}{2\left( \frac{C}{\nu}\right)}\leq \frac{\nu}{2C},$$
et donc la série $\ds \sum_{n=1}^{+\infty}\U_{n} $ est normalement convergente dans $E$,  elle est une solution de l'équation (\ref{Equation_Serie}) et elle définit un élément de $E$ dans la boule $B_{E}(0,  \frac{\nu}{2C})$.\\
\\
Finalement, comme l'équation (\ref{prob_approch_eps}) s'agit de la même  (\ref{N-S-stationnaire}) par l'unicité des solutions de ces équations dans la boule $\Vert \U \Vert_E\leq \frac{\nu}{2C}$ nous avons $\U=\sum_{n=1}^{+\infty}\U_n$. \finpv\\
La décomposition $\ds{\U=\sum_{n=1}^{+\infty}}\U_n$ obtenue dans la Proposition \ref{Prop:descomp_U}  nous permettra d'étudier de façon  précise  la décroissance fréquentielle  de la solution $\U$ et pour cela nous étudions maintenant  l'information sur le support fréquentielle de chaque terme $\U_n$ de cette décomposition. \\
\\
\begin{Proposition}\label{Prop:support-freq} 
Dans le cadre de la Proposition \ref{Prop:descomp_U}, soit $\fe\in \dot{H}^{-1}(\Rt)$ une force qui vérifie la localisation fréquentielle $supp\left( \widehat{\fe}\right)\subset \left\lbrace \xi \in \Rt : \frac{\rho_1}{\ell_0} \leq \vert \xi \vert \leq \frac{\rho_2}{\ell_0}\right\rbrace,$ pour $\ell_0>0$ une échelle d'injection d'énergie et  soit $\U\in E$ la solution des équations des équations de Navier-Stokes stationnaires (\ref{N-S-stationnaire}) que nous assumons  décomposée comme la série $\ds{\U=\sum_{n=1}^{+\infty}\U_n}$ où $\U_1=-\frac{1}{\nu\Delta}\fe$ et pour tout $n\geq 2$ les termes $\U_n$ sont donnés par la formule (\ref{def_U_n}).  Alors les supports des transformées de Fourier des termes $\U_{n}$ vérifient les inclusions:
$$supp\,\left(\widehat{\U_n}\right)\subset  \left\{ \xi \in \Rt: \vert \xi \vert \leq n\frac{\rho_2}{\ell_0}\right\}.$$
\end{Proposition}
\pv  En effet,  par la localisation de la force $\fe$ et comme  $\U_1=-\frac{1}{\nu \Delta}\fe$  nous avons directement que
\begin{equation}\label{Inclusion1}
supp\,\left(\widehat{\U_1}\right)\subset  \left\{  \xi \in \Rt: \vert \xi \vert \leq \frac{\rho_2}{\ell_0}\right\}.
\end{equation}
Pour $n\geq 2$,  par la formule (\ref{def_U_n}) nous avons l'expression suivante pour les termes $\U_{n}$:
$$\U_n= \sum_{k=1}^{n-1}\frac{1}{\nu}B(\U_k,\U_{n-k})=\sum_{k=1}^{n-1}\frac{1}{\nu}\P\left(\frac{1}{\Delta}div(\U_k\otimes \U_{n-k})\right),$$ 
et comme le symbole $\widehat{m}(\xi)$ de l'opérateur $\P\left(\frac{1}{\Delta}div(\cdot )\right)$ est une fonction homogène de degré $-1$, en prenant la transformation de Fourier dans l'expression ci-dessus nous pouvons écrire
\begin{equation}\label{Inclusion2}
\vert \widehat{ \U_n}(\xi)\vert \leq \frac{c}{\nu\vert \xi \vert} \sum_{k=1}^{n-1}\vert \widehat{\U_k}\ast \widehat{\U_{n-k}}(\xi )\vert.
\end{equation} 
Pour le cas $n=2$, nous observons facilement que le support de  $ \widehat{ \U_2}$ dépend du support de $\widehat{\U_1}\ast \widehat{\U_{1}}$ et par (\ref{Inclusion1}) nous avons
$$supp\, \left(\widehat{\U_2}\right)\subset supp\,\left(\widehat{\U_1}\right)+supp\,\left(\widehat{\U_{1}}\right)\subset \left\{  \xi \in \Rt: \vert \xi \vert \leq 2\frac{\rho_2}{\ell_0}\right\}.$$
Donc, en utilisant l'expression (\ref{Inclusion2}), nous obtenons par récurrence les inclusions cherchées:
\begin{equation*}
supp\, \left(\widehat{\U_n}\right)\subset  \left\{ \xi\in \Rt: \vert \xi \vert \leq n\frac{\rho_2}{\ell_0}\right\}. 
\end{equation*}  \finpv\\
Une fois que nous avons étudié les supports fréquentiels des termes de la série $\ds{\U=\sum_{n=1}^{+\infty}\U_n}$ nous avons maintenant tous les outils dont on a besoin pour étudier le décroissance fréquentielle de la solution $\U$ et nous avons donc un tout premier résultat dans le cadre d'un fluide en régime laminaire.\\
\begin{Theoreme}[Première étude de  la décroissance fréquentielle: cadre laminaire]\label{Theo:Ossen_decroissance} Soit $\fe\in \dot{H}^{-1}(\Rt)$ une force à divergence nulle telle que sa transformée de Fourier est bornée et localisée aux fréquences $ \left\lbrace \xi \in \Rt : \frac{\rho_1}{\ell_0} \leq \vert \xi \vert \leq \frac{\rho_2}{\ell_0}\right\rbrace,$ pour $\ell_0>0$ une échelle d'injection d'énergie fixe. Alors, il existe une constante $\eta>0$ (suffisamment petite) telle que si $\Vert \fe \Vert_{\dot{H}^{-1}}<\eta$ alors la solution $\U\in E\subset \dot{H}^{1}(\Rt)$ des équations de Navier-Stokes stationnaires (\ref{N-S-stationnaire}), obtenue par le biais de  la Proposition  \ref{Prop:existence_sol_lam_force_loc}, vérifie la décroissance fréquentielle: 
\begin{equation}\label{dec-frec-lam-theo-1}
\left\vert \widehat{\U}(\xi)\right\vert \leq c_1\,e^{-c_2\vert \xi \vert}, \quad \text{pour tout}\quad \vert \xi \vert \geq \frac{\rho_1}{\ell_0},
\end{equation}  où $c_1=c_1(\nu,\ell_0,\fe)>0$ (avec $\nu>0$ la constante de viscosité du fluide) et $c_2=c_2(\ell_0)>0$.  
\end{Theoreme}  
Nous observons ainsi que cette décroissance fréquentielle est obtenue aux fréquences $\vert \xi \vert\geq \frac{\rho_1}{\ell_0}$ ce qui correspond  au fait que si l'on considère un fluide en régime laminaire (ce qui est caractérisé par la condition de petitesse sur la force) alors  la décroissance exponentielle du spectre d'énergie a également  lieu pour les fréquences de l'ordre de $\frac{1}{\ell_0}$ et non seulement aux hautes fréquences. \\ 
\\

 \dm Par la Proposition \ref{Prop:descomp_U} nous décomposons la solution $\U$ comme la série $\ds{\U=\sum_{n=1}^{+\infty}\U_n}$, où $\U_1=-\frac{1}{\nu\Delta}\fe$ et pour tout $n\geq 2$ les termes $\U_{n}$ sont donnés par l'expression (\ref{def_U_n}). En prenant la transformée de Fourier dans cette série et ensuite, en multipliant par la quantité  $e^{\frac{\ell_0}{\rho_2}\vert \xi \vert}$   nous obtenons   
 \begin{equation}\label{major_espec}
e^{\frac{\ell_0}{\rho_2}\vert \xi \vert} \left\vert \widehat{ \U}(\xi)\right\vert \leq \sum_{n=1}^{+\infty}e^{\frac{\ell_0}{\rho_2}\vert \xi \vert}\left\vert \widehat{\U_n}(\xi)\right\vert= e^{\frac{\ell_0}{\rho_2}\vert \xi \vert}\left\vert \widehat{\U_1}(\xi)\right\vert + \sum_{n=2}^{+\infty}e^{\frac{\ell_0}{\rho_2}\vert \xi \vert}\left\vert \widehat{\U_n}(\xi)\right\vert =(a)+(b),
\end{equation} 
où nous allons montrer que les  termes $(a)$ et $(b)$ ci-dessus sont contrôlés par une constante $c_1=c_1(\nu,\ell_0,\fe)>0$. \\
\\ 
Pour estimer  le terme $(a)$  de (\ref{major_espec}) nous utilisons  la Proposition \ref{Prop:support-freq} d'où l'on sait que  $supp\,\left(\widehat{\U_1}\right)\subset \left\{\xi \in \Rt:  \vert \xi \vert \leq \frac{\rho_2}{\ell_0}\right\}$ et alors nous pouvons écrire
\begin{eqnarray*}
e^{\frac{\ell_0}{\rho_2}\vert \xi \vert}\left\vert \widehat{\U_1}(\xi)\right\vert \mathds{1}_{supp\,\left(\widehat{\U_1}\right)}(\xi)&\leq &e^{\frac{\ell_0}{\rho_2}\vert \xi \vert}\left\vert \widehat{\U_1}(\xi)\right\vert \mathds{1}_{\left\{\xi \in \Rt:  \vert \xi \vert \leq \frac{\rho_2}{\ell_0}\right\} }(\xi)\leq \sup_{\vert \xi\vert <\frac{\rho_2}{\ell_0}}\left( e^{\frac{\ell_0}{\rho_2}\vert \xi \vert} \left\vert \widehat{ \U_1}(\xi)\right\vert \right)\\
&\leq &e \|\widehat{ \U_1}\|_{L^{\infty}},  
\end{eqnarray*} d'où, étant donné que  $\U_1=-\frac{1}{\nu\Delta}\fe$ et de plus, comme $ \widehat{\fe}\in L^{\infty}(\Rt)$ et comme $supp\,\left( \widehat{\fe}\right) \subset  \left\lbrace \xi \in \Rt : \frac{\rho_1}{\ell_0} \leq \vert \xi \vert \leq \frac{\rho_2}{\ell_0}\right\rbrace$ alors nous avons 
$$ e\Vert \widehat{\U_1}\Vert_{L^{\infty}} =\sup_{\frac{\rho_1}{\ell_0}\leq \vert \xi \vert \leq \frac{\rho_2}{\ell_0}}  \frac{e}{\nu \vert \xi \vert^2}\vert \widehat{\fe}(\xi)\vert \leq e\frac{\ell^{2}_{0}}{\nu \rho^{2}_{1}}\Vert \widehat{\fe} \Vert_{L^{\infty}},$$ et nous obtenons ainsi l'estimation 
\begin{equation}\label{estim_U_1}
(a)=e^{\frac{\ell_0}{\rho_2}\vert \xi \vert}\left\vert \widehat{\U_1}(\xi)\right\vert \leq e\frac{\ell^{2}_{0}}{\nu \rho^{2}_{1}}\Vert \widehat{\fe}\ \Vert_{L^{\infty}}.  
\end{equation}
Nous allons maintenant estimer le terme $(b)$  de (\ref{major_espec}). Il s'agit d'estimer tout d'abord  chaque terme $e^{\frac{\ell_0}{\rho_2}\vert \xi \vert}\left\vert \widehat{\U_n}(\xi)\right\vert$ et pour cela nous utiliserons ici  les nombres de Catalan  $A_n$ définis par la formule (\ref{def_A_n}) où nous  allons montrer  qu'il existe deux constantes $C_1=C_1(\nu,\ell_0)>0$ et   $C_2=C_2(\nu,\ell_0)>0$ telle que pour tout entier  $n\geq 2$ nous avons l'estimation suivante: 
\begin{equation}\label{lemme_th_espec}
e^{\frac{\ell_0}{\rho_2}\vert \xi \vert}\left\vert \widehat{\U_n}(\xi)\right\vert\leq C_1\,A_n\left(C_2 \Vert \fe \Vert_{\dot{H}^{-1}} \right)^n, \quad \text{pour tout}\quad \vert \xi \vert \geq \frac{\rho_1}{\ell_0}.
\end{equation} 
 En effet,  par  l'inégalité (\ref{Inclusion2})  nous écrivons
\begin{equation}\label{EstimateUk1}
e^{\frac{\ell_0}{\rho_2}\vert \xi \vert}\left\vert \widehat{\U_n}(\xi)\right\vert\leq e^{\frac{\ell_0}{\rho_2}\vert \xi \vert}\frac{c}{\nu \vert \xi \vert}\sum_{k=1}^{n-1}\left\vert \widehat{\U_n}\ast \widehat{\U_{n-k}}(\xi)\right\vert,
\end{equation}
et comme  $\vert \xi \vert \geq \frac{\rho_1}{\ell_0}$ nous obtenons
$$e^{\frac{\ell_0}{\rho_2}\vert \xi \vert}\left\vert \widehat{\U_n}(\xi)\right\vert
 \leq e^{\frac{\ell_0}{\rho_2}\vert \xi \vert}\frac{c\ell_0}{\rho_1\nu }\sum_{k=1}^{n-1}\left\vert \widehat{\U_n}\ast \widehat{\U_{n-k}}(\xi)\right\vert.$$
D'autre part, par la Proposition  \ref{Prop:support-freq} nous savons que $supp(\widehat{\U_n})\subset \lbrace\xi \in \Rt: \vert \xi \vert \leq \frac{n\rho_2}{\ell_0}\rbrace$ et alors dans l'estimation précédente nous pouvons écrire
\begin{eqnarray*}
e^{\frac{\ell_0}{\rho_2}\vert \xi \vert}\left\vert \widehat{\U_n}(\xi)\right\vert \mathds{1}_{supp(\widehat{\U_n})}(\xi)
 &\leq &e^{\frac{\ell_0}{\rho_2}\vert \xi \vert}\frac{c\ell_0}{\rho_1\nu }\sum_{k=1}^{n-1}\left\vert \widehat{\U_n}\ast \widehat{\U_{n-k}}(\xi)\right\vert\mathds{1}_{\lbrace\xi \in \Rt: \vert \xi \vert \leq \frac{n\rho_2}{\ell_0}\rbrace}(\xi) \\
 & \leq & \underset{|\xi|\geq \frac{\rho_1}{\ell_{0}}}{sup}\left[ e^{\frac{\ell_0}{\rho_2}\vert \xi \vert}\frac{c\ell_0}{\rho_1\nu }\sum_{k=1}^{n-1}\left\vert \widehat{\U_n}\ast \widehat{\U_{n-k}}(\xi)\right\vert\mathds{1}_{\lbrace\xi \in \Rt: \vert \xi \vert \leq \frac{n\rho_2}{\ell_0}\rbrace}(\xi)\right]=(c),
\end{eqnarray*}
d'où, par l'inégalité de H\"older et en utilisant la définition de la norme $\|\cdot\|_{E}$ donnée par la formule (\ref{norme_E})  nous avons
\begin{eqnarray}
 (c) & \leq &e^n \frac{c\ell_0}{\rho_1\nu} \sum_{k=1}^{n-1}\left\Vert \widehat{\U_k}\ast \widehat{\U_{n-k}}\right\Vert_{L^{\infty}}\leq e^n \frac{c\ell_0}{\rho_1\nu} \sum_{k=1}^{n-1}\left\Vert\U_k\otimes\U_{n-k}\right\Vert_{L^{1}} \leq  e^n \frac{c\ell_0}{\rho_1\nu}\sum_{k=1}^{n-1}\Vert \U_k\Vert_{L^2}\Vert  \U_{n-k}\Vert_{L^{2}}\nonumber\\ 
 &\leq  & e^n \frac{c\ell_0}{\rho_1\nu}  \sum_{k=1}^{n-1}\Vert \U_k\Vert_{E}\Vert  \U_{n-k}\Vert_{E}.\label{Controluk}
\end{eqnarray} 

\`A ce stade, nous avons besoin d'estimer  le terme $\Vert \U_k\Vert_{E}$ et pour cela nous revenons à l'inégalité  (\ref{major_catalan}) pour pouvoir écrire 
$$\Vert \U_k \Vert_{E}\leq\frac{\nu}{C} A_k\left(\frac{C}{\nu} \Vert \U_1 \Vert_{E}\right)^k,$$ où $C>0$ est la constante donnée dans (\ref{cont_forme_bilin}) et $A_n$ sont les nombres de Catalan. Ainsi, par l'identité $ \U_1=-\frac{1}{\nu\Delta}\fe$ et  par le Lemme \ref{lemme:controle-force} nous avons alors l'estimation $\Vert \U_1 \Vert_{E}=\Vert-\frac{1}{\nu\Delta}\fe\Vert_{E}\leq  C_{(\nu,\ell_0)} \Vert \fe \Vert_{\dot{H}^{-1}}$ et donc nous obtenons
$$\Vert \U_k \Vert_{E}\leq\frac{\nu}{C} A_k\left(\frac{C}{\nu} C_{(\nu,\ell_0)} \Vert \fe \Vert_{\dot{H}^{-1}} \right)^k\leq \frac{\nu}{C} A_k \left( C^{'}_{(\nu,\ell_0)} \Vert \fe \Vert_{\dot{H}^{-1}} \right)^k. $$
Une fois que nous disposons de cette estimation nous pouvons revenir à l'estimation  (\ref{Controluk}) ci-dessus  pour écrire  
\begin{eqnarray*}
e^n \frac{c\ell_0}{\rho_1\nu}  \sum_{k=1}^{n-1}\Vert \U_k\Vert_{E}\Vert  \U_{n-k}\Vert_{E}&\leq & e^n \frac{c\ell_0}{\rho_1\nu}  \sum_{k=1}^{n-1}\left( \frac{\nu}{C}A_k \left(C^{'}_{(\nu,\ell_0)} \Vert \fe \Vert_{\dot{H}^{-1}} \right)^k\right)\left(\frac{\nu}{C}A_{n-k} \left( C^{'}_{(\nu,\ell_0)} \Vert \fe \Vert_{\dot{H}^{-1}}\right)^{n-k}\right)\\
&\leq &e^n \frac{c\ell_0 \nu}{\rho_1C^2}  \left( C^{'}_{(\nu,\ell_0)} \Vert \fe \Vert_{\dot{H}^{-1}}\right)^{n}\sum_{k=1}^{n-1} A_k A_{n-k}, 
\end{eqnarray*}
et alors par la formule de récurrence $\ds{A_n=\sum_{k=1}^{n-1}A_k\,A_{n-k}}$ donnée dans  (\ref{def_A_n}) nous avons 
\begin{equation*}
e^n \frac{c\ell_0 \nu}{\rho_1C^2} \sum_{k=1}^{n-1}\Vert \U_k\Vert_{E}\Vert  \U_{n-k}\Vert_{E} \leq  \frac{c\ell_0 \nu}{\rho_1C^2} A_n \left( eC^{'}_{(\nu,\ell_0)} \Vert \fe \Vert_{\dot{H}^{-1}}\right)^n,
\end{equation*} et alors si on pose les constantes $C_1=\frac{c\ell_0 \nu}{\rho_1C^2}>0$ et $C_2= eC^{'}_{(\nu,\ell_0)}>0$  et nous obtenons ainsi l'estimation (\ref{lemme_th_espec}). \\
\\
Une fois que l'on a cette estimation nous  pouvons écrire 
\begin{equation*}
\ds{\sum_{n=2}^{+\infty}e^{\frac{\ell_0}{\rho_2}\vert \xi \vert}\left\vert \widehat{\U_n}(\xi)\right\vert\leq C_1 \sum_{n=2}^{+\infty}A_n \left( C_2\Vert \fe \Vert_{\dot{H}^{-1}}\right)^n},
\end{equation*} et par la propriété (\ref{suite_gener}) des nombres de Catalan, si $\Vert \fe \Vert_{\dot{H}^{-1}}<\eta$, avec $\eta>0$ une constante suffisamment petite, alors nous avons   
$$\ds{C_1 \sum_{n=2}^{+\infty}A_n \left( C_2\Vert \fe \Vert_{\dot{H}^{-1}}\right)^n\leq C_1 \sum_{n=1}^{+\infty}A_n \left( C_2\Vert \fe \Vert_{\dot{H}^{-1}}\right)^n=C_1\frac{1-\sqrt{1-4C_2\Vert \fe \Vert_{\dot{H}^{-1}}}}{2}},$$
d'où nous obtenons l'estimation 
\begin{equation}\label{estim_serie}
(b)=\sum_{n=2}^{+\infty}e^{\frac{\ell_0}{\rho_2}\vert \xi \vert}\left\vert \widehat{\U_n}(\xi)\right\vert\leq C_1\frac{1-\sqrt{1-4C_2\Vert \fe \Vert_{\dot{H}^{-1}}}}{2}. 
\end{equation}
Maintenant que nous disposons des estimations (\ref{estim_U_1}) et (\ref{estim_serie}) nous pouvons les remplacer dans l'estimation (\ref{major_espec}) en fixant la constante  $c_1=\max\left(e\frac{\ell^{2}_{0}}{\nu \rho^{2}_{1}}\Vert \widehat{\fe}\ \Vert_{L^{\infty}}, C_1\frac{1-\sqrt{1-4C_2\Vert \fe \Vert_{\dot{H}^{-1}}}}{2} \right)>0$ et la constante $c_2=\frac{\ell_0}{\rho_2}>0$, et    nous avons finalement l'estimation cherchée (\ref{dec-frec-lam-theo-1}).  \finpv\\
\\ 
Comme nous pouvons observer la décroissance fréquentielle exponentielle obtenue dans le Théorème \ref{Theo:Ossen_decroissance} repose essentiellement sur la localisation fréquentielle de la force extérieure  $\fe$ qui entraîne une localisation fréquentielle convenable des termes $\U_n$ lorsqu'on décompose la solution $\U$ comme $\U=\sum_{n=1}^{+\infty}\U_n$.   \\
\\ 
Dans la partie $B)$ ci-dessous, nous étudions cette décroissance fréquentielle  dans un cadre un peu plus général où la force extérieure $\fe$ n'est pas  localisée en fréquence mais elle vérifie une décroissance fréquentielle exponentielle.  
\subsection*{B) Deuxième approche:  décroissance fréquentielle de la force}
Comme annoncé nous supposons ici que la force $\fe\in \dot{H}^{-1}(\Rt)$ vérifie une décroissance fréquentielle $ \vert \widehat{\fe}(\xi)\vert \lesssim e^{-\vert \xi\vert}$ et  nous étudierons un  comportement similaire pour les solutions $\U\in \dot{H}^{1}(\Rt)$ des équations de Navier-Stokes stationnaires (\ref{N-S-stationnaire}). \\
\\
On sait que les propriétés de décroissance de la transformée de Fourier d'une fonction peuvent être regardées comme des propriétés de régularité de telle fonction dans les espaces de Sobolev $\dot{H}^s(\Rt)$ (avec $s>0$), néanmoins, ces espaces ne nous permettent pas d'étudier de façon précise une décroissance fréquentielle exponentielle  et nous avons besoin de passer par le cadre de la classe de Gevrey qui nous introduisons comme suit: la première chose à faire est de définir l'opérateur $e^{\beta  \sqrt{-\Delta}}$ de la façon suivante:
\begin{Definition}\label{Def_noyau_posisson} Soit $\beta>0$. On définit  l'opérateur $e^{\beta \sqrt{-\Delta}}$ au niveau de Fourier par l'identité $\mathcal{F}[ e^{\beta \sqrt{-\Delta}}\varphi](\xi)=e^{\beta \vert \xi \vert}\mathcal{F}[\varphi](\xi)$, pour tout $\varphi$ qui appartient à la classe de Schwartz.
\end{Definition}
Cet opérateur nous sera très  utile pour l'étude de la décroissance fréquentielle exponentielle car le symbole $e^{\beta \vert \xi \vert}$ nous permettra  capturer  ce comportement comme nous le verrons dans la définition suivante:

\begin{Definition}[La classe de Gevrey]\label{Classe de Gevrey} Soit l'espace de Sobolev $\dot{H}^s(\Rt)$ avec $s\in \mathbb{R}$, soit $\beta>0$ et soit l'opérateur $e^{\beta \sqrt{-\Delta}}$ donné dans la Définition \ref{Def_noyau_posisson}. On définit la classe de Gevrey $G^{\beta}_{s}(\Rt)$ par 
$$ G^{\beta}_{s}(\Rt)=\left\lbrace g\in \dot{H}^s(\Rt):  \Vert e^{\beta\sqrt{-\Delta}} g \Vert_{\dot{H}^s}<+\infty \right\rbrace. $$
\end{Definition} 
Grâce  à cette définition nous observons que la quantité $\ds{\Vert e^{\beta\sqrt{-\Delta}} g \Vert_{\dot{H}^s}}$ exprime la décroissance exponentielle de la transformée de Fourier de la fonction $g$. En effet, il suffit d'écrire 
$$\Vert e^{\beta \sqrt{-\Delta}}g \Vert^{2}_{\dot{H}^s}=\int_{\Rt}\vert \xi \vert^{2s} e^{2\beta\vert \xi \vert}\vert \widehat{g}(\xi)\vert^2 d\xi,$$ pour observer que si cette quantité est finie alors il faut que $\widehat{g}(\xi)$ ait une décroissance exponentielle à l'infini. \\
\\
Observons aussi que l'on a $G^{\beta}_{s}(\Rt)\subset  \dot{H}^s(\Rt)$ et en plus, l'espace  $\dot{H}^s(\Rt)$ étant un espace de Banach pour $\vert s \vert <\frac{3}{2}$ alors on a que pour ces valeurs du paramètre $s$ la classe de Gevrey  $G^{\beta}_{s}(\Rt)$ est aussi  un espace de Banach  muni de la norme $\Vert e^{\beta\sqrt{-\Delta}}(\cdot)\Vert_{\dot{H}^S}$.\\
\\
Cet espace fonctionnel  a été introduit par M. Gevrey en 1918 dans \cite{Gevrey} pour l'étude analytique de quelques propriétés de régularité du noyau de la chaleur  et cet espace a été aussi utilisé pour l'étude de l'analycité en temps des solutions mild des équations de Navier-Stokes non stationnaires (voir l'article \cite{FoiasTemam} de C. Foias et R. Temam et l'article \cite{PGLRArt1} de P.G. Lemarié-Rieusset pour plus de références à ce sujet). Dans cette section nous utiliserons la classe de Gevrey pour étudier la décroissance fréquentielle exponentielle  des solutions des équations de Navier-Stokes dans le cadre stationnaire et nous avons ainsi une première caractérisation de ce comportement. \\
\\
En effet,  nous allons montrer que si la force $\fe$ a une décroissance fréquentielle exponentielle, qui est exprimée par le fait que $\fe \in G^{\beta}_{-\frac{3}{2}}(\Rt)$, et si nous supposons en plus un contrôle sur la taille de cette force alors nous pouvons utiliser le principe de contraction de Picard  pour construire une solution $\U$ des équations (\ref{N-S-stationnaire}) dans la classe de Gevrey  $G^{\beta}_{\frac{1}{2}}(\Rt)$. Le choix du paramètre $s=\frac{1}{2}$ dans cette classe de Gevrey (ainsi que le choix du paramètre $s=-\frac{3}{2}$ dans la classe $G^{\beta}_{-\frac{3}{2}}(\Rt)$) est purement technique et nous avons:\\
\begin{Proposition}\label{Prop:dec-freq-gevrey} Soit $\beta>0$ et soit $\fe\in G^{\beta}_{-\frac{3}{2}}(\Rt)$ une force à divergence nulle. Il existe une constante $\eta_1>0$ telle que si $\Vert e^{\beta \sqrt{-\Delta}}\fe \Vert_{\dot{H}^{-\frac{3}{2}}}<\eta_1$ alors il existe $\U\in G^{\beta}_{\frac{1}{2}}(\Rt)$ solution des équations de Navier-Stokes stationnaires (\ref{N-S-stationnaire}) qui est l'unique solution dans la boule $\Vert e^{\beta \sqrt{-\Delta}}\U \Vert_{\dot{H}^{\frac{1}{2}}}<2\eta_1$.  \\
\end{Proposition}   \textbf{Preuve.} Nous écrivons  les équations (\ref{N-S-stationnaire}) comme un problème de point fixe 
$$ \U=\frac{1}{\nu} \P\left( \frac{1}{\Delta}div(\U \otimes \U) \right)-\frac{1}{\nu \Delta}\fe,$$ et nous allons appliquer le principe de contraction Picard dans la classe de Gevrey  $G^{\beta}_{\frac{1}{2}}(\Rt)$. Pour vérifier la continuité de la forme bilinéaire ci-dessus sur cet espace, nous avons besoin  de l'inégalité suivante:
\begin{Lemme}\label{lemme:picard-gevrey} On a $\ds{\left\Vert e^{\beta \sqrt{-\Delta}}\left[\frac{1}{\nu} \P\left( \frac{1}{\Delta}div(\U \otimes \U) \right)\right] \right\Vert_{\dot{H}^{\frac{1}{2}}}\leq \frac{C}{\nu} \left\Vert \left[ e^{\beta\sqrt{-\Delta}}\U\right] \otimes  \left[ e^{\beta\sqrt{-\Delta}}\U\right] \right\Vert_{\dot{H}^{-\frac{1}{2}}}.}$ 
\end{Lemme} 
\pv On commence par écrire
\begin{eqnarray*}
	& & \left\Vert e^{\beta \sqrt{-\Delta}}\left[\frac{1}{\nu} \P\left( \frac{1}{\Delta}div(\U \otimes \U) \right)\right] \right\Vert^{2}_{\dot{H}^{\frac{1}{2}}}=\frac{1}{\nu^2} \int_{\Rt} \vert \xi \vert e^{2\beta \vert \xi \vert} \left\vert \mathcal{F} \left(  \P\left( \frac{1}{\Delta}div(\U \otimes \U) \right)\right)(\xi) \right\vert^2 d\xi \\
	& \leq & \frac{C}{\nu^2} \int_{\Rt} \vert \xi \vert e^{2\beta \vert \xi \vert} \left[ \frac{1}{\vert \xi \vert} \vert \mathcal{F}[\U] \ast \mathcal{F}[\U] (\xi)\right\vert^2 d\xi \leq  \frac{C}{\nu^2} \int_{\Rt} \frac{1}{\vert \xi \vert} \left[  e^{\beta \vert \xi \vert}  \vert \mathcal{F}[\U] \ast \mathcal{F}[\U] (\xi)\vert\right]^2 d\xi,
\end{eqnarray*}  et dans le dernier terme ci-dessus, nous utilisons l'inégalité triangulaire  pour écrire
\begin{eqnarray*}
	& & e^{\beta \vert \xi \vert}  \vert \mathcal{F}[\U] \ast \mathcal{F}[\U] (\xi) \leq   e^{\beta \vert \xi \vert}  \int_{\Rt}\vert \mathcal{F}[\U](\xi-\zeta)\mathcal{F}[\U](\zeta) \vert d\zeta \\
	&\leq & \int_{\Rt}\vert  e^{\beta \vert\xi-\zeta \vert}\mathcal{F}[\U](\xi-\zeta)e^{\beta \vert\zeta \vert}\mathcal{F}[\U](\zeta) \vert d\zeta \leq \left(\vert e^{\beta \vert\xi\vert}\mathcal{F}[\U]\vert \ast  \vert e^{\beta \vert\xi\vert}\mathcal{F}[\U] \vert\right)(\xi),
\end{eqnarray*} et alors, en revenant à l'estimation précédente nous obtenons 
\begin{eqnarray*}
	\frac{C}{\nu^2} \int_{\Rt} \frac{1}{\vert \xi \vert} \left[  e^{\beta \vert \xi \vert}  \vert \mathcal{F}[\U] \ast \mathcal{F}[\U] (\xi)\vert\right]^2 d\xi &\leq &\frac{C}{\nu^2}  \int_{\Rt}\frac{1}{\vert \xi \vert } \left\vert\left(\vert e^{\beta \vert\xi\vert}\mathcal{F}[\U]\vert \ast  \vert e^{\beta \vert\xi\vert}\mathcal{F}[\U] \vert\right)(\xi)\right\vert^2 d\xi\\
	&\leq &\frac{C}{\nu} \Vert e^{\beta \sqrt{-\Delta}} \U \otimes  e^{\beta \sqrt{-\Delta}} \U \Vert^{2}_{\dot{H}^{-\frac{1}{2}}}.
\end{eqnarray*} \finpv 
\\
Maintenant, par  le Lemme \ref{lemme:picard-gevrey},  par les inégalités de Hardy-Littlewood-Sobolev et l'inégalité de H\"older nous pouvons écrire
\begin{eqnarray*}
 & & \left\Vert e^{\beta \sqrt{-\Delta}}\left[\frac{1}{\nu} \P\left( \frac{1}{\Delta}div(\U \otimes \U) \right)\right] \right\Vert_{\dot{H}^{\frac{1}{2}}} \leq  \frac{C}{\nu} \left\Vert \left[ e^{\beta\sqrt{-\Delta}}\U\right] \otimes  \left[ e^{\beta\sqrt{-\Delta}}\U\right] \right\Vert_{\dot{H}^{-\frac{1}{2}}} \\
 &\leq & \frac{C}{\nu} \left\Vert \left[ e^{\beta\sqrt{-\Delta}}\U\right] \otimes  \left[ e^{\beta\sqrt{-\Delta}}\U\right] \right\Vert_{L^{\frac{3}{2}}} \leq \frac{C}{\nu} \left\Vert  e^{\beta\sqrt{-\Delta}}\U \right\Vert_{L^{3}} \left\Vert  e^{\beta\sqrt{-\Delta}}\U \right\Vert_{L^{3}}  \\
 & \leq & \frac{C}{\nu} \left\Vert  e^{\beta\sqrt{-\Delta}}\U \right\Vert_{\dot{H}^{\frac{1}{2}}}  \left\Vert  e^{\beta\sqrt{-\Delta}}\U \right\Vert_{\dot{H}^{\frac{1}{2}}},
\end{eqnarray*} et ainsi nous avons vérifié la continuité de la forme bilinéaire ci-dessus sur la classe de Gevrey $G^{\beta}_{\frac{1}{2}}(\Rt)$.\\
\\
Nous allons maintenant estimer le terme $\left\Vert e^{\beta \sqrt{-\Delta}}\left( \frac{1}{\nu \Delta} \fe \right) \right\Vert_{\dot{H}^{\frac{1}{2}}}$ et  étant donné que  $\fe \in G^{\beta}_{-\frac{3}{2}}(\Rt)$ nous pouvons écrire directement
$\ds{  \left\Vert e^{\beta \sqrt{-\Delta}}\left( \frac{1}{\nu \Delta} \fe \right) \right\Vert_{\dot{H}^{\frac{1}{2}}} \leq \frac{1}{\nu} \Vert e^{\beta\sqrt{-\Delta}} \fe \Vert_{\dot{H}^{-\frac{3}{2}}}}$. \\
\\
Ainsi, il existe une constante $\eta_1>0$ (suffisamment petite) telle que si  $ \Vert e^{\beta\sqrt{-\Delta}} \fe \Vert_{\dot{H}^{-\frac{3}{2}}} <\eta_1 $ alors par le principe de contraction de Picard il existe $\U \in G^{\beta}_{\frac{1}{2}}(\Rt)$ solution des équations (\ref{N-S-stationnaire}) qui est l'unique solution dans la boule $\Vert e^{\beta \sqrt{-\Delta}}\U \Vert_{\dot{H}^{\frac{1}{2}}}<2\eta_1$. \finpv \\
Dans la Proposition \ref{Prop:dec-freq-gevrey}  nous pouvons observer que  comme $\U \in G^{\beta}_{\frac{1}{2}}(\Rt)$ alors nous avons $$\ds{ \int_{\Rt} \vert \xi \vert e^{2\beta \vert \xi \vert} \vert \widehat{\U}(\xi)\vert^2 d\xi <+\infty},$$ ce qui exprime le fait que la solution $\U$ a une décroissance fréquentielle exponentielle  aux hautes fréquences et cette information est mesurée en termes de la norme $L^2$. Néanmoins, comme expliqué à la Section  \ref{Sec:regime-lam-turb}, dans le cadre d'un fluide en régime laminaire la décroissance fréquentielle exponentielle est censée être observée aussi aux basses fréquences et non seulement aux hautes fréquences comme l'on peut regarder dans la quantité ci-dessus. \\
\\
Ainsi, ce résultat n'exprime pas de façon tout à fait satisfaisante la décroissance fréquentielle exponentielle d'un fluide en régime laminaire  mais il nous sera fondamental pour étudier une décroissance fréquentielle plus précise  dans le cadre d'un fluide en régime laminaire (caractérisé par un contrôle sur la force).\\
\begin{Theoreme}[Deuxième étude de la décroissance fréquentielle: cadre laminaire]\label{Theo:gevrey-lam} Soit $\fe \in \dot{H}^{-2}\cap \dot{H}^{-1}(\Rt)$ une force à divergence nulle telle que $\Vert \fe \Vert_{\dot{H}^{-2}} +\Vert \fe \Vert_{\dot{H}^{-1}}<\eta$ pour une constante $\eta>0$ et qui vérifie 
\begin{equation}\label{cond-force-teo2-lam}
\sup_{\xi \in \Rt} \frac{1}{\vert \xi \vert} e^{2\beta \vert \xi \vert} \vert \widehat{\fe}(\xi)\vert \leq c_0 <+\infty,
\end{equation} où $\beta>0$ et $c_0>0$ est une constante suffisamment petite. Alors,  il  existe  $\U \in L^2 \cap \dot{H}^{1}(\Rt)$ solution des équations de Navier-Stokes stationnaires (\ref{N-S-stationnaire}) qui vérifie la décroissance fréquentielle ponctuelle:
\begin{equation}\label{dec-U-teo2-lam}
\vert \widehat{\U}(\xi)\vert \leq c_1 \frac{1}{\vert \xi \vert} e^{-\beta \vert \xi \vert},
\end{equation} pour tout $\vert \xi \vert>0$ et où $c_1>0$ est une constante qui dépend de $\U$. 
\end{Theoreme}   
\dm Nous commençons donc par construire  une solution $\U \in L^2\cap \dot{H}^{1}(\Rt)$.
\begin{Lemme}\label{Lemme:picard-L2-H1} Il existe une constante $\eta>0$ telle que si $\Vert \fe \Vert_{\dot{H}^{-2}} + \Vert \fe \Vert_{\dot{H}^{-1}}<\eta$ alors il existe $\U \in L^2 \cap \dot{H}^{1}(\Rt)$ solution des équations (\ref{N-S-stationnaire}) qui est l'unique solution qui vérifie $\Vert \U \Vert_{\dot{H}^{\frac{1}{2}}}<2\eta$. 
\end{Lemme} \pv On commence par écrire la solution $\U$ comme 
\begin{equation}\label{eq-aux-U}
\U=\frac{1}{\nu}\P \left( \frac{1}{\Delta}div (\U \otimes \U) \right)-\frac{1}{\nu \Delta}\fe,
\end{equation} et nous allons appliquer le principe de contraction de Picard dans  l'espace $E=L^2\cap \dot{H}^1(\Rt)$ (où l'on a $L^{2}\cap \dot{H}^{1}(\Rt)\subset \dot{H}^{\frac{1}{2}}(\Rt)$).  Pour cela on a besoin d'estimer chaque terme de la norme $\Vert \cdot \Vert_{E}=\Vert \cdot \Vert_{L^2}+\Vert \cdot \Vert_{\dot{H}^1}+\Vert \cdot \Vert_{\dot{H}^{\frac{1}{2}}}$.\\
\\
En effet, pour le premier de $\Vert \cdot \Vert_{E}$, par l'estimation (\ref{estim1}) et les inégalités de Hardy-Littlewood-Sobolev nous avons 
\begin{equation}\label{1}
\left\Vert  \frac{1}{\nu}\P \left( \frac{1}{\Delta}div (\U \otimes \U) \right) \right\Vert_{L^2}\leq \frac{c_1}{\nu} \Vert \U \Vert_{L^2}\Vert \U \Vert_{L^3}\leq \frac{c_1}{\nu}\Vert \U \Vert_{L^2} \Vert \U \Vert_{\dot{H}^{\frac{1}{2}}} \leq \frac{c_1}{\nu}\Vert \U \Vert_E \Vert \U \Vert_E.   
\end{equation}
Ensuite, pour le deuxième terme de $\Vert \cdot \Vert_{E}$, par l'estimation (\ref{estim2}) et toujours par les inégalités de Hardy-Littlewood-Sobolev nous avons 
\begin{equation}\label{2}
\left\Vert  \frac{1}{\nu}\P \left( \frac{1}{\Delta}div (\U \otimes \U) \right) \right\Vert_{\dot{H}^1}\leq \frac{c_1}{\nu} \Vert \U \Vert_{\dot{H}^1}\Vert \U \Vert_{L^3}\leq \frac{c_1}{\nu}\Vert \U \Vert_{\dot{H}^1} \Vert \U \Vert_{\dot{H}^{\frac{1}{2}}} \leq \frac{c_1}{\nu}\Vert \U \Vert_E \Vert \U \Vert_E. 
\end{equation}\label{3}
Finalement, par l'estimation (\ref{estim3}) et l'estimation (\ref{estim3}) nous avons 
\begin{equation}
\left\Vert  \frac{1}{\nu}\P \left( \frac{1}{\Delta}div (\U \otimes \U) \right) \right\Vert_{H^{\frac{1}{2}}}  \leq \frac{c_3}{\nu}\Vert \U \Vert_{\dot{H}^{\frac{1}{2}}} \Vert \U \Vert_{\dot{H}^{\frac{1}{2}}}\leq \frac{c_3}{\nu}\Vert \U \Vert_E \Vert \U \Vert_E.
\end{equation}
Ainsi, par les estimations (\ref{1}), (\ref{2}) et (\ref{3}) nous avons la continuité de la forme bilinéaire $\frac{1}{\nu}\P \left( \frac{1}{\Delta}div (\U \otimes \U) \right)$ dans l'espace $E=L^2\cap \dot{H}^1(\Rt)$.\\
\\ 
Nous étudions maintenant le terme $\left\Vert \frac{1}{\nu \Delta}\fe \right\Vert_{E}$. Comme $\fe \in \dot{H}^{-2}\cap \dot{H}^{-1}(\Rt)$ nous pouvons écrire
\begin{eqnarray*}
	\left\Vert \frac{1}{\nu \Delta}\fe \right\Vert_{E} & \leq & \left\Vert \frac{1}{\nu \Delta}\fe \right\Vert_{L^2}+\left\Vert \frac{1}{\nu \Delta}\fe \right\Vert_{\dot{H}^1}+\left\Vert \frac{1}{\nu \Delta}\fe \right\Vert_{\dot{H}^{\frac{1}{2}}} \leq  \frac{c}{\nu} \left[ \Vert \fe \Vert_{\dot{H}^{-2}}+\Vert \fe \Vert_{\dot{H}^{-1}} +\Vert \fe \Vert_{\dot{H}^{-\frac{3}{2}}}\right] \\
	& \leq & \frac{c}{\nu} \left[ \Vert \fe \Vert_{\dot{H}^{-2}}+\Vert \fe \Vert_{\dot{H}^{-1}} +\Vert \fe \Vert^{\frac{1}{2}}_{\dot{H}^{-2}}\Vert \fe \Vert^{\frac{1}{2}}_{\dot{H}^{-1}}\right] \leq \frac{c}{\nu}\left[\Vert \fe \Vert_{\dot{H}^{-2}}+\Vert \fe \Vert_{\dot{H}^{-1}}\right].
\end{eqnarray*}   De cette façon, pour pouvoir appliquer le principe de contraction de Picard il suffit d'avoir le  contrôle  $\Vert \fe \Vert_{\dot{H}^{-2}} +\Vert \fe \Vert_{\dot{H}^{-1}}<\eta$, avec $\eta>0$ une constante suffisamment petite. \finpv\\
 Nous allons maintenant montrer que la solution $\U$ ci-dessus vérifie la décroissance fréquentielle (\ref{dec-U-teo2-lam}). Dans chaque terme de  l'identité (\ref{eq-aux-U}) nous prenons  la transformée de Fourier et nous avons 
$$ \vert \widehat{\U}(\xi)\vert \leq \frac{c}{\vert \xi \vert}\vert \widehat{\U}\ast \widehat{\U}(\xi)\vert +\frac{1}{\nu \vert \xi \vert^2}\vert \widehat{\fe}(\xi)\vert,$$ d'où, en multipliant par $e^{\beta \vert \xi \vert}$ à chaque côté de cette inégalité  on a 
\begin{equation}\label{estim2-gevrey}
e^{\beta \vert \xi \vert}\vert \widehat{\U}(\xi)\vert \leq \frac{c}{\vert \xi \vert}e^{\beta \vert \xi \vert} \vert \widehat{\U}\ast \widehat{\U}(\xi)\vert + \frac{1}{\nu \vert \xi \vert^2}e^{\beta \vert \xi \vert} \vert \widehat{\fe}(\xi)\vert, 
\end{equation}  et l'on cherche à estimer  chaque terme à droite de cette inégalité.\\
\\
Pour estimer le premier terme à droite de (\ref{estim2-gevrey}) nous avons besoin  de vérifier tout d'abord que la solution $\U$ obtenue par le biais du Lemme \ref{Lemme:picard-L2-H1} appartient à la classe  de Gevrey $G^{\beta}_{\frac{1}{2}}(\Rt)$. \\
\\
En effet, comme la force $\fe$ vérifie (\ref{cond-force-teo2-lam}) nous commençons par vérifier  que cette force appartient aussi à la classe de Gevrey $G^{\beta}_{-\frac{3}{2}}(\Rt)$: par (\ref{cond-force-teo2-lam}) nous avons l'estimation 
$$ \vert \widehat{\fe}(\xi)\vert \leq c_0 \vert \xi \vert e^{-2 \beta \vert \xi \vert},$$ et alors nous pouvons écrire 
$$ \Vert e^{\beta \sqrt{-\Delta}}  f \Vert^{2}_{\dot{H}^{-\frac{3}{2}}}=\int_{\Rt} \frac{1}{\vert \xi \vert^3} e^{2\beta \vert \xi \vert}\vert \widehat{\fe}(\xi)\vert^2 d\xi  \leq c^{2}_{0} \int_{\Rt} \frac{1}{\vert \xi \vert^3} e^{2\beta \vert \xi \vert} \left[  \vert \xi \vert^2 e^{-4\beta\vert \xi \vert} \right]  d\xi =c^{2}_{0} C_1<+\infty.$$
Ainsi, pour $\eta_1>0$ la constante de la Proposition \ref{Prop:dec-freq-gevrey} et pour $\eta>0$ la constante du Lemme \ref{Lemme:picard-L2-H1} on fixe maintenant la constante $c_0>0$ de sorte que $c^{2}_{0}C_1<\min(\eta^2,\eta^{2}_1)$, et nous obtenons  $\Vert e^{\beta \sqrt{-\Delta}}\fe \Vert_{\dot{H}^{-\frac{3}{2}}}<\min(\eta,\eta_1)$ et de cette façon,  par la Proposition \ref{Prop:dec-freq-gevrey} il existe  $\V\in G^{\beta}_{\frac{1}{2}}(\Rt)$ solution des équations de Navier-Stokes stationnaires (\ref{N-S-stationnaire}) qui est l'unique solution dans la boule $\Vert e^{\beta \sqrt{-\Delta}}\V \Vert_{\dot{H}^{\frac{1}{2}}}<2\min (\eta_1,\eta)$, mais, nous avons   l'estimation 
$\ds{ \Vert \V \Vert_{\dot{H}^{\frac{1}{2}}}\leq \Vert e^{\beta \sqrt{-\Delta}}\V \Vert_{\dot{H}^{\frac{1}{2}}}}$ d'où nous obtenons 
$\Vert \V \Vert_{\dot{H}^{\frac{1}{2}}}<2\eta$ et ainsi, par l'unicité de la solution $\U$ construite dans le Lemme \ref{Lemme:picard-L2-H1} nous avons $\V=\U$ et donc $\Vert e^{\beta \sqrt{-\Delta}}\U \Vert_{\dot{H}^{\frac{1}{2}}}<+\infty$. \\
\\
Une fois que nous avons $\U\in G^{\beta}_{\frac{1}{2}}(\Rt)$ nous pouvons maintenant estimer le premier terme à droite de (\ref{estim2-gevrey}). En effet, par l'inégalité triangulaire et l'inégalité de Cauchy-Schwarz nous écrivons
\begin{equation}\label{estim3-gevrey}
e^{\beta \vert \xi \vert} \vert \widehat{\U}\ast \widehat{\U}(\xi)\vert \leq \int_{\Rt}\left(e^{\beta\vert \zeta \vert} \vert \widehat{\U}(\zeta)d \zeta \right)\left( e^{\beta \vert \xi -\zeta\vert}\widehat{\U}(\xi-\zeta)\right)d\zeta\leq \Vert e^{\beta\sqrt{-\Delta}} \U \Vert^{2}_{L^2}, 
\end{equation} et nous devons maintenant contrôler le dernier terme à droite ci-dessus.  Par le Lemme \ref{Lemme:picard-L2-H1} nous avons  $\U\in L^2(\Rt)$ et comme $\U\in G^{\beta}_{\frac{1}{2}}(\Rt)$  nous pouvons alors écrire  
\begin{eqnarray*}
\Vert e^{\beta\sqrt{-\Delta}} \U \Vert^{2}_{L^2}&=&\int_{\Rt}e^{2\beta \vert \xi \vert}\vert \widehat{\U}(\xi)\vert^2d\xi = \int_{\vert \xi \vert \leq 1}e^{2\beta \vert \xi \vert}\vert \widehat{\U}(\xi)\vert^2d\xi + \int_{\vert \xi \vert >1}e^{2\beta \vert \xi \vert}\vert \widehat{\U}(\xi)\vert^2d\xi \\
&\leq &\int_{\vert \xi \vert \leq 1}e^{2\beta }\vert \widehat{\U}(\xi)\vert^2d\xi  + \int_{\vert \xi \vert >1} \vert \xi \vert e^{2\beta \vert \xi \vert}\vert \widehat{\U}(\xi)\vert^2d\xi\\
&\leq & e^{2\beta}\Vert \U \Vert^{2}_{L^2}+\Vert e^{\beta \sqrt{-\Delta}}\U \Vert^{2}_{\dot{H}^{\frac{1}{2}}}=C(\beta,\U)<+\infty.
\end{eqnarray*} De cette façon, en revenant à l'estimation (\ref{estim3-gevrey}) nous obtenons 
\begin{equation}\label{estim-bilin-gevrey}
e^{\beta \vert \xi \vert} \vert \widehat{\U}\ast \widehat{\U}(\xi)\vert \leq  C(\beta,\U).  \\
\end{equation} Nous allons maintenant estimer le  deuxième à droite de (\ref{estim2-gevrey}) et pour cela par l'estimation (\ref{cond-force-teo2-lam}) nous écrivons 
\begin{equation}\label{estim4-gevrey}
e^{\beta \vert \xi \vert} \vert \widehat{\fe}(\xi)\vert \leq e^{\beta \vert\xi \vert } \left[ c_0 \vert \xi \vert e^{-2\beta\vert \xi \vert}\right]\leq c_0 \vert \xi \vert. \\
\end{equation} 
Une fois que nous disposons des estimations  (\ref{estim-bilin-gevrey}) et (\ref{estim4-gevrey}) nous les remplaçons  dans l'estimation  (\ref{estim2-gevrey}) pour pouvoir écrire
$$e^{\beta \vert \xi \vert}\vert \widehat{\U}(\xi)\vert \leq \frac{c C(\beta,\U)}{\nu \vert \xi \vert} +\frac{c_0}{\nu \vert \xi \vert}= \frac{C(\beta,\U)+c_0}{\nu} \frac{1}{\vert \xi \vert} =\frac{c_1}{\vert \xi \vert},$$ et ainsi ce théorème est maintenant démontré.  \finpv  \\
Dans la démonstration du Théorème \ref{Theo:gevrey-lam}  que nous venons de faire nous observons que si l'on a un contrôle sur la force  $\fe$ (et donc l'on est dans un régime laminaire) et si cette force a une décroissance fréquentielle exponentielle alors nous pouvons construire une solution $\U$ des équations de Navier-Stokes stationnaires avec la même décroissance fréquentielle.\\
\\
Dans la section qui suit nous verrons comment la classe de Gevrey nous permet aussi d'étudier une décroissance fréquentielle exponentielle des solutions $\U$ mais dans un cadre plus général où l'on ne suppose aucun contrôle sur le force et donc le fluide peut atteindre un état turbulent.  
\subsection{Quelques résultats dans le régime turbulent}\label{Sec:Sec:decroissance-freq-turb}
Dans la section précédente nous avons introduit la classe de Gevrey  $G^{\beta}_{s}(\Rt)$  pour étudier la décroissance fréquentielle  des solutions $\U$ des équations de Navier-Stokes stationnaires (\ref{N-S-stationnaire}); dont les résultats obtenus reposent sur le contrôle qu'on a fait sur la force $\fe$  et donc ces résultats  correspondent au cadre d'un fluide en régime laminaire. \\
\\
Le but de cette section est de généraliser  ces résultats précédents au cadre d'un fluide en régime turbulent et  alors nous ne ferons ici aucun contrôle sur la la taille de la force $\fe$. Dans la Section \ref{Sec:regime-lam-turb} nous avons expliqué que si le fluide est en état turbulent alors la décroissance fréquentielle exponentielle de la solution $\U$ est censée être observée uniquement aux hautes fréquences et nous allons voir comment la classe de Gevrey nous permettra  capturer ce comportement fréquentiel  dans ce cadre turbulent.\\
\\
En effet, pour n'importe quelle force $\fe \in \dot{H}^{-1}(\Rt)$ par  le Théorème \ref{Theo:exstence-sol-stationnaires} nous avons l'existence d'au moins une solution $\U\in \dot{H}^{1}(\Rt)$ des équations (\ref{N-S-stationnaire}) et dans le Théorème \ref{Theo:dec-freq-turb} ci-dessous nous allons  montrer que si la force  vérifie en plus $\fe \in G^{\beta_1}_{-1}(\Rt)$ (où $G^{\beta_1}_{-1}(\Rt) \subset \dot{H}^{-1}(\Rt)$)  alors \emph{toute}  solution $\U\in \dot{H}^{1}(\Rt)$  associée à cette force vérifie  $\U \in G^{\beta_2}_{1}(\Rt)$  pour une certaine constante $\beta_2>0$.  \\
\\
De cette façon, par la définition de la classe de Gevrey nous avons alors que la quantité 
\begin{equation}\label{dec-frec-turb-motiv}
\ds{ \int_{\Rt}\vert \xi \vert^2 e^{2\beta_2 \vert \xi \vert}\vert \widehat{\U}(\xi)\vert^2 d\xi},
\end{equation} est une quantité finie et donc nous pouvons observer que $\vert \widehat{\U}(\xi)\vert$ a une décroissance exponentielle aux hautes fréquences.  \\
\begin{Theoreme}[\'Etude de la décroissance fréquentielle: cadre turbulent]\label{Theo:dec-freq-turb} Soit $\fe \in \dot{H}^{-1}(\Rt)$ une force stationnaire à divergence nulle et soit $\beta_1>0$. Si  $\fe \in G^{\beta_1}_{-1}(\Rt)$  alors toute solution  $\U \in \dot{H}^{1}(\Rt)$ des équations de Navier-Stokes stationnaires (\ref{N-S-stationnaire})  associée à cette force appartient à la classe de Gevrey $G^{\beta_2}_{1}(\Rt)$ où $\beta_2=\beta_2(\beta_1,\U)>0$.  \\
\end{Theoreme} 
Avant de passer à la démonstration de ce théorème nous allons faire les remarques suivantes: si nous comparons ce théorème avec les résultats obtenus  dans la section précédente (le Théorème \ref{Theo:Ossen_decroissance} et le Théorème \ref{Theo:gevrey-lam}) nous pouvons observer que ce théorème est plus générale dans le sens qu'ici nous n'avons fait aucun contrôle sur la taille force $\fe$ et ce résultat est valable aussi dans le cadre d'un fluide en régime turbulent, tandis que  les résultats de la section précédente (où l'on a contrôlé la taille de la force) sont valables uniquement dans le cadre d'un fluide en régime laminaire.\\
\\ 
En revanche, observons maintenant que dans le régime laminaire on obtient une décroissance fréquentielle exponentielle de la solution plus précise que celle obtenue dans le Théorème \ref{Theo:dec-freq-turb} (dans le cadre turbulent). En effet, toujours dans les Théorèmes \ref{Theo:Ossen_decroissance} et Théorème \ref{Theo:gevrey-lam}  démontrés dans le cadre laminaire nous observons que l'on construit une solution  $\U$ qui vérifient  une décroissance fréquentielle ponctuelle: $\vert \widehat{\U}(\xi)\vert \lesssim e^{-\vert \xi \vert}$, tandis que dans le Théorème \ref{Theo:dec-freq-turb} cette décroissance fréquentielle exponentielle en mesurée en termes de la norme  $L^2$ (voir l'expression (\ref{dec-frec-turb-motiv})) et nous n'avons ici une décroissance ponctuelle.\\
\\
\textbf{Démonstration du Théorème \ref{Theo:dec-freq-turb}.} Pour démontrer ce théorème nous avons besoin de considérer pour l'instant le problème d'évolution des équations de Navier-Stokes : 
\begin{equation}\label{N-S-1}
\partial_t \vu+\P(div(\vu\otimes \vu))-\nu\Delta \vu=\vg,\quad div(\vu)=0, \quad  \vu(0,\cdot)=\vu_0,
\end{equation}  où pour $T_0>0$ on a  $\vg \in \mathcal{C}([0,T_0[,\dot{H}^{1}(\Rt))$  une force  à divergence nulle et $\vu_0 \in \dot{H}^{1}(\Rt)$ est une donnée initiale à divergence nulle. L'existence et l'unicité d'une solution (locale en temps) est un résultat classique qui a été développé par H. Fujita et T. Kato dans \cite{FujitaKato}:
\begin{Lemme}\label{Lemme1} Il existe un temps $0<T_1<T_0$ et il existe une fonction  $\vu \in \mathcal{C}([0,T_1[,\dot{H}^{1}(\Rt))$ qui est l'unique  solution des équations (\ref{N-S-1}). 
\end{Lemme} Pour une preuve ce résultat voir le livre \cite{PGLR1}, Théorème $7.1$. \\
\\ 
Une fois que nous avons la solution $\vu$ ci-dessus la première chose à faire est d'étudier la décroissance fréquentielle de cette solution et pour cela nous introduisons l'espace fonctionnel suivant: soit $\beta>0$ et $T>0$, on défini l'espace
 \begin{equation}\label{E}
E_{\beta,T}=\left\lbrace \vu \in \mathcal{C}(]0,T[,\dot{H}^{1}(\Rt)): e^{\beta\sqrt{t}\sqrt{-\Delta}} \vu \in \mathcal{C}(]0,T[,\dot{H}^{1}(\Rt)) \right\rbrace,
\end{equation} muni de la norme $\Vert \cdot \Vert_{\beta,T}= \Vert e^{\beta\sqrt{t}\sqrt{-\Delta}}(\cdot)\Vert_{L^{\infty}([0,T],\dot{H}^{1}(\Rt))}$, où l'opérateur $e^{\beta\sqrt{t}\sqrt{-\Delta}}$ est donné dans la Définition \ref{Def_noyau_posisson}.  Nous avons ainsi  le résultat suivant.
\begin{Proposition}\label{Prop:gevrey-temps} Dans le cadre du Lemme  \ref{Lemme1}, si la force  $\vg\in \mathcal{C}(]0,T_0[,\dot{H}^{1}(\Rt))$ vérifie  $\vg \in  E_{\beta,T_0}$  alors $\vu\in \mathcal{C}(]0,T_1[,\dot{H}^{1}(\Rt))$ l'unique solution des équations (\ref{N-S-1})   vérifie $\vu \in E_{\beta,T_1}$ avec  $0<T_1<T_0$ suffisamment petit.   
\end{Proposition} 
\pv  Nous allons étudier la quantité 
\begin{eqnarray}\label{EstimationPointFixe1}
\Vert \vu_1 \Vert_{\beta,T_1}&=& \left\Vert h_{\nu t}\ast \vu_0+\int_{0}^{t}h_{\nu(t-s)}\ast \vg(s,\cdot)ds-\int_{0}^{t}h_{\nu(t-s)}\ast \P(div(\vu_1 \otimes \vu_1))(s,\cdot)ds\right\Vert_{\beta,T_1}.
\end{eqnarray} 
Les deux premiers termes  de cette expression sont faciles à estimer et nous avons 
\begin{equation}\label{estim-gevrey1} 
\left\Vert h_{\nu t}\ast \vu_0+\int_{0}^{t}h_{\nu(t-s)}\ast \vg(s,\cdot)ds\right\Vert_{\beta,T_1}\leq c(\nu,\beta,T_0)\left(\Vert \vu_0\Vert_{\dot{H}^{1}_{x}}+\Vert e^{\beta\sqrt{t}\sqrt{-\Delta}}\vg \Vert_{L^{\infty}_{t}\dot{H}^{1}_{x}}\right).
\end{equation}
Pour le dernier terme de  (\ref{EstimationPointFixe1}), par la définition de la norme  $\Vert \cdot \Vert_{\beta,T_1}$, 
l'identité de Plancherel et la continuité du Projecteur de Leray nous avons 
\begin{equation*}
\begin{split}
(I)=\left\Vert \int_{0}^{t}h_{\nu(t-s)}\ast \P(div(\vu_1 \otimes \vu_1))ds\right\Vert_{\beta,T_1} = \sup_{0<t<T_1} \left\Vert e^{\beta\sqrt{t}\sqrt{-\Delta}}\left(\int_{0}^{t}h_{\nu(t-s)}\ast \P(div(\vu_1 \otimes \vu_1)) ds \right)\right\Vert_{\dot{H}^{1}_{x}}\\
\leq \sup_{0<t<T_1} c \left\Vert \vert\xi \vert^2 \int_{0}^{t}e^{-\nu(t-s)\vert \xi \vert^2} e^{\beta\sqrt{t}\vert\xi \vert}\left\vert\left(\mathcal{F}[\vu_1]\ast \mathcal{F}[\vu_1]\right)(s,\cdot) \right\vert ds\right\Vert_{L^{2}_{x}},
\end{split}
\end{equation*}
et comme nous avons l'estimation ponctuelle
\begin{equation}\label{eq10}
e^{\beta\sqrt{t}\vert\xi \vert}\left\vert\left(\mathcal{F}[\vu_1]\ast \mathcal{F}[\vu_1]\right)(s,\xi)\right\vert \leq \left[\left( e^{\beta\sqrt{t}\vert\xi \vert}\vert \mathcal{F}[\vu_1]\vert\right)\ast \left( e^{\beta\sqrt{t}\vert\xi \vert}\vert \mathcal{F}[\vu_1]\vert\right)\right](s,\xi),
\end{equation}
dû au fait que  $e^{\beta\sqrt{t}\vert \xi \vert}\leq e^{\beta\sqrt{t}\vert \vert\xi-\eta\vert}e^{\beta\sqrt{t}\vert \eta \vert}$ pour tout $\xi,\,\eta\in \mathbb{R}^3$, alors nous obtenons
\begin{eqnarray*}
(I)&\leq & \sup_{0<t<T_1} c \int_{0}^{t}\left\Vert \vert\xi \vert^\frac{3}{2} e^{-\nu(t-s)\vert \xi \vert^2} \vert\xi \vert^{\frac{1}{2}}\left\vert\left[\left( e^{\beta\sqrt{t}\vert\xi \vert}\vert \mathcal{F}[\vu_1]\vert\right)\ast \left( e^{\beta\sqrt{t}\vert\xi \vert}\vert \mathcal{F}[\vu_1]\vert\right)\right]\right\vert \right\Vert_{L^{2}_{x}}ds.
\end{eqnarray*}
En revenant maintenant à la variable spatiale nous pouvons écrire
$$(I)\leq \sup_{0<t<T_1} c \int_{0}^{t} \left\Vert (-\Delta)^{\frac{3}{4}}h_{\nu(t-s)}\ast  (-\Delta)^{\frac{1}{4}} \left\lbrace\left(\mathcal{F}^{-1}\left[e^{\beta\sqrt{t}\vert\xi \vert}\vert \mathcal{F}[\vu_1]\vert \right]\right) \otimes  \left(\mathcal{F}^{-1}\left[e^{\beta\sqrt{t}\vert\xi \vert}\vert \mathcal{F}[\vu_1]\vert\right]\right)\right\rbrace\right\Vert_{L^{2}_{x}}ds,$$
\begin{eqnarray}
&\leq & \left(c  \int_{0}^{T_1} \left\Vert (-\Delta)^{\frac{3}{4}}h_{\nu(t-s)} \right\Vert_{L^1}ds\right) \left\Vert  \left(\mathcal{F}^{-1}\left[e^{\beta\sqrt{t}\vert\xi \vert}\vert \mathcal{F}[\vu_1]\vert \right]\right) \otimes  \left(\mathcal{F}^{-1}\left[e^{\beta\sqrt{t}\vert\xi \vert}\vert \mathcal{F}[\vu_1]\vert\right]\right) \right\Vert_{L^{\infty}_{t}\dot{H}^{\frac{1}{2}}_{x}}\nonumber\\
&\leq& c\frac{T^{\frac{1}{4}}}{\nu^{\frac{3}{4}}} \left\Vert \mathcal{F}^{-1}\left[e^{\beta\sqrt{t}\vert\xi \vert}\vert \mathcal{F}[\vu_1]\vert \right]\right\Vert_{L^{\infty}_{t}\dot{H}^{1}_{x}} \left\Vert \mathcal{F}^{-1}\left[e^{\beta\sqrt{t}\vert\xi \vert}\vert \mathcal{F}[\vu_1]\vert \right]\right\Vert_{L^{\infty}_{t}\dot{H}^{1}_{x}}\nonumber \\ 
&=&c\frac{T^{\frac{1}{4}}_{1}}{\nu^{\frac{3}{4}}}\Vert \vu_1 \Vert_{\beta,T_1} \Vert \vu_1 \Vert_{\beta,T_1}.\label{estim-gevrey2}
\end{eqnarray} 
Une fois que nous disposons des estimations (\ref{estim-gevrey1}) et (\ref{estim-gevrey2}) on fixe maintenant  $T_1$ suffisamment petit de sorte qu'en appliquant le principe de contraction de Picard nous obtenons une solution $\vu_1\in E_{\beta,T_1}$  des équations (\ref{N-S-1}). Mais, étant donné que  $E_{\beta,T_1}\subset \mathcal{C}(]0,T_1[,\dot{H}^{1}(\mathbb{R}^3))$ on a  $\vu_1\in \mathcal{C}(]0,T_1[,\dot{H}^{1}(\mathbb{R}^3))$ et par l'unicité de la solution  $\vu$ on a  $\vu_1=\vu$ et donc  $\vu \in E_{\beta,T_1}$. \finpv\\
Revenons maintenant aux équations de Navier-Stokes stationnaires (\ref{N-S-stationnaire}) pour monter que toute solution $\U \in \dot{H}^{1}(\Rt)$ vérifie $\U \in G^{\beta_2}_{1}(\Rt)$ avec une constante $\beta_2>0$ que l'on fixera plus tard.\\
\\
 En effet, dans le système d'évolution (\ref{N-S-1}) nous prenons la donnée initiale $\vu_0=\U$  et la force $\vg = e^{-\beta\sqrt{t}\sqrt{-\Delta}}(e^{\beta \sqrt{t}\sqrt{-\Delta}}\fe)$, avec $t\in [0,1]$; et nous allons maintenant vérifier que cette force satisfait les hypothèses du Lemme \ref{Lemme1} et de la Proposition \ref{Prop:gevrey-temps}, c'est à dire, nous devons vérifier que $\vg \in \mathcal{C}(]0,1[,\dot{H}^{1}(\Rt))$ et $e^{\beta\sqrt{t}\sqrt{-\Delta}}\vg \in \mathcal{C}(]0,1[,\dot{H}^{1}(\Rt))$; et pour cela nous allons tout d'abord vérifier que l'on a   $e^{\beta\sqrt{t}\sqrt{-\Delta}}\fe\in \mathcal{C}(]0,1[,\dot{H}^{1}(\Rt))$. On commence donc par écrire
\begin{eqnarray*}
\left\Vert e^{\beta\sqrt{t}\sqrt{-\Delta}}\fe \right\Vert^{2}_{L^{\infty}_{t}\dot{H}^{1}_{x}} &=&\sup_{0<t<1} \int_{\Rt} \vert \xi \vert^2e^{2\beta\sqrt{t} \vert \xi \vert}\left\vert \widehat{\fe}(\xi)\right\vert^2d\xi \leq   \int_{\Rt} \vert \xi \vert^2e^{2\beta \vert \xi \vert}\left\vert \widehat{\fe}(\xi)\right\vert^2d\xi\\
&=& \int_{\Rt} \vert \xi \vert^4e^{2\beta \vert \xi \vert}\left\vert \widehat{\fe}(\xi)\right\vert^2 \frac{d\xi}{\vert \xi \vert^2} = \frac{1}{\beta^4}\int_{\Rt} \left(\beta\vert \xi \vert\right)^4 e^{2\beta\vert \xi \vert}\left\vert \widehat{\fe}(\xi)\right\vert^2 \frac{d\xi}{\vert \xi \vert^2}\\
&\leq & \frac{1}{\beta^4} \int_{\Rt} e^{\beta\vert \xi \vert +2\beta \vert \xi \vert}\left\vert \widehat{\fe}(\xi)\right\vert^2\frac{d\xi}{\vert \xi \vert^2}= \frac{1}{\beta^4} \int_{\Rt} e^{3\beta \vert \xi \vert} \vert \widehat{\fe}(\xi)\vert^2 \frac{d\xi}{\vert \xi \vert^2}. 
\end{eqnarray*}
Comme $\fe \in G^{\beta_1}_{-1}(\Rt)$ par la définition de cet espace fonctionnel nous avons $e^{\beta_1 \sqrt{-\Delta}}\fe \in \dot{H}^{-1}(\Rt)$ et ainsi, en fixant le paramètre  $\beta>0$ tel que $3\beta =2 \beta_1$ nous pouvons alors  écrire 
$$ \frac{1}{\beta^4} \int_{\Rt} e^{3\beta \vert \xi \vert} \vert \widehat{\fe}(\xi)\vert^2 \frac{d\xi}{\vert \xi \vert^2} = \left(\frac{3}{2 \beta_1}\right)^2 \int_{\Rt} e^{2\beta_1\vert \xi \vert}\left\vert \widehat{\fe}(\xi)\right\vert^2\frac{d\xi}{\vert \xi \vert^2}=C_{\beta_1}\Vert e^{\beta_1 \sqrt{-\Delta}}\fe \Vert_{\dot{H}^{-1}} <+\infty,$$ 
d'où nous avons  $e^{\beta\sqrt{t}\sqrt{-\Delta}}\fe\in \mathcal{C}(]0,1[,\dot{H}^{1}(\Rt))$ et alors $\vg \in \mathcal{C}(]0,1[,\dot{H}^{1}(\Rt))$ et $e^{\beta\sqrt{t}\sqrt{-\Delta}}\vg \in \mathcal{C}(]0,1[,\dot{H}^{1}(\Rt))$.\\
\\
Ainsi, par le Lemme \ref{Lemme1} il existe un temps $0<T_1<1$ et une unique solution $\vu \in \mathcal{C}(]0,T_1[,\dot{H}^{1}(\Rt))$ des équations  (\ref{N-S-1}). De plus, comme $e^{\beta\sqrt{t}\sqrt{-\Delta}}\fe\in \mathcal{C}(]0,1[,\dot{H}^{1}(\Rt))$ alors par la Proposition \ref{Prop:gevrey-temps} nous avons  $ e^{\beta\sqrt{t}\sqrt{-\Delta}}\vu \in \mathcal{C}(]0,T_1[,\dot{H}^{1}(\Rt))$. \\
\\ 
D'autre part, comme $\U\in \dot{H}^{1}(\Rt)$  est  une fonction constante en variable de temps nous avons  $\U\in \mathcal{C}(]0,T_1[,\dot{H}^{1}(\Rt))$ et de plus,  comme $\partial_t\U\equiv 0$ et $\vg =e^{-\beta\sqrt{t}\sqrt{-\Delta}}(e^{\beta\sqrt{t}\sqrt{-\Delta}}\fe) =\fe$ nous pouvons alors observer que  $\U \in \mathcal{C}(]0,T_1[,\dot{H}^{1}(\Rt))$ est aussi une solution de l'équation (\ref{N-S-1}) ci-dessus et donc, par l'unicité de la solution $\vu \in \mathcal{C}(]0,T_1[\dot{H}^{1}(\Rt))$ nous obtenons que $\U=\vu$. Ensuite, comme  $e^{\beta\sqrt{t}\sqrt{-\Delta}}\vu \in \mathcal{C}(]0,T_1[,\dot{H}^{1}(\Rt))$ alors nous avons   $e^{\beta\sqrt{t}\sqrt{-\Delta}}\U \in \mathcal{C}(]0,T_1[,\dot{H}^{1}(\Rt))$ pour tout $t\in [0,T_1[$. Ainsi, on fixe $\beta_2=\beta \sqrt{\frac{T_1}{2}}>0$ et alors on a  
\begin{equation*}
 \int_{\mathbb{R}^3}\vert \xi \vert^2 e^{2\beta_2\vert \xi \vert} \vert \U(\xi)\vert^2  d\xi =\big\| e^{\beta\sqrt{\frac{T_1}{2}}\sqrt{-\Delta}}\U\big\|^{2}_{\dot{H}^{1}_{x}}\leq \sup_{0<t<T_1}\Vert e^{\beta\sqrt{t}\sqrt{-\Delta}}\U\Vert^{2}_{\dot{H}^{1}_{x}}<+\infty,
\end{equation*}  et donc $\U\in G^{\beta_2}_{1}(\Rt)$. Le Théorème \ref{Theo:dec-freq-turb} est maintenant démontré.  \finpv

	\chapter{Des théorèmes de type Liouville  }
Dans le chapitre précédent  nous avons étudié la décroissance fréquentielle  des solutions des équations de  Navier-Stokes stationnaires:   
\begin{equation}\label{N-S-intro}
-\nu \Delta \U +(\U \cdot \vec{\nabla})\U+\vec{\nabla}P=\fe, \quad div(\U)=0,
\end{equation}  mais, avant de faire cet étude,  nous avons tout d'abord énoncé un résultat  sur l'existence des solutions $\U$ dans le cadre de l'espace de Sobolev $\dot{H}^{1}(\Rt)$. Plus précisément, par le Théorème \ref{Theo:exstence-sol-stationnaires} page \pageref{Theo:exstence-sol-stationnaires} nous observons que si l'on considère une force  $\fe \in \dot{H}^{-1}(\Rt)$ alors il existe $\U \in \dot{H}^{1}(\Rt)$ une solution des équations (\ref{N-S-intro}) qui vérifie l'estimation
\begin{equation}\label{estim-a-priori}
\nu \Vert \U \Vert_{\dot{H}^{1}}\leq \Vert \fe \Vert_{\dot{H}^{-1}}.
\end{equation} 
Dans ce chapitre nous allons étudier un tout autre problème relié à l'estimation (\ref{estim-a-priori}):  dans cette estimation nous observons que si nous prenons maintenant une force nulle ($\fe =0$) alors la solution des équations 
\begin{equation}\label{N-S-f-nulle-intro}
-\nu \Delta \U +(\U \cdot \vec{\nabla})\U+\vec{\nabla}P=0, \quad div(\U)=0,
\end{equation} obtenue par le  Théorème \ref{Theo:exstence-sol-stationnaires} est forcément la solution triviale $\U=0$; et nous nous posons alors la question de savoir si cette solution $\U=0$ est l'unique solution des équations (\ref{N-S-f-nulle-intro}) dans l'espace $\dot{H}^{1}(\Rt)$. Néanmoins, comme nous l'expliquerons plus en détail  ci-après, cette question est encore une question ouverte
et le problème repose essentiellement sur le fait que, dans l'état actuel de nos connaissances,  on ne sait  pas montrer que toute solution $\U \in \dot{H}^1(\Rt)$ des équations (\ref{N-S-f-nulle-intro})  vérifie l'estimation  (\ref{estim-a-priori}). \\
\\
Le but de ce chapitre est alors d'étudier d'autres espaces fonctionnels $E$  tels que si $\U$ est une solution  des équations de Navier-Stokes stationnaires (\ref{N-S-f-nulle-intro}) et si $\U \in E$ alors on  ait l'identité $\U=0$.  Autrement dit, il s'agit de trouver des espaces fonctionnels $E$ dans lesquelles les équations (\ref{N-S-f-nulle-intro}) possèdent  comme unique solution la solution triviale $\U=0$. Ce problème est également appelé le problème de Liouville pour les équations de Navier-Stokes stationnaires. \\
\\
Dans la Section \ref{sec:intro} ci-dessous nous expliquons plus en détail les enjeux du problème de Liouville pour les équations de Navier-Stokes stationnaires (\ref{N-S-f-nulle-intro}) et nous exposerons quelques résultats connus sur ce problème.  Ensuite, dans la Section \ref{sec:resultats} nous démontrerons  quelques résultats nouveaux que nous avons obtenu  sur  ce sujet.

\section{Introduction}\label{sec:intro} 
Comme annoncé nous allons rappeler ici quelques résultats connus sur le problème de Liouville pour les équations de Navier-Stokes stationnaires (\ref{N-S-f-nulle-intro}) et pour cela nous allons commencer par expliquer un peu plus en détail ce problème.
\subsection{Le problème de Liouville pour les équations de Navier-Stokes stationnaires}\label{sec:intro-Liouville}

Dans cette section nous  introduisons  le problème de Liouville pour les équations de Navier-Stokes stationnaires (\ref{N-S-f-nulle-intro})  que nous allons  étudier tout au long de ce chapitre mais, avant d'entrer dans les détails de ce problème,  il convient tout d'abord de préciser un peu plus notre cadre de travail et on commence donc  par rappeler la définition de solution faible des  ces équations.  

\begin{Definition}[Solution faible]\label{def:Sol-faible}  Une fonction $\U:\Rt \longrightarrow \Rt$ est une solution faible des équations de Navier-Stokes stationnaires (\ref{N-S-f-nulle-intro}) si $\U \in L^{2}_{loc}(\Rt)$ et s'il existe $P\in \mathcal{D}^{'}(\Rt)$ telle que le couple  $(\U,P)$ vérifie ces équations au sens des distributions.
	\end{Definition}
Observons que dans cette définition nous avons $\U \in L^{2}_{loc}(\Rt)$ et cette propriété de la solution nous permet d'assurer que le terme non linéaire des équations (\ref{N-S-f-nulle-intro}) est bien défini dans $\mathcal{D}^{'}(\Rt)$. En effet, comme $\U \in L^{2}_{loc}(\Rt)$ et comme $div(\U)=0$ alors nous pouvons écrire l'identité (toujours au sens des distributions) $(\U \cdot \vec{\nabla})\U=div(\U \otimes \U)$, où  nous observons que si la solution $\U$ appartient à l'espace $L^{2}_{loc}(\Rt)$ alors par l'inégalité de Cauchy-Schwarz nous avons que  le terne non linéaire  $\U \otimes \U$ appartient à l'espace $L^{1}_{loc}(\Rt)$ et donc le terme $div(\U \otimes \U)$ est bien défini dans $\mathcal{D}^{'}(\Rt)$. Ainsi, par l'identité précédente nous avons que le terme $(\U \cdot \vec{\nabla})\U $ est ainsi bien défini au sens des distributions. \\
\\
Maintenant que l'on a précisé notre cadre de travail qui est donné par l'espace $L^{2}_{loc}(\Rt)$, nous cherchons à trouver des espaces fonctionnels $E\subset L^{2}_{loc}(\Rt)$ où l'on puisse résoudre  le problème  de Liouville suivant:  
\vspace{1mm}
\begin{center}
	\fbox{
		\begin{minipage}[l]{157mm}
			Si  $\U \in L^{2}_{loc}(\Rt)$ est une solution des équations de Navier-Stokes stationnaires (\ref{N-S-f-nulle-intro}) et si de plus cette solution vérifie  $\U \in E$ alors on a $\U=0$.
				\end{minipage}} 
\end{center}  
\vspace{5mm} Avant d'expliquer les espaces $E$ que nous allons considérer il est important de remarquer tout d'abord que la condition supplémentaire $\U \in E$ est  vraiment nécessaire pour résoudre ce problème de Liouville.\\
\\
En effet, dans l'exemple ci-dessous nous allons observer que dans le cadre général de l'espace $L^{2}_{loc}(\Rt)$ on ne peut pas s'attendre à ce que la solution triviale $\U=0$ soit l'unique solution des équations (\ref{N-S-f-nulle-intro}). Considérons  la fonction $\psi:\Rt \longrightarrow \mathbb{R}$ définie comme suit: pour $x=(x_1,x_2,x_3)\in \Rt$, 
\begin{equation}\label{def-psi}
\psi(x_1,x_2,x_3)=\frac{x^{2}_{1}}{2}+\frac{x^{2}_{2}}{2}-x^{2}_{3},
\end{equation} et à  partir de cette fonction nous définissons maintenant les fonctions $\U$ et $P$ par les identités
\begin{equation}\label{def-U-ex}
U(x_1,x_2,x_3)=\vec{\nabla}\psi(x_1,x_2,x_3),
\end{equation} et 
\begin{equation}\label{def-P-ex}
P(x_1,x_2,x_3)=-\frac{1}{2}\vert \U(x_1,x_2,x_3)\vert^2.
\end{equation} 
\`A ce stade, observons que la fonction $\U$ appartient à l'espace $L^{2}_{loc}(\Rt)$. En effet,  par la définition de la fonction $\psi$ ci-dessus et par l'identité (\ref{def-U-ex}) nous pouvons écrire
\begin{equation}\label{ident-U-ex}
U(x_1,x_2,x_3)=\vec{\nabla}\psi(x_1,x_2,x_3)=(x_1,x_2,-2x_3),
\end{equation}
 pour observer que l'on a l'estimation $\vert \U(x)\vert \approx \vert x \vert$ d'où nous pouvons en tirer $\U\in L^{2}_{loc}(\Rt)$; et nous avons en plus le résultat suivant:

\begin{Lemme}\label{lemme-exemple} Le couple $(\U,P)$ défini pas les expressions (\ref{def-U-ex}) et (\ref{def-P-ex}) respectivement vérifie les équations de Navier-Stokes stationnaires (\ref{N-S-f-nulle-intro}). \\
\end{Lemme} 
\pv  Observons tout d'abord que par l'identité (\ref{ident-U-ex}) nous avons $\U \in \mathcal{C}^{\infty}(\Rt)$ et alors par l'identité (\ref{def-P-ex}) nous obtenons $P\in \mathcal{C}^{\infty}(\Rt)$. Nous allons maintenant montrer que le couple $(\U,P)$  vérifie les équations  (\ref{N-S-f-nulle-intro}). \\
\\
Par l'identité (\ref{def-psi}) nous avons $\Delta \psi =0$ et ainsi nous pouvons écrire 
$$\nu \Delta \U =\nu \Delta (\vec{\nabla} \psi)=\nu  \vec{\nabla} (\Delta \psi)=0.$$
Ensuite, par le calcul vectoriel nous avons l'identité 
$$ (\U \cdot \vec{\nabla})\U=\frac{1}{2}\vert U \vert^2 + (\vec{\nabla}\wedge \U) \wedge \U,$$ où, comme $\U=\vec{\nabla}\psi$ alors nous avons $ \vec{\nabla}\wedge \U= \vec{\nabla}\wedge (\vec{\nabla}\psi)=0$, et donc nous pouvons écrire 
$$ (\U \cdot \vec{\nabla})\U=  \frac{1}{2}\vert U \vert^2.$$
Avec ces identités et comme $P=-\frac{1}{2}\vert U \vert^2$ nous obtenons l'identité 
$$ -\nu \Delta \U +(\U \cdot \vec{\nabla})\U+\vec{\nabla}P=\frac{1}{2}\vert U \vert^2-\frac{1}{2}\vert U \vert^2=0.$$
De plus, toujours  par l'identité (\ref{ident-U-ex}) nous pouvons en tirer  $div(\U=0)$. \finpv 
\\ 
Nous observons de cette façon que dans le cadre  de l'espace $L^{2}_{loc}(\Rt)$ les équations de Navier-Stokes stationnaires  (\ref{N-S-f-nulle-intro}) ont au moins deux solutions: la solution triviale $\U=0$ et la solution non nulle $\U$ donnée par l'expression (\ref{def-U-ex}), et alors on a besoin de chercher des espaces fonctionnels adéquats  $E\subset L^{2}_{loc}(\Rt)$  où l'on puisse résoudre le problème de Liouville ci-dessus. \\ 
\\
 Nous avons mentionné précédemment que l'espace $\dot{H}^{1}(\Rt)$ semblerait être un espace naturel pour résoudre ce problème de Liouville  et ce fait repose essentiellement sur quelques estimations \emph{formelles} que nous avons fait dans le chapitre précédent (voir l'estimation  (\ref{estimation-a-priori-sol-estationnaires}) page \pageref{estimation-a-priori-sol-estationnaires}) et que nous avons besoin de rappeler  rapidement. On commence par rappeler le résultat suivant:
 \begin{Lemme}\label{Lemme-reg-sol} Soit $(\U,P) \in L^{2}_{loc}(\Rt) \times \mathcal{D}^{'}(\Rt)$ une solution des équations de Navier-Stokes stationnaires (\ref{N-S-f-nulle-intro}). Si $\U \in L^{3}_{loc}(\Rt)$ alors  $\U \in \mathcal{C}^{\infty}(\Rt)$ et $P\in \mathcal{C}^{\infty}(\Rt)$. 
 \end{Lemme}	
 Pour une preuve de ce résultat voir le Théorème $X.1.1.$ du livre \cite{Galdi} page $658$.\\
 \\
De cette façon, pour $\U \in \dot{H}^{1}(\Rt)$ une solution des équations (\ref{N-S-f-nulle-intro}), par les inégalités de Hardy-Littlewood-Sobolev nous avons $\U \in L^6(\Rt)$ et alors $\U \in L^{3}_{loc}(\Rt)$, d'où, par le Lemme \ref{Lemme-reg-sol}  nous obtenons $\U \in \mathcal{C}^{\infty}(\Rt)$ et $P\in \mathcal{C}^{\infty}(\Rt)$ et nous pouvons multiplier chaque terme de l'équation (\ref{N-S-f-nulle-intro}) par la solution $\U$. Ensuite,  nous intégrons en variable d'espace sur $\Rt$ pour obtenir l'identité:
\begin{equation}\label{ipp-estm-a-priori}
-\nu \int_{\Rt} \Delta \U \cdot \U dx+ \int_{\Rt} \left[ (\U\cdot \vec{\nabla}) \U \right] \cdot \U dx+ \int_{\Rt} \vec{\nabla} P \cdot \U dx=0, 
\end{equation} et si nous supposons pour l'instant que la solution $\U$ est  \emph{suffisamment  intégrable} alors, pour le deuxième terme et le troisième terme à gauche de l'identité ci-dessus, comme $div(\U)=0$ et  par une intégration par parties nous pouvons formellement écrire 
\begin{equation*}
\int_{\Rt} [(\U \cdot \vec{\nabla})\U]\cdot \U dx=0,
\end{equation*}
et 
\begin{equation*}
\int_{\Rt} \vec{\nabla}P \cdot \U dx=0.
\end{equation*}
Ainsi,  par l'identité (\ref{ipp-estm-a-priori}) nous pouvons en tirer l'identité  suivante:
\begin{equation}\label{Liouville-H1}
\ds{\int_{\Rt}\vert \vec{\nabla}\otimes \U \vert^2dx=0}, 
\end{equation} et par cette identité on s'attend à ce que l'unique solution $\U \in \dot{H}^{1}(\Rt)$ des équations (\ref{N-S-f-nulle-intro}) soit la solution triviale $\U=0$. \\
\\
Néanmoins, pour  $\U\in \dot{H}^{1}(\Rt)$  l'identité (\ref{ipp-estm-a-priori}) est purement formelle et l'on ne dispose d'aucune information supplémentaire pour assurer que les termes
\begin{equation}\label{int}
\int_{\Rt} [(\U \cdot \vec{\nabla})\U]\cdot \U dx \quad \text{et}\quad \int_{\Rt} \vec{\nabla}P \cdot \U dx,
\end{equation} sont bien définies. En effet, si $\U \in \dot{H}^{1}(\Rt)$ alors  nous avons  $\U \in L^6(\Rt)$ ce qui exprime une certaine décroissance à l'infini de cette fonction, mais pour donner un sens rigoureux aux quantités donnés dans (\ref{int})  nous avons besoin que la solution $\U$ décroisse plus rapidement à l'infini: si nous supposons pour l'instant que la solution vérifie en plus $\U \in L^4(\Rt)$ alors nous pouvons vérifier le résultat suivant 
\begin{Lemme}\label{lemme-tech} Si $\U \in L^4(\Rt)$ alors on a $(\U\cdot \vec{\nabla})\cdot \U \in \dot{H}^{-1}(\Rt)$ et $\vec{\nabla}P \in \dot{H}^{-1}(\Rt)$.
\end{Lemme} Ce résultat est purement technique et sera prouvé  à la fin du chapitre page \pageref{preuve:lemme-tech}. Nous pouvons ainsi observer que si la solution $\U \in \dot{H}^1(\Rt)$ vérifie une décroissance  plus rapide  à l'infini que celle donnée par la propriété $\U \in L^6(\Rt)$ et qui est maintenant mesurée en termes de la norme de l'espace $L^4(\Rt)$, alors  les quantités données dans la formule  (\ref{int}) sont  bien définies  et donc cette fois-ci nous pouvons écrire  l'identité (\ref{Liouville-H1})  qui finalement entraîne l'identité cherchée  $\U=0$. Mais,  à partir de l'information $\U \in \dot{H}^{1}(\Rt)$  malheureusement  nous n'avons aucun argument supplémentaire pour en déduire que $\U$ appartient à l'espace $L^4(\Rt)$. \\
\\
Nous observons ainsi que  le problème de Liouville dans l'espace $\dot{H}^{1}(\Rt)$ s'agit  d'un problème délicat et ce fait   repose essentiellement sur  la décroissance à l'infini de la solution $\U$. En revanche, nous allons observer que  les espaces de Lebesgue  nous fournissent un cadre convenable pour étudier ce problème. \\
\\
En effet, grâce au Lemme \ref{lemme-tech} nous pouvons montrer  un premier résultat sur le problème de Liouville pour les équations de Navier-Stokes stationnaires:\\
\begin{Proposition}\label{Prop:L4} Soit $\U\in L^{2}_{loc}(\Rt)$ une solution des équations de Navier-Stokes stationnaires  (\ref{N-S-f-nulle-intro}). Si $\U\in L^4(\Rt)$ alors $\U=0$.
	\end{Proposition}
\pv Par le Lemme \ref{lemme-tech} nous savons que  $(\U\cdot \vec{\nabla})\cdot \U \in \dot{H}^{-1}(\Rt)$ et $\vec{\nabla}P \in \dot{H}^{-1}(\Rt)$ et comme $\U$ vérifie les équations (\ref{N-S-f-nulle-intro}) nous écrivons $\nu \Delta \U= - (\U\cdot \vec{\nabla})\cdot \U-\vec{\nabla}P \in \dot{H}^{-1}(\Rt)$ pour obtenir $\U \in \dot{H}^{1}(\Rt)$. Comme nous avons l'information $div(\U)=0$ nous disposons alors des identités $\ds{\int_{\Rt} [(\U \cdot \vec{\nabla})\U]\cdot \U dx=0}$ et $\ds{\int_{\Rt} \vec{\nabla}P \cdot \U dx=0}$ qui nous permettent d'écrire l'identité (\ref{Liouville-H1}) et d'en déduire  $\U=0$. \finpv
\\
Motivés par ce premier résultat, nous allons tout d'abord étudier le problème de Liouville pour les équations (\ref{N-S-f-nulle-intro})  dans le cadre des espaces de Lebesgue $L^p(\Rt)$. Plus précisément, il s'agit d'étudier  d'autres valeurs possibles du paramètre d'intégration $p \in [1,+\infty]$ avec lesquels on puisse obtenir des résultats similaires à celui obtenu dans la Proposition \ref{Prop:L4}. Ensuite, nous allons voir qu'il est aussi possible d'étudier ce  problème de Liouville dans des espaces plus générales que les espaces de Lebesgue et dans ce cadre  les espaces de Morrey (qui seront définis dans la section qui suit) nous permettent  d'obtenir des résultats intéressants sur ce problème. \\ 
\\
Mais, pour expliquer  les résultats que nous avons obtenu  nous avons besoin de faire tout d'abord  un très court  rappel des résultats démontrés précédemment.  
\subsection{Des résultats connus}\label{sec:resultats-connus}
Dans cette section nous exposons quelques résultats connus sur le problème de Liouville pour les équations (\ref{N-S-f-nulle-intro})  à partir desquels nous allons établir des résultats nouveaux dans la Section  \ref{sec:resultats} ci-après. \\
\\
Comme annoncé, nous nous intéressons à étudier ce problème de Liouville dans le cadre des espaces de Lebesgue et de Morrey et on commence donc par étudier ce qui à été démontré dans le cadre des espaces de Lebesgue. 

\subsubsection{A) Quelques résultats dans les espaces de Lebesgue}
  Un des résultats le plus connu est celui  obtenue par G. Galdi en $1994$ dans son livre \cite{Galdi}:
\begin{Theoreme}[\cite{Galdi}, Théorème $X.9.5$]\label{Theo:Galdi} Soit $\U \in L^{2}_{loc}(\Rt)$ une solution  des équations de Navier-Stokes stationnaires (\ref{N-S-f-nulle-intro}). Si $\U \in L^{\frac{9}{2}}(\Rt)$ alors $\U=0$.
	\end{Theoreme} 
Si nous comparons ce résultat avec le résultat obtenu dans la Proposition \ref{Prop:L4} nous pouvons observer que  dans le Théorème \ref{Theo:Galdi}  nous prouvons l'unicité de la solution  $\U=0$ des équations (\ref{N-S-f-nulle-intro}) avec une condition de décroissance à l'infini moins forte que celle de la Proposition \ref{Prop:L4}. En effet, dans le Théorème \ref{Theo:Galdi} nous observons que la solution $\U$ appartient à l'espace $L^{\frac{9}{2}}(\Rt)$ tandis que dans la Proposition \ref{Prop:L4} la solution appartient à l'espace $L^4(\Rt)$ et alors dans ce cas cette solution décroît plus rapidement à l'infini. \\
\\
Expliquons maintenant les grandes lignes de la preuve du Théorème \ref{Theo:Galdi}. Ce résultat repose essentiellement sur le fait que si la solution $\U$ appartient à l'espace  $L^{\frac{9}{2}}(\Rt)$ alors nous pouvons vérifier rigoureusement l'estimation (\ref{Liouville-H1}) ce qui entraîne tout d'abord  que cette solution $\U$ appartienne en plus à l'espace $\dot{H}^{1}(\Rt)$ et ensuite que l'on ait l'identité  $\U=0$. Pour vérifier l'estimation (\ref{Liouville-H1}) on montre alors l'estimation technique suivante (voir toujours le Théorème $X.9.5$,   page $729$ du livre \cite{Galdi} pour tous les détails des calculs): pour $R>1$ soit la boule $B_{\frac{R}{2}}=\{ x \in \Rt:  \vert x \vert <\frac{R}{2}\}$, et l'on a  
\begin{equation}\label{estim-galdi}
\int_{B_{\frac{R}{2}}} \vert \vec{\nabla}\otimes \U \vert^2 dx\leq c_1 \left(  \int_{R\leq \vert x \vert \leq 2R}\vert \U(x)\vert^{\frac{9}{2}}dx \right)^{\frac{2}{9}},
\end{equation} où $c_1=c_1(\U,P)>0$ est une constante qui dépend du champ de vitesse $\U$ et de la pression $P$ mais qui ne dépend pas de $R>1$. Ainsi, comme $\U\in L^{\frac{9}{2}}(\Rt)$ nous avons  
\begin{equation}\label{dec-galdi}
\lim_{R\longrightarrow +\infty} \left(  \int_{R\leq \vert x \vert \leq 2R}\vert \U(x)\vert^{\frac{9}{2}}dx \right)^{\frac{2}{9}}=0,
\end{equation}  et donc en revenant à l'estimation (\ref{estim-galdi}) nous obtenons  l'identité  $ \ds{\int_{\Rt}  \vert \vec{\nabla}\otimes \U \vert^2 dx=0}$, qui est l'identité cherchée  (\ref{Liouville-H1}). \\
\\
Il est ainsi  intéressant  d'observer que la décroissance à l'infini de la solution $\U\in L^{\frac{9}{2}}(\Rt)$ exprimée de façon précise par l'identité (\ref{dec-galdi}) entraîne l'identité (\ref{Liouville-H1}) qui, comme mentionné précédemment,  est hors de portée lorsqu'on considère en toute généralité une solution $\U\in \dot{H}^{1}(\Rt)$. \\
\\
\'Etudions maintenant un  deuxième résultat intéressant sur le problème de Liouville dans les espaces de Lebesgue. Nous savons que l'espace $\dot{H}^{1}(\Rt)$ est inclut dans l'espace $L^6(\Rt)$ (par les inégalités de Hardy-Littlewood-Sobolev) et alors nous nous posons la question de savoir si l'on peut résoudre le problème de Liouville pour les équations  (\ref{N-S-f-nulle-intro}) dans cet espace de Lebesgue. Néanmoins, avec les outils que nous avons à notre disposition nous ne sommes pas en mesure de donner une réponse à cette question. En effet, dans la section précédente nous avons déjà souligné le fait que l'information $\U\in L^6(\Rt)$ ne suffit pas pour obtenir rigoureusement l'identité $\U=0$. \\
\\
Dans ce cadre,  le résultat   obtenu par G. Seregin en $2015$ dans son article \cite{Seregin1}  montre que  si l'on ajoute une hypothèse supplémentaire sur la solution: $\U \in L^6\cap BMO^{-1}(\Rt)$, alors on obtient l'identité cherchée $\U=0$.  \\
\\
Pour  définir l'espace $BMO^{-1}(\Rt)$ nous avons besoin  d'introduire rapidement l'espace  $BMO(\Rt)$ (voir le Chapitre $7$ du livre \cite{Grafakos2}). Soit donc $f\in L^{1}_{loc}(\Rt)$ et la boule $B(x_0,R)$ dans $\Rt$. On définit tout d'abord la quantité $m_{B(x_0,R}[f] \in \mathbb{R}$ comme la moyenne de la fonction $f$ sur la boule $B(x_0,R)$: $$m_{B(x_0,R)} [f] = \frac{1}{R^3} \int_{B(x_0,R)} f(x)dx,$$ et 
ensuite on définit l'oscillation moyenne de cette fonction $f$ par la quantité 
$$  \Vert f \Vert_{BMO} = \sup_{x_0\in \Rt, R>0} \frac{1}{R^3} \int_{B(x_0,r_0)} \left\vert f(x) - m_{B(x_0,R)}[f] \right\vert dx. $$ 
Ainsi, l'espace $BMO(\Rt)$ est défini de la façon suivante:
$$ BMO(\Rt)= \{ f\in L^{1}_{loc}(\Rt): \Vert f \Vert_{BMO}<+\infty \},$$ c'est à dire, l'espace $BMO(\Rt)$ est l'espace des fonctions localement intégrables  à oscillations moyennes bornées. \\
\\
Ensuite, en suivant \cite{Seregin1}, nous dirons  que le champ de vecteurs   $\U=(U_1,U_2,U_3)$ appartient à  l'espace $BMO^{-1}(\Rt)$ s'il existe un champ de vecteurs $\vec{V}=(V_1,V_2,V_3)$ où $V_i \in BMO(\Rt)$ pour $1\leq i \leq 3$; tel que l'on a 
\begin{equation}\label{def-BMO-1}
\U=\vec{\nabla}\wedge \V.
\end{equation}  Nous avons ainsi  le résultat suivant.
\begin{Theoreme}[\cite{Seregin1}, Théorème $1.1$]\label{Theo:Seregin1} Soit $\U \in L^{2}_{loc}(\Rt)$ une solution des équations de Navier-Stokes stationnaires (\ref{N-S-f-nulle-intro}). Si $\U \in L^6 \cap BMO^{-1}(\Rt)$ alors $\U=0$.\\
\end{Theoreme}
Nous allons maintenant expliquer  comment cette hypothèse nous permet d'obtenir l'identité $\U=0$.\\
\\
La preuve du Théorème \ref{Theo:Seregin1} donnée dans l'article \cite{Seregin1} suit essentiellement les grandes lignes  de la preuve du Théorème \ref{Theo:Galdi}; et cette preuve est divisée principalement en deux étapes: \\
\begin{enumerate}
	\item[$a)$]  La premier chose à faire est de montrer que  si la solution $\U$ des équations (\ref{N-S-f-nulle-intro})  appartient à l'espace $L^{6}\cap BMO^{-1}(\Rt)$ alors  cette solution appartient aussi à l'espace $\dot{H}^{1}(\Rt)$. Comme $\U \in BMO^{-1}(\Rt)$ alors  par la définition (\ref{def-BMO-1}) on peut écrire  $ \U=\vec{\nabla}\wedge \V$ avec $\V \in BMO(\Rt)$  et  l'on montre qu'il existe une constante $c_2=c_2(\Vert \U \Vert_{L^6}, \Vert \vec{V} \Vert_{BMO})>0$, telle que pour tout $R>0$ on a l'estimation
	\begin{equation}\label{estim-Seregin1}
	\int_{B_{\frac{R}{2}}} \vert \vec{\nabla}\otimes \U \vert^2dx \leq c_2.
	\end{equation}  Cette constante ne dépend pas de $R>0$ et alors si nous prenons la limite lorsque $R\longrightarrow +\infty$ dans l'estimation (\ref{estim-Seregin1}) nous obtenons $\U \in \dot{H}^{1}(\Rt)$. \\
	\\
	\item[$b)$] Une fois que l'on dispose de l'information $\U \in \dot{H}^{1}(\Rt)$, il s'agit de montrer l'identité $\U=0$. La preuve de cette identité donnée dans \cite{Seregin1} est  technique mais nous allons maintenant observer qu'en   passant  par le cadre des espaces  de Besov  on peut  vérifier plus facilement l'identité $\U=0$.\\
	\\
	Rappelons rapidement que pour un paramètre $s>0$ et une fonction $f$  on définit la quantité 
	\begin{equation}
	\Vert f \Vert_{\dot{B}^{-s}_{\infty,\infty}}=\sup_{t>0} t^{\frac{s}{2}} \Vert h_{t}\ast f \Vert_{L^{\infty}},
	\end{equation} où $h_t$ dénote toujours le noyau de la chaleur, et alors l'espace de Besov homogène  $\dot{B}^{-s}_{\infty,\infty}(\Rt)$ est défini comme 
	\begin{equation}\label{def-besov}
	\dot{B}^{-s}_{\infty,\infty}(\Rt)=\{  f \in \mathcal{S}^{'}(\Rt): \Vert f \Vert_{\dot{B}^{-s}_{\infty,\infty}}<+\infty\}.
	\end{equation}
	Une propriété des espaces de Besov qui nous sera très utile par la suite est la suivante:\\
	\begin{Lemme}\label{Lemme-cle-Besov} Soit $E\subset \mathcal{S}^{'}(\Rt)$ un espace de Banach tel que:
	\begin{enumerate}
		\item[1)] pour tout $u\in E$ et pour tout $x_0 \in \Rt$ on a $\Vert u(\cdot - x_0)\Vert_E=\Vert u \Vert_E$, et 
		\item[2)] pour tout $u\in E$ et pour tout $\lambda>0$ on a $\Vert u(\lambda \cdot)\Vert = \lambda^{-1}\Vert u \Vert_{E}$. \\
		\end{enumerate}	Alors $E \subset \dot{B}^{-1}_{\infty,\infty}(\Rt)$.\\
	\end{Lemme}
	Pour une preuve de ce résultat voir le Chapitre $4$ du livre \cite{PGLR2}.\\
	\\
	Maintenant que l'on a introduit les espaces de Besov que nous allons utiliser nous pouvons en déduire l'identité $\U=0$.  En effet, par le Lemme \ref{Lemme-cle-Besov} nous avons $BMO^{-1}(\Rt)\subset \dot{B}^{-1}_{\infty,\infty}(\Rt)$ et comme  $\U \in BMO^{-1}(\Rt)$  alors nous obtenons  $\U\in \dot{B}^{-1}_{\infty,\infty}(\Rt)$ et donc $\U \in \dot{H}^{1} \cap \dot{B}^{-1}_{\infty,\infty}(\Rt)$. Ainsi, par les inégalités de Sobolev précisées (voir l'article \cite{GerardMeyerOru} de P. Gérard, Y. Meyer et F.Oru) nous pouvons écrire 
	\begin{equation}\label{interpolation}
	\Vert \U \Vert_{L^4}\leq \Vert \U \Vert^{\frac{1}{2}}_{\dot{H}^1}\Vert \U \Vert^{\frac{1}{2}}_{\dot{B}^{-1}_{\infty,\infty}},
	\end{equation}  pour obtenir  $\U \in L^4(\Rt)$ et alors, par la Proposition \ref{Prop:L4}, nous avons l'identité $\U=0$. \\  
	\end{enumerate}
Nous pouvons ainsi observer que l'information $\U\in L^6 \cap BMO^{-1}(\Rt)$, est bien exploitée à chaque étape $a)$ et $b)$ ci-dessus et  cette information entraîne essentiellement  $\U \in L^4(\Rt)$ (voir les estimations (\ref{estim-Seregin1}) et \ref{interpolation}) ci-dessus) d'où  on obtient l'identité $\U=0$. \\
\\
Observons aussi  que nous $\U\in L^6 \cap BMO^{-1}(\Rt)$ alors cette solution appartient à deux espaces fonctionnels avec des homogénéités différentes:  l'espace $L^6(\Rt)$ est une espace homogène de degré $-\frac{1}{2}$ tandis que l'espace $BMO^{-1}(\Rt)$ est un espace homogène de degré $-1$. Dans ce cadre, nous allons maintenant observer que l'on peut remplacer l'espace  $L^6(\Rt)$ et l'espace  $BMO^{-1}(\Rt)$ par des espaces de Morrey ( chacun avec la même homogénéité de $L^6$ et $BMO^1$) pour obtenir un résultat similaire à celui obtenu dans le Théorème \ref{Theo:Seregin1}.\\   
\subsubsection{B) Un résultat dans les espaces de Morrey}
On commence par la définition suivante:
\begin{Definition}[Espace de Morrey homogène]\label{Def-Morrey} Soient $1<p\leq q <+\infty$. On définit l'espace de Morrey homogène $\dot{M}^{p,q}(\Rt)$ comme 
	$$ \dot{M}^{p,q}(\Rt)= \{  f\in L^{p}_{loc}(\Rt): \Vert f \Vert_{\dot{M}^{p,q}} <+\infty \}, $$ où 
	\begin{equation}\label{def-morrey}
	\Vert f \Vert_{\dot{M}^{p,q}} =\sup_{x_0 \in \Rt, R>0}
	R^{\frac{3}{q}}  \left(\frac{1}{R^3} \int_{B(x_0,R)}\vert f(x)\vert^p dx \right)^{\frac{1}{p}}.
	\end{equation}
	\end{Definition} 
Les espaces de Morrey sont des espaces fonctionnels plus généraux que les espaces de Lebesgue et  qui mesurent  la décroissance (en moyenne) à l'infini  d'une fonction avec les paramètres $1<p\leq q <+\infty$. En effet, observons tout d'abord que  par l'expression (\ref{def-morrey}) nous pouvons en tirer  l'inclusion $L^q(\Rt)\subset \dot{M}^{p,q}(\Rt)$ et le fait que $\dot{M}^{p,q}(\Rt)$ est un espace homogène de degré $-\frac{3}{q}$. \\
\\
Ensuite, toujours par  l'expression (\ref{def-morrey}), nous observons que si $f\in \dot{M}^{p,q}(\Rt)$ alors  pour tout $x_0\in \Rt$ et $R>0$ nous pouvons écrire
$$ \left(\frac{1}{R^3} \int_{B(x_0,R)}\vert f(x)\vert^p dx \right)^{\frac{1}{p}} \leq \Vert f \Vert_{\dot{M}^{p,q}}  R^{-\frac{3}{q}},$$ ce qui exprime le fait que lorsque  $R>0$ est suffisamment grand alors la moyenne de la fonction $f$ (en termes de la norme $L^p$) sur la boule $B(x_0,R)$   décroît à l'infini comme  $\ds{R^{-\frac{3}{q}}}$ où  l'exposant  $-\frac{3}{q}$ correspond  à homogénéité de l'espace $\dot{M}^{p,q}(\Rt)$.\\
\\
Les espaces de Morrey ont été largement étudiés en connexion avec la régularité des solutions faibles de certaines équations aux dérivées partielles (voir par exemple le Chapitre $13$ du livre \cite{PGLR1} pour des applications à l'étude de la théorie de la régularité des équations de Navier-Stokes); et nous allons maintenant voir comment ces espaces permettent aussi d'étudier le problème de Liouville pour les équations de Navier-Stokes stationnaires. Nous avons ainsi le résultat suivant obtenu par G. Seregin en $2016$ dans l'article \cite{Seregin2}. \\
\begin{Theoreme}[\cite{Seregin2}, Théorème $1.1$]\label{Theo:Seregin2}  Soit $\U\in L^{2}_{loc}(\Rt)$ une solution des équations de Navier-Stokes stationnaires (\ref{N-S-f-nulle-intro}). Si $\U \in \dot{M}^{2,6}\cap \dot{M}^{\frac{3}{2},3}(\Rt)$ alors $\U=0$. 
	\end{Theoreme}
Si nous comparons ce résultat avec le résultat obtenu dans le Théorème \ref{Theo:Seregin1}, où l'on a comme hypothèse $\U \in L^6\cap BMO^{-1}(\Rt)$, nous pouvons observer qu'ici on a  remplacé l'espace $L^6(\Rt)$ par l'espace de Morrey $\dot{M}^{2,6}(\Rt)$  qui a la même homogénéité $-\frac{1}{2}$ (plus précisément on  a $L^6\subset \dot{M}^{2,6}$), tandis que l'espace $BMO^{-1}$ a été remplacé par l'espace de Morrey $\dot{M}^{\frac{3}{2},3}$ qui a la même homogénéité $-1$. \\
\\
La preuve du Théorème \ref{Theo:Seregin2} suit les grandes lignes de la preuve du Théorème \ref{Theo:Seregin1} et cette fois-ci on a que si la solution $\U$ appartient à l'espace $\dot{M}^{2,6}\cap \dot{M}^{\frac{3}{2},3}(\Rt)$ alors il existe une constante $c_3=c_3(\Vert \U \Vert_{\dot{M}^{2,6}},\Vert \U \Vert_{\dot{M}^{\frac{3}{2},3}})>0$ telle que pour tout $R>0$ on a une même estimation que celle donnée dans (\ref{estim-Seregin1}):  
\begin{equation}
\int_{B_{\frac{R}{2}}}\vert \vec{\nabla} \otimes \U \vert^2 dx \leq c_3, 
\end{equation} et alors en prenant la limite lorsque $R \longrightarrow +\infty$ on obtient $\U \in \dot{H}^1(\Rt)$. Ensuite comme $\U \in \dot{M}^{\frac{3}{2},3}(\Rt)$ et comme par le Lemme \ref{Lemme-cle-Besov} l'on a l'inclusion $\dot{M}^{\frac{3}{2},3}(\Rt)\subset \dot{B}^{-1}_{\infty,\infty}(\Rt)$  alors nous avons $\U \in \dot{H}^{1} \cap \dot{B}^{-1}_{\infty,\infty}(\Rt)$ et donc par les inégalité de Sobolev précisées (\ref{interpolation}) nous avons $\U\in L^4(\Rt)$ et finalement grâce à la Proposition \ref{Prop:L4} nous pouvons écrire l'identité $\U=0$. \\
\\
\\
En résumé, nous observons que le problème de Liouville pour les équations de Navier-Stokes stationnaires (\ref{N-S-f-nulle-intro}) peut être résolu dans certains espaces  de Lebesgue et de Morrey; et que l'essentiel de la preuve de ces résultats repose sur le fait que l'on a besoin d'une certaine décroissance à l'infini de la solution $\U$ qui est caractérisée en termes des normes des ces  espaces fonctionnels.\\
\\
Avec cette idée en tête, dans la section qui suit nous démontrons quelques résultats nouveaux sur le problème de Liouville pour les équations   (\ref{N-S-f-nulle-intro}).

\section{Quelques résultats sur le problème de Liouville }\label{sec:resultats} 
En suivant les idées de la section précédente dans la Section \ref{Sec:result-lebesgue} ci-dessous on commence par traiter le cadre des espaces de Lebesgue pour ensuite traiter le cadre des espaces de Morrey dans la Section \ref{sec:cadre-morrey}.
 \subsection{Le cadre des  espaces de Lebesgue}\label{Sec:result-lebesgue} 
 Dans la section précédente nous avons énoncé trois  résultats connus sur le problème de Liouville pour les équations  (\ref{N-S-f-nulle-intro})  dans les espaces de Lebesgue: dans la Proposition \ref{Prop:L4} et dans le Théorème \ref{Theo:Galdi} nous avons observé que ce problème est résolu dans les l'espaces $L^4(\Rt)$ et $L^{\frac{9}{2}}(\Rt)$  respectivement, tandis que dans le Théorème \ref{Theo:Seregin1} nous observons que ce problème est aussi résolu dans un sous-espace particulier de $L^6(\Rt)$ donné par $L^6\cap BMO^{-1}(\Rt)$.  Nous nous posons alors la question s'il est possible d'étudier le problème de Liouville dans des espaces $L^p(\Rt)$ avec d'autres valeurs du paramètre d'intégration $p\in [1,+\infty]$. \\
 \\
 Dans ce cadre nous allons démontrer que l'on peut résoudre le problème de Liouville pour les équations (\ref{N-S-f-nulle-intro}) dans l'espace $L^p(\Rt)$ avec $p\in [3, \frac{9}{2}[$ et pour ces valeurs du paramètre $p$  n'avons pas besoin d'aucune hypothèse additionnelle sur la solution $\U$. En revanche, si nous considérons l'espace $L^p(\Rt)$  avec le paramètre $p\in ]\frac{9}{2},6[$ alors l'hypothèse $\U \in L^p(\Rt)$  ne suffit pas pour pouvoir obtenir l'identité cherchée $\U=0$ et, en suivant les idées du Théorème \ref{Theo:Seregin1},  nous avons besoin de considérer  un sous-espace de $L^p(\Rt)$ qui est donné par l'espace $L^p\cap\dot{B}^{\frac{3}{p}-\frac{3}{2}}_{\infty,\infty}(\Rt)$, où l'espace de Besov $\dot{B}^{\frac{3}{p}-\frac{3}{2}}_{\infty,\infty}(\Rt)$ est défini dans l'expression (\ref{def-besov}). \\
 \\
Nous avons ainsi   le résultat suivant: 
\begin{Theoreme}\label{Theo:Lebesgue} Soit $\U \in L^{2}_{loc}(\Rt)$ une solution  des équations de Navier-Stokes stationnaires   (\ref{N-S-f-nulle-intro}). 
\begin{enumerate}
\item[1)] Si $\U \in L^p(\Rt)$ avec $3 \leq p \leq \frac{9}{2}$ alors on a $\U=0$.
	
\item[2)] Si $\U \in L^p(\Rt)\cap \dot{B}^{\frac{3}{p}-\frac{3}{2}}_{\infty,\infty}(\Rt)$ avec  $\frac{9}{2}<p<6$ alors on a $\U=0$. \\
\end{enumerate}
\end{Theoreme} 
Dans ce résultat nous pouvons observer  que l'espace $L^{\frac{9}{2}}(\Rt)$ (considéré dans le Théorème \ref{Theo:Galdi} de G. Galdi) semble être un espace \emph{limite} pour résoudre le problème de Liouville dans le sens où pour  les valeurs du  paramètre d'intégration $p\in ]\frac{9}{2}, 6[$ on a besoin d'ajouter une hypothèse sur la solution $\U$ comme nous l'observons dans le point $2)$ ci)dessus.\\
\\
Soulignons aussi que le problème de Liouville dans les espaces de Lebesgue pour les valeurs du paramètre d'intégration $p\in [1,3[$ ou $p\in ]6,+\infty[$ est encore (dans nos connaissances actuelles)  une question ouverte. \\ 
\\
\textbf{Démonstration du Théorème \ref{Theo:Lebesgue}}.

\begin{enumerate}
\item[1)]  Nous supposons que la solution $\U$ vérifie $\U \in L^p(\Rt)$ avec $p \in [3,\frac{9}{2}]$.  Nous allons montrer l'identité $\U=0$ et pour cela nous allons suivre les grandes lignes de la preuve du Théorème \ref{Theo:Galdi} donnée dans le livre \cite{Galdi}. \\
\\
On commence par introduire la fonction de troncature suivante: soit $\theta \in \mathcal{C}^{\infty}_{0}(\Rt)$ telle que $0\leq \theta \leq 1$, $\theta(x)=1$ si $\vert x \vert <\frac{1}{2}$ et $\theta(x)=0$ si $\vert x \vert \geq 1$. Soit $R>1$, on définit la fonction $\theta_R(x)$ par $\theta_R(x)=\theta\left( \frac{x}{R}\right)$ où nous pouvons observer que l'on a $\theta_R(x)=1$ si $\vert x \vert < \frac{R}{2}$ et $\theta_R(x)=0$ si $\vert x \vert \geq R$.\\
\\
Nous multiplions maintenant  l'équation  (\ref{N-S-f-nulle-intro}) par la fonction $\theta_R \U$, ensuite nous intégrons sur la boule $B_R=\{x \in \Rt: \vert x \vert <R \}$ et nous obtenons l'identité 
\begin{equation*}
 \int_{B_R} \left(-\nu  \Delta \U +(\U \cdot \vec{\nabla})\U +\vec{\nabla}P \right)\cdot \theta_R \U dx=0.
\end{equation*}  Observons maintenant que comme $\U \in L^p(\Rt)$ avec $p\in [3,\frac{9}{3}]$ alors  $\U \in L^{3}_{loc}(\Rt)$ et par le Lemme \ref{Lemme-reg-sol} nous obtenons $\U \in \mathcal{C}^{\infty}(\Rt)$ et $P\in \mathcal{C}^{\infty}(\Rt)$. Ainsi, chaque terme de l'identité ci-dessus est suffisamment régulier et par des intégrations par parties nous pouvons écrire 
\begin{eqnarray*} 
\nu \int_{B_R} \theta_R \vert \vec{\nabla} \otimes \U \vert^2 dx & =&  \nu  \int_{B_R}  \theta_R \Delta \left(\frac{\vert \U \vert^2}{2}\right)  dx - \sum_{i,j=1}^{3}\int_{B_R} \partial_j \left( U_j \partial_j \left( \frac{U^{2}_{i}}{2} \right) \right) \theta_R dx \\
 & & - \sum_{i=1}^{3}\int_{B_R}(\partial_i P) \theta_R U_i dx \\
 &=&    \nu \int_{B_R}  \theta_R \Delta \left(\frac{\vert \U \vert^2}{2}\right)  dx - \int_{B_R} \theta_R \left[ div \left( \left( \frac{\vert \U \vert^2}{2}+P \right) \U \right) \right] dx \\
 &=& \nu \int_{B_R} \Delta \theta_R  \frac{\vert \U \vert^2}{2} dx + \int_{B_R} \vec{\nabla} \theta_R \cdot \left( \left( \frac{\vert \U \vert^2}{2}+P \right) \U \right) dx. 
\end{eqnarray*} D'autre part, comme $\theta_R (x)=1$ si $\vert x \vert < \frac{R}{2}$ alors nous avons  
$$\ds{ \nu \int_{B_{\frac{R}{2}}} \vert \vec{\nabla}\otimes \U \vert^2 dx \leq  \nu \int_{B_R} \theta_R \vert \vec{\nabla} \otimes \U \vert^2 dx},$$ et par l'identité précédente nous obtenons l'estimation 
\begin{eqnarray}\label{estim-aux-Lioville-1}\nonumber
\nu \int_{B_{\frac{R}{2}}} \vert \vec{\nabla}\otimes \U \vert^2 dx & \leq & \nu \int_{B_R} \Delta \theta_R  \frac{\vert \U \vert^2}{2} dx + \int_{B_R} \vec{\nabla} \theta_R \cdot \left( \left( \frac{\vert \U \vert^2}{2}+P \right) \U \right) dx \\
&= & I_1(R) +I_2(R), 
\end{eqnarray}  et nous allons maintenant montrer que $\ds{\lim_{R \longrightarrow +\infty}I_i(R)=0}$ pour $i=1,2$. \\
\\
En effet, pour le terme $I_1(R)$, par les inégalités de H\"older (avec $\frac{1}{q}+\frac{2}{p}=1$) nous avons 
$$ I_1(R)\leq \left( \int_{B_R} \vert \Delta \theta_R \vert^q dx\right)^{\frac{1}{q}} \left( \int_{B_R} \vert \U \vert^p dx\right)^{\frac{2}{p}} \leq \left( \int_{B_R} \vert \Delta \theta_R \vert^q dx\right)^{\frac{1}{q}} \Vert \U \Vert^{2}_{L^p}.$$ De plus, comme $\theta_R(x)=\theta\left(\frac{x}{R}\right)$ nous avons 
$$ \left( \int_{B_R} \vert \Delta \theta_R \vert^q dx\right)^{\frac{1}{q}}= R^{\frac{3}{q}-2} \Vert \Delta \theta \Vert_{L^q(B_1)},$$ et comme $\frac{1}{q}+\frac{2}{p}=1$ nous pouvons  écrire 
$$ I_1(R)\leq R^{1-\frac{6}{p}} \Vert \Delta \theta \Vert_{L^q(B_1)} \Vert \U \Vert^{2}_{L^p}.$$ Dans cette estimation, observons que comme $p\in [3,\frac{9}{2}]$ alors on a $1-\frac{6}{p}\in [-1,-\frac{1}{3}]$ d'où nous obtenons $\ds{\lim_{R\longrightarrow +\infty} I_1(R)=0}$. \\
\\
\'Etudions maintenant le terme $I_2(R)$.  On commence par observer que comme $\theta_R (x)=1$ si $\vert x \vert <\frac{R}{2}$ et $\theta_R(x)=0$ si $\vert x \vert \geq R$ alors on a $supp \left( \vec{\nabla}\theta_R \right) \subset \{ x\in \Rt: \frac{R}{2} < \vert x \vert <R \}$, et  nous pouvons  écrire  
\begin{eqnarray*}
I_2(R) &=& \int_{B_R} \vec{\nabla} \theta_R \cdot \left( \left( \frac{\vert \U \vert^2}{2}+P \right) \U \right) dx=\int_{\frac{R}{2} < \vert x \vert <R} \vec{\nabla} \theta_R \cdot \left( \left( \frac{\vert \U \vert^2}{2}+P \right) \U \right) dx \\
&\leq & \int_{\frac{R}{2} < \vert x \vert <R} \vert \vec{\nabla} \theta_R \vert \vert \U\vert^3 dx + \int_{\frac{R}{2} < \vert x \vert <R} \vert \vec{\nabla} \theta_R \vert \vert P \vert \vert \U\vert dx\\
&=& (I_2)_a(R)+(I_2)_b(R),
\end{eqnarray*} où, pour étudier ces deux termes ci-dessus nous allons suivre  essentiellement les mêmes lignes dans l'étude du terme $I_1(R)$. En effet, pour le terme $(I_2)_a(R)$, toujours par les inégalités de H\"older avec $ \ds{\frac{1}{r}+\frac{3}{p}=1}$, et par la définition de la fonction $\theta_R$  nous avons 
\begin{eqnarray*} \vspace{2mm}
(I_2)_a(R) &\leq & \left( \int_{\frac{R}{2}<\vert x \vert < R}\vert \vec{\nabla} \theta_R \vert^{r} dx \right)^{\frac{1}{r}} \left( \int_{\frac{R}{2}<\vert x \vert < R} \vert \U \vert^p dx\right)^{\frac{3}{p}} \\ \vspace{4mm}
& \leq & R^{\frac{3}{r}-1} \Vert \vec{\nabla} \theta \Vert_{L^r \left( \frac{1}{2}<\vert x \vert <1 \right)} \Vert \U \Vert^{3}_{L^p \left( \frac{R}{2}<\vert x \vert < R\right)}  \leq   R^{2-\frac{9}{p}} \Vert \vec{\nabla} \theta \Vert_{L^r}  \Vert \U \Vert^{3}_{L^p \left(\frac{R}{2}<\vert x \vert < R \right)}.
\end{eqnarray*}  Ainsi, comme $p\in [3,\frac{9}{2}]$ alors on a $2-\frac{9}{p}\in [-1,0]$ et comme $R>1$  on obtient  $R^{2-\frac{9}{p}}\leq 1$ et par l'estimation précédente nous pouvons écrire  
$$ (I_2)_a(R) \leq  \Vert \vec{\nabla} \theta \Vert_{L^r}  \Vert \U \Vert^{3}_{L^p(\frac{R}{2}<\vert x \vert < R)},$$ d'où, étant donné que $\U \in L^p(\Rt)$ alors nous obtenons $\ds{\lim_{R\longrightarrow+\infty} (I_2)_a(R)=0}$. \\
\\
Pour étudier le terme $ (I_2)_b(R)$, rappelons tout d'abord que comme la vitesse $\U$  appartient à l'espace $L^p(\Rt)$ alors par les inégalités de H\"older et par la continuité du projecteur de Leray  nous avons que la pression $P$ appartient à l'espace $L^{\frac{p}{2}}(\Rt)$ et ainsi, toujours par les inégalités de H\"older nous pouvons écrire 
$$ 	(I_2)_b(R) \leq  \left( \int_{\frac{R}{2}<\vert x \vert < R}\vert \vec{\nabla} \theta_R \vert^{r} dx \right)^{\frac{1}{r}} \left( \int_{\frac{R}{2}<\vert x \vert < R} \vert P \vert^{\frac{p}{2}} dx \right)^{\frac{2}{p}} \left( \int_{\frac{R}{2}<\vert x \vert < R} \vert \U \vert^{p}dx \right)^{\frac{1}{p}},$$ et en suivant les mêmes estimations ci-dessus nous obtenons $\ds{\lim_{R\longrightarrow+\infty} (I_2)_a(R)=0}$. De cette façons nous avons  $\ds{\lim_{R\longrightarrow+\infty} I_2(R)=0}$. \\
\\
Maintenant que l'on dispose de l'information $\ds{\lim_{R\longrightarrow+\infty} I_i(R)=0}$ pour $i=1,2$; nous revenons à l'estimation  (\ref{estim-aux-Lioville-1}) où, en prenant la limite lorsque $R\longrightarrow +\infty$, nous obtenons $ \ds{\nu \int_{\Rt} \vert \vec{\nabla} \otimes \U \vert^2 dx=0}$, et de cette façon nous avons  l'identité cherchée $\U=0$. \\

\item[2)] Nous supposons maintenant $\U \in L^p(\Rt)\cap \dot{B}^{\frac{3}{p}-\frac{3}{2}}_{\infty,\infty}(\Rt)$ avec  $\frac{9}{2}<p<6$ et nous allons montrer que $\U=0$.  La preuve de ce résultat suit les grandes lignes de la preuve du Théorème \ref{Theo:Seregin1} donnée dans l'article \cite{Seregin1} et la première chose à faire est vérifier  la proposition suivante.\\
\begin{Proposition}\label{Prop:Lp-Besov} Soit $\U \in L^p\cap \dot{B}^{\frac{3}{p}-\frac{3}{2}}_{\infty,\infty}(\Rt)$ avec $\frac{9}{2}<p<6$  une solution des équations (\ref{N-S-f-nulle-intro}). Alors  $\U\in \dot{H}^1(\Rt)$. 
\end{Proposition}
\textbf{Preuve.} La preuve de cette proposition repose essentiellement sur l'estimation suivante: soit $R>1$ et la boule $B_R=\lbrace x\in \Rt: \vert x \vert <R \rbrace$, alors on a   
\begin{equation}\label{Cacciopoli}
\int_{B_{\frac{R}{2}}}\vert \vec{\nabla}\otimes \U (x)\vert^2dx \leq C(R)\left(\int_{B_R}\vert \U(x)\vert^pdx\right)^{\frac{2}{p}},
\end{equation} où $C(R)=c_(\U)[R^{1-\frac{6}{p}}+1]$   avec une constante $c_(\U)>0$ qui dépend seulement de la  solution $\U$ et qui ne dépend  pas de $R>1$. \\
\\
Pour vérifier l'estimation (\ref{Cacciopoli}) nous allons suivre quelques idées des articles \cite{Seregin1} et \cite{Seregin2} de G. Seregin et nous introduisons  les fonctions de test $\varphi_R$ et $\W_R$ de la façon suivante: pour $R>1$ fixe,  nous définissons tout d'abord la fonction $\varphi_R$ en considérant  $\varphi_R\in \mathcal{C}^{\infty}_{0}(\Rt)$ telle que vérifie: $0\leq \varphi_R\leq 1$,  pour $\frac{R}{2}\leq \rho<r<R$ fixes on a  $\varphi_R(x)=1$ si $\vert x \vert<\rho$,  $ \varphi_R(x)=0$  si $\vert x \vert\geq r$; et $\Vert \vec{\nabla}\varphi_R\Vert_{L^{\infty}}\leq \frac{c}{r-\rho}$.\\
\\
Ensuite, nous définissons la fonction $\W_R$ comme la solution du problème 
\begin{equation}\label{eq_W_R}
div(\W_R)=\U\cdot \vec{\nabla}\varphi_R, \quad \text{sur}\,\, B_r, \quad \text{et}\quad \W_R=0 \,\, \text{sur}\,\,  \partial B_r;
\end{equation} où par le Lemme $III. 3.1$ page 162 du livre 
\cite{Galdi} il existe une solution $\W_R\in W^{1,p}(B_r)$ qui vérifie en plus l'estimation
\begin{equation}
\Vert \vec{\nabla}\otimes \W_R\Vert_{L^p(B_r)}\leq c \Vert \U\cdot  \vec{\nabla}\varphi_R \Vert_{L^p(B_r)}.
\end{equation} 
Avec les fonctions $\varphi_R$ et $\W_R$ ci-dessus, et avec la solution $\U$ des équations (\ref{N-S-f-nulle-intro}) nous considérons maintenant la fonction $\varphi_R \U-\W_R$ et nous écrivons  
$$ \int_{B_r} \left( -\Delta \U +(\U \cdot \vec{\nabla})\U +\vec{\nabla}p\right)\cdot \left( \varphi_R \U-\W_R \right)dx=0,$$ d'où nous avons 
\begin{equation}\label{eq_ipp1} \nonumber
\int_{B_r} \left( -\Delta \U +(\U \cdot \vec{\nabla})\U \right)\cdot \left( \varphi_R \U-\W_R \right)dx
+\int_{B_r} \left(\vec{\nabla}p\right)\cdot \left( \varphi_R \U-\W_R \right)dx=0.
\end{equation}  \`A ce stade, soulignons que comme $\U \in L^p(\Rt)$ avec $p\in ]\frac{9}{2},6[$ alors $\U \in L^{3}_{loc}(\Rt)$ et donc par le Lemme \ref{Lemme-reg-sol} nous avons $\U \in \mathcal{C}^{\infty}(\Rt)$ et $P\in \mathcal{C}^{\infty}(\Rt)$. De cette façon, chaque terme de l'identité ci-dessus est bien défini. \\
\\
Dans cette identité nous étudions  le deuxième terme à gauche où, étant donné que la fonction $\W_R$ vérifie l'équations  (\ref{eq_W_R}) et que $div(\U)=0$, par une intégration par parties nous avons  
$$ \int_{B_r} \vec{\nabla}p \cdot \left( \varphi_R \U-\W_R \right)dx=-\int_{B_R}p\left( \vec{\nabla}\varphi_R-div(\W_R)+\varphi_R div(\U)\right)dx=0,$$ et nous pouvons ainsi écrire 
$$ \int_{B_r} \left( -\Delta \U +(\U \cdot \vec{\nabla})\U \right)\cdot \left( \varphi_R \U-\W_R \right)dx=0,$$ et toujours par des intégrations par parties  nous obtenons  l'identité suivante  
\begin{eqnarray}\label{identity_ipp}\nonumber
\int_{B_R}\varphi_R \vert \vec{\nabla}\otimes \U\vert^2dx &=& -\int_{B_R}\vec{\nabla}\otimes \U:(\vec{\nabla}\varphi_R\otimes \U)dx+ \int_{B_R}\vec{\nabla}\otimes \W_R : \vec{\nabla}\otimes \U dx \\
& &-\int_{B_R}\left( (\U\cdot \vec{\nabla})\cdot \U \right)\cdot (\varphi_R \U-\W_R)dx=I_1+I_2+I_3. 
\end{eqnarray}
Nous devons estimer les termes $I_1,I_2$ et $I_3$ ci-dessus et nous avons ainsi le lemme technique suivant:  
\begin{Lemme}\label{Lemme1} Soit $\U \in L^p \cap  \dot{B}^{\frac{3}{p}-\frac{3}{2}}_{\infty,\infty}(\Rt)$ une solution des équations (\ref{N-S-f-nulle-intro}). On a les estimations suivantes: 
	\begin{enumerate}
		\item[1)] Comme $\U \in L^p(\Rt)$ alors il existe une constante numérique $c>0$ telle que l'on a  $\ds{I_1+I_2\leq   c\frac{R^{3\left( \frac{1}{2}-\frac{1}{p}\right)}}{r-\rho}\left(\int_{B_r}\vert \vec{\nabla}\otimes \U\vert^2dx\right)^{\frac{1}{2}} \left(\int_{B_r}\vert  \U\vert^pdx\right)^{\frac{1}{p}}.}$
		\item[2)] Comme  $\U \in \dot{B}^{\frac{3}{p}-\frac{3}{2}}_{\infty,\infty}(\Rt) $ alors il existe une constante $c_1=c_1(\U)>0$, qui ne dépend que de la solution $\U$, telle que l'on a $\ds{I_3\leq   c_1\frac{R}{r-\rho}\left(\int_{B_R}\vert \vec{\nabla}\otimes \U\vert^2dx\right)^{\frac{1}{2}} \left(\int_{B_R}\vert  \U\vert^pdx\right)^{\frac{1}{p}}.
		}$ \vspace{3mm}
	\end{enumerate}	
\end{Lemme}
L'estimation donnée dans le point $1)$ repose essentiellement sur les inégalités de H\"oder, tandis que l'estimation donnée dans le point $2)$ repose sur quelques propriétés des espaces de Besov et nous allons vérifier ces estimations en détail a la fin du chapitre page \pageref{preuve:Lemme1}. Ainsi, avec  les estimations ci-dessus on pose tout d'abord la constante $c(\U)=\max(c,c_1(\U))>0$ et on définit la quantité  
\begin{equation}\label{c(R)}
c(R)=c(\U)[R^{3\left( \frac{1}{2}-\frac{1}{P}\right)}+R].
\end{equation} Avec cette quantité 
nous revenons maintenant à l'identité (\ref{identity_ipp})  pour écrire  
$$ \int_{B_r}\varphi_R\vert \vec{\nabla}\otimes \U\vert^2dx\leq \frac{c(R)}{(r-\rho)} \left(\int_{B_R}\vert \vec{\nabla}\otimes \U\vert^2dx\right)^{\frac{1}{2}} \Vert U \Vert_{L^p}.$$ 
D'autre part, comme $\varphi_R(x)=1$ si $\vert x \vert <\rho$  alors nous avons 
$\ds{\int_{B_{\rho}}\vert \vec{\nabla}\otimes \U \vert^2 dx \leq \int_{B_r}\varphi_R\vert \vec{\nabla}\otimes \U\vert^2dx}$, et par l'estimation précédente nous écrivons 
$$ \int_{B_{\rho}}\vert \vec{\nabla}\otimes \U \vert^2 \leq \frac{c(R)}{(r-\rho)} \left(\int_{B_R}\vert \vec{\nabla}\otimes \U\vert^2dx\right)^{\frac{1}{2}} \Vert U \Vert_{L^p},$$ où, en appliquant les inégalités de Young au terme de droite  nous obtenons l'estimation 
\begin{equation}\label{Young}
\int_{B_{\rho}}\vert \vec{\nabla}\otimes \U \vert^2 dx \leq  \frac{1}{4} \int_{B_r}\vert \vec{\nabla}\otimes \U\vert^2dx +4\frac{[c(R)]^2 \Vert \U \Vert^{2}_{L^p}}{(r-\rho)^2},
\end{equation} et  nous allons choisir les paramètres $\rho>0$ et $r>0$ (avec $\frac{R}{2}<\rho<r<R$) convenablement pour obtenir l'estimation cherchée (\ref{Cacciopoli}). \\
\\
En effet, pour tout $k \in \mathbb{N}$ positif on pose  $ \ds{\rho_k=\frac{R}{2^{\frac{1}{k}}}}$, et dans l'estimation (\ref{Young}) on pose les paramètres $\rho=\rho_k$ et $r=\rho_{k+1}$ (où l'on a $\frac{R}{2} \leq \rho_k < \rho_{k+1} <R$) pour écrire 
$$ \int_{B_{\rho_k}}\vert \vec{\nabla}\otimes \U \vert^2 \leq \frac{1}{4}\int_{B_{\rho_{k+1}}}\vert \vec{\nabla}\otimes \U\vert^2dx+ 4\frac{[c(R)]^2 \Vert \U \Vert^{2}_{L^p}}{(\rho_{k+1}-\rho_k)^2}.$$ Nous devons maintenant étudier le deuxième terme à droite ci-dessus. \'Etant donné que  $ \ds{\rho_k=\frac{R}{2^{\frac{1}{k}}}}$ alors nous avons  $\ds{(\rho_{k+1}-\rho_k)^2 =R^2(2^{-\frac{1}{k+1}} -2^{-\frac{1}{k}})^2 \geq c \frac{  R^2}{k^2}}$, et ainsi  nous écrivons 
$$ 4\frac{[c(R)]^2 \Vert \U \Vert^{2}_{L^p}}{(\rho_{k+1}-\rho_k)^2} \leq  4c\,k^2 \frac{[c(R)]^2 \Vert \U \Vert^{2}_{L^p}}{R^2},$$ d'on nous obtenons l'estimation 
$$  \int_{B_{\rho_k}}\vert \vec{\nabla}\otimes \U \vert^2 \leq \frac{1}{4}\int_{B_{\rho_{k+1}}}\vert \vec{\nabla}\otimes \U\vert^2dx +4c\,k^2 \frac{[c(R)]^2 \Vert \U \Vert^{2}_{L^p}}{R^2},$$ 
et comme et $\rho_{k+1} < R$ nous écrivons la formule de récurrence suivante
$$ \int_{B_{\rho_k}}\vert \vec{\nabla}\otimes \U \vert^2 \leq  \frac{1}{4}\int_{B_{R}}\vert \vec{\nabla}\otimes \U\vert^2dx +4c\,k^2 \frac{[c(R)]^2 \Vert \U \Vert^{2}_{L^p}}{R^2}.$$
\`A cette stade, on itère cette formule de récurrence pour $k=1,\cdots, n$ et comme $\rho_{1}=\frac{R}{2}$ nous obtenons 
$$ \int_{B_{\frac{R}{2}}}\vert \vec{\nabla}\otimes \U \vert^2dx\leq \frac{1}{4^n}\int_{B_{R}}\vert \vec{\nabla}\otimes \U\vert^2dx +4c \frac{[c(R)]^2 \Vert \U \Vert^{2}_{L^p}}{R^2} \left[ \sum_{k=1}^{n} \frac{k^2}{4^k}\right],$$ d'où, 
en prenant maintenant la limite lorsque $n\longrightarrow +\infty$ nous avons  l'estimation 
$$\ds{\int_{B_{\frac{R}{2}}}\vert \vec{\nabla}\otimes \U \vert^2dx\leq c\frac{[c(R)]^2 \Vert \U \Vert_{L^p}}{R^2}}.$$ Finalement,  par l'identité (\ref{c(R)}) nous savons que $c(R)=c(\U)[R^{3\left( \frac{1}{2}-\frac{1}{p}\right)}+R]$ d'où nous obtenons   $\ds{\frac{[c(R)]^2}{R^2}\leq  c(\U)[R^{1-\frac{6}{p}} +1]}$. Ainsi, on pose  la constante $$\ds{ C(R)= c(\U)[R^{1-\frac{6}{p}} +1]},$$ et  nous obtenons finalement l'estimation  (\ref{Cacciopoli}). \\
\\
Maintenant que l'on dispose de cette estimation, observons que comme   $C(R)=c(\U)[R^{1-\frac{6}{p}}+1]$, et comme l'on a $p<6$ alors l'exposant $1-\frac{6}{p}$ est  une quantité négative et donc on a $\ds{\lim_{R\longrightarrow +\infty} C(R)\leq c(\U)<+\infty}$. Ainsi, en prenant la limite lorsque $R\longrightarrow +\infty$ dans l'estimation (\ref{Cacciopoli}) et de plus, étant donné que $\U \in L^p(\Rt)$ nous obtenons $ \Vert \U \Vert^{2}_{\dot{H}^1}\leq c(\U) \Vert \U \Vert^{2}_{L^p}<+\infty$. Nous avons de cette façon $\U \in \dot{H}^1(\Rt)$ et la Proposition \ref{Prop:Lp-Besov} est maintenant vérifiée.  \finpv  
\\
Par  la Proposition \ref{Prop:Lp-Besov}   nous avons  que la solution $\U$ appartient  à l'espace $\dot{H}^1(\Rt)$ et comme l'on a en plus $\U \in \dot{B}^{\frac{3}{p}-\frac{3}{2}}(\Rt)$ nous pouvons alors vérifier que $\U \in L^q(\Rt)$, pour un certain $q\in \left]3,\frac{9}{2}\right[$.    En effet, si l'on pose $\beta=\frac{3}{2}-\frac{3}{p}$ (où comme  $\frac{9}{2}<p<6$ alors on a $\frac{5}{6}<\beta <1$) alors 
par les inégalités de Sobolev précisées (voir toujours l'article \cite{GerardMeyerOru}) nous pouvons écrire 
$$ \Vert \U  \Vert_{L^q} \leq c \Vert  \U \Vert^{\theta}_{\dot{H}^{1}} \Vert \U \Vert^{1-\theta}_{\mathcal{B}^{-\beta}_{\infty,\infty}},$$ avec  $\theta=\frac{2}{q}$ et  $\beta=\frac{\theta}{1-\theta}$;  et par ces identités nous avons la relation   $q=\frac{2}{\beta}+2$ où, comme $\frac{5}{6}<\beta <1$ alors nous obtenons  $3<q<\frac{9}{2}$.\\ 
\\
Une que nous disposons de la l'information $\U \in L^q(\Rt)$ avec  $q\in \left]3,\frac{9}{2}\right[$ alors  par la Proposition \ref{Prop:unicite-Lp-intermediaire} nous avons l'identité cherchée $\U=0$. Le Théorème \ref{Theo:Lebesgue} est maintenant démontré. \finpv 
\\
\end{enumerate}
Dans la section qui suit nous nous intéressons à étudier le problème de Liouville pour les équations  (\ref{N-S-f-nulle-intro}) lorsque la solution $\U$ appartient à certains espaces de Morrey.
  
\subsection{Le cadre des espaces de Morrey}\label{sec:cadre-morrey}  
Comme annoncé nous considérons dans cette section les espaces de Morrey $\dot{M}^{p,q}(\Rt)$ (avec $1<p\leq q <+\infty$)  qui ont été  introduits dans la Définition \ref{Def-Morrey}; et en suivant les idées de la section précédente nous étudions pour quelles valeurs des paramètres $p$ et $q$ ci-dessus nous pouvons résoudre le problème de Liouville suivant:  si la solution  $\U$  des équations  de Navier-Stokes stationnaires (\ref{N-S-f-nulle-intro}) vérifie $\U \in \dot{M}^{p,q}(\Rt)$ alors $\U=0$. \\
\\
Rappelons rapidement que par le Théorème  \ref{Theo:Seregin2} (obtenue par G. Seregin dans l'article \cite{Seregin2}) nous avons que ce problème de Liouville est résolu dans l'espace $\dot{M}^{\frac{3}{2},3}\cap \dot{M}^{2,6}(\Rt)$ et l'essentiel des idées de la preuve de ce résultat (voir  la Section \ref{sec:resultats-connus} pour plus de détails) repose sur le fait que pour obtenir l'identité $\U=0$ nous avons besoin de considérer deux espaces de Morrey avec homogénéités différentes: l'espace $\dot{M}^{\frac{3}{2},3}(\Rt)$ est espace homogène de degré $-1$, tandis que l'espace $\dot{M}^{2,6}(\Rt)$ est un espace homogène de degré $-\frac{1}{2}$. \\
\\
Dans ce cadre, nous allons observer que l'on peut  généraliser ce résultat de la façon suivante: nous allons supposer  maintenant  que la solution  vérifie $\U \in \dot{M}^{2,3}(\Rt)\cap \dot{M}^{2,q}(\Rt)$, avec $3<q<+\infty$, où nous observons que l'on a conservé un espace de Morrey homogène de degré $-1$: $\dot{M}^{2,3}(\Rt)$, mais l'espace $\dot{M}^{2,6}(\Rt)$ a été remplacé par n'importe quel espace $\dot{M}^{2,q}(\Rt)$ qui a une homogénéité $-\frac{3}{q}\in ]-1,0[$, et c'est dans ce sens que le résultat suivant est une généralisation du Théorème \ref{Theo:Seregin2}.  \\ 
 \begin{Theoreme}\label{Theo:Seregin-generalise} Soit $\U \in L^{2}_{loc}(\Rt)$ une solution  des équations de Navier-Stokes stationnaires (\ref{N-S-f-nulle-intro}). Si $\U\in \dot{M}^{2,3}(\Rt)\cap \dot{M}^{2,q}(\Rt)$ avec $3<q<+\infty$ alors $\U=0$. \\
 \end{Theoreme} 
\dm  Nous considérons tout d'abord la solution stationnaire $\U$ comme la donnée initiale du problème de Cauchy pour les équations de Navier-Stokes:
\begin{equation}\label{N-S-aux2}
\partial_t \vu + (\vu \cdot \vec{\nabla})\vu -\nu \Delta \vu +\vec{\nabla} p =0, \quad div(\vu)=0, \quad \vu(0,\cdot)=\U.
\end{equation} 
Ensuite, par le  Théorème $8.2$ du livre \cite{PGLR1}   il existe un  temps $T_0>0$, et une fonction  $\vu\in \mathcal{C}([0,T_0[,\dot{M}^{2,q}(\Rt) )$ qui est une solution du problème de Cauchy (\ref{N-S-aux2}) et qui vérifie en plus l'estimation 
\begin{equation}\label{estim:Linf}
\sup_{0<t<T_0} t^{\frac{3}{2 q}}\Vert \vu(t,\cdot)\Vert_{L^{\infty}}<+\infty. 
\end{equation}
De plus, par le Théorème $8.4$ du livre \cite{PGLR1}  nous avons que la solution $\vu$ est l'unique solution de (\ref{N-S-aux2}) et comme $\U \in \mathcal{C}([0,T_0[,\dot{M}^{2,q}(\Rt))$ est aussi une solution du problème (\ref{N-S-aux2})  nous avons alors l'identité $\vu=\U$. Ainsi, par l'estimation (\ref{estim:Linf}) nous pouvons écrire 
\begin{equation}\label{estim-U-Linf}
\left(\frac{T_0}{2}\right)^{\frac{3}{2 q}} \Vert \U \Vert_{L^{\infty}}\leq \sup_{0<t<T_0} t^{\frac{1}{2}}\Vert \U\Vert_{L^{\infty}}<+\infty, 
\end{equation} pour obtenir  $\U \in L^{\infty}(\Rt)$. \\
\\
Maintenant que l'on dispose de l'information $\U \in L^{\infty}(\Rt)$, nous utilisons l'information additionnelle  $\U \in \dot{M}^{2,3}(\Rt)$ pour montrer que $\U =0$. La première chose à faire est de vérifier la proposition suivante: 
\begin{Proposition}\label{Prop:morrey-sobolev} Si $\U \in \dot{M}^{2,3}\cap L^{\infty}(\Rt)$ est une solution des équations de Navier-Stokes stationnaires (\ref{N-S-f-nulle-intro}) alors $\U \in \dot{H}^{1}(\Rt)$.
\end{Proposition}
\pv  En suivant quelques idées des articles \cite{Seregin1} et \cite{Seregin2}  nous considérons la fonction de troncature suivante: pour $R>1$ fixe, la fonction $\varphi_R \in \mathcal{C}^{\infty}_{0}(\Rt)$ est telle que $0\leq \varphi_R \leq 1$, $\varphi_R(x)=1$ si $\vert  x \vert < \frac{R}{2}$, $\varphi_R(x)=0$ si $\vert x \vert >R$ et l'on a $\Vert \vec{\nabla} \varphi_R \Vert_{L^{\infty}}\leq \frac{c}{R}$ et $\Vert \Delta \varphi_R \Vert_{L^{\infty}}\leq \frac{c}{R^2}$ où $c>0$ est une constante qui ne dépend pas de $R>1$. \\
\\
Une fois que l'on a introduit la fonction $\varphi_R$ ci-dessus, pour  la solution $\U $ nous considérons maintenant la fonction  $\varphi_R \U$ et comme $\U$ et $P$ vérifie les équations (\ref{N-S-f-nulle-intro})  nous pouvons écrire 
$$ \int_{B_R} \left( -\nu \Delta \U +(\U \cdot \vec{\nabla})\U+\vec{\nabla}P\right)\cdot \varphi_R \U dx=0,$$ où $B_R$ dénote toujours la boule  $\{ x \in \Rt:  \vert x \vert <R\}$. Nous allons maintenant étudier cette identité et pour cela nous avons besoin du résultat technique suivant:
\begin{Lemme}\label{Lemme:morrey-L-inf-interp} Soit $\U \in \dot{M}^{2,3}\cap L^{\infty}(\Rt)$. Alors on a $\Vert \U \Vert_{\dot{M}^{3,\frac{9}{2}}}\leq c \Vert U \Vert^{\frac{2}{3}}_{\dot{M}^{2,3}}\Vert \U \Vert^{\frac{1}{3}}_{L^{\infty}}$. \end{Lemme}
La preuve de ce lemme repose sur les inégalités de interpolation et elle sera faite en détail  à la fin du chapitre page \pageref{preuve:lemme-morrey-L-inf-interp}. Avec l'information $\U \in \dot{M}^{3,\frac{9}{2}}(\Rt)$, observons tout d'abord que l'on a $\U \in L^{3}_{loc}(\Rt)$ et alors par le Lemme \ref{Lemme-reg-sol} on obtient $\U \in \mathcal{C}^{\infty}(\Rt)$ et $P\in \mathcal{C}^{\infty}(\Rt)$. De cette façon, nous pouvons faire une intégration par parties dans l'identité précédente et nous obtenons   
\begin{equation}\label{iden1-aux}
\int_{B_R} \varphi_R \vert \vec{\nabla}\otimes \U \vert^2 dx= \int_{B_R} \frac{\vert \U \vert^2}{2}\Delta \varphi_R dx +\int_{B_R} \left( \frac{\vert \U \vert^2}{2}+P \right) (\U \cdot \vec{\nabla}\varphi_R)dx= I_1(R)+I_2(R),
\end{equation} où nous cherchons à estimer les termes $I_1(R)$ et $I_2(R)$.\\
\\
Pour  le terme $I_1(R)$ de l'identité (\ref{iden1-aux}), comme  $\Vert \Delta \varphi_R \Vert_{L^{\infty}}\leq \frac{c}{R^2}$ nous avons  
$$ \vert I_1(R) \vert \leq  \frac{c}{R}\int_{B_R}\vert \U \vert^2 dx \leq \frac{c}{R^2}R^{6\left(\frac{1}{2}-\frac{1}{3} \right)} \left( \int_{B_R}\vert \U \vert^3dx\right)^{\frac{2}{3}}\leq \frac{c}{R}\left( \int_{B_R}\vert \U \vert^3dx\right)^{\frac{2}{3}},$$ et alors, comme  $\U \in \dot{M}^{3,\frac{9}{2}}(\Rt)$ nous pouvons écrire $\ds{\left( \int_{B_R}\vert \U \vert^3dx\right)^{\frac{2}{3}} \leq \Vert \U \Vert^{2}_{\dot{M}^{3,\frac{9}{2}}}R^{6\left( \frac{1}{3}-\frac{2}{9}\right)}}$, d'où nous obtenons l'estimation
\begin{equation}\label{estim_I_1_Prop_Th_2}
\vert I_1(R) \vert \leq \frac{c}{R}\Vert \U \Vert^{2}_{\dot{M}^{3,\frac{9}{2}}}R^{6\left( \frac{1}{3}-\frac{2}{9}\right)} \leq \frac{c}{R^{\frac{1}{3}}}  \Vert \U \Vert^{2}_{\dot{M}^{3,\frac{9}{2}}}. 
\end{equation} 
Pour le terme $I_2(R)$ de l'identité (\ref{iden1-aux}), étant donné que  $\Vert \vec{\nabla} \varphi_R \Vert_{L^{\infty}}\leq \frac{c}{R}$ alors nous avons 
\begin{equation}\label{estim-aux-morrey}
 \vert I_2(R) \vert \leq \frac{c}{R}\int_{B_R}\vert \U \vert^3 dx+ \frac{c}{R} \int_{B_R} \vert p \vert \vert \U \vert dx=(I_2)_a+(I_2)_b,
 \end{equation} et nous avons encore besoin d'estimer les termes $(I_2)_a$ et $(I_2)_b$ ci-dessus. Pour le terme $(I_2)_a$, toujours comme $\U \in \dot{M}^{3,\frac{9}{2}}(\Rt)$ nous pouvons écrire $\int_{B_R}\vert \U \vert^3dx\leq \Vert \U\Vert^{3}_{\dot{M}^{3,\frac{9}{2}}}R^{9\left(\frac{1}{3}-\frac{2}{9}\right)}$ et alors nous obtenons 
 \begin{equation}\label{estim-I2a}
 (I_2)_a\leq\ds{\frac{c}{R}\int_{B_R}\vert \U \vert^3 dx \leq  \frac{c}{R}\Vert \U\Vert^{3}_{\dot{M}^{3,\frac{9}{2}}}R^{9\left(\frac{1}{3}-\frac{2}{9}\right)}\leq c \Vert \U\Vert^{3}_{\dot{M}^{3,\frac{9}{2}}}}.
 \end{equation}
Nous allons maintenant étudier le  terme $(I_2)_b$ de l'identité (\ref{estim-aux-morrey}) et pour cela nous avons besoin du lemme technique suivant.
\begin{Lemme}\label{Lemme:pression} Soient $\U,P$ une solution des équations de Navier-Stokes stationnaires (\ref{N-S-f-nulle-intro}). Si $\U \in  \dot{M}^{p,q}(\Rt)$ avec $p\geq 2$ et $q\geq 3$ alors $P \in \dot{M}^{\frac{p}{2},\frac{q}{2}}(\Rt)$. 
\end{Lemme} La preuve de ce lemme  repose essentiellement sur la continuité du projecteur de Leray sur les espaces de Morrey et cette preuve sera faite à la fin du chapitre page \pageref{preuve:lemme-pression}. Comme $\U \in \dot{M}^{3,\frac{9}{2}}(\Rt)$ alors  par le Lemme  \ref{Lemme:pression} nous avons  $P\in \dot{M}^{\frac{3}{2},\frac{9}{4}}(\Rt)$; et par les inégalités de  H\"older  nous écrivons
\begin{eqnarray}\label{estim_I2b} \nonumber
	(I_2)_b &\leq &\frac{c}{R}\left(\int_{B_R}\vert p \vert^{\frac{3}{2}} \right)^{\frac{2}{3}} \left(\int_{B_R}\vert \U \vert^3 \right)^{\frac{1}{3}}\leq \frac{c}{R}\left[ \Vert p \Vert_{\dot{M}^{\frac{3}{2},\frac{9}{4}}}R^{3\left(\frac{2}{3}-\frac{4}{9} \right)} \right]\left[ \Vert \U \Vert_{\dot{M}^{3,\frac{9}{2}}}R^{3\left(\frac{1}{3}-\frac{2}{9} \right)} \right]\\
	&\leq& c \Vert p \Vert_{\dot{M}^{\frac{3}{2},\frac{9}{4}}} \Vert \U \Vert_{\dot{M}^{3,\frac{9}{2}}}\leq c \Vert \U \Vert^{3}_{\dot{M}^{3,\frac{9}{2}}}.
\end{eqnarray}
Ainsi, avec les estimations (\ref{estim-I2a}) et (\ref{estim_I2b}) nous revenons à l'inégalité (\ref{estim-aux-morrey}) pour écrire
\begin{equation}\label{estim_term_I_2_Th_2}
\vert I_2(R)\vert\leq c \Vert \U \Vert^{3}_{\dot{M}^{3,\frac{9}{2}}}.
\end{equation} Revenons maintenant à l'identité  (\ref{iden1-aux}) où, étant donné que $\varphi_R(x)=1$ si $\vert x \vert <\frac{R}{2}$, et de plus, par les estimations  (\ref{estim_I_1_Prop_Th_2}) et (\ref{estim_term_I_2_Th_2}), nous pouvons alors écrire
$$\int_{B_{\frac{R}{2}}}\vert \vec{\nabla}\otimes \U\vert^2dx \leq \int_{B_{R}}\psi_R \vert \vec{\nabla}\otimes \U\vert^2dx \leq \frac{c}{R^{\frac{1}{3}}}  \Vert \U \Vert^{2}_{\dot{M}^{3,\frac{9}{2}}} +c \Vert \U \Vert^{3}_{\dot{M}^{3,\frac{9}{2}}},$$ d'où, nous prenons  la limite lorsque  $R \longrightarrow+\infty$ et nous avons $\U \in \dot{H}^{1}(\Rt)$. La Proposition \ref{Prop:morrey-sobolev} est démontrée. \finpv
Maintenant nous avons l'information dont on a besoin pour montrer l'identité $\U=0$. En effet, par le Lemme \ref{Lemme-cle-Besov} nous avons  $\dot{M}^{2,3}(\Rt) \subset \dot{B}^{-1}_{\infty,\infty}(\Rt)$ et comme  $\U \in \dot{M}^{2,3}(\Rt)$ alors $ \U \in \dot{B}^{-1}_{\infty,\infty}(\Rt)$. Ensuite, para la Proposition \ref{Prop:morrey-sobolev} nous avons ainsi $\U \in \dot{H}^{1}(\Rt)$ et alors par les inégalités de Sobolev précisées (\ref{interpolation}) nous avons $\U \in L^4(\Rt)$. Ainsi, par le point $1)$ du Théorème \ref{Theo:Lebesgue} nous avons finalement $\U =0$. Le théorème \ref{Theo:Seregin-generalise}  est maintenant démontré \finpv 
\\
Dans la démonstration du Théorème  \ref{Theo:Seregin-generalise} que nous venons de faire nous pouvons observer que nous avons besoin de l'information $\U \in \dot{M}^{2,3}(\Rt)$ et $\U \in \dot{M}^{2,q}(\Rt)$ (avec $3<q<+\infty$) pour obtenir l'identité  $\U=0$. Néanmoins,  dans le  théorème  ci-dessous nous allons observer que si l'on considère maintenant un sous-espace particulier de $\dot{M}^{2,3}(\Rt)$ alors nous n'avons pas besoin de considérer un deuxième espace de Morrey (avec une homogénéité différente) pour obtenir $\U=0$.  \\
\\
Nous commençons donc par introduire, en toute généralité, les espaces fonctionnels  qui nous utiliserons ci-après. 
 \begin{Definition}[L'espace $\overline{M}^{p,q}(\Rt)$]\label{Def:sous-espace-morrey} Soient $1<p\leq q <+\infty$. On définit l'espace $\overline{M}^{p,q}(\Rt)$ comme l'adhérence de l'espace des fonction de test $\mathcal{C}^{\infty}_{0}(\Rt)$ dans l'espace de Morrey  $\dot{M}^{p,q}(\Rt)$ donné dans la Définition \ref{def-morrey}.
\end{Definition}
Une fois que l'on a défini ces espaces fonctionnels, on fixe maintenant le paramètre $q=3$ et pour $p\in ]1,3]$ on considère l'espace  $\overline{M}^{p,3}(\Rt) \subset\dot{M}^{p,3}(\Rt)$. Observons tout d'abord que par la définition ci-dessus nous avons $\overline{M}^{p,3}(\Rt) \subset\dot{M}^{p,3}(\Rt)$  et comme  l'on a aussi l'inclusion $\dot{M}^{p,3}(\Rt) \subset \dot{M}^{2,3}(\Rt)$ alors  l'espace $\overline{M}^{p,3}(\Rt)$ est un sous-espace de $\dot{M}^{2,3}(\Rt)$ et il est  un espace homogène de degré $-1$. \\
\\
Nous allons démontrer que l'on peut résoudre un problème de Liouville pour les équations  (\ref{N-S-f-nulle-intro}) lorsque la solution $\U$ appartient à l'espace $\overline{M}^{p,3}(\Rt)$ pour des valeurs convenables du paramètre $p$.\\
 \begin{Theoreme}\label{Theo:Morrey-2} Soit l'espace $\overline{M}^{p,3}(\Rt)$ avec $p \in ]2,3]$, donné la Définition \ref{Def:sous-espace-morrey}. Soit $\U \in L^{2}_{loc}(\Rt)$ une solution  des équations de Navier-Stokes stationnaires (\ref{N-S-f-nulle-intro}). Si $\U\in \overline{M}^{p,3}(\Rt)$  alors $\U=0$. \\
 \end{Theoreme}
Si nous comparons ce résultat avec le Théorème \ref{Theo:Seregin-generalise} nous observons que dans le cadre particulier de l'espace $\overline{M}^{p,3}(\Rt)$ on n'a pas besoin de faire aucune hypothèse additionnelle sur la solution $\U$ pour obtenir l'identité $\U=0$. \\
\\ 
\dm Pour démontrer ce théorème nous allons suivre les grandes lignes de la démonstration du Théorème  \ref{Theo:Seregin-generalise} et on commence par considérer la solution $\U$ comme la donnée initiale du problème de Cauchy (\ref{N-S-aux2}). Ainsi, toujours  par   Théorème $8.2$ du livre \cite{PGLR1} 
il existe une fonction $\vu \in \mathcal{C}([0,T_0[,\overline{M}^{p,3}(\Rt) )$  qui est l'unique solution du problème de Cauchy (\ref{N-S-aux2}) (grâce au Théorème $8.4$ du livre \cite{PGLR1}) et qui vérifie en plus l'estimation 
 \begin{equation}\label{estim} 
 \sup_{0<t<T_0} t^{\frac{1}{2}}\Vert \vu(t,\cdot)\Vert_{L^{\infty}}<+\infty. 
 \end{equation} Mais, comme la solution stationnaire $\U$ est aussi une solution du problème (\ref{N-S-aux2}) alors on a $\vu=\U$ et donc la solution $\U$ vérifie l'estimation (\ref{estim}) d'où nous pouvons en tirer $\U\in L^{\infty}(\Rt)$. \\
 \\
 Maintenant que l'on dispose de l'information $\U \in \overline{M}^{p,3} \cap L^{\infty}(\Rt)$ nous allons montrer l'identité $\U=0$. En effet, rappelons tout d'abord que l'on a l'inclusion $\overline{M}^{p,3} \subset \dot{M}^{2,3}(\Rt)$ et nous avons ainsi $\U \in \dot{M}^{2,3}(\Rt) \cap L^{\infty}(\Rt)$. Ensuite, par la Proposition \ref{Prop:morrey-sobolev} nous avons $\U \in \dot{H}^{1}(\Rt)$ et alors 
 $\U \in \dot{H}^{1}(\Rt) \cap \dot{M}^{2,3}(\Rt) \subset \dot{H}^{1}(\Rt) \cap \dot{B}^{-1}_{\infty,\infty}(\Rt)$, d'où, toujours par les inégalités de Sobolev précisées (\ref{interpolation}) nous avons $\U\in L^4(\Rt)$, et par le point $1)$ du Théorème  \ref{Theo:Lebesgue} nous pouvons finalement écrire $\U=0$. \finpv

\section{Preuve des lemmes techniques}

\subsubsection{Preuve du Lemme \ref{lemme-tech} page \pageref{lemme-tech}}\label{preuve:lemme-tech}
On commence par montrer que si $\U \in L^4(\Rt)$ alors on a $\vec{\nabla} P \in \dot{H}^{-1}(\Rt)$. En effet, pour le champ de vitesse $\U=(U_1,U_2,U_3)$ nous pouvons écrire  la pression $P$ comme $\ds{P=\sum_{i,j}^{3}\frac{1}{-\Delta}\partial_i \partial_j (U_i U_j)}$,  où nous observons que comme $\U \in L^4(\Rt)$ alors par les inégalités de H\"older nous avons $U_i \,U_j \in L^2(\Rt)$ (pour $1\leq i,j\leq 3$)  et donc $P \in L^2(\Rt)$. De cette façon nous obtenons $\vec{\nabla} P \in \dot{H}^{-1}(\Rt)$. \\
\\
Vérifions maintenant que $(\U\cdot \vec{\nabla})\cdot \U \in \dot{H}^{-1}(\Rt)$. Nous écrivons la solution $\U$ comme
\begin{equation}\label{N-S-pf-aux}
\U=\P \left( \frac{1}{\nu \Delta} div(\U \otimes \U) \right),
\end{equation} où nous observons que comme $\U \in L^4(\Rt)$ alors toujours par les inégalités de H\"older nous avons $\U \otimes \U\in L^2(\Rt)$ et donc $\P \left( \frac{1}{\nu \Delta} div(\U \otimes \U) \right) \in \dot{H}^1(\Rt)$. Ainsi, par l'identité (\ref{N-S-pf-aux}) nous avons $\U \in \dot{H}^1(\Rt)$ et alors $\Delta \U \in \dot{H}^{-1}(\Rt)$. \\
\\
Nous écrivons finalement $(\U\cdot \vec{\nabla})\cdot \U=\Delta \U -\vec{\nabla}P $ pour obtenir   $(\U\cdot \vec{\nabla})\cdot \U \in \dot{H}^{-1}(\Rt)$. \finpv
\subsubsection{Preuve du Lemme \ref{Lemme1} page \pageref{Lemme1}} 
\begin{enumerate}\label{preuve:Lemme1}
	\item[1)] On commence par estimer le terme $I_1$. En utilisant les  inégalités de H\"older  (avec $\frac{1}{2}=\frac{1}{q}+\frac{1}{p}$) et comme $\Vert \vec{\nabla}\varphi_R\Vert_{L^{\infty}}\leq \frac{c}{r-\rho}$ nous pouvons écrire
	\begin{eqnarray}\label{estimate_I_1}\nonumber
	\vert I_1 \vert &\leq&  \left(\int_{B_r}\vert \vec{\nabla}\otimes \U\vert^2dx\right)^{\frac{1}{2}}\left(\int_{B_r}\vert \vec{\nabla}\varphi_R\otimes \U\vert^2dx\right)^{\frac{1}{2}}\leq  \left(\int_{B_r}\vert \vec{\nabla}\otimes \U\vert^2dx\right)^{\frac{1}{2}} \left(\int_{B_r}\vert \vec{\nabla}\varphi_R\vert^qdx\right)^{\frac{1}{q}} \times \\ \nonumber
	& & \times  \left(\int_{B_R}\vert  \U\vert^pdx\right)^{\frac{1}{p}}  \leq c\frac{r^{\frac{3}{q}}}{r-\rho} \left(\int_{B_r}\vert \vec{\nabla}\otimes \U\vert^2dx\right)^{\frac{1}{2}} \left(\int_{B_r}\vert  \U\vert^pdx\right)^{\frac{1}{p}} \\
	& \leq & c\frac{R^{3\left( \frac{1}{2}-\frac{1}{p}\right)}}{r-\rho}\left(\int_{B_r}\vert \vec{\nabla}\otimes \U\vert^2dx\right)^{\frac{1}{2}} \left(\int_{B_r}\vert  \U\vert^pdx\right)^{\frac{1}{p}}. 
	\end{eqnarray}
	Estimons maintenant le terme $I_2$. En utilisant l'inégalité de  Cauchy-Schwarz et comme la fonction $\W_R$ vérifie $\Vert \vec{\nabla}\otimes \W_R\Vert_{L^p(B_r)}\leq c \Vert \U\cdot \vec{\nabla}\varphi_R \Vert_{L^p(B_r)}$  nous avons 
	\begin{eqnarray}\label{estmate_I_2}\nonumber
	 \vert I_2 \vert &\leq &  \left(\int_{B_r}\vert \vec{\nabla}\otimes \U\vert^2dx\right)^{\frac{1}{2}}\left(\int_{B_r}\vert \vec{\nabla}\otimes \W_R\vert^2dx\right)^{\frac{1}{2}}\leq  \left(\int_{B_r}\vert \vec{\nabla}\otimes \U\vert^2dx\right)^{\frac{1}{2}} r^{3\left(\frac{1}{2}-\frac{1}{p}\right)} \times \\ \nonumber
	 & & \times  \left(\int_{B_r}\vert \vec{\nabla}\otimes \W_R\vert^pdx\right)^{\frac{1}{p}} 
	\leq  \left(\int_{B_r}\vert \vec{\nabla}\otimes \U\vert^2dx\right)^{\frac{1}{2}} r^{3\left(\frac{1}{2}-\frac{1}{p}\right)}\left(\int_{B_r}\vert \vec{\nabla}\varphi\cdot \U \vert^pdx\right)^{\frac{1}{p}} \\
	&\leq & c\frac{R^{3\left( \frac{1}{2}-\frac{1}{p}\right)}}{r-\rho}\left(\int_{B_r}\vert \vec{\nabla}\otimes \U\vert^2dx\right)^{\frac{1}{2}} \left(\int_{B_r}\vert  \U\vert^pdx\right)^{\frac{1}{p}}.
	\end{eqnarray}
	\item[2)] \'Etudions maintenant le terme $I_3$ défini par
	 
\begin{equation}\label{I_3}
I_3= -\int_{B_r}\left(U_1\partial_1\U+U_2\partial_2\U+U_3\partial_3\U\right) \cdot (\varphi_R \U-\W_R)dx.
\end{equation} On commence par observer que l'on peut toujours écrire $\U=\vec{\nabla}\wedge \V$. En effet, il suffit de définir 
\begin{equation}\label{def-V}
\V= \vec{\nabla}\wedge \left( \frac{1}{-\Delta} \U \right),
\end{equation} et comme $div(\U)=0$ on a l'identité 
$$ \vec{\nabla}\wedge \V=\vec{\nabla}\wedge \left( \vec{\nabla}\wedge \left( \frac{1}{-\Delta} \U \right) \right)=\vec{\nabla} \left( div \left( \frac{1}{-\Delta} \U  \right)\right) - \Delta \left( \frac{1}{-\Delta} \U  \right) =\U.$$ Observons aussi le fait que comme $\U \in \mathcal{B}^{\frac{3}{p}-\frac{3}{2}}_{\infty,\infty}(\Rt)$ alors par l'identité (\ref{def-V}) on a $\V\in \mathcal{B}^{\frac{3}{p}-\frac{1}{2}}_{\infty,\infty}(\Rt)$.  De plus, pour mener à bien les estimations dont on aura besoin nous allons poser la fonction $\V^{*}=\V-\V(0)$ et comme  l'on a $\vec{\nabla}\wedge \V =\vec{\nabla}\wedge \V^{*}$ (car $\V(0)\in \Rt$ est un vecteur constant) et nous allons écrire $\U=\vec{\nabla}\wedge  \V^{*}$, c'est à dire, nous avons les identités $U_i=\partial_j V^{*}_{k}-\partial_{k}V^{*}_{j}$ où les indices $i,j,k=1,2,3$ sont ordonnés par la règle de la main droite.\\
\\
De cette façon, dans l'expression (\ref{I_3}) nous remplaçons $U_i$ par  $\partial_j V^{*}_{k}-\partial_{k}V^{*}_{j}$  et nous obtenons l'identité suivante:
\begin{eqnarray*}
	I_3 &=&-\int_{B_r}\left((\partial_2 V^{*}_3-\partial_3 V^{*}_2)\partial_1\U+(\partial_3 V^{*}_1-\partial_1 V^{*}_3)\partial_2\U+(\partial_1 V^{*}_2-\partial_2 V^{*}_1)\partial_3\U\right)\cdot(\varphi_R \U-\W_R)dx\\
	&=&-\int_{B_r} \sum_{i=1}^{3} \left((\partial_jV^{*}_{k}-\partial_k V^{*}_{j})\partial_i \U \right)\cdot (\varphi_R \U-\W_R)dx\\
	&=& -\int_{B_r}  \sum_{i=1}^{3} \left[\partial_jV^{*}_{k} (\partial_i \U)\cdot (\varphi_R \U-\W_R)-\partial_k V^{*}_{j}(\partial_i \U)\cdot (\varphi_R \U-\W_R)\right]dx,
\end{eqnarray*} et en faisant une intégration par parties nous avons
\begin{eqnarray}\label{estim:I3} \nonumber
	I_3&=& \int_{B_r} \sum_{i=1}^{3}\left[  V^{*}_{k}(\partial_j \partial_i \U )\cdot (\varphi_R \U-\W_R) +V^{*}_{k}(\partial_i\U)\cdot   \partial_j(\varphi_R \U-\W_R)-V^{*}_{j}(\partial_k\partial_i\U)\cdot(\varphi_R \U-\W_R) \right.\\ \nonumber
	& &\left.-V^{*}_{j}(\partial_i\U)\cdot\partial_k(\varphi_R \U-\W_R) \right]dx \\ \nonumber
	&=& \int_{B_r} \sum_{i=1}^{3} \left[ \left( V^{*}_{k}(\partial_j \partial_i \U )-V^{*}_{j}(\partial_k\partial_i\U)\right)\cdot (\varphi_R \U-\W_R)\right]dx  + \int_{B_r} \sum_{i=1}^{3} \left[V^{*}_{k}(\partial_i\U)\cdot   \partial_j(\varphi_R \U-\W_R)\right.\\
	& & \left.-V^{*}_{j}(\partial_i\U)\cdot\partial_k(\varphi_R \U-\W_R)\right]dx =(I_3)_a+(I_3)_b.
\end{eqnarray}
Nous devons étudier encore les termes  $(I_3)_a$ et $(I_3)_b$ ci-dessus. Pour le terme $(I_3)_a$, comme les indices $i,j,k=1,2,3$ sont ordonnées par la règle de la main droite nous avons alors  $\sum_{i=1}^{3} \left[ \left( V^{*}_{k}(\partial_j \partial_i \U )-V^{*}_{j}(\partial_k\partial_i\U)\right)\right]=(0,0,0)$ et  donc nous avons 
\begin{equation}\label{estim:Ia}
(I_3)_a= \int_{B_r} \sum_{i=1}^{3} \left[ \left( V^{*}_{k}(\partial_j \partial_i \U )-V^{*}_{j}(\partial_k\partial_i\U)\right)\cdot (\varphi_R \U-\W_R)\right]dx =0. 
\end{equation}

Pour le terme  $(I_3)_b$ nos écrivons
\begin{eqnarray*}
	(I_3)_b &=& \int_{B_r}\sum_{i=1}^{3}\left[ V^{*}_{k}\partial_i \U\cdot (\partial_j \varphi_R \U+\varphi_R\partial_j \U-\partial_j\W_R)  -V^{*}_{j}\partial_i\U \cdot(\partial_k \varphi_R \U+\varphi_R\partial_k \U-\partial_k \W_R)\right]dx\\
	&=&\int_{B_r}\underbrace{\sum_{i=1}^{3} \left[ V^{*}_{k}\partial_i \U \cdot (\varphi_R \partial_j\U)-V^{*}_{j}\partial_i\U \cdot(\varphi_R\partial_k\U) \right]}_{(a)} dx+\int_{B_r} \sum_{i=1}^{3} \left[ V^{*}_{k}\partial_i\U\cdot (\partial_I\varphi_R \U-\partial_j\W_R)\right. \\
	&=&\left. - V^{*}_{j}\partial_i \U \cdot(\partial_k \varphi_R \U-\partial_k \W_R)\right]dx, 
\end{eqnarray*} 
où, toujours par le fait que les indices $i,j,k=1,2,3$ sont ordonnés par la règle de la main droite nous avons que le terme $(a)$ ci-dessus est égale à zéro et ainsi nous avons 
\begin{equation}\label{estim:Ib}
(I_3)_b=\int_{B_r} \sum_{i=1}^{3} \left[ V^{*}_{k}\partial_i\U\cdot (\partial_I\varphi_R \U-\partial_j\W_R)- V^{*}_{j}\partial_i \U \cdot(\partial_k \varphi_R \U-\partial_k \W_R)\right]dx.
\end{equation}
Une fois que nous disposons  des identités (\ref{estim:Ia}) et (\ref{estim:Ib}) nous revenons au terme $I_3$ donné dans (\ref{estim:I3}) pour écrire 
$$ I_3= \int_{B_r} \sum_{i=1}^{3} \left[ V^{*}_{k}\partial_i\U\cdot (\partial_I\varphi_R \U-\partial_j\W_R)- V^{*}_{j}\partial_i \U \cdot(\partial_k \varphi_R \U-\partial_k \W_R)\right]dx,$$ d'où nous avons
$$ \vert I_3 \vert \leq \int_{B_r} \vert \V^{*}\vert \vert \vec{\nabla}\otimes \U\vert \vert \vec{\nabla}\varphi_R\otimes \U\vert dx + \int_{B_r} \vert \V^{*}\vert \vert \vec{\nabla}\otimes \U\vert \vert \vec{\nabla}\otimes \W_R\vert dx.$$ Dans les deux termes à droite de cette estimation nous appliquons tout d'abord l'inégalité de Cauchy-Schwarz, ensuite l'inégalité de  H\"older  (avec $\frac{1}{2}=\frac{1}{q}+\frac{1}{p}$) et comme $\Vert \vec{\nabla}\varphi_R \Vert_{L^{\infty}}\leq \frac{c}{r-\rho}$ et $\Vert \vec{\nabla}\otimes \W_R\Vert_{L^{p}(B_r)}\leq c \Vert \vec{\nabla}\varphi_R \otimes \U\Vert_{L^p{B_r}}$, alors nous obtenons
$$ \vert I_3 \vert \leq \frac{c}{r-\rho}\left(\int_{B_r}\vert \V^{*}\vert^{q} \right)^{\frac{1}{q}} \left(\int_{B_R}\vert \vec{\nabla}\otimes \U\vert^2dx\right)^{\frac{1}{2}} \left(\int_{B_R}\vert  \U\vert^pdx\right)^{\frac{1}{p}}.$$  
Dans cette dernière estimation nous avons besoin d'étudier le terme  $\ds{\left(\int_{B_r}\vert \V^{*}\vert^{q} \right)^{\frac{1}{q}}}$. Rappelons que $\V\in \mathcal{B}^{\frac{3}{p}-\frac{1}{2}}_{\infty,\infty}(\Rt)$, et comme $\frac{9}{2}<p<6$ alors $0<\frac{3}{p}-\frac{1}{2}<\frac{1}{6}$ et donc $\V$ est une fonction $\alpha-$H\"olderienne  avec  $\alpha = \frac{3}{p}-\frac{1}{2}<\frac{1}{6}$ et nous avons ainsi $ \ds{\sup_{0<\vert x \vert<r}\frac{\vert \V(x)-\V(0)}{\vert x \vert^{\frac{3}{p}-\frac{1}{2}<\frac{1}{6}}}\leq \Vert \V \Vert_{\mathcal{B}^{\frac{3}{p}-\frac{1}{2}}_{\infty,\infty}}}$, d'où nous obtenons $\Vert V-\V(0)\Vert_{L^{\infty}(B_r)}\leq c_0 r^{\frac{3}{p}-\frac{1}{2}}$, avec $c_0=\Vert \V \Vert_{\mathcal{B}^{\frac{3}{p}-\frac{1}{2}}_{\infty,\infty}}$. \\
\\
\'A ce stade, rappelons aussi que  $\V^{*}(x)=\V(x)-\V(0)$, et comme $\frac{1}{2}=\frac{1}{q}+\frac{1}{p}$ alors nous obtenons
$$ \left(\int_{B_r}\vert \V^{*}\vert^{q} \right)^{\frac{1}{q}}\leq \Vert \V^{*}\Vert_{L^{\infty}(B_r)} r^{\frac{3}{q}}\leq c_0 r^{\frac{3}{p}-\frac{1}{2}+\frac{3}{q}} \leq c_0\,r\leq c_0 \, R ,$$  d'où, en posant la constante $c_1=\max(c,c_1)$  nous avons estimation  cherchée:
\begin{equation*}
\vert I_3 \vert \leq  c_1\frac{R}{r-\rho}\left(\int_{B_R}\vert \vec{\nabla}\otimes \U\vert^2dx\right)^{\frac{1}{2}} \left(\int_{B_R}\vert  \U\vert^pdx\right)^{\frac{1}{p}}.
\end{equation*} Observons finalement que la constante $c_1>0$ dépend de la fonction $\V$ qui dépend finalement de solution $\U$ par l'identité (\ref{def-V}) et nous écrivons alors $c_1=c_1(\U)$. \finpv
\end{enumerate} 

\subsubsection{Preuve du Lemme \ref{Lemme:morrey-L-inf-interp} page \pageref{Lemme:morrey-L-inf-interp}}\label{preuve:lemme-morrey-L-inf-interp}
Nous allons montrer l'estimation 
\begin{equation}\label{estim-aux}
\Vert \U \Vert_{\dot{M}^{3,\frac{9}{2}}}\leq c \Vert U \Vert^{\frac{2}{3}}_{\dot{M}^{2,3}}\Vert \U \Vert^{\frac{1}{3}}_{L^{\infty}}.
\end{equation}
Soient $x_0\in \Rt$ et $R>0$ et soir la boule $B(x_0,R)\subset \Rt$.  Par les inégalités d'interpolation nous avons l'estimation 
$$ \left( \int_{B(x_0,R)} \vert \U \vert^3 dx\right)^{\frac{1}{3}} \leq c \left[\left(\int_{B(x_0,R)} \vert \U \vert^2 dx \right)^{\frac{1}{2}} \right]^{\frac{2}{3}} \Vert \U \Vert^{\frac{1}{3}}_{L^{\infty}}, $$ et alors  nous pouvons écrire 
$$ R^{-\frac{1}{3}} \left( \int_{B(x_0,R)} \vert \U \vert^3 dx\right)^{\frac{1}{3}} \leq c \left[R^{-\frac{1}{2}} \left(\int_{B(x_0,R)} \vert \U \vert^2 dx \right)^{\frac{1}{2}} \right]^{\frac{2}{3}} \Vert \U \Vert^{\frac{1}{3}}_{L^{\infty}},$$ d'où, dans l'estimation à gauche ci-dessus  nous écrivons $\ds{R^{-\frac{1}{3}}= R^{\frac{3}{\frac{9}{2}}-\frac{3}{3} }}$,  dans l'estimation à droite ci-dessus nous écrivons $\ds{R^{-\frac{1}{2}}=R^{\frac{3}{3} -\frac{3}{2}}} $ pour obtenir
$$  R^{\frac{3}{\frac{9}{2}}-\frac{3}{3} } \left( \int_{B(x_0,R)} \vert \U \vert^3 dx\right)^{\frac{1}{3}} \leq c  \left[R^{\frac{3}{3} -\frac{3}{2}} \left(\int_{B(x_0,R)} \vert \U \vert^2 dx \right)^{\frac{1}{2}} \right]^{\frac{2}{3}} \Vert \U \Vert^{\frac{1}{3}}_{L^{\infty}},$$ et donc nous avons 
$$ \sup_{x_0 \in \Rt, R>0} R^{\frac{3}{\frac{9}{2}}-\frac{3}{3} } \left( \int_{B(x_0,R)} \vert \U \vert^3 dx\right)^{\frac{1}{3}} \leq c \left[ \sup_{x_0 \in \Rt, R>0} R^{\frac{3}{3} -\frac{3}{2}} \left(\int_{B(x_0,R)} \vert \U \vert^2 dx \right)^{\frac{1}{2}} \right]^{\frac{2}{3}} \Vert \U \Vert^{\frac{1}{3}}_{L^{\infty}}.$$ Ainsi, par la définition des quantités $\Vert \U \Vert_{\dot{M}^{3,\frac{9}{2}}}$ et $\Vert \U \Vert_{\dot{M}^{2,3}}$ données dans la Définition \ref{Def-Morrey} nous écrivons  l'estimation cherchée (\ref{estim-aux}). \finpv  
\subsubsection{Preuve du Lemme \ref{Lemme:pression} page \pageref{Lemme:pression}}\label{preuve:lemme-pression}
Pour le champ de vitesse  $\U=(U_i)_{i=1,2,3}$ nous écrivons  la pression $P$ comme   $$\ds{P=\sum_{i,j}^{3}\frac{1}{-\Delta}\partial_i \partial_j (U_i U_j)=\sum_{i,j}^{3}\mathcal{R}_i\mathcal{R}_j (U_i U_j)},$$ où $\mathcal{R}_i=\frac{\partial_i}{\sqrt{-\Delta}}$
est la transformée de Riesz. Ainsi,  étant donné que l'opérateur  $\mathcal{R}_i\mathcal{R}_j$ est borné dans les espaces de Morrey $\dot{M}^{p,q}(\Rt)$ avec $p\geq 2$ et $q \geq 3$ (voir le livre \cite{PGLR1}, page $171$) et par les inégalités de H\"older nous avons l'estimation
$$ \Vert P \Vert_{\dot{M}^{\frac{p}{2},\frac{q}{2}}}\leq c \sum_{i,j=1}^{3}\Vert \mathcal{R}_i \mathcal{R}_j(U_i U_j)\Vert_{\dot{M}^{\frac{p}{2},\frac{q}{2}}}\leq c \Vert \U \otimes \U \Vert_{\dot{M}^{\frac{p}{2},\frac{q}{2}}}\leq c \Vert \U \Vert^{2}_{\dot{M}^{p,q}}.$$ \finpv

	\newpage
	\thispagestyle{empty}
	
	\hbox{\includegraphics[width=8.6cm]{ED_EDMH-h.jpg}}

	\bigskip
	\noindent\fbox{\parbox{\textwidth}{
	{\bf Titre : }  Descriptions déterministes de la turbulence dans les équations de Navier-Stokes 
	
	\medskip
	{\bf Mots Clefs : }  Équations de Navier-Stokes; système stationnaire; théorie K41
		
	\medskip
	{\bf Résumé : } Cette thèse est consacrée à l'étude déterministe de la turbulence dans les équations de Navier-Stokes; et elle est divisée en quatre chapitres indépendants. 
	Le premier chapitre s'agit d'une discussion rigoureuse sur  l'étude la loi de dissipation d'énergie,  proposée par théorie de la turbulence K41, dans le cadre déterministe des équations de Navier-Stokes homogènes et incompressibles, avec une force externe stationnaire (la force ne dépende que de la variable spatiale) et posées sur l'espace $\Rt$ tout entier.  Le but de ce chapitre est de mettre en évidence le fait que si nous considérons les équations de Navier-Stokes posées sur $\Rt$ alors certains quantités physiques, nécessaires pour l'étude de la loi de dissipation de Kolmogorov, n'ont pas une définition rigoureuse et alors pour donner un sens à ces quantités on propose  de considérer  les équations de Navier-Stokes mais avec un terme additionnel d'amortissement .  Dans le cadre de ces équations de Navier-Stokes amorties, on obtient des estimations du taux de dissipation d'énergie selon la loi de dissipation de Kolmogorov. \\
	Dans le deuxième chapitre  on s'intéresse à l'étude des solutions stationnaires des équations de Navier-Stokes amorties introduites dans le chapitre précédent. Ces solutions stationnaires correspondent à un type particulier des solutions qui ne dépendent que de la variable d'espace:   la motivation pour étudier ces solutions stationnaires étant donné que  la force externe que nous considérons tout au long de cette thèse est une fonction stationnaire. Dans ce chapitre on étudie essentiellement deux propriétés des solutions stationnaires: la première propriété correspond à la stabilité de ces solutions où on montre que si l'on contrôle la force externe des équations de Navier-Stokes amorties  alors toute solution non stationnaire (qui dépend de la variable d'espace et aussi de la variable de temps) converge vers une solution stationnaire lorsque le temps tend à l'infini. La deuxième propriété porte sur l'étude de la décroissance en variable spatiale des ces solutions stationnaires. \\
	Dans le troisième chapitre on continue à étudier les solutions stationnaires des équations de Navier-Stokes, mais cette fois-ci on considère les équations de Navier-Stokes classiques (sans aucun terme d'amortissement) .  Le but de ce chapitre est d'étudier un tout autre problème relié à l'étude déterministe de la turbulence et qui porte sur la décroissance de la transformée de Fourier des solutions stationnaires. En effet, selon la théorie de la turbulence K41, si le fluide est en régime laminaire on s'attend à observer une décroissance exponentielle de la transformée de Fourier des solutions stationnaires et cette décroissance à lieu dès les bases fréquences, tandis que si le fluide est en régime turbulent alors on s'attend à observer cette même décroissance exponentielle mais seulement aux hautes fréquences. Ainsi, à l'aide des outils de l'analyse de Fourier, dans ce chapitre on donne  des descriptions précises sur cette décroissance  exponentielle fréquentiel (dans le régime laminaire et dans le régime turbulent) des solutions stationnaires. \\
	Dans le quatrième et dernier chapitre on revient aux solutions stationnaires des équations de Navier-Stokes  (on considère toujours les équations  classiques) et on étude l'unicité de ces solutions dans le cas particulier où la  force externe est nulle. En suivant essentiellement  quelques idées des travaux précédents de G. Seregin,  on étudie l'unicité des ces solutions tout d'abord  dans les cadres des espaces de Lebesgue et ensuite dans le cadre  plus général  des espaces de Morrey.
	}}
\newpage
\thispagestyle{empty}
	
	\bigskip
	\noindent\fbox{\parbox{\textwidth}{
	{\bf Title : }  Deterministic descriptions of the turbulence in the Navier-Stokes equations
	
	\medskip
	{\bf Keys words : }  Navier-Stokes equations; stationary system; K41 theory
	
	\medskip
	{\bf Abstract : } 
	This PhD thesis is devoted to deterministic study of the turbulence in the Navier-Stokes equations. The thesis is divided in four independent chapters.\\
	The first chapter involves a rigorous discussion about the energy's dissipation law, proposed by theory of the turbulence K41, in the deterministic setting  of the homogeneous and incompressible Navier-Stokes equations, with a stationary external force  (the force only depends of the spatial variable) and on the whole space $\Rt$. The energy's dissipation law,  also called the Kolmogorov's dissipation law, characterizes the energy's dissipation rate  (in the form of heat) of a turbulent fluid and this law was developed by A.N. Kolmogorov in 1941. However, its  deduction  (which uses mainly tools of statistics) is not fully understood   until our days and then an active research area   consists in studying this law in the rigorous framework  of the Navier-Stokes equations which describe in a mathematical way   the fluids motion   and in particular the movement of  turbulent fluids.  In this setting,  the purpose of this chapter is to highlight the fact that if we consider the Navier-Stokes equations  on $\Rt$ then certain physical quantities, necessary  for the study of the  Kolmogorov's dissipation law, have no a rigorous definition and then to give a sense  to these quantities we suggest to consider  the  Navier-Stokes equations with an additional damping term. In the framework of these damped equations, we obtain some estimates for the energy's dissipation rate according to the Kolmogorov's dissipation law.\\
	In the second chapter we are interested in  study the stationary solutions of the damped Navier-Stokes introduced in the previous chapter. These stationary solutions are a particular type of solutions which do not depend of the temporal variable and their study is motivated by the fact that we always  consider the Navier-Stokes equations with a stationary external force. In this chapter we study two properties of the stationary solutions: the first property concerns the stability of these solutions where we prove that if we have a control on the external force then all non stationary solution (with  depends of both spatial and temporal variables) converges toward a stationary solution. The second property concerns the decay in spatial variable of the stationary solutions. These properties of stationary solutions are a consequence of the damping term introduced in the Navier-Stokes equations. \\
	In the third chapter we still study the stationary solutions of  Navier-Stokes equations  but now we consider the classical equations  (without any additional damping term). The purpose of this chapter is to  study an other problem related to the  deterministic description of the turbulence: the frequency decay of the stationary solutions. Indeed, according to the K41 theory, if the fluid is in a laminar setting then the stationary solutions of the Navier-Stokes equations must exhibit a exponential frequency decay which starts at lows  frequencies. But, if the fluid is in  a turbulent setting then this exponential frequency decay must be observed only at highs frequencies. In this chapter, using some Fourier analysis tools, we give a precise description of this exponential frequency decay in the laminar and in the turbulent setting.  \\
	In the fourth and last chapter we return to the stationary solutions of the classical Navier-Stokes equations  and we study the uniqueness of these solutions in the particular case without any external force.  Following  some ideas of G. Seregin, we study the uniqueness of these solutions first in the framework  of  Lebesgue spaces of and then in the a general  framework of Morrey spaces.
	}}
	
	\vfill
	\hfill \includegraphics[width=1cm]{pictoParis-Saclay}

\end{document}